\documentclass[12pt]{article}

\usepackage{amssymb}

\usepackage{amsmath}

\usepackage[mathcal]{euscript}

\usepackage{amsthm}

\usepackage{graphicx}

\usepackage{physics}

\usepackage{float}

\usepackage{tikz}
\usepackage{xcolor}

\newcommand{\RR}{\mathbb{R}}

\newcommand{\XX}{\mathbb{X}}
\newcommand{\YY}{\mathbb{Y}}
\newcommand{\ZZ}{\mathbb{Z}}
\newcommand{\DD}{\mathbb{D}}
\newcommand{\CC}{\mathbb{C}}
\newcommand{\e}{\varepsilon}
\newcommand{\EE}{\mathsf{E}}
\newcommand{\PP}{\mathsf{P}}

\begin{document}

\begin{center}
{\bf \large  Algorithms of  Phase Space Reduction \\ and Asymptotics of Hitting Times for \vspace{1mm} \\  Perturbed Semi-Markov Processes}
\end{center}

\vspace{1mm}

\begin{center}
Dmitrii Silvestrov\footnote{
Dmitrii Silvestrov, Department of Mathematics, Stockholm University, Sweden \\
E-mail: silvestrov@math.su.se} \\
\end{center}

{\bf Abstract}. The paper presents new asymptotic recurrent algorithms of phase space reduction for regularly and singularly perturbed semi-Markov processes. These algorithms  give effective conditions of weak convergence for distributions and convergence of expectations for hitting times as well as recurrent  formulas for computing  the corresponding normalisation functions, Laplace transforms for limiting distributions and limits for expectations.  \\

{\bf 1. Introduction} \\

The main goals of this paper are: to  present new asymptotic recurrent  algorithms of phase space reduction for regularly and singularly perturbed finite semi-Markov processes,  to give effective conditions for weak convergence of distributions and convergence of expectations  
for hitting times for regularly and singularly perturbed finite semi-Markov processes,  and to  construct  effective  recurrent algorithms  for funding the corresponding normalisation functions, Laplace transforms for 
limiting distributions and limits for expectations  of hitting times.

Random functionals similar with hitting times  are known under different names such as  first passage times, absorption times, and  first-rare-event times in theoretical studies,  and  as  lifetimes, first failure times, extinction times, etc., in applications. Limit theorems for such functionals for Markov type processes  are objects of long term research studies.

The cases of Markov chains and semi-Markov processes with finite phase spaces are the most deeply investigated. We refer here to  selected works, which contain related results,  [1 - 7, 13, 16 - 34, 36, 38, 39,  41 - 43, 47 - 51, 53 - 55, 57, 58, 60, 62,  66 - 71, 73 - 89].

There exists also a huge  bibliography of works, which contain limit theorems for hitting times and related functionals for Markov chains and semi-Markov processes with countable and arbitrary phase spaces. 

Here, we  would like only to mention books, where one can find materials on perturbed Markov chains, semi-Markov processes with finite, countable and arbitrary phase spaces and related problems. These are, [5, 6, 9, 10, 12, 14, 15, 31, 35 - 37, 39, 40, 44 - 46, 48, 49, 51, 53, 59,  61, 64, 65, 70, 71, 79, 82 - 84, 86, 87]. 
% \cite{Ani11}, \cite{Ani12},  \cite{BLM}, \cite{Bor2}, \cite{GySi4}, \cite{Kal3}, \cite{Kal13}, 
% \cite{Ka8}, \cite{Kie2}, \cite{Kij}, \cite{KK2}, \cite{KoLi7}, \cite{KSw}, \cite{KoTu2}, \cite{KT4}, \cite{Kov2}, 
% \cite{KoKuSh}, \cite{MeT5}, \cite{ObB}, \cite{Sen4}, \cite{Sen10}, \cite{Sil1}, \cite{Sil25}, \cite{SiSi2}, 
% \cite{St9}, \cite{St10}, \cite{SS1}, \cite{YZ2}, \cite{YZ4}.  

%Seneta (1973, 2006), Silvestrov (1974, 2004), Kovalenko (1975),  Korolyuk and Turbin  (1976, 1978), Courtois %(1977), Kalashnikov  (1978, 1997), Anisimov (1988, 2008), Stewart and  Sun (1990), Korolyuk and Swishchuk %(1992, 1995), Meyn and  Tweedie (1993, 2009), Kartashov (1996),   Kovalenko,  Kuznetsov and Shurenkov %(1996), Kijima (1997), Borovkov (1998), Stewart (1998, 2001),  Yin and Zhang (1998, 2005, 2013), Korolyuk, %V.S. and Korolyuk, V.V. (1999),    Bini, Latouche and Meini (2005), Koroliuk and  Limnios (2005), Gyllenberg %and  Silvestrov (2008),   Avrachenkov, Filar and Howlett (2013) and Obzherin and Boyko (2015),  and %Silvestrov, D. and Silvestrov, S. (2017a). 

The aggregation/disaggregation  is one  of the most effective and widely used approaches in studies of asymptotics for hitting times and related functionals, especially for singularly perturbed models of Markov chains and semi-Markov processes.  Works  [3 - 5, 11, 17, 26, 32, 33, 44 - 48, 54, 63,  75, 78, 79, 80, 82 - 84, 85 - 89] represent related  asymptotic results. 
% \cite{SiAn1}, \cite{Han4}, \cite{KoTu1}, \cite{KoTu2}, \cite{Ani3}, \cite{Ani4}, \cite{Ani11}}, \cite{Turb}, \cite{CaSt}, \cite{La3}, 
% \cite{Sc5}, \cite{SS1}, \cite{Hassin}, \cite{KSw}, \cite{KK2}, \cite{AFH}, \cite{GKP}, \cite{St9}, \cite{St10}, \cite{KoLi7}, \cite{YZ2}, \cite{YZ4}, % \cite{SiS1}, \cite{SiSi2}, \cite{SiSi4D1}, \cite{Sil10}
%Simon and Ando (1961),  Hanen (1963),  Korolyuk and Turbin (1970, 1976), Anisimov (1971a, 1971b, 1988, 2008), Turbin (1971),  
%Cao and Stewart  (1985), Stewart and Sun (1990),  Latouche (1991), Schweitzer (1991), Hasin and Haviv (1992), Korolyuk and 
%Swishchuk (1992, 1995), Korolyuk, V.S. and Korolyuk, V.V. (1999),  Gambin,  Krzy\.{z}anowski and Pokarowski (2008), Stewart %(1998, 2001),  Koroliuk and  Limnios (2005),  Yin and Zhang (1998, 2005, 2013), Silvestrov, D. and Silvestrov, S. (2016b, 2017a, %2017b) and Silvestrov (2018)  represents the corresponding asymptotic results. 

We refer to books [31, 53, 72, 79]
%Silvestrov (2004), Gyllenberg and Silvestrov (2008), Silvestrov, D. and Silvestrov. %S (2017a), 
%Kovalenko, Kuznetsov and Shurenkov (1996). 
and papers [52, 63, 78],
% Schweitzer (1991), Kovalenko (1994) and Silvestrov D. and Silvestrov S. (2016),
%\cite{GySi4}, \cite{Sil2}, \cite{SiSi4D1},  \cite{KoKuSh}, \cite{Kov6}, \cite{SiS1}, \cite{Sc5}, \cite{KoKuSh}
  where one can find comprehensive  bibliographies of works in the area,  supplemented by the corresponding bibliographical remarks.

In the present paper, the above mentioned approach is used for constructing new asymptotic recurrent algorithms of phase space reduction  and applying them to study of asymptotics of hitting times  for regularly and singularly  perturbed ergodic type finite semi-Markov processes. 

We consider semi-Markov processes $\eta_{\e}(t)$ with a finite phase space $\XX$ and assume that  transition characteristics of these processes (transition probabilities $p_{\e, ij}$ of the corresponding embedded Markov chains  
$\eta_{\e,n}$  and distribution functions $F_{\e, ij}(\cdot)$ of transition  times) depend on some perturbation parameter $\e \in (0, 1]$ and converge  in some natural sense to the corresponding characteristics of the limiting (unperturbed) semi-Markov process $\eta_{0}(t)$,  as $\e \to 0$. 

The singular character of perturbation model means that:  \vspace{1mm}

{\bf (a)} the transition probabilities $p_{\e, ij}$ converge to the corresponding limiting transition probabilities $p_{0, ij}$, as $\e \to 0$, \vspace{1mm}

{\bf (b)} the phase space $\XX$  split in one or several closed classes of communicative states and possibly a class of transient states, for the limiting Markov chain $\eta_{0, n}$. \vspace{1mm}

The regular perturbation model is the particular case, where the phase space $\XX$  consists from one class of communicative states, for the Markov chain $\eta_{0, n}$. 

The ergodic type of the perturbation model means that, together with  condition {\bf (a)}, the following perturbation conditions hold: \vspace{1mm}

{\bf (c)} the distribution functions  $F_{\e, ij}(\cdot \, v_{\e, i})$ of transition  times,  normalised by some  ``local'' normalisation functions $v_{\e, i} \in [1,  \infty)$, weakly converge, as $\e \to0$,  to some limiting distribution functions $F_{0, ij}(\cdot)$, which are not concentrated at zero,  \vspace{1mm}

{\bf (d)} the first  moments $e_{\e, ij} = \int_0^\infty t F_{\e, ij}(dt)$, are finite for $\e \in (0, 1]$ and, being normalised by the same normalisation functions $v_{\e, i}$,  converge, as $\e \to 0$,  to the first moments of the corresponding limiting distribution functions $e_{0, ij} 
= \int_0^\infty t F_{0, ij}(dt)$,  also  assumed to be finite. \vspace{1mm}

The normalisation functions $v_{\e, i}$ realise some kind of initial ``space-depen\-dent''  asymptotic compression of time. Due to convergence condition {\bf (c)}, this compression of time  prevents to normalised transition times of semi-Markov processes $\eta_\e(t)$ be asymptotically stochastically unbounded or vanishing to zero random variables.

The object of our interest are the first hitting times $\tau_{\e, \DD}$ to some domain $\DD \subset \XX$, for the semi-Markov processes $\eta_\e(t)$. It is assumed that: \vspace{1mm}

{\bf (e)} $\PP_i \{ \tau_{\e, \DD} < \infty \} = 1, i \in \XX$, for $\e \in (0, 1]$. \vspace{1mm}

While,  it is possible 
that $\tau_{\e, \DD} \stackrel{\PP}{\longrightarrow} \infty$ as $\e \to 0$. 

We are interested to describe asymptotics for distributions $G_{\e, \DD, ij}(\cdot) =  \PP_i \{\tau_{\e, \DD}$ $\leq \cdot,  \eta_\e(\tau_{\e, \DD}) = j \}$ and expectations $E_{\e, \DD, ij} = \EE_i  \tau_{\e, \DD} {\rm I}( \eta_\e(\tau_{\e, \DD})= j)$ that is to  find  appropriate  normalisation functions $\check{v}_{\e, i} \in [1,  \infty)$ and $\bar{v}_{\e, i} \in [1,  \infty)$, which provide weak convergence  of normalised distributions $G_{\e, \DD, ij}(\cdot \, \check{v}_{\e, i})$  and  convergence of normalised expectations  $\bar{v}_{\e, i}^{-1} E_{\e, \DD, ij}$,   as $\e \to 0$. 

We especially are interested to find conditions,  under which one can choose the normalisation functions 
$\bar{v}_{\e, i} = \check{v}_{\e, i}$ and get convergence of normalised expectations  to the first moments of the corresponding limiting distributions for hitting times. The interest to such conditions  is caused by our intention to apply  results concerned asymptotics of hitting times for getting  ergodic and quasi-ergodic theorems for singularly perturbed semi-Markov type processes. 

We present effective asymptotic recurrent algorithms, which let us find appropriate normalisation functions, get convergence relations for distributions and expectations of hitting time, compute limiting Laplace transforms for distributions and limits for expectations of hitting times. 

The first two algorithms  let one solve the above  asymptotic problems for the main case, where an initial state of the semi-Markov processes $\eta_{\e}(t)$ belongs to domain $\overline{\DD}$. In this case,  distributions $G_{\e, \DD, ij}(\cdot), i \in \overline{\DD}, j \in \DD$ are completely determined by transition characteristics $p_{\e, ij}$ and $F_{\e, ij}(\cdot)$ of the semi-Markov process $\eta_{\e}(t)$,  for $i \in  \overline{\DD}, j \in \XX$. This makes it possible to weaken perturbation conditions   {\bf (a)},  {\bf (c)}, {\bf (d)}, and to assume that they hold only for $i \in  \overline{\DD}, j \in \XX$.

The first algorithm is based on recurrent alternating application of  procedures of two types. The first one is the procedure  of removing virtual transitions from trajectories of  perturbed  semi-Markov processes. The second one is the procedure of one-state reduction of phase space for  perturbed  semi-Markov processes. 

Let us shortly describe the first cycle of application for  these procedures.

The $1$st type procedure removes virtual transitions (of the form $i \to i $, for states $i \in \overline{\DD}$) from trajectories of the semi-Markov process $\eta_{\e}(t)$. It replaces process $\eta_{\e}(t)$ by the new semi-Markov process $\tilde{\eta}_{\e}(t)$, which has  the same phase space $\XX$, but new transition probabilities $\tilde{p}_{\e, ij}$ and new distribution functions $\tilde{F}_{\e, ij}(\cdot)$ and expectations   $\tilde{e}_{\e, ij}$ of transition times. These transition probabilities, Laplace transforms   $\tilde{\phi}_{\e, ij}(s) = \int_0^\infty e^{-st} \tilde{F}_{\e, ij}(dt)$ and expectations of transition times are expressed as simple rational transformations of  analogous characteristics for the  semi-Markov process $\eta_\e(t)$. The important property of this procedure is that the distribution $G_{\e, \DD, ij}(\cdot)$ is invariant with respect to it (i.e., this distribution coincides with the analogous distribution 
$\tilde{G}_{\e, \DD, ij}(\cdot)$ for process $\tilde{\eta}_\e(t)$),  in the case where the initial state $i \in \overline{\DD}$. 

In order to compensate the aggregation of transition times provided by this procedure,  the initial local normalisation functions for transition times $v_{\e, i}$  (used in the convergence conditions of type {\bf (c)})  are  replaced by new ones,  $\tilde{v}_{\e, i} = (1 - p_{\e, ii} )^{-1} v_{\e, i}$, for $i \in \overline{\DD}$. 

We find  some effective conditions of asymptotic comparability for the initial transition probabilities $p_{\e, ij}$  and normalisation functions $v_{\e, i}$, which provide holding of the basic convergence conditions of types {\bf (a)} -- {\bf (e)}  for semi-Markov processes $\tilde{\eta}_\e(t)$ and give  explicit rational formulas for computing the limiting transition probabilities $\tilde{p}_{0,ij}$ and expectations $\tilde{e}_{0, ij}$   of transition times  and rational formulas combined with scaling of argument for the limiting Laplace transforms $\tilde{\phi}_{0, ij}(s)$,   via analogous characteristics for the  semi-Markov process $\eta_0(t)$. 

The $2$nd type procedure realises exclusion of specially chosen state $k$ from domain 
$\overline{\DD}$. It replaces process $\tilde{\eta}_{\e}(t)$ by the new semi-Markov process 
$_k\eta_{\e}(t)$, which has  the new reduced phase space $_k\XX = \XX \setminus \{ k \}$, new transition probabilities $_kp_{\e, ij}$ and new distribution functions $_kF_{\e, ij}(\cdot)$ and expectations   $_ke_{\e, ij}$ of transition times. These transition probabilities, Laplace transforms   $_k\phi_{\e, ij}(s) = \int_0^\infty e^{-st} \, _kF_{\e, ij}(dt)$ and expectations of transition times are expressed as simple rational transformations of   analogous characteristics for the semi-Markov process $\tilde{\eta}_\e(t)$. The important property of this procedure is that the distribution $\tilde{G}_{\e, \DD, ij}(\cdot)$ is invariant with respect to it (i.e.,  this distribution coincides with the analogous distribution $_kG_{\e, \DD, ij}(\cdot)$ for process $_k\eta_\e(t)$),  in the case where the initial state $i \in \overline{\DD}_k =  \overline{\DD} \setminus \{ k \}$. 

The local normalisation functions  $_kv_{\e, i} = \tilde{v}_{\e, i}, i \in \, _k\overline{\DD}$ (used in the convergence conditions of type {\bf (c)}) are the same for the semi-Markov processes $_k\eta_\e(t)$ and $\tilde{\eta}_\e(t)$. In order to preserve the property of non-concentration the limiting distributions for transition times at zero, state $k$ is chosen among the least absorbing states in domain $\overline{\DD}$,  for which:  \vspace{1mm}
 
 {\bf (f)} quotients $\tilde{v}_{\e, k}/ \tilde{v}_{\e, i} \to \tilde{w}_{0, ki} \in [0, \infty)$ as $\e \to 0$, for $i \in \overline{\DD}$.  \vspace{1mm}
 
 We find  some effective conditions of asymptotic comparability for the initial transition probabilities $p_{\e, ij}$  and normalisation functions $v_{\e, i}$, which provide holding of the basic convergence conditions of types {\bf (a)} -- {\bf (f)}  for semi-Markov processes $_k\eta_\e(t)$ and give  explicit rational formulas for computing the limiting transition probabilities $_kp_{0,ij}$ and expectations $_ke_{0, ij}$   of transition times  and rational formulas combined with scaling of argument for the limiting Laplace transforms $_k\phi_{0, ij}(s)$,   via analogous characteristics for the  semi-Markov process $\tilde{\eta}_0(t)$. 

The procedures of removing virtual transitions and one-state reduction of phase space
are repeated recurrently according alternating sequence of types  $\langle 1, 2,  \ldots, 1, 2, 1 \rangle$, with sequential exclusion of states  $k_1, k_2, \ldots$, $k_{\bar{m} -1}$ (chosen at every step according the corresponding variant of condition {\bf(f)}) from domain $\overline{\DD} = \{ k_1, \ldots, k_{\bar{m}} \}$ ($\bar{m}$ is the number of states in domain $\overline{\DD}$), until this domain  will be reduced to the one-state set $\{ k_{\bar{m}} \}$.  

The state $k_{\bar{m}} $ is, in  some sense, the  most absorbing state in domain $\overline{\DD}$. The  above recurrent algorithm let us  construct  the semi-Markov process, which we denote as  $_{\bar{k}_{\bar{m}-1}}\tilde{\eta}_{\e}(t)$ (it depends on sequence of removed states $\bar{k}_{\bar{m}-1} = \langle k_1, k_2, \ldots, k_{\bar{m} -1} \rangle$), with  the phase space $_{\bar{k}_{\bar{m}-1}}\XX = \DD \cup \{ k_{\bar{m}} \}$,  compute the corresponding normalisation function $\check{v}_{\e, k_{\bar{m}}} $, transition  probabilities $_{\bar{k}_{\bar{m}-1}} \tilde{p}_{\e, ij}$ of the embedded  Markov chain, Laplace transforms for distribution functions  $_{\bar{k}_{\bar{m}-1}}\tilde{F}_{\e, ij}(\cdot)$ and expectations $_{\bar{k}_{\bar{m}-1}}\tilde{e}_{\e, ij}$ of transition times, for $\e \in (0, 1]$ and for the limiting case $\e = 0$, as well as to prove the corresponding convergence relations given in conditions of type {\bf (a)}, {\bf (c)}, and {\bf (d)}.

The semi-Markov process $_{\bar{k}_{\bar{m}-1}}\tilde{\eta}_{\e}(t)$ does not makes virtual transitions.  If the initial state $\eta_\e(0) = \, _{\bar{k}_{\bar{m}-1}}\tilde{\eta}_{\e}(0) = k_{\bar{m}}$, process $_{\bar{k}_{\bar{m}-1}}\tilde{\eta}_{\e}(t)$  hits domain $\DD$ just after the first jump. Therefore,  distribution $_{\bar{k}_{\bar{m}-1}}G_{\e, \DD, k_{\bar{m}} j} (\cdot) = \, _{\bar{k}_{\bar{m}-1}}\tilde{F}_{\e, k_{\bar{m}} j}(\cdot)$ $\times \, _{\bar{k}_{\bar{m}-1}} \tilde{p}_{\e, k_{\bar{m}} j}$.  
  At the same time,  the mentioned above invariance properties of distributions 
$G_{\e, \DD, ij}(\cdot)$ with respect  to the procedures of removing virtual transitions and phase space reduction  imply that distribution $G_{\e, \DD, k_{\bar{m}} j}(\cdot) = \, _{\bar{k}_{\bar{m}-1}}G_{\e, \DD, k_{\bar{m}} j} (\cdot)$ and expectation $E_{\e, \DD, k_{\bar{m}} j} = \, _{\bar{k}_{\bar{m}-1}}E_{\e, \DD, k_{\bar{m}} j}$. 

Using the above mentioned asymptotic recurrent  relations for distributions of transition times for semi-Markov processes with reduced phase spaces  and the above remarks concerned invariance properties for these distributions, we find the normalisation functions $\check{v}_{\e, k_{\bar{m}}} $,  prove  weak convergence of distributions  $G_{\e, \DD, k_{\bar{m}} j}(\cdot \, \check{v}_{\e, k_{\bar{m}}})$ and convergence of expectations $\check{v}_{\e, k_{\bar{m}}}^{-1} E_{\e, \DD, k_{\bar{m}} j}$, respectively,  to the corresponding 
limiting distributions $G_{0, \DD, k_{\bar{m}} j}(\cdot)$ and their first moments $E_{0, \DD, k_{\bar{m}} j}$.

The second  algorithm is based on backward recurrent relations, which express Laplace transforms for distributions $G_{\e, \DD, k_{n} j}(\cdot)$ and expectations $E_{\e, \DD, k_n j}$  as linear expressions, respectively, of Laplace transforms for distributions $G_{\e, \DD, k_{r} j}(\cdot), r = \bar{m}, \ldots, n +1$ and expectations $E_{\e, \DD, k_r j}, r = \bar{m}, \ldots, n +1$, for $n = \bar{m} -1, \ldots, 1$. These relations  make it possible to find  the corresponding normalisation functions, to get  convergence relations  for distributions   and expectations  of hitting times,  and to compute Laplace transforms for the corresponding limiting distributions and limits for expectations, for initial states $k_n,  n = \bar{m} -1, \ldots, 1$.

The third algorithm is based on relations, which express Laplace transforms for distributions $G_{\e, \DD, i j}(\cdot)$ and expectations $E_{\e, \DD, i j}$, for $i \in \DD$,   as linear expressions of, respectively, Laplace transforms for distributions $G_{\e, \DD, k_{r} j}(\cdot), k_r$ $\in \overline{\DD}$ and expectations $E_{\e, \DD, k_r j}$, $k_r \in \overline{\DD}$. These relations  make it possible to  find the corresponding normalisation functions,  to get convergence relations  for distributions  and expectations of hitting times, and to compute   Laplace transforms for the corresponding limiting distributions and limits for expectations, for initial states $i \in \DD$. Here, perturbation conditions  mentioned in {\bf (a)},  {\bf (c),} and {\bf (d)}  are involved also for states $i \in  \DD, j \in \XX$.

It is useful to note that, the corresponding limiting distributions $G_{0,  \DD, ij}(\cdot)$ are not concentrated at zero, for $i \in \XX, j \in \DD$ such that the corresponding limiting hitting probabilities $P_{0, \DD, ij} = G_{0, \DD, ij}(\infty) > 0$. These distributions and the 
corresponding normalisation functions  depend on the  sequence of states  $\bar{k}_{\bar{m}-1} = \langle k_1, k_2, \ldots, k_{\bar{m} -1} \rangle$ chosen for sequential exclusion from domain $\overline{\DD}$, but only up to some computable scaling factors.  

It worth also to mention some interesting phenomena, which take place for singularly perturbed models and make the difference between singularly and simpler regularly  perturbed models.

First, it is possible that the normalisation functions $\check{v}_{\e, i}$ for hitting times can tend to  $\infty$  as $\e \to 0$ with different rates,  for  initial states with different asymptotic absorbing grades.

Second, it is possible that the normalisation functions $\check{v}_{\e, i}$, used in the weak convergence relations for hitting times, guarantee convergence of  normalised expectations of hitting times to the first moment of the corresponding limiting distribution 
$G_{0, \DD, ij}(\cdot)$,  only in the case,  where the initial state $i  = k_{\bar{m}}$ is the most absorbing state  in domain $\overline{\DD}$.  However, it may be that, for some  other less absorbing initial  states, different normalisation functions 
$\bar{v}_{\e, i}$ should be used for expectations of hitting times, in order they would converge to some finite limits  $\bar{E}_{0, \DD, i j}$.  Some additional conditions should be required, in order expectations of hitting times  normalised by functions $\check{v}_{\e, i}$ would converge to the first moments of the corresponding  limiting distributions for hitting times.

It worth also to additionally comment  the question about mentioned above asymptotic comparability conditions for the initial transition probabilities $p_{\e, ij}$ and the initial normalisation functions $v_{\e, i}$, which should be imposed on them, in order for the described above 
asymptotic recurrent algorithms to work properly.  

The most effective variant of such conditions is based on the notion of a complete family of asymptotically comparable functions. Some family
${\cal H} = \{ h(\cdot) \}$ of  positive functions defined on interval $(0, 1]$ is a complete family of asymptotically comparable functions if it possesses two properties: {\bf (1)} ${\cal H}$ is closed with respect to summation, multiplication and division operations, {\bf (2)} there exists limit $\lim_{\e \to 0} h(\e) = a[h(\cdot)] \in [0, \infty]$, for any function $h(\cdot)  \in {\cal H}$.  

The simplest example is the family ${\cal H}_1 = \{ h(\cdot) \}$, for which there exist, for any function $h(\cdot)  \in {\cal H}_1$,  constants  $a_h > 0, b_h \in (- \infty, \infty)$ such that $h(\e)/a_h\e^{b_h} \to 1$ as $\e \to 0$. 

Other examples of complete families of asymptotically comparable functions, based on natural combinations of polynomial, exponential and logarithmic functions,  are given in Appendix A. 

The asymptotic comparability condition mentioned above is based on the assumptions that:  \vspace{1mm}

{\bf (g)}  there exist sets $\YY_i \subseteq \XX, \, i \in \XX$ such that the transition probabilities $p_{\cdot, ij},  j \in \YY_i, \,  i \in \XX$ belong to some complete family of asymptotically comparable functions ${\cal H}$, while  $p_{\e, ij} = 0, \e \in (0, 1]$, for $j \in \overline{\YY}_i, i \in \XX$,  \vspace{1mm}

{\bf (h)}   the initial normalisation functions \, $v_{\cdot, i}, \, j \in \XX$ \, also  belong to the family ${\cal H}$.  

\vspace{1mm}

The described above asymptotic recurrent algorithms, based on conditions of types {\bf (a)} -- {\bf (h)}, provide existence of all intermediate limits involved in the asymptotic recurrent  relations used for  computing  limiting Laplace transforms for distributions and limits for expectation of hitting times. Computing of the corresponding normalisation functions do require only recurrent application of some rational formulas to the initial normalisation functions $v_{\e, i}$ and the initial transition probabilities $p_{\e, ij}$.  

The above algorithms became more effective in the cases,  where the asymptotic comparability conditions  {\bf (g)}   and 
{\bf (h)} are based  on some concrete  complete family of asymptotically comparable functions ${\cal H}$ admitting effective computation of limits related to summation, multiplication and division operations. For example,  this relates to the above mentioned family ${\cal H}_1$, and families  ${\cal H}_2$ and  ${\cal H}_3$ described in Appendix A.
 
In this case, computing of the limiting Laplace transforms for distributions  and limits for expectations of hitting times requires only recurrent application of some rational transformations to limiting transition probabilities and expectations  and rational transformations combined with scaling of argument to limiting Laplace transforms used in the described above asymptotic recurrent algorithms. Coefficients of these rational transformations and scaling parameters can be computed with the use of recurrent  formulas based on operational rules (for computing limits for sums, products and quotients) for functions from the corresponding complete family of asymptotically comparable functions. 

This paper includes ten sections. In Section 2, we formulate the basic perturbation conditions. In Section 3, we describe the procedure of removing virtual transitions for regularly and  singularly perturbed semi-Markov processes and describe weak asymptotics for distributions of hitting times for the simplest case, where an initial state for the corresponding semi-Markov processes belongs to  a one-state domain $\overline{\DD}$ (Theorem 1).  In Section 4, we describe the procedure  of one-state reduction of phase space for regularly and  singularly perturbed semi-Markov processes. In Section 5, we shortly repeat the description of procedure of removing virtual transitions for  regularly and  singularly perturbed semi-Markov processes resulted by one-state reduction of phase space and  describe weak asymptotics for distributions  of hitting times for the case, where an initial state belongs  to a two-states domain 
$\overline{\DD}$ (Theorems 2 and 3). In Section 6, we describe the multi-step algorithm of recurrent removing virtual transitions and reduction of phase space  for regularly and  singularly perturbed semi-Markov processes. In Section 7, we present weak limit theorems for hitting times and recurrent formulas for computing normalisation functions and Laplace transforms for limiting distributions of hitting times, for the general case, where an initial state belongs to a  multi-states domain $\overline{\DD}$ (Theorems 4 -- 7).  In Section 8, we  present weak limit theorems for hitting times and recurrent formulas for computing normalisation functions, and Laplace transforms for limiting distributions of hitting times, for the case where an initial state belongs to domain $\DD$ (Theorem 8). In Section 9, we present the corresponding limit theorems  for expectations of hitting times (Theorems 9 -- 12). In Section 10, we discuss some generalisations of results presented in previous sections, in particular, those connected with conditions of non-concentration of distributions of transition and hitting times at zero, reward interpretations of  hitting times and related multivariate and  real-valued reward functionals of hitting type.  Also, numerical examples illustrating the asymptotic recurrent algorithms of phase space reduction and asymptotic results presented in the paper are given in this section. In Appendix A, a definition and some basic properties  of complete families of asymptotically comparable functions and examples of such families are given, and  their connection with asymptotic recurrent  algorithms for perturbed semi-Markov processes  is described. \\

{\bf  2. Perturbed Semi-Markov  Processes and Hitting Times} \\

In this section, we define perturbed semi-Markov processes and formulate basic perturbation and convergence conditions. We also define random functionals known as hitting times, which are the main objects of our interest.  We shall see that asymptotic analysis of hitting times to some domain $\DD$ is essentially different for the cases, where an initial state belongs to domain $\overline{\DD}$ or to domain $\DD$. In Sections 2 -- 7 we consider the main first case, where an initial state belongs to domain $\overline{\DD}$.

\vspace{1mm}

{\bf 2.1 Perturbed semi-Markov processes and hitting times}. Let  $\XX = \{1, \ldots, M \}$ be a finite  set.

Let also, $(\eta_{\e, n}, \kappa_{\e, n}), n = 0, 1, \ldots$ be, for every $\e \in (0, 1]$  a Markov renewal process, i.e., a homogeneous Markov chain,  with a phase space $\XX \times [0, \infty)$ and  transition probabilities, for $(i, s), (j, t) \in \XX \times [0, \infty)$,
\begin{equation}\label{trans}
Q_{\e, ij}(t) =  \PP \{ \eta_{\e, 1} = j,  \kappa_{\e, 1} \leq t / \eta_{\e, 0} = i,  \kappa_{\e, 0} = s \}.
\end{equation}

The first component of the Markov renewal process $\eta_{\e, n}, n = 0, 1, \ldots$ is itself  a homogeneous, so-called embedded, Markov chain,  with the phase space $\XX$ and transition probabilities, for $i, j \in \XX$,
\begin{equation}\label{transa}
p_{\e, ij} =  Q_{\e, ij}(\infty) = \PP \{ \eta_{\e, 1} = j,   / \eta_{\e, 0} = i \}.
\end{equation}

The above Markov renewal process is used to define a semi-Markov process, 
\begin{equation}\label{semi}
\eta_{\e}(t) = \eta_{\e, \nu_{\e}(t)}, t \geq 0,
\end{equation}
where $\zeta_{\e, n} = \kappa_{\e, 1} + \cdots + \kappa_{\e, n}, n = 1, 2, \ldots, \zeta_{\e, 0} = 0$, are the corresponding instantss of jumps, and $\nu_\e(t) = \max(n \geq 1: \zeta_{\e, n}  \leq t)$ is the number of jumps in an interval $[0, t], t \geq 0$,  for the above semi-Markov process.

Let us introduce distributions functions of sojourn times for the semi-Markov process $\eta_\e(t)$ that are, for $i \in \XX$, 
 \begin{equation}\label{sojou}
 F_{\e, i}(t) = \PP_i \{\kappa_{\e, 1} \leq t \} = \sum_{j \in \XX}Q_{\e, ij}(t), \ t \geq 0. 
 \end{equation}
 
 Here and henceforth, notations $\PP_i$ and $\EE_i$ are used for conditional probabilities and expectations, under condition $\eta_\e(0) = i$. 

We  assume  that the following regularity condition holds:
\begin{itemize}
\item [${\bf A}$:]  $F_{\e, i}(0)  < 1, i \in \XX$, for $\e \in (0, 1]$.
\end{itemize}

Condition ${\bf A}$   excludes a.s. instant sequential jumps and  guarantees that $\PP_i \{ \lim_{n \to \infty} \zeta_{\e, n} = \infty \} = 1$, for any $i \in \XX$ and $\e \in (0, 1]$. Therefore, process $\eta_{\e}(t), t \geq 0$ is well defined at the interval $[0, \infty)$, for every $\e \in (0, 1]$.

Let $\DD$ be a nonempty subset of $\XX$. 

The object of our interest are random functionals, which are known as hitting times,
\begin{equation}
\tau_{\e, \DD} = \sum_{n = 1}^{\nu_{\e, \DD}} \kappa_{\e, n}, \ {\rm where} \ \nu_{\e, \DD} = \min(n \geq 1: \eta_{\e, n} \in \DD).
\end{equation}

Let us assume that the following condition holds:
\begin{itemize}
\item [${\bf B}$:] {\bf (a)} $p_{\e, ij} > 0, \e \in (0, 1]$ or $p_{\e, ij} = 0, \e \in (0, 1]$, for every $i \in \overline{\DD}, j \in \XX$, {\bf (b)} or any $i \in \overline{\DD}$, there exists a chain of states $i = j_0, j_1 \in \overline{\DD}, \ldots, j_{n_{i}- 1} \in \overline{\DD},  j_{n_{i}} \in \DD$ such that $\prod_{1 \leq l \leq n_{i}} p_{1, j_{l-1} j_l} > 0$.
\end{itemize}

Due to condition ${\bf B}$ {\bf (a)}, product $\prod_{1 \leq l \leq n_{ij}} p_{\e, j_{l-1} j_l}$ is a positive number or equals zero  simultaneously for all $\e \in (0, 1]$. 

Condition ${\bf B}$  {\bf (b)}  implies that, for $\e = 1$, and, due the above remark, for any $\e \in (0, 1]$
and $i \in \overline{\DD}$,  
\begin{equation}\label{hitta}
\PP_i \{ \nu_{\e, \DD} < \infty \} = \PP_i \{ \tau_{\e, \DD} < \infty \}  = 1.
\end{equation} 

Moreover, in this case relation (\ref{hitta}) also holds for  any $\e \in (0, 1]$
and $i \in \DD$. Indeed, in this case, 
\begin{align}\label{hittany}
\PP_i \{ \nu_{\e, \DD} < \infty \} & = \PP_i \{ \tau_{\e, \DD} < \infty \}  \vspace{1mm} \nonumber \\
& = \PP_i \{ \eta_{\e, 1} \in \DD \} + \sum_{k \in \overline{\DD}}  \PP_k \{ \nu_{\e, \DD} < \infty \} p_{\e, ik} 
\vspace{1mm} \nonumber \\
& =  \sum_{r \in \DD} p_{\e, ir} + \sum_{k \in \overline{\DD}}  p_{\e, ik}  = 1. 
\end{align} 
 
Let us introduce, for $i \in \XX, j \in \DD$ and $\e \in (0, 1]$, distribution functions,
\begin{equation}\label{wurew}
G_{\e,\DD, ij}(t) = \PP \{ \tau_{\e, \DD} \leq t, \eta_{\e}( \tau_{\e, \DD}) = j \}, \, t \geq 0,   
\end{equation}
Laplace transforms,
\begin{equation}\label{wurewop}
\Psi_{\e. \DD, ij}(s) = \int_0^\infty e^{-st} G_{\e, \DD, ij}(dt), s \geq 0, 
\end{equation}
expectations,
\begin{equation}\label{wurewlo}
E_{\e,\DD, ij} = \EE \tau_{\e, \DD} {\rm I}( \eta_{\e}( \tau_{\e, \DD}) = j),  
\end{equation}
and hitting probabilities,  
\begin{equation}\label{wurewlom}
P_{\e,\DD, ij} = \PP \{ \eta_{\e}( \tau_{\e, \DD}) = j \} = G_{\e,\DD, ij}(\infty).  
\end{equation}

We are interested to find conditions, which would imply that there exist normalisation functions $\check{v}_{\e, i} > 0, i \in \XX$ such that the following weak convergence relation holds, for $i \in \XX, j \in \DD$,
\begin{equation}\label{convera}
G_{\e, \DD, ij}(\cdot \, \check{v}_{\e, i}) \Rightarrow  G_{0, \DD, ij}(\cdot)  = F_{0, \DD, ij}(\cdot)P_{0, \DD, ij} 
\  {\rm as} \ \e \to 0,
\end{equation}
where: (a) $F_{0, \DD, ij}(\cdot)$ are proper distribution functions such that $F_{0, \DD, ij}(0) < 1$, for $i \in \XX, j \in \DD$,   (b)  $P_{0, \DD, ij}, j \in \DD$ is a discrete distribution, i.e.,  $P_{0, \DD, ij} \geq 0, j \in \DD$ and $\sum_{j \in \DD}P_{0, \DD, ij} = 1$, for $i \in \XX$.

We also would like to find conditions, which  would imply that there exist normalisation functions $\bar{v}_{\e, i} > 0, i \in \XX$ such that the following convergence relation holds, for $i \in \XX, j \in \DD$,
\begin{equation}\label{converan}
\bar{v}_{\e, i}^{-1}E_{\e, \DD, ij} \to  \bar{E}_{0, \DD, ij} \  {\rm as} \e \to 0,
\end{equation}
where limits  $\bar{E}_{0, \DD, ij} \in  (0, \infty)$, if  $P_{0, \DD, ij} > 0$. 

We are also interested to find additional conditions, which  would make it possible to take the normalisation functions $\bar{v}_{\e, i}  = \check{v}_{\e, i}$ and get, for $i \in \XX, j \in \DD$, limits for expectations,  
\begin{equation}\label{converanuk}
  \bar{E}_{0, \DD, ij}  = E_{0, \DD, ij} = \int_0^\infty t G_{0, \DD, ij}(dt). 
\end{equation}

We are going to present such conditions and to  construct effective recurrent algorithms, which would let us compute the above normalisation functions, Laplace transforms $\Psi_{0.\DD, ij}(s) =
\int_0^\infty e^{-st} G_{0, \DD, ij}(dt), s \geq 0$, limits for expectations  $\bar{E}_{0, \DD, ij}$ and 
$E_{0, \DD, ij}$, and the limting hitting probabilities  $P_{0, \DD, ij}$, for $i \in \XX, j \in \DD$. 
\vspace{1mm}

The case, where $\DD = \XX$ is trivial, since $\nu_{\e, \XX} = 1$ and  $\tau_{\e, \XX} = \kappa_{\e, 1}$.  In this case, the asymptotic relations given in conditions ${\bf C}$  -- ${\bf E}$ formulated below solve the asymptotic problems formulated above,  

We, therefore, assume that the number of states in the phase space $M \geq 2$ and domain $\DD \subset \XX$, i.e., the number of states in this domain $m < M$. Respectively, the number of states $\bar{m} = M - m$ in domain 
$\overline{\DD}$ is at least $1$.

It is obvious that distributions  $G_{\e, \DD, ij}(t), t \geq 0, j \in \DD,  i \in \overline{\DD}$ are completely determined by transition probabilities $Q_{\e, kr} (t), t \geq 0,   r \in \XX, k \in \overline{\DD}$. This means that distributions $G_{\e, \DD, ij}(t), t \geq 0,  j \in \DD, i \in \overline{\DD}$  coincide for semi-Markov processes with any form of transition probabilities $Q_{\e, kr} (t), t \geq 0,  r \in \XX, k \in \DD$.

This makes it possible to essentially simplify the model in the case, where we are interested to investigate asymptotics of hitting times for some fixed domain $\DD$ and only for initial states $i \in \overline{\DD}$. This case in considered in Sections 2 -- 7. For example, we can assume  that the transition probabilities  of semi-Markov processes $\eta_\e(t)$ satisfy the following condition:
\begin{itemize}
\item [${\bf \hat{B}}$:] $Q_{\e, kr}(t) =  \frac{1}{m} {\rm I}( r \in \DD){\rm I}(t \geq 1)$, for $t \geq 0,  r \in \XX, k \in \DD$ and $\e \in (0, 1]$.
\end{itemize}

Note that, in this case,  relation $F_{\e, i}(0)  < 1$, given in condition ${\bf A}$,  remains to hold, for $i \in \DD$ and $\e \in (0, 1]$. Also, condition ${\bf B}$ {\bf (a)} remains to hold, for  $i \in \DD$.  Condition ${\bf B}$ {\bf (b)} is not influenced by 
condition ${\bf \hat{B}}$.
\vspace{1mm}

 {\bf 2.2. Perturbation conditions}. We assume that the following condition holds:
\begin{itemize}
\item [${\bf C}$:]  $p_{\e, ij} \to p_{0, ij}$ as $\e \to 0$, for $i \in \overline{\DD},  j \in \XX$. 
\end{itemize}

Here and henceforth, symbol $\e \to 0$ means that $0 < \e \to 0$.

In the case, where condition  ${\bf \hat{B}}$  is assumed to hold, the convergence relation given in condition ${\bf C}$ also holds, for $i \in \DD,  j \in \XX$.

Since matrix ${\mathbf P}_\e = \| p_{\e, ij}  \|$ is stochastic, conditions ${\bf \hat{B}}$ and ${\bf C}$   imply that  
matrix ${\mathbf P}_0 = \| p_{0, ij}  \|$ is also stochastic. Let $\eta_{0, n}, n = 0, 1, \ldots$ be a Markov chain with the phase space $\XX$ and the matrix of transition probabilities ${\mathbf P}_0$.  Condition ${\bf C}$  makes it possible to interpret the Markov chains $\eta_{\e, n}$, for $\e \in (0, 1]$, as a perturbed version of the Markov chain $\eta_{0, n}$.

In the case, where condition  ${\bf \hat{B}}$  is assumed to hold, conditions ${\bf B}$ and ${\bf C}$ also imply that  domain $\DD$ is a closed class of communicative states, while domain $\overline{\DD}$ is, either, a class of transient states or a closed class of states, which can possess an arbitrary communicative structure, for the limiting Markov chain $\eta_{0, n}$. 

The case, where  domain  $\overline{\DD}$ is one class of communicative states relates to the model with regular perturbations. The case, where some states in domain $\overline{\DD}$ do not communicates, relates to the model with singular perturbations. 

Let us introduce, for $i \in \XX$ and $\e  \in [0, 1]$,  sets,
\begin{equation}\label{hore}
\YY_{\e, i} = \{ j \in \XX: p_{\e, ij} > 0 \}.
\end{equation}

Condition ${\bf B}$  implies that,  for $i \in \overline{\DD}$ and $\e \in (0, 1]$,
\begin{equation}\label{bywerta}
\YY_{\e, i} = \YY_{1, i}.
\end{equation} 

If, also,  condition ${\bf C}$ holds, then, for $i \in \overline{\DD}$,
\begin{equation}\label{bywertak}
\YY_{0, i} \subseteq \YY_{1, i}. 
\end{equation} 

In the case, where condition  ${\bf \hat{B}}$  is assumed to hold, sets $\YY_{\e, i} = \YY_{1, i} = \DD, \e \in [0, 1]$, for $i \in \DD$. 

The transition probabilities $Q_{\e, ij}(t)$ can, for every $\e \in (0, 1]$,  be represented in the following  form,  for $t \geq 0, i, j \in \XX$,
\begin{equation}\label{transasma}
Q_{\e, ij}(t) =   F_{\e, ij}(t) p_{\e, ij},       
\end{equation}
where 
\begin{equation}\label{transabasa}
F_{\e, ij}(t) = \PP \{ \kappa_{\e, 1} \leq t /   \eta_{\e, 0} = i, \eta_{\e, 1} = j \}.
\end{equation}

If $j \in \YY_{1, i}, i \in \XX$, then,
\begin{equation}\label{herew}
F_{\e, ij}(t) = p_{\e, ij}^{-1}Q_{\e, ij}(t), t \geq 0. 
\end{equation}

Let  us recall  the distribution function of the sojourn time in state $i \in \XX$ for the semi-Markov process  $\eta_{\e}(t)$, 
\begin{align}\label{diotr}
F_{\e, i}(t)  & = \PP \{  \kappa_{\e, 1} \leq t / \eta_{\e, 0} = i \}  \vspace{2mm}  \nonumber \\ 
& = \sum_{j \in \XX}Q_{\e, ij}(t) \vspace{2mm}    = \sum_{j \in \YY_{1, i}}Q_{\e, ij}(t) = \sum_{j \in \YY_{1, i}}F_{\e, ij}(t)p_{\e, ij}.
\end{align}

If  $j \in \overline{\YY}_{1, i}, i \in \XX$, an arbitrary distribution function  concentrated on $[0, \infty)$ can play the role of $F_{\e, ij}(t)$. 
We use the  standard variant and choose,
\begin{equation}\label{diotra}
 F_{\e, ij}(t) =  F_{\e, i}(t),\  t \geq 0.
 \end{equation}
 
Let  $v_{\e, i}, \e \in (0, 1]$  be, for every $i \in \overline{\DD}$, a function taking values in interval 
$[1, \infty)$. We shall refer to functions $v_{\e,  i}$   as to ``initial local''  normalisation functions.

We  assume that the following condition holds: 
\begin{itemize}
\item [${\bf D}$:]   {\bf (a)} $F_{\e, ij}(\cdot \, v_{\e,  i}) = \PP \{ \kappa_{\e, 1}/ v_{\e,  i}   \leq \cdot /  \eta_{\e, 0} = i, \eta_{\e, 1} = j  \} \Rightarrow F_{0, ij}(\cdot)$ as $\e \to 0$, for $j \in \YY_{1, i}, i \in \overline{\DD}$, {\bf (b)} $F_{0, ij}(\cdot), j \in \YY_{1, i}, i \in \overline{\DD}$ are proper distributions functions such 
that $F_{0, ij}(0) < 1, j \in \YY_{1, i}, i \in \overline{\DD}$,   {\bf (c)} $v_{\e, i} \to v_{0, i} \in [1, \infty]$ as $\e \to 0$, for $i \in \overline{\DD}$. 
\end{itemize}

Here and henceforth, symbol $\Rightarrow$ is used to  denote  weak convergence of distribution functions. 

Condition ${\bf B}$, ${\bf C}$ and ${\bf D}$ imply that the following relation of weak convergence holds, for $i \in \overline{\DD}$,
\begin{align}\label{alsok}
F_{\e, i}( \cdot v_{\e,  i}) & = \sum_{j  \in \YY_{1, i}}  F_{\e, ij}(\cdot v_{\e,  i} ) p_{\e, ij} \makebox[35mm]{} \vspace{2mm} \nonumber \\
& \Rightarrow \sum_{j  \in \YY_{1, i}}  F_{0, ij}( \cdot) p_{0, ij} 
\vspace{2mm} \nonumber \\
%\end{align*}
%\begin{align}
&  = \sum_{j  \in \YY_{0, i}}  F_{0, ij}( \cdot) p_{0, ij} =  F_{0, i}(\cdot) \ {\rm as} \ \e \to 0.
\end{align}

Obviously, $F_{0, i}(t)$ is a proper distribution function such that  $F_{0, i}(0) < 1$, for $i \in  \overline{\DD}$.

Condition ${\bf D}$  {\bf (b)} guarantees  that the initial local normalisation functions $v_{\e, i}$ do not deliver excessive levels of local normalisation. 

Let us introduce Laplace transforms, for $i, j \in \XX$ and $\e \in [0, 1]$,  
\begin{equation}\label{lapla}
\phi_{\e, ij}(s) = \int_0^\infty e^{-s t} F_{\e, ij}(dt), \ s \geq 0.  
\end{equation}

Condition ${\bf D}$   can be also re-formulated  in the following  equivalent form:
\begin{itemize}
\item [${\bf D}'$:]   {\bf (a)} $\phi_{\e, ij}(s /v_{\e,  i})  \to \phi_{0, ij}(s)$ as $\e \to 0$, for $s \geq 0$ and $j \in \YY_{1, i}, i \in \XX$, {\bf (b)} $\phi_{0, ij}(s) = \int_0^\infty e^{-s t} F_{0, ij}(dt), s \geq 0, j \in \YY_{1, i}, i \in \overline{\DD}$ are Laplace transforms of proper distributions functions such that $F_{0, ij}(0) < 1, j \in \YY_{1, i}, i \in \overline{\DD}$,  {\bf (c)}  $v_{\e,  i} \to v_{0,  i} \in [0, \infty]$ as $\e \to 0$, for $i \in \overline{\DD}$. 
\end{itemize} 

It is useful to note that relation $F_{0, ij}(0) < 1$ holds if and only if $\phi_{0, ij}(s) < 1, s > 0$. Note, also,  that these  inequalities hold  if $\phi_{0, ij}(s_0) < 1$, for some $s_0 > 0$. 

Note also that, if condition  ${\bf \hat{B}}$  is assumed to hold, the asymptotic relations given in conditions ${\bf D}$ and ${\bf D}'$ also hold for $i \in \DD$. In this case, the initial normalisation functions $v_{\e, i} = 1, \e \in [0, 1]$, the limiting distribution functions $F_{0, ij}(t)= F_{0, i}(t) = {\rm I}(t \geq 1), t \geq 0$, and the limiting Laplace transforms $\phi_{0, ij}(s) = \phi_{0, i}(s)  = e^{-s}, s \geq 0$, for $i, j \in \DD$.
 
Let now introduce the semi-Markov transition probabilities, for $t \geq 0$, $i, j \in \XX$,
\begin{equation}\label{transad}
Q_{0, ij}(t) = F_{0, ij}(t) p_{0, ij}. 
\end{equation} 

Probabilities $Q_{0, ij}(t), t \geq 0, i, j \in \XX$ can serve as transition probabilities for a Markov renewal process  $(\eta_{0, n}, \kappa_{0, n}), n = 0, 1, \ldots$, with the phase space $\XX \times [0, \infty)$. 
We also can define the corresponding semi-Markov process,
$\eta_{0}(t) = \eta_{0, \nu_{0}(t)}, t \geq 0$,
where $\zeta_{0, n} = \kappa_{0, 1} + \cdots + \kappa_{0, n}, n = 1, 2, \ldots, \zeta_{0, 0} = 0$, are the corresponding instants of jumps, and 
$\nu_0(t) = \max(n \geq 1: \zeta_{0, n}  \leq t)$ is the number of jumps in an interval $[0, t], t \geq 0$ for the above semi-Markov process.  

According condition ${\bf D}$, probabilities $F_{0, i}(0) < 1, i \in \XX$. This implies that $\PP_i\{ \lim_{n \to \infty} \zeta_{0, n} = \infty \} = 1$, for $i \in \XX$. Thus,  the process $\eta_0(t)$ is well 
defined on the time interval $[0, \infty)$. 

Conditions ${\bf C}$ and  ${\bf D}$ make it possible to interpret the semi-Markov process $\eta_{\e}(t)$, for $\e \in (0, 1]$, as a perturbed version of the semi-Markov process $\eta_{0}(t)$.

Let us also introduce expectations of inter-jump  times, for $i, j \in \XX$ and $\e \in [0, 1]$, 
\begin{equation}\label{expr}
e_{\e, ij} = \int_0^\infty t F_{\e, ij}(dt) 
\end{equation}
and
\begin{equation}\label{expromt}
e_{\e, i} = \int_0^\infty t F_{\e, i}(dt)  = \sum_{j  \in \YY_{\e, i}} e_{\e, ij} p_{\e, ij}. 
\end{equation}

Finally, we assume that the following condition holds:
\begin{itemize}
\item [${\bf E}$:] {\bf (a)} $e_{\e, ij}  < \infty, j \in \YY_{1, i}, i \in \XX$, for every $\e \in (0, 1]$,  {\bf (b)} $e_{\e, ij}/ v_{\e, i}   \to  e_{0, ij} < \infty$ as $\e \to 0$, for $j \in \YY_{1, i}, i \in \XX$.
\end{itemize}

Thus, it is assumed not only weak convergence of distribution functions $F_{\e, ij}(\cdot v_{\e,  i})$ to distribution functions $F_{0, ij}(\cdot )$, 
as $\e \to 0$,  but also convergence of their first moments $e_{\e, 0, ij} / v_{\e,  i}$ as $\e \to 0$  to the corresponding first moments $e_{0, ij}$ of the distribution function $F_{0, ij}(\cdot )$, for $j \in \YY_{1, i}, i \in \XX$. 

Note also that conditions ${\bf D}$  and  ${\bf E}$   imply that 
$e_{0, ij} \in (0, \infty)$, for  $j \in \YY_{1, i}, i \in \XX$.  

Conditions ${\bf B}$, ${\bf C}$, ${\bf D}$,  and ${\bf E}$  also imply the following asymptotic relation  holds, for $i \in \XX$,
\begin{align}\label{alsokala}
\frac{e_{\e, i}}{v_{\e,  i}}  & = \sum_{j  \in \YY_{1, i}}  \frac{e_{\e, ij}}{v_{\e,  i}} p_{\e, ij}  \vspace{2mm} \nonumber \\
& \to \sum_{j  \in \YY_{1, i}}  e_{0, ij} p_{0, ij} =   \sum_{j  \in \YY_{0, i}}  e_{0, ij} p_{0, ij}
\vspace{2mm} \nonumber \\
%\end{align}
%\begin{align}
& = e_{0, i} \ {\rm as} \ \e \to 0.
\end{align}

It is useful to note that expectation $e_{0, i} \in (0, \infty)$, for $i \in \overline{\DD}$.

Note also that, if  ${\bf \hat{B}}$  is assumed to hold, the asymptotic relations given in condition ${\bf E}$ also hold for $i \in \DD$. In this case, the  limiting expectations $e_{0, ij} = e_{0, i} = 1$, for $i, j \in \DD$.
\\

{\bf 3. Removing of Virtual Transitions  for  Perturbed \\ \makebox[10mm]{} Semi-Markov Processes} \\ 

In this section, we describe the procedure of removing virtual transitions of the form $i \to i$ from  trajectories of perturbed semi-Markov processes. We show that hitting times are asymptotically invariant with respect to this procedure. We also formulate conditions which guarantee that basic perturbation conditions imposed on the initial semi-Markov processes also hold for the semi-Markov processes with removed virtual 
transitions. Also, we  give explicit formulas for re-calculating  normalisation functions, limiting distributions and expectations  
in the corresponding perturbation conditions for  the semi-Markov processes with removed virtual 
transitions. \vspace{1mm}

{\bf 3.1 Procedure of removing virtual transitions for perturbed semi-Markov processes}. Let us assume that $\e \in (0, 1]$
and  conditions ${\bf A}$ and ${\bf B}$  hold. 

Let us define stopping times for Markov chain $\eta_{\e, n}$ that are, for  $r = 0, 1, \ldots$,
\begin{equation}\label{stop}
\theta_{\e}[r] = {\rm I}(\eta_{\e, r} \in \overline{\DD}) \min(n > r: \eta_{\e, n} \neq \eta_{\e, r}) + {\rm I}(\eta_{\e, r} \in \DD) (r+1).
\end{equation}

By the definition,  $\theta_{\e}[r]$ is, either the first after $r$ moment of change of state $\eta_{\e, r}$ by the Markov chain 
$\eta_{\e, n}$, if  $\eta_{\e, r} \in \overline{\DD}$,  or $r +1$,  if $\eta_{\e, r} \in \DD$.

Let us also  define sequential stopping times,
\begin{equation}\label{opasd}
\mu_{\e, n} =\theta_{\e}[\mu_{\e, n-1}], n = 1, 2, \ldots, \ {\rm where} \ \mu_{\e, 0} = 0.
\end{equation}

Let us now construct a new Markov renewal process $(\tilde{\eta}_{\e,  n}, \tilde{\kappa}_{\e,  n}), n = 0, 1, \ldots$ with the phase space $\XX \times [0, \infty)$ using the following recurrent relations,
\begin{equation}\label{recur}
(\tilde{\eta}_{\e,  n}, \tilde{\kappa}_{\e, n}) = \left\{ 
\begin{array}{cll}
(\eta_{\e, 0}, 0), & \text{for} \ n = 0,  \vspace{2mm} \\
(\eta_{\e,  \mu_{\e, n}}, \sum_{l =  \mu_{\e,  n-1} + 1}^{\mu_{\e, n}} \kappa_{\e, l}),  & \text{for} \ n = 1, 2, \ldots.  
\end{array}
\right.
\end{equation}
 
We also can define the corresponding semi-Markov process,
\begin{equation}\label{semin}
\tilde{\eta}_{\e}(t) = \tilde{\eta}_{\e,  \tilde{\nu}_{\e}(t)}, t \geq 0,
\end{equation}
where $\tilde{\zeta}_{\e,  n} = \tilde{\kappa}_{\e,  1} + \cdots + \tilde{\kappa}_{\e,   n}, n = 1, 2, \ldots, \tilde{\zeta}_{\e,   0} = 0$, are the corresponding instants of jumps, and $\tilde{\nu}_{\e}(t) = \max(n \geq 1: \tilde{\zeta}_{\e,  n}  \leq t)$ is the number of jumps in an interval $[0, t], t \geq 0$ for the above semi-Markov process. 
  
Below, we use symbols $*$ and the upper index $^{(*n)}$ to indicate the corresponding convolution operation and $n$ times self convolution for  proper or improper distribution functions.  

The definition of stopping times $\mu_{\e, n}$ implies that the transition probabilities for the above 
Markov renewal processes are determined by the following relation,  
\begin{equation*}
\tilde{Q}_{\e, ij}(t) =  \PP \{ \tilde{\eta}_{\e, 1} = j,  \tilde{\kappa}_{\e, 1} \leq t  / \tilde{\eta}_{\e, 0} = i \} \makebox[46mm]{}
\end{equation*}
\begin{equation}\label{transaka}
\makebox[15mm]{} = \left\{
\begin{array}{cll}
Q_{\e, ij}(t)  & \text{for} \  t \geq 0, j \in  \XX,  i \in \DD, \vspace{2mm} \\
0 & \text{for} \  t \geq 0, j = i, i \in \overline{\DD}, \vspace{2mm} \\
\sum_{n = 0}^\infty Q^{(*n)}_{\e, ii}(t) * Q_{\e, ij}(t), & \text{for} \  t \geq 0, j \neq i, i \in  \overline{\DD}, 
\end{array}
\right.   
\end{equation}
where, as usual, $Q^{(*0)}_{\e, ii}(t) = {\rm I}(t \geq 0), t \geq 0$.

Respectively, the  transition probabilities  for the embedded Markov chain $\tilde{\eta}_{\e, n} $ are given by the following relation, 
\begin{equation}\label{gopet}
 \tilde{p}_{\e, ij}  = \PP \{ \tilde{\eta}_{\e,  1} = j / \tilde{\eta}_{\e, 0} = i \}  
=  \left\{
\begin{array}{cll}
p_{\e, ij} & \text{for} \ j \in \XX, i \in \DD,  \vspace{2mm} \\ 
0 & \text{for} \ j = i, i \in \overline{\DD}, \vspace{2mm} \\ 
 \frac{p_{\e, ij}}{1 - p_{\e, ii}} & \text{for} \ j \neq i, i \in \overline{\DD}. 
\end{array}
\right.
\end{equation}

Note that condition ${\bf B}$ implies that probabilities $p_{\e, ii} < 1, i \in \overline{\DD}$, for every $\e \in (0, 1]$.

We are going to describe assumptions, under which  basic conditions  ${\bf A}$ -- ${\bf E}$ hold for the semi-Markov processes $\tilde{\eta}_{\e}(t)$ constructed with the use of the described above procedure of removing virtual transitions for the semi-Markov processes $\eta_\e(t)$. 

\vspace{1mm}

{\bf 3.2 Condition ${\bf A}$}. Condition  ${\bf A}$ and ${\bf B}$, assumed to hold for the Markov processes  $\eta_{\e}(t)$, imply that this condition ${\bf A}$ holds for the Markov processes $\tilde{\eta}_{\e}(t)$. 

Let us  introduce, for $i \in \XX$ and $\e \in (0, 1]$, distributions functions of sojourn times for the semi-Markov process $\tilde{\eta}_{\e}(t)$,   
 \begin{align}\label{sojouko}
 \tilde{F}_{\e, i}(t) & = \PP_i \{\tilde{\kappa}_{\e,  i} \leq t \} \vspace{2mm} \nonumber \\ 
 & = \sum_{j \in \XX} \tilde{Q}_{\e,  ij}(t), \ t \geq 0. 
 \end{align}

Since $ \tilde{F}_{\e, i}(t)  =  F_{\e, i}(t), t \geq 0$, for $i \in \DD$, probabilities $\tilde{F}_{\e, i}(0) < 1, i \in \DD$, for $\e \in (0, 1]$.   

Conditions ${\bf A}$ and ${\bf B}$,  for the semi-Markov processes $\eta_\e(t)$, and relation (\ref{transaka}) implies that, for $i \in \overline{\DD}$ and $\e \in (0, 1]$,
\begin{align}\label{ertuk}
\tilde{F}_{\e, i}(0) & =  \sum_{n = 0}^\infty Q^{n}_{\e, ii}(0)  \sum_{j \neq i}Q_{\e, ij}(0)  = \frac{F_{\e, i}(0) - Q_{\e, ii}(0)}{1 - Q_{\e, ii}(0)}  < 1.
 \end{align}
 
 Relation (\ref{transaka}) obviously implies that inequality $\tilde{F}_{\e, i}(0) < 1$ also holds, 
 for $i \in \DD$.
 
 Thus, the regularity condition  ${\bf A}$ holds for the semi-Markov processes $\tilde{\eta}_{\e}(t)$.  
  
 Therefore,  process   $\tilde{\eta}_{\e}(t)$ is well defined on the interval $[0, \infty)$, for every $\e \in (0, 1]$. 
 
 In what follows, we  denote by ${\bf \tilde{A}}$  condition ${\bf A}$ for the semi-Markov processes $\tilde{\eta}_{\e}(t)$. 
 \vspace{1mm}

{\bf 3.3 Conditions ${\bf B}$ and ${\bf \hat{B}}$}. Condition  ${\bf B}$, assumed to hold for the Markov processes  $\eta_{\e}(t)$, implies that  condition ${\bf B}$ holds for the semi-Markov processes $\tilde{\eta}_{\e}(t)$.

It is obviously follows from relation (\ref{gopet}),  for condition  ${\bf B}$ {\bf (a)}. 

As far as condition ${\bf B}$ {\bf (b)} is concerned, one can always assume that a chain of states $j_0, \dots, j_{n_i}$, appearing in condition ${\bf B}$ {\bf (b)}, satisfies an additional assumption that $j_{l-1} \neq j_{l}, l = 1, \ldots, n_{i}$. For any such chain, $\prod_{l = 1}^{n_i} \tilde{p}_{\e,  j_{l-1} j_{l}} \geq \prod_{l = 1}^{n_i} p_{\e, j_{l-1} j_{l}} > 0$. Therefore, condition  ${\bf B}$ {\bf (b)} also holds for the semi-Markov processes $\tilde{\eta}_{\e, \bar{k}_0}(t)$.

 In what follows, we  denote by ${\bf \tilde{B}}$  condition ${\bf B}$ for the semi-Markov processes $\tilde{\eta}_{\e}(t)$. 
 
 Note also that, if condition  ${\bf \hat{B}}$  is assumed to hold for the semi-Markov processes 
$\eta_{\e}(t)$, it also holds for the semi-Markov processes $\tilde{\eta}_{\e}(t)$, since, 
according relation (\ref{transaka}), the transition probabilities $\tilde{Q}_{\e, kr}(t) =  
Q_{\e, kr}(t) =  \frac{1}{m} {\rm I}( r \in \DD){\rm I}(t \geq 1)$, for $t \geq 0,  r \in \XX, k \in \DD$ and $\e \in (0, 1]$.  \vspace{1mm}

{\bf 3.4 Condition ${\bf C}$}.  Let us introduce the set, 
\begin{equation}\label{setra}
\ZZ_{0} = \{i \in \XX: p_{0, ii} = 1 \}.
\end{equation}

According to condition ${\bf C}$,  we can interpret states $i \in \XX \setminus \ZZ_0$ and states $j \in  \ZZ_{0}$, respectively,  as asymptotically non-absorbing  and absorbing states. Respectively, the cases, where 
set  $\ZZ_{0} = \emptyset$ or $\ZZ_{0} \neq \emptyset$ relate,  respectively, to the models of regularly and singularly perturbed semi-Markov processes. 

The following condition, additional  to the basic perturbation conditions  ${\bf A}$, ${\bf B}$,  and ${\bf C}$, plays an important role in the procedure of removing of virtual  
transitions for the semi-Markov processes $\eta_\e(t)$:  
\begin{itemize}
\item [${\bf \tilde{C}}$:]  $\tilde{p}_{\e, ij}  = {\rm I}(j \neq i)\frac{p_{\e, ij}}{1 - p_{\e, ii}}  \to \tilde{p}_{0, ij} \in [0, 1]$ as $\e \to 0$, for $j \in \XX, i \in \overline{\DD}$.
\end{itemize}

Recall that condition ${\bf B}$ implies that probabilities $p_{\e, ii} < 1, i \in \overline{\DD}$, for every $\e \in (0, 1]$.

Note also that probabilities $p_{\e, ij} = 0, \e \in (0, 1]$, and, thus, probabilities $\tilde{p}_{0,  ij} = 0$, for $j \in \overline{\YY}_{1, i}, i \in \overline{\DD}$.

Condition ${\bf \tilde{C}}$ plays the role of condition ${\bf C}$ for semi-Markov 
processes $\tilde{\eta}_{\e}(t)$.

If condition  ${\bf \hat{B}}$  is assumed to hold, the asymptotic relation given in condition ${\bf \tilde{C}}$ also holds, for $j \in \XX, i \in \DD$. In this case, probabilities $\tilde{p}_{\e, ij}   = \frac{1}{m}
{\rm I}(j \in \DD), j \in \XX, i \in \DD$, for all $\e \in (0, 1]$, and, thus, the above asymptotic relation holds and the corresponding limiting probabilities $\tilde{p}_{0, ij}   = \frac{1}{m}{\rm I}(j \in \DD)$, for $j \in \XX, i \in \DD$.

Since, in this case, matrix ${\tilde{\mathbf P}}_{\e} = \| \tilde{p}_{\e, ij} \|$ is stochastic, 
for $\e \in (0, 1]$, conditions  ${\bf \hat{B}}$  and ${\bf \tilde{C}}$  imply that  matrix ${\tilde {\mathbf P}}_{0} = \| \tilde{p}_{0,  ij} \|$ is also stochastic. 

Let $\tilde{\eta}_{0, n}, n = 0, 1, \ldots$ be a Markov chain with the phase space $\XX$ and the matrix of transition probabilities 
$\tilde{{\mathbf P}}_{0}$. Condition ${\bf \tilde{C}}$  makes it possible to interpret the Markov chains $\tilde{\eta}_{\e, n}$, for $\e \in (0, 1]$, as perturbed version of the Markov chain $\tilde{\eta}_{0, n}$. 

Note that condition ${\bf C}$  implies that, for $i \in \XX \setminus \ZZ_0, j \in \overline{\DD}$,
\begin{equation}\label{noter}
\tilde{p}_{\e,  ij}  = {\rm I}(j \neq i)\frac{p_{\e, ij}}{1 - p_{\e, ii}}  \to \tilde{p}_{0, ij} =  {\rm I}(j \neq i)\frac{p_{0, ij}}{1 - p_{0, ii}} \ {\rm as} \ \e \to 0.
\end{equation}

If the set of asymptotically absorbing states $\ZZ_{0} = \emptyset$, then relation (\ref{noter}) implies that condition 
${\bf \tilde{C}}$ holds.

If set $\ZZ_{0} \neq \emptyset$, condition ${\bf \tilde{C}}$  is partly a new one. Indeed, $p_{\e, ij}, 1 - p_{\e, ii} \to 0$ as $\e \to 0$, for  $i \in \ZZ_{0}, j \in \XX$ and, thus, condition ${\bf C}$  does not imply holding of the convergence relation given in condition ${\bf \tilde{C}}$,  for $i \in \ZZ_{0}, j \in \XX$. 

Since, probabilities $\tilde{p}_{\e, ii} = 0$, for  $i \in \overline{\DD}$ and $\e \in [0, 1]$, condition ${\bf \tilde{C}}$ is, in fact, equivalent to the following condition:
\begin{itemize}
\item [${\bf \tilde{C}}'$:]   $\frac{\tilde{p}_{\e, ij}}{1 - \tilde{p}_{\e, ii}} = \tilde{p}_{\e, ij}  \to \frac{\tilde{p}_{0, ij}}{1 - \tilde{p}_{0, ii}}  = \tilde{p}_{0, ij}$  as $\e \to 0$, for $j \in \XX, i \in \overline{\DD}$. 
\end{itemize}

Condition  ${\bf \tilde{C}}'$ plays the role of condition ${\bf \tilde{C}}$ for the semi-Markov processes $\tilde{\eta}_\e(t)$.
 
Let us introduce the following condition:
\begin{itemize}
\item [${\bf C}_{0}$:]  $q_{\e}[i i' i''] = \frac{p_{\e, i i'}}{p_{\e, i i''}} \to  q_{0}[i i' i''] \in [0, \infty]$ as $\e \to 0$, 
for $i' \in \XX, i'' \in \YY_{1, i}$, $i \in \overline{\DD}$.
\end{itemize}

Note that probabilities $p_{\e, i i'} = 0, \e \in (0, 1]$, for $i' \in \overline{\YY}_{1, i},  i \in \overline{\DD}$, and, thus, the asymptotic relation given in condition  ${\bf C}_{0}$ automatically holds, with limits  $q_{0, i i' i''} = 0$, for $i' \in  \overline{\YY}_{1, i}, i'' \in \YY_{1, i}, i \in \overline{\DD}$. \vspace{1mm}

{\bf Lemma 1}. {\em Condition ${\bf C}_{0}$ is sufficient for holding of condition ${\bf C}$.} \vspace{1mm}

{\bf Proof}. Condition ${\bf C}_{0}$ implies that, for $j \in \YY_{1, i}, i \in \overline{\DD}$,
\begin{align}\label{quponrew}
p_{\e, ij}  & = \frac{p_{\e, ij}}{\sum_{r \in \XX} p_{\e, ir}}  = \big(  \sum_{r \in \XX} \frac{p_{\e, ir}}{p_{\e, ij}} \big)^{-1}   
\vspace{2mm} \nonumber \\
& \to    \big(  \sum_{r \in \XX} q_{0}[irj]  \big)^{-1}  = p_{0, ij} \ {\rm as} \ \e \to 0.
\end{align}

Note that limits in relation (\ref{quponrew})  satisfy relations,  $p_{0, ij} \geq 0, j \in \YY_{1, i}$, \, 
$\sum_{j \in \YY_{1, i}} p_{0, ij} = 1$, for $i \in \overline{\DD}$. $\Box$
 \vspace{1mm} 

{\bf Lemma 2}. {\em Condition ${\bf C}_{0}$ is sufficient for holding of conditions  ${\bf \tilde{C}}$ and ${\bf \tilde{C}}'$.} \vspace{1mm}

{\bf Proof}. Condition ${\bf C}_{0}$ implies that, for $j \neq i, j \in \YY_{1, i}, i \in \overline{\DD}$,
\begin{align}\label{quponta}
\tilde{p}_{\e, ij}  & = \frac{p_{\e, ij}}{1 - p_{\e, ii}}  = \frac{p_{\e, ij}}{\sum_{r \neq i} p_{\e, ir}}  
 = \big(  \sum_{r \neq i} \frac{p_{\e, ir}}{p_{\e, ij}} \big)^{-1}   \vspace{2mm} \nonumber \\
% \end{align*}
 %\begin{align}
& \to    \big(  \sum_{r \neq i} q_{0}[irj]  \big)^{-1}  = \tilde{p}_{0, ij} \ {\rm as} \ \e \to 0.
\end{align}

Note that limits in relation  (\ref{quponta}) satisfy relations  $\tilde{p}_{0,  ij}\geq 0, j \neq i, j \in \YY_{1, i}$, \, $\sum_{j \neq i, j \in \YY_{1, i}} \tilde{p}_{0,  ij}  = 1$, for $i \in \XX$. $\Box$
\vspace{1mm}

If condition  ${\bf \hat{B}}$  is assumed to hold, the asymptotic relation given in condition ${\bf C}_0$ also holds, for 
$i' \in \XX, i'' \in \YY_{1, i}$, $i \in \DD$. In this case, probabilities $\tilde{p}_{\e, ij}   = \frac{1}{m}
{\rm I}(j \in \DD), j \in \XX, i \in \DD$, for $\e \in [0, 1]$, and sets $\YY_{1, i} =\DD, i \in \DD$. Thus, the above asymptotic relation holds and the corresponding limiting quantities $q_0[ii'i'']   = {\rm I}(i' \in \DD)$, for $i'', i \in \DD$.
\vspace{1mm}
 
 Let us   define sets, for $i \in \overline{\DD}$ and  $\e \in [0, 1]$, 
\begin{equation}\label{setana}
 \tilde{\YY}_{\e,  i}  = \{ j \in \XX: \tilde{p}_{\e,  ij} > 0 \}. 
\end{equation}

Condition ${\bf B}$ implies that condition ${\bf \tilde{B}}$ holds, and, thus, for $i \in \overline{\DD}$ and $\e \in (0, 1]$,
 \begin{equation}\label{setaser}
 \tilde{\YY}_{\e, i}  =  \tilde{\YY}_{1, i}. 
\end{equation}

Also, conditions ${\bf B}$ and  ${\bf C}$ imply that conditions  ${\bf \tilde{B}}$ and  ${\bf \tilde{C}}$ hold, and, thus, 
for $i \in \overline{\DD}$,
 \begin{equation}\label{setaserba}
 \tilde{\YY}_{0, i} \subseteq \tilde{\YY}_{1, i}. 
\end{equation}

It is readily seen that sets $\YY_{1, i}$ and $\tilde{\YY}_{1, i}$ are connected by the following relation, for $i \in \overline{\DD}$,
 \begin{equation}\label{setasernu}
 \tilde{\YY}_{1, i}   = \YY_{1, i} \setminus \{ i \}. 
\end{equation}

Note also that, under condition ${\bf \hat{B}}$, sets $\tilde{\YY}_{\e, i} = \DD, \e \in [0, 1]$,  for $i \in \DD$. \vspace{1mm} 

{\bf 3.5 Condition ${\bf D}$}.  The corresponding distribution function $\tilde{F}_{\e,  ij}(t), t \geq 0$ is defined, for $\e \in (0, 1]$ by the  following relation, for $j \in  \tilde{\YY}_{1,   i}, i \in \overline{\DD}$,  
\begin{align}\label{transakai}
\tilde{F}_{\e,   ij}(t) & =  \PP \{ \tilde{\kappa}_{\e,  1} \leq t  / \tilde{\eta}_{\e, 0} = i, \tilde{\eta}_{\e, 1} = j \}  \vspace{2mm} \nonumber \\
& =  \frac{1}{\tilde{p}_{\e,  ij}} \sum_{n = 0}^\infty F^{(*n)}_{\e, ii}(t) * F_{\e, ij}(t) p_{\e, ii}^n p_{\e, ij},  \  t \geq 0,
\end{align}
and,  for $j \notin  \tilde{\YY}_{1, i} , i \in \overline{\DD}$,
\begin{equation}\label{taswet}
\tilde{F}_{\e,  ij}(t) = \tilde{F}_{\e,  i}(t), \ t \geq 0, 
\end{equation}
where, for $i  \in \overline{\DD}$,
\begin{align}\label{taswetas}
\tilde{F}_{\e,  i}(t) & =  \PP \{ \tilde{\kappa}_{\e,  1} \leq t  / \tilde{\eta}_{\e,  0} = i \}  \vspace{2mm} \nonumber \\
&= \sum_{j \in \XX} \tilde{Q}_{\e, ij}(t)  =  \sum_{j \in \tilde{\YY}_{1, i}} \tilde{F}_{\e, ij}(t) \tilde{p}_{\e,  ij}, t \geq 0, 
\end{align}

The corresponding Laplace transform $\tilde{\phi}_{\e, ij}(s)$ takes, for 
$\e \in (0, 1]$,  the following form,  for $j \in  \tilde{\YY}_{1, i} , i \in \overline{\DD}$, 
\begin{align}\label{trawet}
\tilde{\phi}_{\e, ij}(s) & =  \EE \{ e^{- s \tilde{\kappa}_{\e,  1}}  /   \tilde{\eta}_{\e, 0} = i, \tilde{\eta}_{\e, 1} = j \} 
= \int_0^\infty e^{-st} \tilde{F}_{\e, ij}(dt)  \vspace{2mm} \nonumber \\
& = \frac{1}{\tilde{p}_{\e,  ij}} \cdot   \frac{\phi_{\e, ij}(s)  p_{\e, ij}}{1 - \phi_{\e, ii}(s)p_{\e, ii}} 
=  \frac{\phi_{\e, ij}(s)(1 - p_{\e, ii})}{1 - \phi_{\e, ii}(s)p_{\e, ii}} \vspace{2mm}\nonumber \\
& =    \frac{\phi_{\e, ij}(s)}{1 + p_{\e, ii} (1 - p_{\e, ii})^{-1}(1 - \phi_{\e, ii}(s))}, \  s \geq 0.  
\end{align}
and,   for $j \notin  \tilde{\YY}_{1,  i} , i \in  \overline{\DD}$, 
\begin{equation}\label{trawetata}
\tilde{\phi}_{\e,  ij}(s)   = \tilde{\phi}_{\e,  i}(s), s \geq 0,  
\end{equation}
where
\begin{align}\label{taswetop}
\tilde{\phi}_{\e,  i}(s)  = \EE \{ e^{- s \tilde{\kappa}_{\e,  1}}  /   \tilde{\eta}_{\e, 0} = i \}   = \sum_{j \in \tilde{\YY}_{1,   i}} \tilde{\phi}_{\e, ij}(s) \tilde{p}_{\e,  ij}, \ s \geq 0.
\end{align}

Let us now assume that, additionally to ${\bf A}$ -- ${\bf C}$, conditions ${\bf D}$ (or, equivalently,  
${\bf D}'$) and ${\bf E}$  hold.

In order to compensate aggregation of  inter-jump times used in relation (\ref{recur}),  let us introduce  new local normalisation functions 
$\tilde{v}_{\e, i} \in [1, \infty), \e \in (0, 1]$, defined by the following relation, for $i \in \overline{\DD}$,
\begin{equation}\label{compreg}
 \tilde{v}_{\e,  i} = (1 - p_{\e, ii})^{-1} v_{\e,  i}.
 \end{equation} 
 
 Note that, according condition ${\bf B}$ {\bf (b)}, probability $1 - p_{\e, ii} \in (0, 1], \e \in (0, 1]$, for $i \in \overline{\DD}$. 

 Conditions ${\bf C}$ and ${\bf D}$  imply that, for $i \in \overline{\DD}$, 
\begin{equation}\label{vana}
 \tilde{v}_{\e,  i} \to \tilde{v}_{0, i} \ {\rm as} \ \e \to 0,   
 \end{equation}
 where
 \begin{equation}\label{vanaba}
 \tilde{v}_{0,  i} =  \left\{
 \begin{array}{cll}
 (1 - p_{0, ii})^{-1}v_{0, i}    & \text{if} \  p_{0, ii} < 1 \ \text{and} \  v_{0, i} < \infty,  \vspace{1mm} \\
 \infty & \text{if} \  p_{0, ii} = 1 \ \text{or} \   v_{0, i} = \infty. 
 \end{array}
 \right.
\end{equation}

Using conditions ${\bf C}$, ${\bf D}'$  and relation (\ref{trawet}), we get the following relation, 
for $j \in \tilde{\YY}_{1,  i}, i \in  \overline{\DD} \setminus \ZZ_{0}$, 
\begin{equation*}
\tilde{\phi}_{\e,  ij}(s /   \tilde{v}_{\e, i} ) 
= \left\{
\begin{array}{lll}
   \frac{\phi_{\e, ij}( (1 - p_{\e, ii}) s / v_{\e,  i})(1 - p_{\e, ii})}{1 - 
p_{\e, ii}\phi_{\e, ii}( (1 - p_{\e, ii})  s /  v_{\e, i}))} & \  \text{if} \ i \in \YY_{1, i}  \makebox[1mm]{}
\vspace{2mm} \\
\phi_{\e, ij}(s / v_{\e,  i}) & \  \text{if} \ i \notin \YY_{1, i}
\end{array}
\right.
\end{equation*}
\begin{equation}\label{trawetba}
\to  \tilde{\phi}_{0,  ij}(s) \ {\rm as} \ \e \to 0, \ {\rm for} \ s \geq 0, 
\end{equation}
where
\begin{equation}\label{adii}
\tilde{\phi}_{0,  ij}(s) = \left\{
\begin{array}{lll}
\frac{\phi_{0, ij}((1 - p_{0, ii})s)(1 - p_{0, ii})}{1 - p_{0, ii} \phi_{0, ii}((1 - p_{0, ii})s ))} & \ \text{if} \ i \in \YY_{1, i},  
\vspace{2mm} \\
\phi_{0, ij}(s) & \ \text{if} \ i \notin \YY_{1, i}. 
\end{array}
\right.
\end{equation}

Obviously, $\tilde{\phi}_{0, ij}(s) = \int_0^\infty e^{-st} \tilde{F}_{0,  ij}(dt) = \EE e^{-s \tilde{\kappa}_{0, ij}}, s \geq 0$,  is the Laplace transform for the distribution function $\tilde{F}_{0, ij}(\cdot)$ of some non-negative random variable 
$\tilde{\kappa}_{0,  ij}$. If $i \in \YY_{1, i}$, i.e., $p_{0, ii} > 0$,  the random variable $\tilde{\kappa}_{0,  ij}$ is a geometric type sum of random variables, i.e., $\tilde{\kappa}_{0, ij} = \sum_{n = 1}^{\mu_{0, i} -1} (1 - p_{0, ii})  \kappa_{0, ii, n} +  (1 - p_{0, ii}) \kappa_{0, ij}$, where: (a) random variables $\mu_{0, i}, \kappa_{0, ij}, \kappa_{0, ii, n}, n = 1, 2, \ldots$ are independent; (b) $\mu_{0, i}$ is a geometrically distributed random variable with parameter $1 - p_{0, ii}$ (which takes value $n$ with probability $(1- p_{0, ii})p_{0, ii}^{n -1}$,  for $n = 1, 2, \ldots$);  (c) $\kappa_{0, ii, n}, n = 1, 2, \ldots$ are   i.i.d. random variables with the Laplace transform  $\phi_{0, ii}(s)$; (d) $\kappa_{0, ij}$ is a random variable with the Laplace transform $\phi_{0, ij}(s)$.  If $i \notin \YY_{1, i}$, i.e., $p_{0, ii} = 0$, the random variable $\tilde{\kappa}_{0,  ij} = 
\kappa_{0, ij}$. 

Obviously, $\tilde{F}_{0, ij}(0) < 1$, for $j \in \tilde{\YY}_{1, i}, i \in \overline{\DD} \setminus \ZZ_{0}$.

Let $\hat{\kappa}_{\e, ij, n}, n = 1, 2, \ldots$ be, for every $j \in \YY_{1,  i}, i \in \overline{\DD}$ and $\e \in (0, 1]$,  i.i.d.  random variables with distribution function $\hat{F}_{\e, ij}(\cdot) = F_{\e, ij}(\cdot \, v_{\e,  i})$. Let, also,   $\hat{\kappa}_{0, ij}$ be, for every $j \in \YY_{1,i}, i \in \overline{\DD}$, a random variable with the distribution function $\hat{F}_{0, ij}(\cdot) = 
F_{0, ij}(\cdot)$. Let also $\hat{e}_{\e, ij} = \int_0^\infty t \hat{F}_{\e, ij}(dt)$, for $j \in \YY_{1, i}, i \in \overline{\DD}$ and $\e \in [0, 1]$. 

By the central criterium of convergence for sums of independent random variables (see, for example, [56]) 
conditions  ${\bf D}$ and  ${\bf E}$ imply that the following relation of weak LLN (law of large numbers) type holds, for any $0 < u_\e \to \infty$ as $\e \to 0$ and 
$j \in \YY_{1, i}, i \in \overline{\DD}$, 
\begin{equation}\label{weak}
\frac{\sum_{n \leq u_\e} \hat{\kappa}_{\e, ij, n}}{u_{\e}}   \stackrel{{\rm d}}{\longrightarrow} \hat{e}_{0, ij} \ {\rm as} \ \e \to 0.
\end{equation}

Indeed, the above conditions imply that (a) $\hat{F}_{\e, ij}(\cdot) \Rightarrow \hat{F}_{0, ij}(\cdot)$ as $\e \to 0$, and (b) $\hat{e}_{\e, ij} 
  \to \hat{e}_{0, ij}  < \infty$ as $\e \to 0$, for every $j \in \YY_{1, i}, i \in \overline{\DD}$. 

Let $0 < s_k \to \infty$ as $k \to \infty$ be a sequence of continuity points for  the distribution function $\hat{F}_{0, ij}(\cdot)$. The above asymptotic relations (a) and (b) obviously imply  that, for any $t >  0$, 
\begin{align}\label{central}
\varlimsup_{\e \to 0} \int_{t u_\e}^\infty s \hat{F}_{\e, ij}(ds) & \leq \varlimsup_{\e \to 0}  \int_{s_k}^\infty s \hat{F}_{\e, ij}(ds)   
\makebox[50mm]{} \vspace{3mm} \nonumber \\
%\end{align*} 
%\begin{align}
& \makebox[5mm]{} = \varlimsup_{\e \to 0}(\hat{e}_{\e, ij} - 
\int_0^{s_k} s \hat{F}_{\e, ij}(ds)) 
 \vspace{3mm} \nonumber \\
&  \makebox[5mm]{}  = \hat{e}_{0, ij} - \int_0^{s_k}s \hat{F}_{0, ij}(ds) \to 0 \ {\rm  as} \ k \to \infty,
\end{align} 
and, thus,  the following relation holds, for any $t >0$, 
\begin{equation}\label{sdertak}
\varlimsup_{\e \to 0} \int_{t u_\e}^\infty s \hat{F}_{\e, ij}(ds) = 0. 
\end{equation}  

Relation (\ref{sdertak}) implies that, for any $t > 0$,  
\begin{align}\label{sdertaka}
u_\e \PP_1 \{u_\e^{-1} \hat{\kappa}_{\e, ij, 1} >  t \} & = u_\e  (1 - \hat{F}_{\e, ij}(tu_\e)) \vspace{1mm} \nonumber \\
& \leq t^{-1}\int_{t u_\e}^\infty s \hat{F}_{\e, ij}(ds) \to 0 \ {\rm as} \  \e \to 0. 
\end{align}  

Also, relation  (\ref{sdertak}) implies that,  for any $t > 0$,   
\begin{equation}\label{sdertakta}
u_\e \EE_1 u_\e^{-1} \hat{\kappa}_{\e, ij, 1} {\rm I}( u_\e^{-1} \hat{\kappa}_{\e, ij, 1} \leq t) = \int^{tu_\e}_0 s \hat{F}_{\e, ij}(ds) \to \hat{e}_{0, ij} \ {\rm  as} \ \e \to 0.
\end{equation}  

Relations (\ref{sdertaka}) and  (\ref{sdertakta}) imply, by the criterion of central convergence (see, for example,  [56]), that relation (\ref{weak})  holds.

Let us introduce Laplace transforms, for $j \in \YY_{1, i}, i \in \overline{\DD}$ and $\e \in [0, 1]$,
\begin{align}\label{treas}
\hat{\phi}_{\e, ij}(s) & = \int_0^\infty e^{- st} \hat{F}_{\e, ij}(dt) \vspace{2mm} \nonumber \\ 
& = \phi_{\e, ij}(s/v_{\e, i}), s \geq 0.
\end{align}

As well known, relation (\ref{weak}) is equivalent to the following relation,  
\begin{align}\label{centra} 
u_\e ( 1 - \hat{\phi}_{\e, ij}(s /  u_\e)) & = u_\e ( 1 - \phi_{\e, ij}(s /  u_\e v_{\e, i})) \vspace{1mm} \nonumber \\ 
& \to  e_{0, ij}s \ {\rm as} \ \e \to 0, \ {\rm for} \ s \geq 0. 
\end{align}

Since, $p_{\e, ii} \to 1$ as $\e \to 0$, for $i \in \overline{\DD} \cap \ZZ_{0}$, we can, in this case,  choose $u_\e = (1 - p_{\e, ii})^{-1}$. 

Note also that state $i \in \YY_{1, i}$, for $i \in  \overline{\DD} \cap \ZZ_{0}$. 

Conditions ${\bf C}$, ${\bf D}'$, ${\bf E}$  and relations  (\ref{trawet}) and (\ref{centra}) imply that the following relation holds, for $j \in \tilde{\YY}_{1, i}, i \in \overline{\DD} \cap \ZZ_{0}$,   
\begin{align}\label{trawetbama}
\tilde{\phi}_{\e, ij}(s /  \tilde{v}_{\e, i} ) & =   \frac{\phi_{\e, ij}(  s / (1 - p_{\e, ii})^{-1} v_{\e,  i})}{1 
+ p_{\e, ii} (1 - p_{\e, ii})^{-1} (1 - \phi_{\e, ii}(  s / (1 - p_{\e, ii})^{-1}  v_{\e, i}))} \vspace{2mm} \nonumber \\ 
& \to   \frac{1}{1 +  e_{0, ii} s} = \tilde{\phi}_{0,  ij}(s) \ {\rm as} \ 
\e \to 0, \ \, {\rm for} \ s \geq 0.   
\end{align}

In this case,  $\tilde{\phi}_{0,  ij}(s) = \tilde{\phi}_{0,  i}(s)  =  \int_0^\infty e^{-st} \tilde{F}_{0,  i}(dt)$  is the Laplace transform of a 
exponentially distributed random variable $\tilde{\kappa}_{0, ij} =  \kappa_{0, i}$, with parameter $e_{0, ii}^{-1}$.

Also, condition  ${\bf \tilde{C}}$ and relations (\ref{taswetop}), (\ref{trawetba}) and (\ref{trawetbama}) imply that the following relation takes place, for $i \in \overline{\DD}$,
\begin{align}\label{tatop}
\tilde{\phi}_{\e,  i}(s)  & = \sum_{j \in \tilde{\YY}_{1,   i}} \tilde{\phi}_{\e, ij}(s) \tilde{p}_{\e,  ij} \vspace{2mm} \nonumber \\ 
& \to \sum_{j \in \tilde{\YY}_{1,   i}} \tilde{\phi}_{0, ij}(s) \tilde{p}_{0,  ij}  
= \sum_{j \in \tilde{\YY}_{0,   i}} \tilde{\phi}_{0, ij}(s) \tilde{p}_{0,  ij}\vspace{2mm} \nonumber \\ 
& = \tilde{\phi}_{0,  i}(s) \ {\rm as} \ \e \to 0, {\rm for} \ s \geq 0.
\end{align}   

Relations (\ref{trawetba}) and (\ref{trawetbama}) imply that condition ${\bf D}'$ and, thus, also condition ${\bf D}$,  hold for semi-Markov processes $\tilde{\eta}_{\e}(t)$, with the Laplace transforms of the corresponding limiting distribution functions given in the above relations, and the new normalisation functions  $\tilde{v}_{\e, i}, i \in \overline{\DD}$. 

 In what follows, we can also denote by ${\bf \tilde{D}}$  and ${\bf \tilde{D}}'$, respectively,   conditions ${\bf D}$ and ${\bf D}'$ for the semi-Markov processes $\tilde{\eta}_{\e}(t)$ (expressed, for condition ${\bf D}'$,   in the form 
 of relations (\ref{trawetba}) and (\ref{trawetbama})). 
 
 Note also that, if condition  ${\bf \hat{B}}$  is assumed to hold, the asymptotic relations given in conditions ${\bf \tilde{D}}$ and ${\bf \tilde{D}}'$ also hold for $i \in \DD$. In this case, the normalisation functions $\tilde{v}_{\e, i} = 
v_{\e, i} = 1, \e \in (0, 1]$, the limiting distribution functions $\tilde{F}_{0, ij}(t) = \tilde{F}_{0, i}(t) = {\rm I}(t \geq 1), t \geq 0$, and the limiting Laplace transforms $\tilde{\phi}_{0, ij}(s) = \tilde{\phi}_{0, i}(s)  = e^{-s}, s \geq 0$, for $i, j \in \DD$.
 \vspace{1mm}

{\bf 3.6 Condition ${\bf E}$}.  Let us assume that conditions ${\bf A}$ --  ${\bf E}$   hold.

Expectation $\tilde{e}_{\e,  ij}$ takes 
the following form, for $j \in \tilde{\YY}_{1, i}, i \in \overline{\DD}$,
\begin{align}\label{twet}
\tilde{e}_{\e,  ij} & =  \EE \{  \tilde{\kappa}_{\e,  1}  /   \tilde{\eta}_{\e,  0} = i, \tilde{\eta}_{\e,  1} = j \}  \vspace{1mm} \nonumber \\   
& = e_{\e, ij} + p_{\e, ii} (1 - p_{\e, ii})^{-1} e_{\e, ii}.
 \end{align}
 
Also, for $i \in  \overline{\DD}$,
\begin{equation}\label{diotrasavom}
\tilde{e}_{\e,  i} = \EE \{ \tilde{\kappa}_{\e,  1} /   \tilde{\eta}_{\e,  0} = i \}   = 
\sum_{j \in \tilde{\YY}_{1,  i}} \tilde{e}_{\e,   ij}\tilde{p}_{\e,    ij}.
\end{equation}
and,  for $j \notin \tilde{\YY}_{1, i}, i \in \XX$,
\begin{equation}\label{twetok}
\tilde{e}_{\e,  ij}  =  \tilde{e}_{\e,  i}.
 \end{equation}

Conditions ${\bf B}$ -- ${\bf E}$ and relation  (\ref{twet}) imply that, for $j \in \tilde{\YY}_{1, i}, i \in \overline{\DD}$,
\begin{equation*}
\frac{\tilde{e}_{\e, ij}}{\tilde{v}_{\e,  i}}  = \left\{
\begin{array}{lll}
(1 - p_{\e, ii}) \frac{e_{\e, ij}}{v_{\e,  i}} 
+ p_{\e, ii}  \frac{e_{\e, ii}}{v_{\e,  i}} & \ \text{if} \  i \in \YY_{1, i} \makebox[10mm]{} \vspace{2mm} \\
\frac{e_{\e, ij}}{v_{\e,  i}} & \ \text{if} \  i \notin \YY_{1, i} 
\end{array}
\right.
\end{equation*}
\begin{equation}\label{expofava}
\to  \tilde{e}_{0,   ij} = \int_0^\infty t \tilde{F}_{0, ij}(dt) \ {\rm as} \ \e \to 0, \makebox[20mm]{}
\end{equation}
where
\begin{equation}
\tilde{e}_{0,   ij}  = \left\{
\begin{array}{lll}
(1 -  p_{0, ii}) e_{0, ij} + p_{0, ii}  e_{0, ii} & \ \text{if} \  i \in \YY_{1, i}, \\
e_{0, ij}& \ \text{if} \  i \notin \YY_{1, i}. 
\end{array}
\right.
\end{equation}

Note that the above formula works for both cases, where $j \in \tilde{\YY}_{1,    i}, i \in \overline{\DD} \setminus \ZZ_{0}$  and $j \in \tilde{\YY}_{1, i}, i \in  \overline{\DD} \cap \ZZ_{0}$. 

Also, the following relation takes place, for $i \in \overline{\DD}$, 
\begin{align}\label{expofavanu}
\frac{\tilde{e}_{\e, i}}{\tilde{v}_{\e, i}}  &  = \sum_{j \in \tilde{\YY}_{1,  i}}  \frac{\tilde{e}_{\e, ij}}{v_{\e,  i}} \tilde{p}_{\e, ij}
 \vspace{2mm} \nonumber \\
&  \to \sum_{j \in \tilde{\YY}_{1,  i}} \tilde{e}_{0, ij} \tilde{p}_{0, ij} =
\sum_{j \in \tilde{\YY}_{0,  i}} \tilde{e}_{0, ij} \tilde{p}_{0, ij}
 \vspace{2mm} \nonumber \\
& = \tilde{e}_{0, i} = \int_0^\infty t \tilde{F}_{0, i}(dt) \ {\rm as} \ \e \to 0.
\end{align}

Relation (\ref{expofavanu})  implies that condition ${\bf E}$  holds for the semi-Markov  processes $\tilde{\eta}_{\e}(t)$, with the first moments for the corresponding limiting distribution functions given in the above relation, and the local normalisation functions $v_{\e, i}, i \in \overline{\DD}$. 

 In what follows, we can also denote by ${\bf \tilde{E}}$ conditions ${\bf E}$  for the semi-Markov processes $\tilde{\eta}_{\e}(t)$ (expressed in the form 
 of relation (\ref{expofava})). 
 
 Note also that, if  ${\bf \hat{B}}$  is assumed to hold, the asymptotic relations given in condition 
 ${\bf \tilde{E}}$ also hold for $i \in \DD$. In this case, the  limiting expectations $\tilde{e}_{0, ij} = 
 \tilde{e}_{0, i}  = 1$, for $i, j \in \DD$.
 \vspace{1mm}
 
{\bf 3.7 Summary}. The following lemma summarise the above remarks made in Subsections 3.2 -- 3.6.  \vspace{1mm}

{\bf Lemma 3}. {\em Let conditions ${\bf A}$ -- ${\bf E}$  and ${\bf \tilde{C}}$ hold for the semi-Markov processes 
$\eta_\e(t)$. Then, conditions ${\bf A}$ -- ${\bf E}$ and ${\bf \tilde{C}}$  also hold for  the  semi-Markov processes 
$\tilde{\eta}_{\e}(t)$, respectively, in the form of conditions  ${\bf \tilde{A}}$ -- ${\bf \tilde{E}}$ and ${\bf \tilde{C}}'$.}
\vspace{1mm}

{\bf Remark 1}.  Condition ${\bf \tilde{C}}$ plays the role of condition ${\bf C}$ and condition  ${\bf \tilde{C}}'$, which is equivalent to condition  ${\bf \tilde{C}}$,  plays the role of condition ${\bf \tilde{C}}$,   for the semi-Markov processes $\tilde{\eta}_{\e}(t)$. 
\vspace{1mm}

{\bf Remark 2}. Condition ${\bf C}_{0}$ is sufficient for  holding of conditions ${\bf C}$, ${\bf \tilde{C}}$ and ${\bf \tilde{C}}'$. \vspace{1mm}

{\bf 3.8 Hitting times for semi-Markov processes with removed virtual transitions}. Let us introduce hitting times for semi-Markov processes $\tilde{\eta}_{\e}(t)$, 
\begin{equation}
\tilde{\tau}_{\e, \DD} = \sum_{n = 1}^{\tilde{\nu}_{\e,  \DD}} \kappa_{\e,  n}, \ {\rm where} \ 
\tilde{\nu}_{\e,  \DD} = \min(n \geq 1: \tilde{\eta}_{\e, n} \in \DD).
\end{equation}

The following lemma present an important invariance property  of hitting times $\tau_{\e, \DD}$.
\vspace{1mm}

{\bf Lemma 4}, {\em Let condition ${\bf A}$ and ${\bf B}$   hold, and, in sequel, conditions ${\bf \tilde{A}}$ and ${\bf \tilde{B}}$ hold. Then, the following relation takes place, for $\e \in (0, 1]$,  
\begin{equation}\label{iderek}
\PP_i \{ \tau_{\e, \DD} = \tilde{\tau}_{\e, \DD}, \,  \eta_{\e}(\tau_{\e, \DD}) = 
\tilde{\eta}_{\e}(\tilde{\tau}_{\e, \DD}) \} = 1, \, i \in \XX.   
\end{equation}}
\makebox[3mm]{} {\bf Proof}. The following obvious relation connects 
semi-Markov processes $\eta_\e(t), t \geq 0$ and $\tilde{\eta}_{\e}(t), t \geq 0$,
for $i \in \XX$ and  $\e \in (0, 1]$,
\begin{equation}\label{ident}
\PP_i \{ \eta_\e(t) = \tilde{\eta}_{\e}(t), \, t \geq 0 \} = 1.
\end{equation}

Let us define the following variant of hitting time, for $s  \geq 0$,
\begin{equation}\label{intra}
\tau_\e[s] = \inf(t \geq s: \eta_\e(t) \in \DD),
\end{equation}
and
\begin{equation}\label{intramop}
\tilde{\tau}_\e[s] = \inf(t \geq s: \tilde{\eta}_\e(t) \in \DD).
\end{equation}

Relation (\ref{ident}) obviously implies that, for $s \geq 0, i \in \overline{\DD}$ and $\e \in (0, 1]$, 
\begin{equation}\label{ideremik}
\PP_i \{\tau_\e[s] = \tilde{\tau}_\e[s], \, 
\eta_{\e}(\tau_{\e}[s]) = \tilde{\eta}_\e(\tilde{\tau}_\e[s]) \} = 1.
\end{equation}

Also, relations  (\ref{ident})  and (\ref{ideremik}) imply that, for $i \in \overline{\DD}$ and $\e \in (0, 1]$, 
\begin{align}\label{iderem}
& \PP_i \{\tau_{\e, \DD} = \tau_\e[0] = \tilde{\tau}_\e[0] = \tilde{\tau}_{\e, \DD} , 
\vspace{1mm} \nonumber \\
& \quad \ \ \eta_{\e}(\tau_{\e, \DD})  = \eta_{\e}(\tau_\e[0]) =  \tilde{\eta}_{\e}(\tilde{\tau}_\e[0]) = \tilde{\eta}_{\e}(\tilde{\tau}_{\e, \DD}) \} = 1.
\end{align}

Relations (\ref{stop}) -- (\ref{recur}) imply that: (a) if $\eta_{\e, 0} = i \in \DD$ and 
$\eta_{\e, 1}  \in \DD$, then the random functional $\theta_\e[0] = \nu_{\e, \DD} = \tilde{\nu}_{\e, \DD} = 1$, and, thus, $\tau_{\e, \DD} 
= \tilde{\tau}_{\e, \DD} = \kappa_{\e, 1}$, and $\eta_\e(\tau_{\e, \DD})  = \tilde{\eta}_\e(\tilde{\tau}_{\e, \DD}) = \eta_{\e, 1}$; (b)  if $\eta_{\e, 0} = i \in \DD$ and $\eta_{\e, 1}  \in \overline{\DD}$, then 
$\tau_{\e, \DD}  = \kappa_{\e, 1} + \tau_{\e}[\kappa_{\e, 1}]$ and  $\tilde{\tau}_{\e, \DD}  = \kappa_{\e, 1} + 
\tilde{\tau}_{\e}[\kappa_{\e, 1}]$.

Using the above remarks and the Markov property of the Markov renewal process $(\eta_{\e, n}, \kappa_{\e, n})$ and relations (\ref{ideremik}), (\ref{iderem})  we get, for $i \in \DD$,  and $\e \in (0, 1]$,
\begin{align}\label{ideremop}
& \PP_i \{ \tau_{\e, \DD} = \tilde{\tau}_{\e, \DD}, \eta_{\e}(\tau_{\e, \DD}) = 
\tilde{\eta}_{\e}(\tilde{\tau}_{\e, \DD}) \}     
\vspace{2mm}  \nonumber \\
&  \quad = \PP_i \{ \tau_{\e, \DD} = \tilde{\tau}_{\e, \DD}, \eta_{\e}(\tau_{\e, \DD}) = 
\tilde{\eta}_{\e}(\tilde{\tau}_{\e, \DD}), \eta_{\e, 1} \in \DD \} 
\vspace{2mm}  \nonumber \\
&  \quad \quad +  \PP_i \{ \tau_{\e, \DD} = \tilde{\tau}_{\e, \DD}, \eta_{\e}(\tau_{\e, \DD}) = 
\tilde{\eta}_{\e}(\tilde{\tau}_{\e, \DD}), \eta_{\e, 1} \in \overline{\DD} \} \vspace{2mm}  \nonumber \\
&  \quad = \PP_i\{ \eta_{\e, 1} \in \DD \} \vspace{2mm}  \nonumber \\
&  \quad \quad + \sum_{k \in \overline{\DD}} \int_0^\infty
\PP_k \{ s + \tau_{\e}[s] =  s + \tilde{\tau}_{\e}[s], \eta_{\e}(\tau_{\e}[s]) = 
\tilde{\eta}_{\e}(\tilde{\tau}_{\e}[s]) \} Q_{\e, ik}(ds) \makebox[3mm]{}
\vspace{2mm}  \nonumber \\
&  \quad = \sum_{r \in \DD} p_{\e, ir} 
+ \sum_{k \in \overline{\DD}} \int_0^\infty Q_{\e, ik}(ds) = 1.  
\end{align}

The prof is completed. $\Box$

Let us also introduce distributions,  for $i \in \XX$ and $\e \in (0, 1]$,  
\begin{equation}\label{trewq}
\tilde{G}_{\e, \DD, ij}(t) = \PP_i \{ \tilde{\tau}_{\e,  \DD} \leq t, \tilde{\eta}_{\e}(\tilde{\tau}_{\e, \DD}) = j) \}, \, t \geq 0, j \in \DD.
\end{equation}

The following lemma is an obvious corollary of Lemma 4. \vspace{1mm}

{\bf Lemma 5}. {\em Let condition ${\bf A}$ and ${\bf B}$   hold, and, in sequel, conditions ${\bf \tilde{A}}$ and ${\bf \tilde{B}}$ hold. Then, the following relation takes place, for  $\e \in (0, 1]$,  
\begin{equation}\label{iderenasba}
G_{\e, \DD, ij}(t)   = \tilde{G}_{\e,  \DD, ij}(t), \ t \geq 0,  j \in \DD, \, i \in \XX.  
\end{equation}}
\makebox[3mm]{} Let us also introduce Laplace transforms, for $i \in \XX$ and $\e \in (0, 1]$,
\begin{equation}\label{trewqol}
\Psi_{\e, \DD, ij}(s) = \EE_i \exp\{- s \tau_{\e, \DD}\} {\rm I}(\eta_{\e}(\tau_{\e, \DD}) = j), \, s \geq 0, j \in \DD.
\end{equation}
and  
\begin{equation}\label{trewqoki}
\tilde{\Psi}_{\e,  \DD, ij}(s) = \EE_i \exp\{- s \tilde{\tau}_{\e, \DD}\} {\rm I}(\tilde{\eta}_{\e}(\tilde{\tau}_{\e, \DD}) = j) \}, 
\, s \geq 0, j \in \DD.
\end{equation}

The following lemma re-formulates propositions of Lemma 5 in the equivalent form,  in terms of  Laplace transforms for hitting times. 

\vspace{1mm}

{\bf Lemma 6}. {\em  Let condition ${\bf A}$ and ${\bf B}$   hold, and, in sequel, conditions ${\bf \tilde{A}}$ and ${\bf \tilde{B}}$ hold. Then, the following relation takes place, for $\e \in (0, 1]$,
\begin{equation}\label{iderenasnop}
\Psi_{\e, \DD, ij}(s)   = \tilde{\Psi}_{\e, ij}(s), \ s \geq 0,  j \in \DD, \, i \in \XX.    
\end{equation}}
\makebox[4mm]{} Lemmas 4 -- 6 let us reduce study of asymptotics for distributions of hitting times  for semi-Markov processes $\eta_\e(t)$ to the case of more simple  semi-Markov processes $\tilde{\eta}_{\e}(t)$.  

The distributions $\tilde{G}_{\e,  \DD, ij}(t), t \geq 0, j \in \DD, i \in \overline{\DD}$ are determined by transition probabilities $\tilde{Q}_{\e, kr}(t), t \geq 0,   r \in \XX, k \in \overline{\DD}$, which themselves
are determined by transition probabilities $Q_{\e, kr}(t), t \geq 0,   r \in \XX, k \in \overline{\DD}$.

This makes it possible to essentially simplify the model in the case, where we are interested to investigate asymptotics of hitting times for some fixed domain $\DD$ and only for initial states $i \in \overline{\DD}$. This case in considered in Sections 2 -- 7. For example, we can assume in what follows in these sections that the transition probabilities  of semi-Markov processes $\eta_\e(t)$ and, in sequel, the transition probabilities of semi-Markov processes  
$\tilde{\eta}_\e(t)$, satisfy condition ${\bf \hat{B}}$. 

\vspace{1mm}

{\bf 3.9 Weak asymptotics of hitting times for the case of one-state domain $\overline{\DD}$}. The above remarks and Lemmas 4 -- 6 let us describe asymptotics for distributions of hitting times for the simplest case, where domain 
$\overline{\DD} = \{ i \}$  is a one-state set.  

In this case,  the hitting time $\tau_{\e, \DD} = \tilde{\tau}_{\e, \DD}  =  \tilde{\kappa}_{\e, 1}$ and $\eta_{\e}(\tau_{\e, \DD}) = \tilde{\eta}_{\e}(\tilde{\tau}_{\e, \DD}) = \tilde{\eta}_{\e, 1}$, if $\eta_\e(0) = i$.

Thus, the following relation takes place,
\begin{equation}\label{cotyre}
G_{\e, \DD, ij}(t) = \tilde{F}_{\e, ij}(t) \tilde{p}_{\e, ij}, t \geq 0, j \in \DD.
\end{equation}

According to this relation, one can expect that the corresponding limiting distribution should take the form $G_{0, \DD, ij}(\cdot) = 
\tilde{F}_{0, ij}(\cdot) \tilde{p}_{0, ij}$,  where  the limiting probabilities $\tilde{p}_{0, ij}$  and Laplace transforms of distribution functions $\tilde{F}_{0,  ij}(\cdot)$  are given: (a) by relations (\ref{noter}), (\ref{trawetba}) and (\ref{adii}), if $p_{0, ii} < 1$, or (b) by relation given in condition ${\bf \tilde{C}}$ and relation (\ref{trawetbama}), if $p_{0, ii} =1$, while  (c) the normalising function $\tilde{v}_{\e, i}$ is given by relation {\rm (\ref{compreg})}. 

The hitting time is in the case of one-state domain $\overline{\DD}$ is a geometric 
random sum. The  following simple theorem is, in fact,  a slight modification of the well known  weak convergence limit
theorem for geometric random sums. We refer to  books [36, 72], where the related comments and references can be found.
 \vspace{1mm}

{\bf Theorem 1}. {\em Let domain $\overline{\DD} = \{ i \}$  be a one-state set and conditions ${\bf A}$ -- ${\bf E}$ and ${\bf \tilde{C}}$ hold for 
semi-Markov processes $\eta_\e(t)$.  Then,  the following asymptotic relation takes place, for $j \in \DD$, 
\begin{equation}\label{cotyreva}
G_{\e, \DD, ij}(\cdot \, \tilde{v}_{\e, i}) \Rightarrow  G_{0, \DD, ij}(\cdot) = \tilde{F}_{0, ij}(\cdot) \tilde{p}_{0, ij} \ {\rm as} \ \e \to 0. 
\end{equation}} 
\makebox[4mm]{} {\bf Remark 3}. The distribution functions $\tilde{F}_{0, ij}(\cdot), j \in \DD$ are not concentrated in $0$, i.e., $\tilde{F}_{0, ij}(0) < 1$, for $j \in \DD$. 

{\bf Remark 4}. Conditions ${\bf C}$ and ${\bf \tilde{C}}$  can be replaced in Lemma 3 and Theorem 1  by condition ${\bf C}_{0}$. \vspace{1mm}

{\bf Remark 5}. Conditions ${\bf E}$  and  
${\bf \tilde{C}}$ can be omitted in Theorem 1, for the case (a) $p_{0, ii} < 1$. \\

{\bf 4. One-State Reduction of Phase Space  for Perturbed \\ \makebox[11mm]{} Semi-Markov Processes} \\ 

In this section, we describe an asymptotic  one-state phase space reduction procedure for perturbed semi-Markov processes.  We show that hitting times are asymptotically invariant with respect to this procedure. We also formulate conditions, which guarantees that basic perturbation conditions imposed on the initial semi-Markov processes also holds for the semi-Markov processes with reduced phase space.  Also, we give explicit formulas for re-calculating  normalisation functions, limiting distributions and expectations  
in the corresponding perturbation conditions for  the semi-Markov processes with reduced phase space.  

\vspace{1mm}

{\bf 4.1 An asymptotic one-step procedure of  phase space reduction for perturbed semi-Markov processes}. We assume now that the number of states $\bar{m}$ in the domain $\overline{\DD}$ is larger than  $1$. 

We also assume that $\e \in (0, 1]$ and conditions ${\bf A}$ and ${\bf B}$ hold for the semi-Markov processes $\eta_\e(t)$ and, in sequel, for   the semi-Markov processes $\tilde{\eta}_{\e}(t)$. 

Let us chose some state $k \in \overline{\DD}$ and introduce the reduced phase space $_k\XX = \XX \setminus \{ k \}$ and the reduced domain $_k\overline{\DD} = \overline{\DD} \setminus \{ k \}$. 

Let us define the stopping times for the Markov chain $\tilde{\eta}_{\e, n}$ that are, for  $r = 0, 1, \ldots$,   
\begin{equation}\label{fotyr}
_k\alpha_{\e}[r] = \min(n > r: \tilde{\eta}_{\e, n} \in \, _k\XX ).
\end{equation} 

By the definition, $_k\alpha_{\e}[r] $ is the first after $r$ time of hitting into the reduced phase space  
$_k\XX$ by the Markov chain $\tilde{\eta}_{\e, n}$. 

Since the Markov chain $\tilde{\eta}_{\e,  n}$ does not  make virtual transitions in domain $\overline{\DD}$, the following relation  takes place, for $r = 0, 1, \ldots$, 
\begin{equation}\label{relas}
_k \alpha_{\e}[r]  = \left\{
 \begin{array}{cll}
 r + 1 & \text{if} \  \tilde{\eta}_{\e,   r+1} \in \, _k\XX, \\
 r + 2 & \text{if} \  \tilde{\eta}_{\e,   r+1} = k.
 \end{array}
 \right.
 \end{equation} 

Let us also define sequential stopping times,
\begin{equation}\label{serwer}
_k\beta_{\e, n} = \, _k\alpha_{\e}[_k\beta_{\e,  n-1}], n = 1, 2, \ldots, \ {\rm where} \ \, _k\beta_{\e,  0} 
= {\rm I}(\tilde{\eta}_{\e,  0} = k).
\end{equation}

Let us now construct a new Markov renewal process $(_k\eta_{\e, n}, \, _k\kappa_{\e, n}), n = 0, 1, \ldots$, with the phase space $_k\XX \times [0, \infty)$ using the following recurrent relations,
\begin{equation}\label{recurask}
(_k\eta_{\e,  n}, \, _k \kappa_{\e,  n}) = \left\{ 
\begin{array}{cll}
(\tilde{\eta}_{\e, \, _k\beta_{\e,  0}}, 0) & \text{for} \ n = 0,  \vspace{2mm} \\
(\tilde{\eta}_{\e,   \, _k\beta_{\e,   n}}, \sum_{l =  \, _k\beta_{\e,   n-1} + 1}^{_k\beta_{\e,   n}} \tilde{\kappa}_{\e,    l}) & \text{for} \ n = 1, 2, \ldots.  
\end{array}
\right.
\end{equation}

We also can define the corresponding reduced semi-Markov process,
\begin{equation}\label{semink}
_k\eta_{\e}(t) = \eta_{\e, \, _k\nu_{\e}(t)}, t \geq 0,
\end{equation}
where $_k\zeta_{\e, n} = \, _k\kappa_{\e, 1} + \cdots + \, _k\kappa_{\e, n}, n = 1, 2,
 \ldots, \, _k\zeta_{\e, 0} = 0$, are the corresponding instants of jumps, and 
$_k\nu_{\e}(t) = \max(n \geq 1: \, _k\zeta_{\e, n}  \leq t)$ is the number of jumps in an interval $[0, t], t \geq 0$ for the above semi-Markov process.  

The definition of stopping times $_k\beta_{\e, n}$ implies that the transition probabilities for the above 
Markov renewal processes are determined by the following relation,  for $ t \geq 0, i, j \in \, _k\XX$,
 \begin{equation*}
_kQ_{\e, ij}(t)  =  \PP \{ _k\eta_{\e,  1} = j, \, _k\kappa_{\e, 1} \leq t  /  \, _k\eta_{\e, 0} = i \} \makebox[46mm]{}
\end{equation*}
\begin{equation}\label{transakak} 
\makebox[3mm]{} = \left\{
\begin{array}{lll} 
Q_{\e, ij}(t)  + Q_{\e, ik}(t) * \tilde{Q}_{\e, kj}(t) & \text{for} \ j \in \, _k\XX,  i \in \DD, \vspace{2mm} \\
\tilde{Q}_{\e,  ij}(t)  + \tilde{Q}_{\e, ik}(t) * \tilde{Q}_{\e, kj}(t) & \text{for} \ j \neq i,  
j \in \, _k\XX,  i \in \, _k\overline{\DD}, \vspace{2mm} \\
\tilde{Q}_{\e, ik}(t) * \tilde{Q}_{\e,  ki}(t) & \text{for} \ j = i,   i \in \, _k\overline{\DD}. 
 \end{array}
 \right.
\end{equation}

Respectively, the  transition probabilities  for the embedded Markov chain $_k\eta_{\e, n} $ are given by the following relation, for 
$i, j \in \, _k\XX$,
\begin{equation*}
\makebox[14mm]{}  _kp_{\e, ij}  = \PP \{ _k\eta_{\e, 1} = j / \, _k\eta_{\e, 0} = i \} \makebox[65mm]{}
\end{equation*}
\begin{equation}\label{gopetk}
= \left\{
\begin{array}{lll} 
p_{\e, ij} +  p_{\e, ik} \, \tilde{p}_{\e, kj} & \text{for} \ j \in \, _k\XX,   i \in \DD,   \vspace{2mm}  \\ 
\tilde{p}_{\e,  ij} +  \tilde{p}_{\e,  ik} \, \tilde{p}_{\e, kj} & \text{for} \ j \neq i,  
j \in \, _k\XX,  i \in \, _k\overline{\DD},  \vspace{2mm}  \\
 \tilde{p}_{\e, ik} \, \tilde{p}_{\e, ki} & \text{for} \ j = i, i \in \, _k\overline{\DD}.
\end{array}
\right.
\end{equation}

It is worth to note that the semi-Markov process $_k\eta_{\e}(t)$ can be considered as result of the two-stages procedure applied to the initial semi-Markov process  $\eta_\e(t)$. At the first stage, the semi-Markov process $\eta_\e(t)$ is transformed in the semi-Markov process $\tilde{\eta}_{\e}(t)$, with the use of the  procedure of removing virtual transitions. At the second stage, the semi-Markov process 
$\tilde{\eta}_{\e}(t)$ is transformed in the semi-Markov process $_k\eta_{\e}(t)$, with the  use of  the described above procedure of the phase space reduction.

We are going to describe assumptions, under which  basic conditions ${\bf A}$ --${\bf E}$ and condition ${\bf \tilde{C}}$ hold for the semi-Markov processes 
$_k\eta_{\e}(t)$, constructed with the use of the described above procedure for the semi-Markov processes 
$\eta_{\e}(t)$. Note that the phase space for the semi-Markov processes $_k\eta_{\e}(t)$ is $_k\XX$ and domain $_k\overline{\DD}$ plays the role of 
domain $\overline{\DD}$.  \vspace{1mm}

{\bf 4.2 Condition ${\bf A}$}.  Let conditions  ${\bf A}$ and ${\bf B}$ hold, and, in sequel, conditions ${\bf \tilde{A}}$ and ${\bf \tilde{B}}$ hold. 

Let us introduce distribution functions of sojourn times, for $t \geq 0, i \in \; _k\XX$ and $\e \in (0, 1]$, 
\begin{equation}\label{hokiu}
_kF_{\e, i}(t)  =  \PP \{ _k\kappa_{\e,  1} \leq t  / \, _k\eta_{\e, 0} = i \} 
= \sum_{j \in \, _k\XX} \, _kQ_{\e, ij}(t). 
\end{equation}

Relations  (\ref{transakak}) and (\ref{hokiu}) imply that, for $t \geq 0, i \in \DD$ and $\e \in (0, 1]$,
\begin{equation}\label{hokiuba}
_kF_{\e,  i}(t)  =  \sum_{j \neq k }Q_{\e, ij}(t) + Q_{\e, ik}(t) * \tilde{F}_{\e,  k}(t), 
\end{equation}
and, thus, for $i \in \DD$ and $\e \in (0, 1]$,
\begin{align}\label{hokiubana}
_kF_{\e,  i}(0)  & =  \sum_{j \neq k}Q_{\e, ij}(0) + Q_{\e, ik}(0)  \tilde{F}_{\e, k}(0)  \vspace{1mm} \nonumber \\
& \leq \sum_{j \neq k}Q_{\e, ij}(0) + Q_{\e, ik}(0) = F_{\e, i}(0) < 1.  
\end{align}

Also, the above relations imply that, for $ t \geq 0, i \in \, _k\overline{\DD}$ and $\e \in (0, 1]$, 
\begin{equation}\label{hokiubata}
_kF_{\e,  i}(t)  =  \sum_{j \neq i, k}\tilde{Q}_{\e,  ij}(t) 
+ \tilde{Q}_{\e, ik}(t) * \tilde{F}_{\e, k}(t), 
\end{equation}
and, thus, for $i \in \, _k\overline{\DD}$ and $\e \in (0, 1]$,
\begin{align}\label{hokiubar}
_kF_{\e,  i}(0)  & =  \sum_{j \neq i, k}\tilde{Q}_{\e,  ij}(0) 
+ \tilde{Q}_{\e, ik}(0)  \tilde{F}_{\e,  k}(0)  \vspace{1mm} \nonumber \\
& \leq \sum_{j \neq i, k}\tilde{Q}_{\e,  ij}(0) + \tilde{Q}_{\e, ik}(0) = \tilde{F}_{\e, i}(0) < 1.  
\end{align}

 Thus, the regularity condition  ${\bf A}$ holds for the semi-Markov processes $_k\eta_{\e}(t)$. 
 
 Therefore,  the semi-Markov process   $_k\eta_{\e}(t)$ is well defined on the interval $[0, \infty)$, 
 for every $\e \in (0, 1]$. 
 
 In what follows, we denote by $_k{\bf A}$  condition ${\bf A}$ for the semi-Markov processes $_k\eta_{\e}(t)$ (expressed in the form  of relations (\ref{hokiuba}) and (\ref{hokiubar})). 
 \vspace{1mm}

{\bf 4.3 Conditions ${\bf B}$ and ${\bf \hat{B}}$}.  Condition ${\bf B}$  assumed to hold for semi-Markov processes $\eta_\e(t)$ and, in sequel, for  semi-Markov processes   $\tilde{\eta}_{\e}(t)$, imply that  condition ${\bf B}$ holds for  semi-Markov processes  $_k\eta_{\e}(t)$.

Relation (\ref{gopetk}) obviously implies this for condition ${\bf B}$ {\bf (a)}. 

Condition ${\bf B}$ for the semi-Markov processes $\tilde{\eta}_{\e}(t)$, is, under condition  ${\bf B}$ {\bf (a)},  equivalent to the assumption that, for any  state $i \in \overline{\DD}$,  there exist an integer number $n = n_{i}$ and a chain of states $i = j_0, j_1 \in \overline{\DD}, \ldots, j_{n_i - 1}  \in \overline{\DD}, j_{n_i} \in \DD$ such that $j_{l-1} \neq j_{l}, l = 1, \ldots, n_i $ and $\prod_{l = 1}^{n_i} \tilde{p}_{1,  j_{l-1} i_{l}} > 0$.  
Let us assume that states $j_0, j_{n_i} \in \, _k\XX$. If states $j_{l-1}, j_{l} \in \, _k\XX$, for some $1 \leq l \leq n_i$, then probability 
$_kp_{1, j_{l-1} j_{l}} \geq \tilde{p}_{1,  j_{l-1} j_{l}}  > 0$. If state $j_l = k$, for some $1 \leq l \leq n$, then states $j_{l-1}, j_{l +1} \in \, _k\XX$ and probability $_kp_{1, j_{l-1} j_{l+1}} \geq 
\tilde{p}_{1,  j_{l-1} k} \tilde{p}_{\e,  k j_{l+1}}  > 0$.  Let us now $i = j_0 = j'_0, j'_1, \ldots, j'_{n'_i} = j_{n_i}$  be a new chain of states constructing by exclusion from the chain $i = j_0, j_1, \ldots, j_{n_i} $ all states  $j_l, 1 \leq l \leq n_i - 1$ such that $j_l = k$. The above inequalities  imply that for this new chain 
$\prod_{l = 1}^{n'_i} \, _kp_{1, j'_{l-1} j'_{l}} > 0$. Therefore, condition  ${\bf B}$  {\bf (b)}  also holds for the semi-Markov processes $_k\eta_{\e}(t)$. 

In what follows, we denote by $_k{\bf B}$  condition ${\bf B}$ for the semi-Markov processes $_k\eta_{\e}(t)$. \vspace{1mm}

Condition ${\bf \hat{B}}$,  assumed to hold for  semi-Markov processes $\eta_\e(t)$ and, in sequel,  for  semi-Markov processes  $\tilde{\eta}_{\e}(t)$, implies that condition ${\bf \hat{B}}$ also holds for semi-Markov processes  $_k\eta_{\e}(t)$.  It follows from relation (\ref{transakak})  that, in this case, the transition probabilities $Q_{\e, ij}(t) = \tilde{Q}_{\e, ij}(t) = \, _kQ_{\e, ij}(t) = \frac{1}{m} {\rm I}(r \in \DD) {\rm I}(t \geq 1),  t \geq 0, i \in \DD, j \in \, _k\XX$. 
\vspace{1mm}

{\bf 4.4 Conditions ${\bf C}$ and ${\bf \tilde{C}}$}.  Let us assume that conditions ${\bf A}$ -- ${\bf C}$,  and  ${\bf \tilde{C}}$  hold, and, in sequel, conditions 
${\bf \tilde{A}}$ -- ${\bf \tilde{B}}$ hold. 

Condition ${\bf \tilde{C}}$ implies the following condition holds:
\begin{itemize}
\item [$_k{\bf C}$:]  $_kp_{\e, ij} =  \tilde{p}_{\e, ij} +  \tilde{p}_{\e, ik} \tilde{p}_{\e, kj}   \to \, _kp_{0, ij}  = \tilde{p}_{0, ij} +  \tilde{p}_{0, ik} \tilde{p}_{0,  kj}$ 
as $\e \to 0$, for $j \in \, _k\XX, i \in \, _k\overline{\DD}$.
\end{itemize}

Condition $_k{\bf C}$,  in fact, plays the role of condition ${\bf C}$ for  the semi-Markov 
processes $_k\eta_{\e}(t)$.

If also condition ${\bf \hat{B}}$ is assumed to hold for the semi-Markov processes  $\eta_\e(t)$, and, in sequel, for  the semi-Markov processes  $\tilde{\eta}_\e(t)$, then transition probabilities $_kp_{\e, ij}  = p_{\e, ij} = \frac{1}{m} {\rm I}(j \in \DD), j \in \, _k\XX, i \in \DD$, for $\e \in [0, 1]$. Thus,
the asymptotic relations given in condition $_k{\bf C}$ also hold, for $j \in \, _k\XX,  i \in \DD$. 

Since, matrix $_k\mathbf{P}_{\e} = \|\, _kp_{\e, ij} \|$ is stochastic, for $\e \in (0, 1]$, conditions ${\bf \hat{B}}$ and 
$_k{\bf C}$  imply that  matrix $_k\mathbf{P}_{0} = \| \, _kp_{0, ij} \|$ is also stochastic. 

Let $_k\eta_{0, n}, n = 0, 1, \ldots$ be a Markov chain with the phase space $_k\XX$ and the matrix of transition probabilities $_k\mathbf{P}_{0}$. Conditions conditions ${\bf \hat{B}}$ and $_k{\bf C}$ make it possible to interpret the Markov chains $_k\eta_{\e, n}$, for 
$\e \in (0, 1]$, as perturbed version of the Markov chain $_k\eta_{0, n}$. 

Also, let us  consider the following condition for the semi-Markov processes $\eta_\e(t)$:
\begin{itemize}
\item [$_k{\bf \tilde{C}}$:]  $_k\tilde{p}_{\e, ij} = {\rm I}(j \neq i)\frac{_kp_{\e, ij}}{1 - \, _kp_{\e, ii}}  \to \, _k\tilde{p}_{0, ij} \in [0, 1]$ as $\e \to 0$, for 
$j \in \, _k\XX, i \in \, _k\overline{\DD}$.
\end{itemize}

Note that conditions ${\bf B}$ holds for the semi-Markov processes $_k\eta_\e(t)$, and, thus,  probabilities $_kp_{\e, ii} < 1, \e \in (0, 1]$, for $i \in \, _k\overline{\DD}$.

Probabilities $_k\tilde{p}_{\e, ii}  = 0, \e \in (0, 1]$, and, thus, probabilities $_k\tilde{p}_{0, ii}  = 0$, for $i \in \, _k\overline{\DD}$. Also, probabilities 
$_k\tilde{p}_{\e, ij}  = 0, \e \in (0, 1]$,  and, thus, probabilities $_k\tilde{p}_{0, ij}  = 0$, for $j \notin \, _k \YY_{1, i}, i \in \, _k\overline{\DD}$.

Condition $_k{\bf \tilde{C}}$  plays the role of condition ${\bf \tilde{C}}$ for semi-Markov 
processes $_k\eta_{\e}(t)$.

Let us also  define, sets, for $i \in \XX$ and $\e \in [0, 1]$, 
\begin{equation}\label{seta}
 _k\YY_{\e,  i}  = \{ j \in \, _k\XX: \, _kp_{\e,  ij} > 0 \}.
\end{equation}

Condition ${\bf B}$ implies that conditions ${\bf \tilde{B}}$ and  $_k{\bf B}$ hold, and, thus, for $i, k \in \overline{\DD}$ and 
$\e \in (0, 1]$, 
\begin{equation}\label{noity}
_k\YY_{\e,  i}  =  \, _k\YY_{1, i}. 
\end{equation}

Also, conditions ${\bf B}$, ${\bf C}$ and ${\bf \tilde{C}}$ imply that conditions   $_k{\bf B}$ and $_k{\bf C}$ hold, and, thus,
for $i, k \in \overline{\DD}$,
\begin{equation}\label{noityny}
_k\YY_{0,  i}  \subseteq  \, _k\YY_{1, i}. 
\end{equation}

Sets $\tilde{\YY}_{1,  i}$ and  $_k\YY_{1,  i}$, are connected by the following 
relations, for $i, k \in \overline{\DD}$,
\begin{equation}\label{noitynybut}
_k\YY_{1,  i}  = \left\{
\begin{array}{cl}
(\tilde{\YY}_{1, i} \cup  \tilde{\YY}_{1, k}) \setminus \{ k \} & \ \text{if} \ k \in \tilde{\YY}_{1, i}, \\
\tilde{\YY}_{1, i} & \ \text{if} \ k \notin \tilde{\YY}_{1, i}.
\end{array}
\right.
\end{equation}

Note also that, under condition ${\bf \hat{B}}$, sets $_k\YY_{\e, i} = \DD, \e \in [0, 1]$, for $i \in \DD$.  
 \vspace{1mm}

{\bf 4.5 Conditions of asymptotic comparability for transition probabilities of reduced embedded Markov chains}. It is natural to try to find condition which  would imply holding of conditions  $_k{\bf C}$ and $_k{\bf \tilde{C}}$, for any $k \in \overline{\DD}$, and  would be expressed more explicitly in  terms of initial transition probabilities $p_{\e, ij}$.  

Let us introduce the following condition:
\begin{itemize}
\item [$_k{\bf C}_{0}$:]  $_kq_{\e}[i i' i''] = \frac{_kp_{\e, i i'}}{_kp_{\e, i i''}} \to  \, _kq_{0}[i i' i''] \in [0, \infty]$ as $\e \to 0$, 
for $i' \in \, _k\XX, i'' \in \, _k\YY_{1, i}$, $i \in \, _k\overline{\DD}$.
\end{itemize}

Note that probabilities $_kp_{\e, i i'} = 0, \e \in (0, 1]$, for $i' \in \, _k\overline{\YY}_{1, i},  i \in \, _k\overline{\DD}$, and, thus, the asymptotic 
relation given in condition  $_k{\bf C}_{0}$ automatically holds, with limits  $_kq_{0}[i i' i''] = 0$, for 
$i' \in  \, _k\overline{\YY}_{1, i}, i'' \in \, _k\YY_{1, i}, i \in \, _k\overline{\DD}$. \vspace{1mm}

{\bf Lemma 7}. {\em Condition $_k{\bf C}_{0}$ is sufficient for holding of condition $_k{\bf C}$.} \vspace{1mm}

{\bf Proof}. Lemma 7 is a corollary of Lemma  1, which, just, should be applied to the semi-Markov processes $_k\eta_{\e}(t)$ instead of 
the semi-Markov processes $\eta_{\e}(t)$. In particular, the corresponding analogue of relation (\ref{quponrew}), which express limiting probabilities $p_{0, ij}$ as functions of quantities 
$q_{0, irj}$ takes the following form,  $j \in \, _k\YY_{1, i}, i \in \, _k\overline{\DD}$,
\begin{align}\label{quponrewad}
_kp_{\e, ij}  & = \frac{_kp_{\e, ij}}{\sum_{r \in \, _k\XX} \, _kp_{\e, ir}}  = \big(  \sum_{r \in \, _k\XX} \frac{_kp_{\e, ir}}{_kp_{\e, ij}} \big)^{-1}   
\vspace{2mm} \nonumber \\
& \to    \big(  \sum_{r \in \, _k\XX} \, _kq_{0}[irj]  \big)^{-1}  = \, _kp_{0, ij} \ {\rm as} \ \e \to 0.
\end{align}

Note that limits in relation (\ref{quponrewad})  satisfy relations,  $_kp_{0, ij} \geq 0, j \in \, _k\YY_{1, i}$, \, 
$\sum_{j \in \, _k\YY_{1, i}} \, _kp_{0, ij} = 1$, for $i \in \, _k\overline{\DD}$. $\Box$
 \vspace{1mm} 

{\bf Lemma 8}. {\em Condition $_k{\bf C}_{0}$ is sufficient for holding of condition  $_k{\bf \tilde{C}}$.} \vspace{1mm}

{\bf Proof}. Lemma 8 is a corollary of Lemma  2, which, just, should be applied to the semi-Markov processes $_k\eta_{\e}(t)$ instead of 
the semi-Markov processes $\eta_{\e}(t)$. In particular, the corresponding analogue of relation (\ref{quponta}), which express limiting probabilities $\tilde{p}_{0, ij}$ as functions of quantities 
$q_{0, irj}$ takes the following form,  $j \neq i, j \in \, _k\YY_{1, i}, i \in \, _k\overline{\DD}$,
\begin{align}\label{qupontaad}
_k\tilde{p}_{\e, ij}  & = \frac{_kp_{\e, ij}}{1 - \, _kp_{\e, ii}}  = \frac{_kp_{\e, ij}}{\sum_{r \neq i} \, _kp_{\e, ir}}  
 = \big(  \sum_{r \neq i, k} \frac{_kp_{\e, ir}}{_kp_{\e, ij}} \big)^{-1}   \vspace{2mm} \nonumber \\
& \to    \big(  \sum_{r \neq i, k} \, _kq_{0}[irj]  \big)^{-1}  =  \, _k\tilde{p}_{0, ij} \ {\rm as} \ \e \to 0.
\end{align}

Note that limits in relation  (\ref{quponta}) satisfy relations  $_k\tilde{p}_{0,  ij} \geq 0, j \neq i, j \in \, _k\YY_{1, i}$, \, 
$\sum_{j \neq i, j \in \, _k\YY_{1, i}} \, _k\tilde{p}_{0,  ij}  = 1$, for $i \in \, _k\overline{\DD}$. $\Box$ \vspace{1mm}
 
Let us introduce the following condition:
\begin{itemize}
\item [${\bf C}_{1}$:]  $q_{\e}[i i' i'', j j' j''] = \frac{p_{\e, i i'}p_{\e, j j'}}{p_{\e, i i''}p_{\e, j j''}}  \to  q_{0}[i i' i'', j j' j'']   \in [0, \infty]$ as $\e \to 0$, for $i', j' \in \XX, i'' \in \YY_{1, i}, j'' \in \YY_{1, j},  i, j \in \overline{\DD}$. 
\end{itemize}

Note that product $p_{\e, i i'}p_{\e, j j'} = 0, \e \in (0, 1]$, if  $i' \in \overline{\YY}_{1, i}$ or $j' \in \overline{\YY}_{1, j}$, for 
any $i, j \in \overline{\DD}$. In such cases, asymptotic relation given in condition ${\bf C}_{2}$ automatically holds and the corresponding limits 
$q_{0}[i i' i'', j j' j'']  = 0$.

Condition ${\bf C}_{1}$ is stronger than condition ${\bf C}_{0}$. Indeed, if to choose $j' = j''$ in the asymptotic relations given in condition ${\bf C}_{1}$, then these relations reduce to the corresponding asymptotic relations given in condition ${\bf C}_{0}$.

However, it is possible that condition ${\bf C}_{0}$ holds, while 
condition ${\bf C}_{1}$ does not. \vspace{1mm}

{\bf Lemma 9}. {\em Condition ${\bf C}_{1}$ is sufficient for holding of condition $_k{\bf C}_0$, and, thus, conditions $_k{\bf C}$ and $_k{\bf \tilde{C}}$, for any $k \in \overline{\DD}$.} \vspace{1mm}

{\bf Proof}. Condition ${\bf B}$ holds for the semi-Markov processes $\eta_\e(t)$, $\tilde{\eta}_\e(t)$, and $_k\eta_\e(t)$, and, thus, probabilities 
$_kp_{\e, ii''},  1 - p_{\e, ii} = \sum_{l \neq i} p_{\e, il}, 1 - p_{\e, kk} = \sum_{r \neq k} p_{\e, kr} \in (0, 1], \e \in (0,1]$, and, thus, 
$\sum_{r \neq k} p_{\e, ii''}p_{\e, kr} + p_{\e, ik}p_{\e, ki''} > 0$,  for $i'' \in \, _k\YY_{1, i}, i \in \, _k\overline{\DD}, k \in \overline{\DD}$. Therefore, the following relation takes place, for $i' \in \, _k\XX, i'' \in \, _k\YY_{1, i}, i \in \, _k\overline{\DD}, k \in \overline{\DD}$,
\begin{align}\label{qupona}
_kq_{\e}[ii'i'']  & = \frac{_kp_{\e, ii'}}{_kp_{\e, ii''}}  = \frac{\tilde{p}_{\e, ii'} + \tilde{p}_{\e, ik} \tilde{p}_{\e, ki'}}{ 
\tilde{p}_{\e, ii''} + \tilde{p}_{\e, ik}\tilde{p}_{\e, ki''}}  \vspace{2mm} \nonumber \\
& = \frac{\sum_{ l \neq k}p_{\e, ii'} p_{\e, kl} + p_{\e, ik}p_{\e, ki'}}{ 
\sum_{r \neq k} p_{\e, ii''}p_{\e, kr} + p_{\e, ik}p_{\e, ki''}}. 
\end{align}

Also, since condition ${\bf B}$ holds for the semi-Markov processes $\eta_\e(t)$, $\tilde{\eta}_\e(t)$ and $_k\eta_\e(t)$, 
probabilities $_kp_{\e, ii'} >0,  \e \in (0, 1]$, and, thus, $\sum_{ l \neq k}p_{\e, ii'} p_{\e, kl} + p_{\e, ik}p_{\e, ki'} > 0,  \e \in (0, 1]$, 
for $i' \in \, _k\YY_{1, i}, i \in \, _k\overline{\DD},  k \in \overline{\DD}$.  Also, every product $p_{\e, ii'} p_{\e, kl},  l \neq k$ and 
product $p_{\e, ik}p_{\e, ki'}$, either take positive value,  for every  $\e \in (0, 1]$, or equals $0$, for every  $\e \in (0, 1]$. Let us 
introduce indicators  ${\rm I}_{ii'kl} =  {\rm I}(p_{1, ii'} p_{1, kl} > 0),  l \neq k$ and ${\rm I}_{ikki'} = {\rm I}(p_{1, ik}  p_{1, ki'} > 0)$. 
By the above remarks, at least one of these indicators take value $1$, for every $i' \in \, _k\YY_{1, i}, i \in \, _k\overline{\DD}, 
k \in \overline{\DD}$.

Using the above remarks, conditions ${\bf C}_{0}$ and ${\bf C}_{1}$,   and relation  (\ref{qupona}), we get the following relation,  
for $i', i'' \in \, _k\YY_{1, i}, i \in \, _k\overline{\DD}, k \in \overline{\DD}$,
\begin{align}\label{optrew}
_kq_{\e}[ii'i'']  & = \sum_{l \neq k} {\rm I}_{ii'kl} 
\big(\sum_{r \neq k} \frac{p_{\e, ii''} p_{\e, kr}}{p_{\e, ii'} p_{\e, kl}} 
+ \frac{p_{\e, ik}p_{\e, ki''}}{p_{\e, ii'}p_{\e, kl}}) \big)^{-1} \vspace{2mm} \nonumber \\
&  \quad + {\rm I}_{ikki'} \big( \sum_{r \neq k} \frac{p_{\e, ir} p_{\e, ki''}}{p_{\e, ik} p_{\e, ki'} } + \frac{p_{\e, ki''}}{p_{\e, ki'} }) \big)^{-1} \vspace{2mm} \nonumber \\
& \to \sum_{l \neq k} {\rm I}_{ii'kl} \big( \sum_{r \neq k}   q_{0}[ii''i', krl] + q_{0}[iki', ki''l] ) \big)^{-1} 
\vspace{2mm} \nonumber \\
%\end{align*}
%\begin{align}\label{optrew}
& \quad + {\rm I}_{ikkj} \big( \sum_{r \neq k}  q_{0}[irk, ki''i'] + q_{0}[ki''i']) \big)^{-1} \vspace{2mm} \nonumber \\
& = \, _kq_{0}[ii'i'']  \in [0, \infty] \ {\rm as} \ \e \to 0, 
\end{align}
where every product of the form ${\rm I}_{\mathbf{\cdot}} ( \mathbf{\cdot} )^{-1}$ in the above sums should be counted as $0$ if the corresponding indicator equals $0$. $\Box$

In what follows, the following condition will also be used in computations  of the corresponding limiting Laplace transforms for transition times of perturbed reduced semi-Markov processes $_k\eta_\e(t)$:
\begin{itemize}
\item [$_k{\bf \hat{C}}$:]  $_k\hat{q}_{\e}[ij] =  \frac{\tilde{p}_{\e,  ij}}{_kp_{\e,   ij}}  
\to \, _k\hat{q}_{0}[ij] \in [0, 1]$ as $\e \to 0$, for $j \in \, _k\YY_{1, i}, i \in \, _k\overline{\DD}$.
\end{itemize}

Probabilities $\tilde{p}_{\e, ij} = {\rm I}(j \neq i)\frac{p_{\e, ij}}{1 - p_{\e, ii}}  = 0, \e \in (0, 1]$,  for $j \in \overline{\YY}_{1, i} \cup \{ i \}$, and, thus, limits $_k\hat{q}_{0}[ij]  = 0$, for $j \in ( \overline{\YY}_{1, i} \cup \{ i \}) \cap \, _k\YY_{1, i}, i \in \, _k\overline{\DD}$.\vspace{1mm}

{\bf Lemma 10}. {\em Condition ${\bf C}_{1}$ is sufficient for holding of condition $_k{\bf \hat{C}}$, for any $k \in \overline{\DD}$.} 
\vspace{1mm}

{\bf Proof}. Condition ${\bf B}$ holds for the semi-Markov processes $\eta_\e(t)$, $\tilde{\eta}_\e(t)$, and $_k\eta_\e(t)$, and, thus, probabilities 
$1 - p_{\e, ii} = \sum_{l \neq i} p_{\e, il}, 1 - p_{\e, kk} = \sum_{r \neq k} p_{\e, kr} \in (0, 1], \e \in (0,1]$, and, thus, 
$\sum_{r \neq k}p_{\e, ij}p_{\e, kr} + p_{\e, ik} p_{\e, kj} > 0$,  for $j \neq i, j \in \, _k\YY_{1, i}, i \in \, _k\overline{\DD}, k \in \overline{\DD}$. Therefore, the following relation takes place, for $j \neq i, j \in \, _k\YY_{1, i}, i \in \, _k\overline{\DD}, k \in \overline{\DD}$,
\begin{align}\label{guotras}
_k\hat{q}_{\e}[ij] & =  \frac{\tilde{p}_{\e,  ij}}{\tilde{p}_{\e, ij} + \tilde{p}_{\e, ik}\tilde{p}_{\e, kj}}   \vspace{2mm} \nonumber \\
& = \frac{\sum_{l \neq k}p_{\e, ij}p_{\e, kl}}{\sum_{r \neq k}p_{\e, ij}p_{\e, kr} + p_{\e, ik} p_{\e, kj}}. 
\end{align}

Also, since condition ${\bf B}$ holds for the semi-Markov processes $\eta_\e(t)$, $\tilde{\eta}_\e(t)$ and $_k\eta_\e(t)$, probabilities $\tilde{p}_{\e, ij} >0,  \e \in (0, 1]$, and, thus, $\sum_{ l \neq k}p_{\e, ij} p_{\e, kl}  > 0,  \e \in (0, 1]$, for $j \neq i, j \in \, _k\YY_{1, i}, i \in \, _k\overline{\DD}, k \in \overline{\DD}$.  Also, every product $p_{\e, ij} p_{\e, kl},  l \neq k$, either take positive value,  for every  $\e \in (0, 1]$, or equals $0$, for every  $\e \in (0, 1]$. Let us introduce indicators  ${\rm I}_{ijkl} =  {\rm I}(p_{1, ij} p_{1, kl} > 0),  l \neq k$. By the above remarks, at least one of these indicators take value $1$, for every $j \neq i, j \in \, _k\YY_{1, i}, i \in \, _k\overline{\DD}, k \in \overline{\DD}$.

Using the above remarks, conditions ${\bf C}_{0}$ and ${\bf C}_{1}$,   and relation  (\ref{guotras}), we get the following relation, for $j \neq i, j \in \, _k\YY_{1, i}, i \in \, _k\overline{\DD}, k \in \overline{\DD}$,
\begin{align}\label{guotrasad}
_k\hat{q}_{\e}[ij]  & = \sum_{l \neq k} {\rm I}_{ijkl} \big( \sum_{r \neq k} \frac{p_{\e, kr}}{p_{\e, kl}} 
+ \frac{p_{\e, ik} p_{\e, kj}}{p_{\e, ij} p_{\e, kl}} \big)^{-1}  \vspace{2mm} \nonumber \\ 
& \to \sum_{l \neq k} {\rm I}_{ijkl} \big( \sum_{r \neq k} q_{0}[krl] + q_{0}[ikj, kjl] \big)^{-1}  = \, _kq_{0, ij}    \ {\rm as} \ \e \to 0.
\end{align}

Since, $_k\hat{q}_{\e}[ij]  \in [0, 1], \e \in (0, 1]$, then $_k\hat{q}_{0}[ij]  \in [0, 1]$, for $j \neq 1, j \in \, _k\YY_{1, i}, i \in \, _k\overline{\DD}, k \in \overline{\DD}$. $\Box$ 

If condition  ${\bf \hat{B}}$  is assumed to hold, the asymptotic relation given in condition ${\bf C}_1$ also holds, for 
$i', j' \in \XX, i'' \in \YY_{1, i}, j'' \in \YY_{1, j}$,  $i \in \XX, j \in \DD$ and  $i \in \DD, j \in \XX$. 
In this case, probabilities $\tilde{p}_{\e, ij}   = \frac{1}{m}
{\rm I}(j \in \DD), j \in \XX, i \in \DD$, for $\e \in [0, 1]$, and sets $\YY_{1, i} =\DD, i \in \DD$. Thus, the above asymptotic relation holds and the corresponding limiting quantities $q_0[ii'i'', jj'j'']   = 
{\rm I}(i' \in \DD){\rm I}(j' \in \DD)$, for $i'', j'', i, j \in \DD$, $q_0[ii'i'', jj'j'']   = q_0[ii'i'']{\rm I}(j' \in \DD)$,  for $i'' \in \YY_{1, i}, i \in \overline{\DD} j'', j \in \DD$, and $q_0[ii'i'', jj'j'']   = {\rm I}(i' \in \DD)q_0[jj'j'']$,  for $i, i'', \DD, j' \in \YY_{1, j}, j \in \overline{\DD}$. $\Box$
\vspace{1mm}

{\bf 4.6 Conditions of asymptotic comparability for normalising functions for transition times of reduced  semi-Markov processes}. Let us introduce the following condition of asymptotic comparability for normalising functions:
\begin{itemize}
\item [${\bf \tilde{F}}$:] $\tilde{w}_{\e, ji}  = \frac{\tilde{v}_{\e, j}}{\tilde{v}_{\e, i}} =  
\frac{(1 - p_{\e, jj})^{-1} v_{\e, j}}{(1 - p_{\e, ii})^{-1} v_{\e, i}} \to \tilde{w}_{0, ji}  \in [0, \infty]$  as  $\e \to 0$, 
for $j, i \in \, \overline{\DD}$.
\end{itemize}

Also let us introduce condition:
\begin{itemize}
\item [$_k{\bf \hat{F}}$:] $\tilde{w}_{\e, ki} = \frac{\tilde{v}_{\e, k}}{\tilde{v}_{\e, i}} =  
\frac{(1 - p_{\e, kk})^{-1} v_{\e, k}}{(1 - p_{\e, ii})^{-1} v_{\e, i}} \to \tilde{w}_{0, ki}  
\in [0, \infty)$  as  $\e \to 0$, for $i \in \, \overline{\DD}$.
\end{itemize}

Note that $\tilde{w}_{\e, jj} = 1, \e \in [0, 1]$, for $j\in  \overline{\DD}$. 

In some sense, condition $_k{\bf \hat{F}}$  means that it is assumed that state $k$ is one of the least absorbing 
states in domain $\overline{\DD}$. 

Let denote by $\overline{\DD}^*$ the set of states $k \in \overline{\DD}$, for which condition $_k{\bf \hat{F}}$ holds. 

Condition ${\bf \tilde{F}}$ obviously imply that set $\overline{\DD}^* \neq \emptyset$. 

Obviously, $\tilde{w}_{0, ki} \in (0, \infty)$, for $k, i \in \overline{\DD}^*$, while $\tilde{w}_{0, ki} = 0$, for $k \in \overline{\DD}^*$, $i \in 
\overline{\DD} \setminus \overline{\DD}^*$. 

Also,   $\overline{\DD}^* = \overline{\DD}$, if $\tilde{w}_{0, ki} \in (0, \infty)$, for any $k, i \in \overline{\DD}$. 

Set $\overline{\DD}^*$  can be found with the use of  the following simple algorithm. Let us order by an arbitrary way states in domain 
$\overline{\DD}$, i.e., represent this domain  in the form $\overline{\DD} = \{ i_1, \ldots, i_{\bar{m}} \}$. Let  us now define states $i^*_n, n = 1, \ldots, \bar{m}$ using the following 
recurrent procedure, 
\begin{align}\label{proce}
& i^*_1 = i_1,   \vspace{2mm} \nonumber \\
& i_2^* = i^*_1 {\rm I}(\tilde{w}_{0, i_2 i^*_1} = \infty) +  i_2 {\rm I}(\tilde{w}_{0,  i_2 i^*_1} \in [0, \infty)), \vspace{2mm} \nonumber \\
& \ldots \vspace{2mm} \nonumber \\
& i^*_{\bar{m}} = i^*_{\bar{m} -1} {\rm I}(\tilde{w}_{0,  i_{\bar{m}} i^*_{\bar{m} -1}} = \infty)  
+  i_{\bar{m}} {\rm I}(\tilde{w}_{0,  i_{\bar{m}} i^*_{\bar{m} -1}} \in [0, \infty)). 
\end{align}

The following obvious proposition takes place. \vspace{1mm}

{\bf Lemma 11}. {\em Let condition ${\bf \tilde{F}}$ holds. Then, state $i^*_{\bar{m}} \in \overline{\DD}^*$ and set $\overline{\DD}^* = 
\{ i_n \in  \overline{\DD}:  \tilde{w}_{0,  i_n, i^*_{\bar{m}}} \in (0, \infty) \}$.} \vspace{1mm} 

It is also natural to try to find condition which would imply holding of conditions  ${\bf \tilde{F}}$ and  would be expressed more explicitly in  terms of initial transition probabilities $p_{\e, ij}$. 

Let us introduce condition: 
\begin{itemize}
\item [${\bf F}_0$:] $u_{\e}[i i', j j'']  =   
\frac{p_{\e, i i'} v_{\e, i}^{-1}}{p_{\e, j j''} v_{\e, j}^{-1}} \to u_{0}[i i', j j'']  \in [0, \infty]$  as  $\e \to 0$, for $i' \in \XX, 
j'' \in \YY_{1, j}, i, j \in \overline{\DD}$. 
\end{itemize}

Probabilities $p_{\e, i i'} = 0, \e \in (0, 1]$, for $i' \in \overline{\YY}_{1, i}, i \in  \overline{\DD}$ and, thus, limits $u_{\e, i i', j j''}  = 0$, for $i' \in \overline{\YY}_{1, i}, j'' \in \YY_{1, j}, i, j \in \overline{\DD}$.

\vspace{1mm}

 {\bf Lemma 12}. {\em Condition ${\bf F}_{0}$ is sufficient for holding of condition ${\bf \tilde{F}}$.} 
\vspace{1mm}

{\bf Proof}. Condition ${\bf B}$ holds for the semi-Markov processes $\eta_\e(t)$ and $\tilde{\eta}_\e(t)$, and, thus, probabilities 
$1 - p_{\e, ii} = \sum_{r \neq i} p_{\e, ir}, 1 - p_{\e, jj} = \sum_{s \neq j} p_{\e, js} \in (0, 1], \e \in (0,1]$,  for 
for $i, j \in \overline{\DD}$. Thus, the following relation takes place, for  $i, j \in \overline{\DD}$, 
\begin{align}\label{desa}
\tilde{w}_{\e, ji} & = \frac{(1 - p_{\e, jj})^{-1} v_{\e, j}}{(1 - p_{\e, ii})^{-1} v_{\e, i}} 
\vspace{2mm} \nonumber \\
& = \frac{(1 - p_{\e, ii}) v_{\e, i}^{-1}}{(1 - p_{\e, jj}) v_{\e, j}^{-1}}  \vspace{2mm} \nonumber \\
& = \frac{\sum_{r \neq i} p_{\e, ir}v_{\e, i}^{-1}}{\sum_{s \neq j} p_{\e, js}v_{\e, j}^{-1}}. 
\end{align}

Also, since condition ${\bf B}$ holds for the semi-Markov processes $\eta_\e(t)$ and $\tilde{\eta}_\e(t)$, probabilities $1 - \tilde{p}_{\e, ii} >0,  \e \in (0, 1]$, and, thus, $\sum_{ r \neq i}p_{\e, ir}  > 0,  \e \in (0, 1]$, for $i \in \overline{\DD}$.  Also, every product $p_{\e, ir},  r \neq i$, either take positive value,  for every  $\e \in (0, 1]$, or equals $0$, for every  $\e \in (0, 1]$. Let us introduce indicators  ${\rm I}_{ir} =  {\rm I}(p_{1, ir}  > 0),  
r \neq i$. By the above remarks, at least one of these indicators take value $1$, 
for every  $j, i \in \overline{\DD}$, 
\begin{align}\label{desana}
\tilde{w}_{\e, ji} & = \sum_{r \neq i} {\rm I}_{ir}  \big(\sum_{s \neq j} \frac{p_{\e, js}v_{\e, j}^{-1}}{p_{\e, ir}v_{\e, i}^{-1}} \big)^{-1} \vspace{2mm} \nonumber \\
& \to  \sum_{r \neq i} {\rm I}_{ir}  \big(\sum_{s \neq j} u_{0}[js, ir] \big)^{-1}  \vspace{2mm} \nonumber \\
& = \tilde{w}_{0, ji}  \in [0, \infty] \ {\rm as} \ \e \to 0, 
\end{align}
where every product of the form ${\rm I}_{\mathbf{\cdot}} ( \mathbf{\cdot} )^{-1}$ in the above sums should be counted as $0$ if the corresponding indicator equals $0$. $\Box$

The comparability condition ${\bf F}_{0}$  is expressed in a more explicit form via transition probabilities $p_{\e, ij}$ and normalising functions $v_{\e, i}$, than condition  ${\bf \tilde{F}}$. 

\vspace{1mm}

{\bf 4.7 Condition ${\bf D}$}. Let us assume that conditions ${\bf A}$ -- ${\bf C}$,  ${\bf D}$ (and, thus, ${\bf D}'$),  ${\bf E}$, and ${\bf \tilde{C}}$  hold.  In this case, conditions ${\bf \tilde{A}}$ -- ${\bf \tilde{B}}$,  ${\bf \tilde{D}}$ (and, thus, ${\bf\tilde{D}}'$),  ${\bf \tilde{E}}$, and $_k{\bf A}$ -- $_k{\bf C}$  hold. Let us additionally assume that conditions $_k{\bf \tilde{C}}$, $_k{\bf \hat{C}}$, and $_k{\bf \hat{F}}$  hold. 

Relations (\ref{transakak}) and (\ref{gopetk}) imply that the corresponding distribution functions $_kF_{\e, ij}(t), t \geq 0$ are defined by the  following relation, for $j \in \, _k\YY_{1, i}, i \in \,  _k\overline{\XX}$, 
\begin{equation}\label{traalk}
_kF_{\e, ij}(t)  =  \PP \{ _k\eta_{\e, 1} = j,  \, _k\kappa_{\e, 1} \leq t  / \, _k\eta_{\e, 0} = i, \, _k\eta_{\e,  1} = j \} \makebox[15mm]{}
\end{equation}
\begin{equation*}\label{transakalk}
\quad = \left\{
\begin{array}{lll}
F_{\e, ij}(t)\frac{p_{\e, ij}}{p_{\e, ij} + p_{\e ik} \tilde{p}_{\e,  kj}} & \vspace{2mm} \\ 
+    F_{\e, ik}(t) * \tilde{F}_{\e,  kj}(t) \frac{p_{\e ik} \tilde{p}_{\e,  kj}}{p_{\e, ij} + p_{\e ik} \tilde{p}_{\e,  kj}} 
&  \text{for} \ j \in \, _k\XX,  i \in \DD, \vspace{2mm} \\
 \tilde{F}_{\e, ij}(t)\frac{\tilde{p}_{\e,  ij}}{\tilde{p}_{\e,  ij} + \tilde{p}_{\e,  ik} \tilde{p}_{\e,   kj}} & \vspace{2mm} \\ 
 +   \tilde{F}_{\e,  ik}(t) * \tilde{F}_{\e,  kj}(t) \frac{\tilde{p}_{\e,  ik} \tilde{p}_{\e,  kj}}{\tilde{p}_{\e,  ij} 
 + \tilde{p}_{\e,  ik} \tilde{p}_{\e,   kj}} 
 & \text{for} \ j \neq i,  j \in \, _k\XX,  i \in \, _k\overline{\DD}, \vspace{2mm} \\
 \tilde{F}_{\e, ik}(t) * \tilde{F}_{\e,  ki}(t)    & \text{for} \ j = i,   i \in \, _k\overline{\DD}. 
 \end{array}
 \right.
\end{equation*}

Also, let us define, distribution functions $_kF_{\e,  i}(t), t \geq 0$, for $i \in \, _k\XX$,
\begin{align*}
_kF_{\e,  i}(t)  & = \PP \{ _k\kappa_{\e,  1} \leq t / \, _k\eta_{\e,  0} = i \} \vspace{2mm} \nonumber \\
\end{align*}
\begin{align}\label{diotraske}
& = \sum_{j \in \, _k\YY_{1,  i}} \, _kQ_{\e, ij}(t) = \sum_{j \in \, _k\YY_{1, i}} \, _kF_{\e,  ij}(t) \, _kp_{\e,  ij}.
\end{align}

Relations (\ref{traalk}) and (\ref{diotraske}) imply that the corresponding Laplace transforms $\, _k\phi_{\e,  ij}(s), s \geq 0$ take the following form,  for $j \in \, _k\YY_{1, i}, i \in \,  _k\overline{\DD}$, 
\begin{equation}\label{trawetk}
_k\phi_{\e,  ij}(s)  =  \EE \{ e^{- s \, _k\kappa_{\e,  1}}  /   \, _k\eta_{\e,  0} = i, \, _k\eta_{\e, 1} = j \}  \makebox[39mm]{}
\end{equation}
\begin{equation*}
\makebox[5mm]{} = \left\{
\begin{array}{lll}
\tilde{\phi}_{\e,  ij}(s) \frac{\tilde{p}_{\e,  ij}}{\tilde{p}_{\e,  ij} + \tilde{p}_{\e,  ik} \tilde{p}_{\e,   kj}}    & \vspace{2mm} \\ 
+ \tilde{\phi}_{\e,  ik}(s) \tilde{\phi}_{\e,  kj}(s) \frac{\tilde{p}_{\e,  ik} \tilde{p}_{\e,  kj}}{\tilde{p}_{\e,  ij} 
 + \tilde{p}_{\e,  ik} \tilde{p}_{\e,   kj}}    & \text{for} \ j \neq i,  j \in \, _k\YY_{1, i},  i \in \, _k\overline{\DD}, \vspace{2mm} \\
 \tilde{\phi}_{\e,  ik}(s) \tilde{\phi}_{\e,  ki}(s)  & \text{for} \ j = i \in \, _k\YY_{1, i},   i \in \, _k\overline{\DD},  
  \end{array}
 \right.
\end{equation*}
and the corresponding Laplace transforms $_k\phi_{\e,  i}(s), s \geq 0$ take the following form, for $i \in \, _k\overline{\DD}$,
\begin{align}\label{diotrasavk}
_k\phi_{\e,  i}(s)  & = \EE \{ e^{- s \, _k\kappa_{\e,  1}}  /   \, _k\eta_{\e,  0} = i \}  \nonumber \\
&  = \sum_{j \in \, _k\YY_{1, i}}  \, _k\phi_{\e,  ij}(s) \, _kp_{\e,  ij}, s \geq 0.
\end{align}

Relation (\ref{trawetk})  gives a hint to try to use the local normalisation functions  $_kv_{\e, i}, i \in \, _k\overline{\DD}$  (for getting  the asymptotic relations appearing in condition ${\bf D}$ for the semi-Markov processes $_k\eta_{\e}(t)$) defined by the following relation,  
\begin{equation}\label{normali}
_kv_{\e, i}  = \tilde{v}_{\e, i}, \ i \in \, _k\overline{\DD}.
\end{equation}

Relations  (\ref{trawetba}) and (\ref{trawetbama}) play the role of asymptotic relation given in condition ${\bf D}'$ for  semi-Markov 
processes $\tilde{\eta}_{\e}(t)$. Using these relations, relations (\ref{trawetk}), (\ref{normali}) and 
conditions  $_k{\bf \hat{C}}$ and ${\bf _k\hat{G}}$,  we get the following relation, for $j \in \, _k\YY_{1, i} 
\setminus \{i \}, i \in \, _k\overline{\DD}$,
\begin{equation*}
_k\phi_{\e, ij}(s / \, _kv_{\e, i})  \makebox[101mm]{}
\end{equation*}
\begin{equation*}
= \left\{
\begin{array}{lll}
\tilde{\phi}_{\e, ij}(s / \tilde{v}_{\e, i}) \frac{\tilde{p}_{\e,  ij}}{_kp_{\e,   ij}}
\vspace{2mm}  \\ 
+ \tilde{\phi}_{\e,  ik}(s / \tilde{v}_{\e,  i}) \tilde{\phi}_{\e,  kj}( \frac{\tilde{v}_{\e, k}}{\tilde{v}_{\e,  i}} s  / \tilde{v}_{\e, k}) 
\frac{\tilde{p}_{\e, ik} \tilde{p}_{\e,   kj}}{_kp_{\e,  ij}}  & \text{if} \ j, k \in \tilde{\YY}_{1, i}, j \in \tilde{\YY}_{1, k} \makebox[2mm]{}
\vspace{2mm} \\
\tilde{\phi}_{\e, ij}(s / \tilde{v}_{\e, i}) & \text{if} \ k \notin \tilde{\YY}_{1, i} \, \text{or} \,  j \notin \tilde{\YY}_{1, k}
\vspace{2mm} \\
\tilde{\phi}_{\e,  ik}(s / \tilde{v}_{\e,  i}) \tilde{\phi}_{\e,  kj}( \frac{\tilde{v}_{\e, k}}{\tilde{v}_{\e,  i}} s  / \tilde{v}_{\e, k})
 & \text{if} \ j \notin \tilde{\YY}_{1, i} \vspace{2mm}
\end{array}
\right.
\end{equation*}
\begin{equation}\label{trawetkasopa}
\to \, _k\phi_{0, ij}(s) \ {\rm as} \ \e \to 0, \ {\rm for} \ s \geq 0,  \makebox[51mm]{}
\end{equation} 
where
\begin{equation*}
_k\phi_{0, ij}(s)  \makebox[112mm]{}
\end{equation*}
\begin{equation}\label{trawetkpa}
\makebox[2mm]{}= \left\{
\begin{array}{lll}
\tilde{\phi}_{0, ij}(s) \, _kq_{0}[ij]
\vspace{2mm}  \\ 
+ \,  \tilde{\phi}_{0,  ik}(s) \tilde{\phi}_{0,  kj}(\tilde{w}_{0,  ki} s) (1 - \, _kq_0[ij]) 
  & \text{if} \ j, k \in \tilde{\YY}_{1, i}, j \in \tilde{\YY}_{1, k}, \makebox[2mm]{}
\vspace{2mm} \\
\tilde{\phi}_{0, ij}(s) & \text{if} \ k \notin \tilde{\YY}_{1, i} \, \text{or} \,  j \notin \tilde{\YY}_{1, k},
\vspace{2mm} \\
\tilde{\phi}_{0,  ik}(s) \tilde{\phi}_{0,  kj}( \tilde{w}_{0,  ki}  s)
 & \text{if} \ j \notin \tilde{\YY}_{1, i},  
\end{array}
\right.
\end{equation}

Since, $i \notin \tilde{\YY}_{1, i}$, then $j \in \, _k\YY_{1,  i} \cap \{ i \}$ if and only if  
$k \in \tilde{\YY}_{1, i}, j \in \tilde{\YY}_{1, k}$. This holds for any $i \in \, _k\overline{\DD}$.

That is why, for  $j \in \, _k\YY_{1,  i} \cap \{ i \},  i \in \, _k\overline{\DD}$,
\begin{align}\label{trawetkasf}
_k\phi_{\e,  ii}(s / \, _kv_{\e, i})  & = \tilde{\phi}_{\e, ik}(s / \tilde{v}_{\e, i}) \tilde{\phi}_{\e, kj}( \frac{\tilde{v}_{\e, k}}{\tilde{v}_{\e, i}}  
s / \tilde{v}_{\e, k})   \vspace{2mm} \nonumber \\
& \to  \, _k\phi_{0,  ii}(s)  \ {\rm as} \ \e \to 0, \ {\rm for} \ s \geq 0, 
\end{align}
where
\begin{equation}\label{nuwert}
_k\phi_{0,  ii}(s)  =  \tilde{\phi}_{0,  ik}(s) \tilde{\phi}_{0,  ki}(\tilde{w}_{0,  ki} s).
\end{equation}

Note   that, under condition ${\bf \hat{B}}$, the asymptotic relations (\ref{trawetkasopa}) also holds for 
$j \in \tilde{\YY}_{1,  i} = \DD, i \in \DD$, with limiting distribution functions $_kF_{0, ij}(t)    = {\rm I}(t \geq 1), \ t \geq 0, j, i \in \DD$, the limiting Laplace transforms $_k\phi_{0, ij}(s)    = e^{-s}, \ s \geq 0,  j, i \in \DD$,  and the  local normalising functions 
$_kv_{\e, i}  = \tilde{v}_{\e, i} = v_{\e,  i} \equiv 1, i \in \DD$. 

Also, the following relation holds, for $i \in \, _k\overline{\DD}$,
\begin{align}\label{trawetkasfa}
_k\phi_{\e,  i}(s / \, _kv_{\e, i}) &  \to \sum_{j \in \, _k\YY_{1, i}} \, _k\phi_{0, ij}(s)\, _kp_{0, ij} 
\vspace{2mm} \nonumber \\ 
& = \sum_{j \in \, _k\YY_{0, i}} \, _k\phi_{0, ij}(s)\, _kp_{0, ij}
\vspace{2mm} \nonumber \\ 
& = \, _k\phi_{0,  i}(s)  \ {\rm as} \ \e \to 0, \ {\rm for} \ s \geq 0.
\end{align}

In the case, where  $j \in \, _k\YY_{1,  i} \setminus \{ i \}, i \in \, _k\overline{\DD}$, the corresponding limiting random variable 
$_k\kappa_{0, ij}$, with the Laplace transform $_kf_{0,  ij}(s) = \EE e^{- s \, _k\kappa_{0,  ij}}$ = 
$\int_0^\infty \, e^{-st} \, _kF_{0,  ij}(dt), s \geq 0$,  can be represented in the  form, $_k\kappa_{0, ij} =$ 
$\tilde{\kappa}_{0, ij} \, _k\chi_{0,  ij}$ $+ (\tilde{\kappa}_{0,  ik}  + \tilde{w}_{0,  ki} \tilde{\kappa}_{0,  kj}) (1 - \, _k\chi_{0, ij})$, where: (a) $\tilde{\kappa}_{0, ij}, \tilde{\kappa}_{0,  ik}$, and $\tilde{\kappa}_{0,  kj}$ are random variables with 
Laplace transforms, respectively, $\tilde{\phi}_{0,  ij}(s)$,  $\tilde{\phi}_{0,  ik}(s)$, and $\tilde{\phi}_{0,  kj}(s)$, (b)  
$_k\chi_{0, ij}$ is a random variable taking two values $1$ and $0$ with probabilities, respectively,  
$_kq_{0}[ij]$ and $1 - \, _kq_{0}[ij]$, if $j, k \in \tilde{\YY}_{1, i}, j \in \tilde{\YY}_{1, k}$, or $1$ and $0$, if $k \notin \tilde{\YY}_{1, i} \, \text{or} \,  j \notin \tilde{\YY}_{1, k}$, or $0$ and $1$, if   $j \notin \tilde{\YY}_{1, i}$, (c) the random variables $\tilde{\kappa}_{0, ij}, \tilde{\kappa}_{0, ik}, \tilde{\kappa}_{0, kj}$ and $_k\chi_{0, ij}$ are independent.   

In the case, where $j \in \, _k\YY_{1,  i} \cap \{ i \},  i \in \, _k\overline{\DD}$, the corresponding limiting random variable 
$_k\kappa_{0,  ii}$, with the Laplace transform $_k\phi_{0,  ii}(s) =  \EE e^{- s \, _k\kappa_{0, ii}}$
$= \int_0^\infty \, e^{-st} \, _kF_{0,  ii}(dt), s \geq 0$, can be represented in the following form, $_k\kappa_{0,  ii}$  
$ = \tilde{\kappa}_{0,  ik}  + \tilde{w}_{0, ki} \tilde{\kappa}_{0,  ki}$, where 
(d) \, $\tilde{\kappa}_{0,  ik}$ and $\tilde{\kappa}_{0,  kj}$ are random variables with 
Laplace transforms, respectively, $\tilde{\phi}_{0, ik}(s)$ and  $\tilde{\phi}_{0, ki}(s)$, (e) random variables 
$\tilde{\kappa}_{0,  ik}$ and $\tilde{\kappa}_{0, ki}$ are independent. 

Condition ${\bf \tilde{D}}$  holds for semi-Markov processes $\tilde{\eta}_{\e}(t)$, and, thus, $\tilde{F}_{0,  ij}(0)$ $< 1$, for 
$j \in \YY_{1,  i}, i \in \overline{\DD}$.  This and relations (\ref{trawetkasopa}) and (\ref{trawetkasf}), 
obviously, imply that $_kF_{0,  ij}(0) < 1$, for $i \in \, _k\YY_{1,  i}, i \in \, _k\overline{\DD}$. 

Thus, relations (\ref{trawetkasopa}) and  (\ref{trawetkasf}) imply that condition ${\bf D}'$  and, thus, also condition ${\bf D}$,   holds for  semi-Markov 
processes $_k\eta_{\e}(t)$, with the Laplace transforms of the corresponding limiting distribution functions given in the above relations, and the normalisation functions  $_kv_{\e, i}, i \in \, _k\overline{\DD}$. 

 In what follows, we can also denote by $_k{\bf D}$  and $_k{\bf \tilde{D}}'$, respectively,   conditions ${\bf D}$ and ${\bf D}'$ for the semi-Markov processes $_k\eta_{\e}(t)$ 
 (expressed, for condition ${\bf D}'$,   in the form of relations  (\ref{trawetkasopa}) and  (\ref{trawetkasf})). 
 
  Note also that, if condition  ${\bf \hat{B}}$  is assumed to hold, the asymptotic relations given in conditions 
  $_k{\bf D}$ and $_k{\bf D}'$ also hold for $i \in \DD$. In this case, the normalisation functions $_kv_{\e, i} 
  = 1, \e \in (0, 1]$, the limiting distribution functions $_kF_{0, ij}(t) = \, _kF_{0, i}(t) = {\rm I}(t \geq 1), t \geq 0$, and the limiting Laplace transforms $_k\phi_{0, ij}(s) = \, _k\phi_{0, i}(s)  = e^{-s}, s \geq 0$, for $i, j \in \DD$.
 \vspace{1mm}

{\bf 4.8 Condition ${\bf E}$}. Let us assume that conditions ${\bf A}$ -- ${\bf C}$, ${\bf D}$ (and, thus, ${\bf D}'$), ${\bf E}$, and ${\bf \tilde{C}}$. In this case, conditions ${\bf \tilde{A}}$ -- ${\bf \tilde{C}}$,  ${\bf \tilde{D}}$ (and, thus, ${\bf\tilde{D}}'$),  ${\bf \tilde{E}}$  and $_k{\bf A}$ -- $_k{\bf C}$  hold. Let us additionally assume that conditions $_k{\bf \tilde{C}}$, $_k{\bf \hat{C}}$, and $_k{\bf \hat{F}}$  hold and, thus, condition $_k{\bf D}$
also holds. 

As follows from relation (\ref{traalk}), the corresponding expectation $\, _ke_{\e, ij}$ takes the following form, for 
$j \in \, _k\YY_{1, i}, i \in \,  \, _k\overline{\DD}$,
\begin{equation*}\label{twetk}
_ke_{\e, ij}  =  \EE \{  \, _k\kappa_{\e, 1}  /   \, _k\eta_{\e, 0} = i, \, _k\eta_{\e,  1} = j \} \makebox[66mm]{}
\end{equation*} 
\begin{equation}\label{diotrask} 
\makebox[2mm]{} = \left\{
\begin{array}{lll}
\tilde{e}_{\e,  ij} \frac{\tilde{p}_{\e,  ij}}{_kp_{\e,  ij}}  
+  (\tilde{e}_{\e,  ik} + \tilde{e}_{\e,  kj})\frac{\tilde{p}_{\e,  ik} \tilde{p}_{\e,  kj}}{_kp_{\e, ij}} &  \text{for} \ j \neq i,  j \in \, _k\YY_{1, i},  i \in \, _k\overline{\DD},   \vspace{2mm} \\
\tilde{e}_{\e,  ik} + \tilde{e}_{\e,  ki}  & \text{for} \ j = i \in \, _k\YY_{1, i},   i \in \, _k\overline{\DD}. 
 \end{array}
 \right.
\end{equation}

Also, for $i \in \, _k\overline{\DD}$,
\begin{align}\label{diotrasavomk}
_ke_{\e,  i} & = \EE \{_k\kappa_{\e,   1} /  \, _k\eta_{\e,  0} = i \}  \vspace{2mm} \nonumber \\ 
&  = \sum_{j \in \, _k\YY_{1,  i}} \, _ke_{\e,  ij} \, _kp_{\e,  ij}.
\end{align}
and, for $j \notin \, _k\YY_{1,  i}, i \in \,  _k\overline{\DD}$  
\begin{equation}\label{againab} 
_ke_{\e,  ij} =  \, _ke_{\e,  i}. 
\end{equation} 

Condition  ${\bf E}$  is represented for semi-Markov processes $\tilde{\eta}_{\e}(t)$ by relation 
(\ref{expofava}). Using this relation, relation (\ref{diotrasavomk}) and 
conditions ${\bf D}'$, $_k{\bf \hat{C}}$, and  $_k{\bf \hat{F}}$,   
we get the following relation, for $j \in \, _k\YY_{1, i}, i \in \, _k\overline{\DD}$,
\begin{equation*}
_ke_{\e,  ij} / \, _kv_{\e,  i} \makebox[108mm]{}
\end{equation*}
\begin{equation*}
= \left\{
\begin{array}{lll}
 \frac{\tilde{e}_{\e,  ij}}{\tilde{v}_{\e, i}} \frac{\tilde{p}_{\e,   ij}}{_kp_{\e,  ij}}  
 + \big( \frac{\tilde{e}_{\e,  ik}}{\tilde{v}_{\e,  i}} + \frac{\tilde{e}_{\e,  kj}}{\tilde{v}_{\e, k}} \frac{\tilde{v}_{\e, k}}{\tilde{v}_{\e, i}} \big) \frac{\tilde{p}_{\e,  ik} \tilde{p}_{\e, kj}}{_kp_{\e,  ij}}  & \ \text{if} \ \ j, k \in \tilde{\YY}_{1, i}, j \in \tilde{\YY}_{1, k} 
\vspace{2mm} \\
\frac{\tilde{e}_{\e,  ij}}{\tilde{v}_{\e, i}} \frac{\tilde{p}_{\e,   ij}}{_kp_{\e,  ij}} & \ \text{if} \ k \notin \tilde{\YY}_{1, i} \
\text{or} \  j \notin \tilde{\YY}_{1, k} 
\vspace{2mm} \\
\big( \frac{\tilde{e}_{\e,  ik}}{\tilde{v}_{\e,  i}} + \frac{\tilde{e}_{\e,  kj}}{\tilde{v}_{\e, k}} \frac{\tilde{v}_{\e, k}}{\tilde{v}_{\e, i}} \big) \frac{\tilde{p}_{\e,  ik} \tilde{p}_{\e, kj}}{_kp_{\e,  ij}} & \ \text{if} \ j \notin \tilde{\YY}_{1, i} \vspace{2mm} 
\end{array}
\right.
\end{equation*}
\begin{equation}\label{trawetkasnop}
\to \, _ke_{0,  ij} =  \int_0^\infty t \, _kF_{0, ij}(dt) \ {\rm as} \ \e \to 0, \makebox[45mm]{}  
\end{equation}
where
\begin{equation}\label{etkasnop}
_ke_{0,  ij}  = \left\{
\begin{array}{lll}
\tilde{e}_{0,  ij} \, _kq_{0, ij} & \vspace{2mm} \\
+  (\tilde{e}_{0,  ik} + \tilde{e}_{0,  kj} w_{0, ki}) (1 - \, _kq_{0, ij}) 
 & \ \text{if} \ \ j, k \in \tilde{\YY}_{1, i}, j \in \tilde{\YY}_{1, k},  \vspace{2mm} \\
\tilde{e}_{0,  ij}  & \ \text{if} \ k \notin \tilde{\YY}_{1, i} \
\text{or} \  j \notin \tilde{\YY}_{1, k}, \vspace{2mm} \\
\tilde{e}_{0,  ik} + \tilde{e}_{0,  kj} w_{0, ki}  & \ \text{if} \ j \notin \tilde{\YY}_{1, i}.
\end{array}
\right.
\end{equation}

Also, for  $i \in \,  _k\overline{\DD}$,
\begin{align}\label{trawetkasfanop}
\frac{_ke_{\e,  i}}{ \, _kv_{\e, i}} &  \to  \sum_{j \in \, _k\YY_{1,  i}} \, _ke_{0,  ij} \, _kp_{0, ij} 
\vspace{2mm} \nonumber \\
& = \sum_{j \in \, _k\YY_{0,  i}} \, _ke_{0,  ij} \, _kp_{0, ij}  = \, _ke_{0,  i}   \ {\rm as} \ \e \to 0.
\end{align}

Thus, relation (\ref{trawetkasnop}) implies that condition ${\bf E}$  holds for  semi-Markov 
processes $_k\eta_{\e}(t)$, with  the expectations given in the above relation, and the normalisation functions  $_kv_{\e, i}, i \in \, _k\overline{\DD}$. 

In what follows, we can also denote by $_k{\bf E}$  condition ${\bf E}$  for the semi-Markov processes $_k\eta_{\e}(t)$ 
 (expressed   in the form of relations  (\ref{trawetkasnop})). 
 
 Note also that, if condition  ${\bf \hat{B}}$  is assumed to hold, the asymptotic relations given in condition 
  $_k{\bf E}$ also hold for $i \in \DD$. In this case, the corresponding limiting expectations $_ke_{0, ij} = \, _ke_{0, i} = 1, i, j \in \DD$.

 \vspace{1mm}

{\bf 4.9 Summary}. The following lemma summarises the above remarks.  \vspace{1mm}

{\bf Lemma 13}. {\em Let conditions ${\bf A}$ -- ${\bf E}$, and ${\bf \tilde{C}}$, $_k{\bf \tilde{C}}$, $_k{\bf \hat{C}}$, $_k{\bf \hat{F}}$  hold. Then, 
conditions ${\bf A}$ -- ${\bf E}$ and  ${\bf \tilde{C}}$ also hold for  the  semi-Markov processes 
$_k\eta_{\e}(t)$, respectively, in the form of conditions  $_k{\bf A}$ -- $_k{\bf E}$ and  $_k{\bf \tilde{C}}$.}
\vspace{1mm}

{\bf Remark 6}. Conditions ${\bf C}$ and ${\bf \tilde{C}}$ imply holding of condition $_k{\bf C}$, which plays the role of condition
${\bf C}$ for the semi-Markov processes  $_k\eta_{\e}(t)$. Condition $_k{\bf \tilde{C}}$ plays the role of condition
${\bf \tilde{C}}$ for the semi-Markov processes  $_k\eta_{\e}(t)$.
\vspace{1mm}

{\bf Remark 7}. Condition ${\bf C}_{2}$ implies that conditions ${\bf C}$, ${\bf \tilde{C}}$, ${\bf \tilde{C}}'$, $_k{\bf C}$, $_k{\bf \tilde{C}}$. and $_k{\bf \hat{C}}$ hold.

\vspace{1mm}

{\bf Remark 8}. Condition  ${\bf \tilde{F}}$  implies that set $\overline{\DD}^*$ of states $k \in \overline{\DD}$, for which condition ${\bf _k\hat{G}}$ holds, is not empty. \vspace{1mm}
Condition ${\bf F}_0$ implies holding of condition ${\bf \tilde{F}}$. \vspace{1mm}

{\bf 4.10  Hitting times for semi-Markov processes with reduced phase space}. Let us introduce hitting times for semi-Markov processes $_k\eta_{\e}(t)$, 
\begin{equation}\label{opop}
_k\tau_{\e, \DD} = \sum_{n = 1}^{_k\nu_{\e, \DD}} \, _k\kappa_{\e,  n}, \ {\rm where} \  
_k\nu_{\e,  \DD} = \min(n \geq 1: \, _k\eta_{\e, n} \in \DD).
\end{equation}

The definitions of semi-Markov processes $\eta_{\e}(t)$, $\tilde{\eta}_{\e}(t)$,  and 
$_k\eta_{\e}(t)$ imply that the following lemma takes place. \vspace{1mm}

{\bf Lemma 14}. {\em Let conditions ${\bf A}$, ${\bf B}$ hold, and, in sequel, conditions ${\bf \tilde{A}}$, ${\bf \tilde{B}}$, $_k{\bf A}$,  $_k{\bf B}$ hold. Then, the following relations takes place, for 
 $\e \in (0, 1]$,
\begin{align}\label{idereba}
& \PP_i \{ \tau_{\e, \DD} = \tilde{\tau}_{\e, \DD} =  \, _k\tau_{\e,  \DD}, 
\makebox[69mm]{} \vspace{2mm} \nonumber \\  
& \quad \quad  \eta_{\e}(\tau_{\e \DD}) = \tilde{\eta}_{\e}(\tilde{\tau}_{\e \DD}) = \,   _k\eta_{\e}(_k\tau_{\e, \DD}) \} = 1, \ i \in \, _k\XX,  
\end{align}
and
\begin{align}\label{iderebas}
& \PP_{k} \{ \tau_{\e, \DD} =  \tilde{\tau}_{\e, \DD} =  \tilde{\kappa}_{\e, 1} {\rm I}(\tilde{\eta}_{\e, 1} \in \DD) +
(\tilde{\kappa}_{\e, 1}  + \,  _k\tau_{\e, \DD}) {\rm I}(\tilde{\eta}_{\e, 1} \in \, _k\overline{\DD}),
\vspace{2mm} \nonumber \\
& \quad \quad  \quad  \eta_{\e}(\tau_{\e, \DD})  = \tilde{\eta}_{\e}(\tilde{\tau}_{\e, \DD})  
= \tilde{\eta}_{\e, 1} {\rm I}(\tilde{\eta}_{\e, 1} \in \DD) \vspace{2mm} \nonumber \\
& \quad \quad  \quad \quad \quad \quad  + \, _k\eta_{\e}(_k\tau_{\e, D}) {\rm I}(\tilde{\eta}_{\e, 1} \in \, _k\overline{\DD}) \} = 1.  
\end{align}}

{\bf Proof}. Recall that state $k \in \overline{\DD}$, $_k\overline{\DD} = \overline{\DD} \setminus \{k \}$, and $_k\XX = \XX \setminus \{k \} = \DD \cup \, _k\overline{\DD}$.  
 
Relation (\ref{iderek}) given in Lemma 4 implies that, for  $\e \in (0,1]$,
\begin{equation}\label{nersa}
\PP_i \{\tau_{\e, \DD} =  \tilde{\tau}_{\e, \DD}, 
\eta_{\e}(\tau_{\e \DD}) = \tilde{\eta}_{\e}(\tilde{\tau}_{\e \DD}) \} = 1, \ i \in \XX.
\end{equation}

Relations (\ref{fotyr}) -- (\ref{recurask}) imply that $\beta_{\e, 0} = 0$, if $\tilde{\eta}_{\e, 0} = i \in \, _k\XX$. In this case, the hitting times $\tilde{\nu}_{\e, \DD}$ and $_k\nu_{\e, \DD}$  are connected by the following relation, for $\e \in (0,1]$,
 \begin{equation}\label{relanaka}
\PP_i\{ \tilde{\nu}_{\e, \DD}  = \, _k\beta_{\e, \, _k\nu_{\e, \DD}}  \} = 1, \ i \in \, _k\XX. 
 \end{equation}

Thus,  the  hitting times $\tilde{\tau}_{\e, \DD}$ and $_k\tau_{\e, \DD}$ are connected by the following relation, for $\e \in (0,1]$,
\begin{align}\label{relanak}
\PP_i \{ \tilde{\tau}_{\e, \DD}  & = \sum_{n = 1}^{\tilde{\nu}_{\e, \DD}} \tilde{\kappa}_{\e, n}  
 = \sum_{n = 1}^{_k\beta_{\e, \, _k\nu_{\e, \DD}}} \tilde{\kappa}_{\e, n}  \vspace{2mm} \nonumber \\
&  =  \sum_{n = 1}^{_k\nu_{\e, \DD}} \, _k\kappa_{\e, n} 
= \, _k\tau_{\e, \DD} \} = 1, \ i \in \, _k\XX,   
\end{align}
and the random variables $\tilde{\eta}_\e(\tilde{\tau}_{\e, \DD})$ and $_k\eta_\e(_k\tau_{\e, \DD})$ are connected by the following relation, for $\e \in (0,1]$, 
\begin{align}\label{relanakboj}
\PP_i \{ \tilde{\eta}_\e(\tilde{\tau}_{\e, \DD})  &  = \tilde{\eta}_{\e, \tilde{\nu}_{\e, \DD}}  =  
\tilde{\eta}_{\e, _k\beta_{\e, \, _k\nu_{\e, \DD}}} \vspace{2mm} \nonumber \\
& = \,   _k\eta_{\e, \, _k\nu_{\e, \DD}}  = \, _k\eta_\e(_k\tau_{\e, \DD}) \} = 1, i \in \, _k\XX. 
\end{align}

Relations (\ref{nersa}), (\ref{relanak}),  and (\ref{relanakboj}) imply that relation (\ref{idereba}) holds.

Relations (\ref{fotyr}) -- (\ref{recurask}) also imply that $\beta_{\e, 0} = 1$ and $\tilde{\eta}_{\e, 1} \in 
\, _k\XX$, if $\tilde{\eta}_{\e, 0} = k$. In this case, relation, analogous to (\ref{relanaka}),  takes the following form, for $\e \in (0, 1]$, 
 \begin{equation}\label{relanakamva}
\PP_k\{ \tilde{\nu}_{\e, \DD}  = {\rm I}(\tilde{\eta}_{\e, 1} \in \DD) + 
\, _k\beta_{\e, \, _k\nu_{\e, \DD}}{\rm I}(\tilde{\eta}_{\e, 1} \in \, _k\overline{\DD})  \} = 1. 
 \end{equation}
 
Thus,  the  hitting times $\tilde{\tau}_{\e, \DD}$ and $_k\tau_{\e, \DD}$ are connected by the following relation, for $\e \in (0,1]$,
\begin{align*}
\PP_k \{ \tilde{\tau}_{\e, \DD}  & = \sum_{n = 1}^{\tilde{\nu}_{\e, \DD}} \tilde{\kappa}_{\e, n}
\makebox[70mm]{} \vspace{2mm} \nonumber \\
& =  \tilde{\kappa}_{\e, 1} {\rm I}(\tilde{\eta}_{\e, 1} \in \DD) + (\tilde{\kappa}_{\e, 1}  + \sum_{n = 2}^{_k\beta_{\e, \, _k\nu_{\e, \DD}}} \tilde{\kappa}_{\e, n}) {\rm I}(\tilde{\eta}_{\e, 1} \in \, _k\overline{\DD})  
\vspace{2mm} \nonumber \\
\end{align*}
\begin{align}\label{relanaktov}
& = \tilde{\kappa}_{\e, 1} {\rm I}(\tilde{\eta}_{\e, 1} \in \DD)  + (\tilde{\kappa}_{\e, 1}  
+ \sum_{n = 1}^{_k\nu_{\e, \DD}} \, _k\kappa_{\e, n}) {\rm I}(\tilde{\eta}_{\e, 1} \in \, _k\overline{\DD}) 
\vspace{2mm} \nonumber \\
& =  \tilde{\kappa}_{\e, 1} {\rm I}(\tilde{\eta}_{\e, 1} \in \DD) + (\tilde{\kappa}_{\e, 1}  + \, _k\tau_{\e, \DD}) {\rm I}(\tilde{\eta}_{\e, 1} \in \, _k\overline{\DD}) \} = 1,   
\end{align}
and the random variables $\tilde{\eta}_\e(\tilde{\tau}_{\e, \DD})$ and $_k\eta_\e(_k\tau_{\e, \DD})$ are connected by the following relation, for $\e \in (0,1]$, 
\begin{align}\label{iderebasburn}
& \PP_{k} \{  \tilde{\eta}_{\e}(\tilde{\tau}_{\e, \DD})  
= \tilde{\eta}_{\e, 1} {\rm I}(\tilde{\eta}_{\e, 1} \in \DD) + \, _k\eta_{\e}(_k\tau_{\e, D}) {\rm I}(\tilde{\eta}_{\e, 1} \in \, _k\overline{\DD}) \} = 1. 
\end{align}

Relations (\ref{nersa}),  (\ref{relanaktov}), and (\ref{iderebasburn}) imply that relation (\ref{iderebas}) holds. $\Box$

\vspace {1mm}

Let us also introduce distributions, for $i \in \, _k\XX, k \in \overline{\DD}$ and $\e \in (0, 1]$,
\begin{equation}\label{moku}
_kG_{\e,  \DD, ij}(t) = \PP_i \{_k\tau_{\e, \DD} \leq t, \, _k\eta_{\e}(_k\tau_{\e \DD}) = j \},  \ t \geq 0, j \in \DD. 
\end{equation}

The following lemma is a corollary of Lemma 4. \vspace{1mm}

{\bf Lemma 15}. {\em Let conditions ${\bf A}$, ${\bf B}$ hold, and, in sequel, conditions ${\bf \tilde{A}}$, ${\bf \tilde{B}}$, $_k{\bf A}$, $_k{\bf B}$ hold. Then, the following relation takes place, for  
$\e \in (0, 1]${\rm :} \vspace{1mm}  
\begin{equation}\label{iderebanoma}
G_{\e, \DD,  ij}(t)  = \tilde{G}_{\e, \DD,  ij}(t) = \, _kG_{\e, \DD, ij}(t), \, t \geq 0,  j \in \DD, \,  i \in \, _k\XX, 
\end{equation}
and
\begin{align}\label{iderebasnom}
& G_{\e, \DD, kj}(t)   = \tilde{G}_{\e, \DD, kj}(t)    \vspace{3mm} \nonumber \\
& \quad \quad = \tilde{F}_{\e, k j}(t)\tilde{p}_{\e,  kj} + \sum_{r \in \, _k\overline{\DD}} (\tilde{F}_{\e, kr}(t) * \,  _kG_{\e, \DD, rj}(t) \tilde{p}_{\e, kr} \vspace{2mm} \nonumber \\
& \quad \quad = \tilde{F}_{\e, k j}(t)\tilde{p}_{\e,  kj}  +  \sum_{r \in \, _k\overline{\DD}} (\tilde{F}_{\e, kr}(t) * G_{\e, \DD, rj}(t) \tilde{p}_{\e, kr}, \, 
t \geq 0,  j \in \DD.  
\end{align}}
\makebox[3mm]{} {\bf Proof}. Equalities  given in relations (\ref{iderebanoma}) are obvious corollaries of relation (\ref{idereba}).

The first equality given in relations (\ref{iderebasnom}) is  follows from relation 
(\ref{iderebas}).

Relations (\ref{fotyr}) -- (\ref{recurask}) imply that, in the case where $\tilde{\eta}_{\e, 0} = k$ the random variable $\beta_{\e, 0} = 1$ and, thus,  the random functional $(_k\tau_{\e, \DD}, \, _k\eta_{\e}(_k\tau_{\e, \DD}))$ is determined  by trajectory of the Markov renewal process $(\tilde{\eta}_{\e, n}, \tilde{\kappa}_{\e, n})$ for $n \geq 1$. 
This makes it possible to use relation (\ref{iderebas}) and  the Markov property of Markov renewal process $(\tilde{\eta}_{\e, n}, \tilde{\kappa}_{\e, n})$,  and  to get  the following relation, for $t \geq 0, j \in \DD$ and $\e \in (0, 1]$,
\begin{align}\label{basd}
\tilde{G}_{\e, \DD, kj}(t)  & = \PP_k \{ \tilde{\kappa}_{\e, 1} \leq t, \tilde{\eta}_{\e, 1} = j \} 
\vspace{2mm} \nonumber \\
& \quad + \sum_{r \in \, _k\overline{\DD}} \PP_k \{ \tilde{\kappa}_{\e, 1} +  \, 
_k\tau_{\e, \DD} \leq t,  \, _k\eta_{\e}(_k\tau_{\e, \DD}) = j, \tilde{\eta}_{\e, 1}  = r \} 
\vspace{2mm} \nonumber \\
& = \tilde{F}_{\e, k j}(t)\tilde{p}_{\e,  kj} + \sum_{r \in \, _k\overline{\DD}} (\tilde{F}_{\e, kr}(t) * \,  _kG_{\e, \DD, rj}(t) \tilde{p}_{\e, kr} 
\end{align}

Thus, the second equality given in relation (\ref{iderebasnom}) holds.

The third equality  given in relation (\ref{iderebasnom})  follows from the second equality given in this relation and relation (\ref{iderebanoma}). $\Box$ \vspace{1mm}

\makebox[3mm]{} Let us also introduce Laplace transforms, for $i \in \, _k\XX$ and $\e \in (0, 1]$,  
\begin{equation}\label{trewqmo}
_k\Psi_{\e, \DD, ij}(s) = \EE_i \exp\{- s \, _k\tau_{\e, \DD}\} {\rm I}(_k\eta_{\e}(_k\tau_{\e, \DD}) = j) \}, 
\, s \geq 0, j \in \DD.
\end{equation}

The following lemma re-formulates propositions of Lemma 15 in the equivalent form of asymptotic relations for Laplace transforms of hitting times. 

\vspace{1mm}

{\bf Lemma 16}. {\em Let conditions ${\bf A}$, ${\bf B}$ hold, and, in sequel, conditions ${\bf \tilde{A}}$, ${\bf \tilde{B}}$, $_k{\bf A}$, $_k{\bf B}$ hold. Then, the following relation takes place, for  
$\e \in (0, 1]$, 
 \begin{equation}\label{mokubada}
\Psi_{\e, \DD, ij}(s)  = \tilde{\Psi}_{\e, \DD, ij}(s) = \, _k\Psi_{\e, \DD, ij}(s), \ s \geq 0,  j \in \DD, \, i \in \, _k\XX \makebox[12mm]{}
\end{equation} 
and 
\begin{align}\label{iderenasfat}
& \Psi_{\e, \DD, kj}(s)  = \tilde{\Psi}_{\e, \DD,  kj}(s)  \vspace{2mm} \nonumber \\ 
& \quad \quad \ =  \tilde{\phi}_{\e, kj}(s)\tilde{p}_{\e, kj}  + \sum_{r \in \, _k\overline{\DD}}  \, 
_k\Psi_{\e, \DD, rj}(s) \tilde{\phi}_{\e, kr}(s) \tilde{p}_{\e,  kr} 
\vspace{2mm} \nonumber \\ 
& \quad \quad \ =  \tilde{\phi}_{\e, kj}(s)\tilde{p}_{\e, kj}  + \sum_{r \in \, _k\overline{\DD}}  \Psi_{\e, \DD, rj}(s) 
\tilde{\phi}_{\e, kr}(s) \tilde{p}_{\e,  kr}, \ s \geq 0,  j \in \DD.
\end{align}}

Lemmas 4 -- 6 and 14 -- 16  let  one reduce study of asymptotics for distributions of hitting times  for 
semi-Markov processes $\eta_\e(t)$ and $\tilde{\eta}_{\e}(t)$  to the case of more simple semi-Markov processes $_k\eta_{\e}(t)$.  

The distributions $_kG_{\e, \DD, ij}(t), t \geq 0, j \in \DD,  i \in \, _k\overline{\DD}$ are determined by transition probabilities $_kQ_{\e, \DD, lr}(t), t \geq 0,   r \in \, _k\XX, l \in \, _k\overline{\DD}$, which itself
are determined by transition probabilities  $\tilde{Q}_{\e, \DD, lr}(t), t \geq 0,   r \in \XX, l \in \, 
\overline{\DD}$. The latter transition probabilities  are determined by transition probabilities 
$Q_{\e, lr}(t), t \geq 0,   r \in \XX, l \in \, \overline{\DD}$.

This makes it possible to essentially simplify the model in the case, where we are interested to investigate asymptotics of hitting times for some fixed domain $\DD$ and only for initial states $i \in \overline{\DD}$. This case in considered in Sections 2 -- 7. For example, we can assume in what follows in these sections that the transition probabilities  of semi-Markov processes $\eta_\e(t)$ and, in sequel, the transition probabilities of semi-Markov processes  $\tilde{\eta}_\e(t)$ and $_k\eta_\e(t)$, satisfy condition 
${\bf \hat{B}}$. 

The above remarks and Theorem 1 let us describe asymptotics for distributions of hitting times for another simplest case, where domain $\overline{\DD} = \{ i, k \}$  is a two-states set.  

Let us  assume that condition $_k{\bf F}$ holds,  i.e., state $k$ is less or equally absorbing with state $i$. 

According Lemma 15, the following relation takes place, 
\begin{equation}\label{rytew}
G_{\e, \DD, ij}(t) =\, _kG_{\e, \DD, ij}(t),  t \geq 0, j \in \DD.
\end{equation}

 It is obvious that, in this case,  domain  $_k\overline{\DD} = \{ i \}$ is a 
one-state set. 

This makes it possible to apply Theorem 1 to the semi-Markov processes 
$_k\eta_\e(t)$ and, due to relation (\ref{rytew}),   to describe asymptotics for distributions $G_{\e, \DD, ij}(t)$. 

However, to do this, we first should apply the procedure of removing virtual transitions to the semi-Markov processes $_k\eta_\e(t)$ and to construct the corresponding semi-Markov processes  
$_k\tilde{\eta}_\e(t)$. \\

{\bf 5. Removing of Virtual Transitions  for  Perturbed \\ \makebox[10mm]{} Reduced Semi-Markov Processes} \\ 

In this section, we shortly describe the procedure of removing virtual transitions of the form $i \to i$ from  trajectories of perturbed reduced semi-Markov processes $_k\eta_{\e}(t)$. This procedure is analogous to the of removing virtual transitions  for the semi-Markov processes
$\eta_{\e}(t)$ described in Section 3. The difference is only that the semi-Markov processes $\eta_{\e}(t)$ should be replaced by the semi-Markov processes $_k\eta_{\e}(t)$. This let us to shorten this description, just, by replacing computations and proofs by 
references  to the corresponding computations and proofs given in Section 3.
\vspace{1mm}

{\bf 5.1 Removing of virtual transition for perturbed reduced semi-Markov processes}. 
Let us assume that $\e \in (0, 1]$
and  conditions ${\bf A}$ and ${\bf B}$  hold for the semi-Markov processes $\eta_\e(t)$ and, thus, for the 
semi-Markov processes $\tilde{\eta}_\e(t)$ and $_k\eta_\e(t)$.  

Let us define stopping times for Markov chain $_k\eta_{\e, n}$ that are, for  $r = 0, 1, \ldots$,
\begin{equation}\label{stopasd}
_k\theta_{\e}[r] = {\rm I}(_k\eta_{\e, r} \in \, _k\overline{\DD}) \min(n > r: \, _k\eta_{\e, n} \neq \, _k\eta_{\e, r}) + {\rm I}(_k\eta_{\e, r} \in \, _k\DD)(r +1).
\end{equation}

By the definition,  $_k\theta_{\e}[r]$ is, either the first after $r$ moment of change of state $_k\eta_{\e, r}$ by the Markov chain $_k\eta_{\e, n}$, if $_k\eta_{\e, r} \in \, _k\overline{\DD}$, or $r +1$,  if $_k\eta_{\e, r} \in \, _k\DD$. One can refer to $_k\theta_{\e}[r]$ as to the conditional first after $r$ moment of change the state $_k\eta_{\e, r}$ by the Markov chain $_k\eta_{\e, n}$.

Let us also  define sequential stopping times,
\begin{equation}\label{opasdasd}
_k\mu_{\e, n} = \, _k\theta_{\e}[_k\mu_{\e, n-1}], n = 1, 2, \ldots, \ {\rm where} \ _k\mu_{\e, 0} = 0.
\end{equation}

Let us now construct a new Markov renewal process $(_k\tilde{\eta}_{\e,  n}, \, _k\tilde{\kappa}_{\e,  n}), n = 0, 1, \ldots$ with the phase space 
$_k\XX \times [0, \infty)$ using the following recurrent relations,
\begin{equation}\label{recurasd}
(_k\tilde{\eta}_{\e,  n}, \, _k\tilde{\kappa}_{\e, n}) = \left\{ 
\begin{array}{cll}
(_k\eta_{\e, 0}, 0) & \text{for} \ n = 0,  \vspace{2mm} \\
(_k\eta_{\e,  \, _k\mu_{\e, n}}, \sum_{l =  \, _k\mu_{\e,  n-1} + 1}^{_k\mu_{\e, n}} \, _k\kappa_{\e, l}) & \text{for} \ n = 1, 2, \ldots.  
\end{array}
\right.
\end{equation}
 
We also can define the corresponding semi-Markov process,
\begin{equation}\label{seminasd}
_k\tilde{\eta}_{\e}(t) = \, _k\tilde{\eta}_{\e,  \, _k\tilde{\nu}_{\e}(t)}, t \geq 0,
\end{equation}
where $_k\tilde{\zeta}_{\e,  n} = \, _k\tilde{\kappa}_{\e,  1} + \cdots + \, _k\tilde{\kappa}_{\e,   n}, n = 1, 2, \ldots, \, _k\tilde{\zeta}_{\e,   0} = 0$, are the corresponding instants of jumps, and $_k\tilde{\nu}_{\e}(t) = \max(n \geq 1: \, _k\tilde{\zeta}_{\e,  n}  \leq t)$ is the number of jumps in an interval $[0, t], t \geq 0$ for the above semi-Markov process. 

The transition probabilities for the above 
Markov renewal processes are determined by the following relation analogous to relation (\ref{transaka}),  
\begin{equation*}\label{transakan}
_k\tilde{Q}_{\e, ij}(t) =  \PP \{ _k\tilde{\eta}_{\e, 1} = j,  \, _k\tilde{\kappa}_{\e, 1} \leq t  / \, _k\tilde{\eta}_{\e, 0} = i \} \makebox[46mm]{}
\end{equation*}
\begin{equation}\label{transakannop}
\makebox[13mm]{} = \left\{
\begin{array}{cll}
_kQ_{\e, ij}(t) & \text{for} \  t \geq 0, j \in \, _k\XX, i \in \, _k\DD, \vspace{2mm} \\
0 & \text{for} \  t \geq 0, j = i, i \in \, _k\overline{\DD}, \vspace{2mm} \\
\sum_{n = 0}^\infty \, _kQ^{(*n)}_{\e, ii}(t) * \, _kQ_{\e, ij}(t), & \text{for} \  t \geq 0, j \neq i, i \in \, _k\overline{\DD}, 
\end{array}
\right.   
\end{equation}

Respectively, the  transition probabilities  for the embedded Markov chain $_k\tilde{\eta}_{\e, n} $ are given by the following relation, 
\begin{equation}\label{gopetasd}
_k\tilde{p}_{\e, ij}  = \PP \{ _k\tilde{\eta}_{\e,  1} = j / \, _k\tilde{\eta}_{\e, 0} = i \}  
=  \left\{
\begin{array}{cll}
_kp_{\e, ij} & \text{if} \ j \in \, _k\XX, i \in \, _k\DD,  \vspace{2mm} \\ 
0 & \text{if} \ j = i, i \in \, _k\overline{\DD}, \vspace{2mm} \\ 
 \frac{_kp_{\e, ij}}{1 - \, _kp_{\e, ii}} & \text{if} \ j \neq i, i \in \, _k\overline{\DD}. 
\end{array}
\right.
\end{equation}

Note that condition  ${\bf B}$ implies that probabilities $_kp_{\e, ii} < 1, i \in \overline{\DD}$, for every $\e \in (0, 1]$. 
\vspace{1mm}

{\bf 5.2 Conditions ${\bf A}$, ${\bf B}$ and ${\bf \hat{B}}$}.  Let conditions ${\bf A}$ and ${\bf B}$ hold. In this case, conditions ${\bf \tilde{A}}$, ${\bf \tilde{B}}$,  $_k{\bf A}$,  and $_ k{\bf B}$ also hold. Thus, by remarks made in 
Subsections 3.2 and 3.3, which should be applied to the semi-Markov processes $_k\eta_\e(t)$ and $_k\tilde{\eta}_\e(t)$ instead of the semi-Markov processes $\eta_\e(t)$ and $\tilde{\eta}_\e(t)$, conditions ${\bf A}$ and ${\bf B}$
also holds for the semi-Markov processes $_k\tilde{\eta}_\e(t)$.

 Therefore,  process   $_k\tilde{\eta}_{\e}(t)$ is well defined on the interval $[0, \infty)$, for every $\e \in (0, 1]$. 
 
 In what follows, we can also denote by $_k{\bf \tilde{A}}$ and $_k{\bf \tilde{B}}$, respectively,  conditions ${\bf A}$ and ${\bf B}$ for the semi-Markov processes $_k\tilde{\eta}_{\e}(t)$. 
 
 It is also useful to note that condition ${\bf \hat{B}}$ assumed to hold for the semi-Markov processes $\eta_\e(t)$, and, thus,  for the semi-Markov processes $\tilde{\eta}_\e(t)$ and $_k\eta_\e(t)$, also holds for the semi-Markov processes  $_k\tilde{\eta}_{\e}(t)$.  It follows from relation (\ref{transakannop}) according to which, the transition probabilities $_kQ_{\e, ij}(t), t \geq 0$   and  
$_k\tilde{Q}_{\e, ij}(t), t \geq 0$  coincide for $i \in \DD, j \in \, _k\XX$. \vspace{1mm}

{\bf 5.3 Conditions ${\bf C}$ and ${\bf \tilde{C}}$}.  Let us assume that conditions ${\bf A}$ -- 
${\bf C}$,  and  ${\bf \tilde{C}}$,   $_k{\bf \tilde{C}}$ hold. In this case condition ${\bf \tilde{A}}$, ${\bf \tilde{B}}$, $_k{\bf A}$, $_k{\bf B}$,  and $_k{\bf C}$ hold.

Condition  $_k{\bf \tilde{C}}$ plays the role of condition ${\bf \tilde{C}}$, for the semi-Markov processes $_k\tilde{\eta}_{\e}(t)$, and, also,  the role of condition ${\bf C}$ for the semi-Markov processes $_k\tilde{\eta}_{\e}(t)$.

Since, probabilities $_k\tilde{p}_{\e, ii} = 0$, for  $i \in \, _k\overline{\DD}$ and $\e \in [0, 1]$, condition $_k{\bf \tilde{C}}$ is, in fact, equivalent to the following condition:
\begin{itemize}
\item [$_k{\bf \tilde{C}}'$:]   $\frac{_k\tilde{p}_{\e, ij}}{1 - \, _k\tilde{p}_{\e, ii}} = \, _k\tilde{p}_{\e, ij}  \to \, \frac{_k\tilde{p}_{0, ij}}{1 - \, _k\tilde{p}_{0, ii}}  
= \, _k \tilde{p}_{0, ij}$  as $\e \to 0$, for $j \neq i, j \in \, _k\XX, i \in \, _k\overline{\DD}$.
\end{itemize}

Condition  $_k{\bf \tilde{C}}'$ plays the role of condition ${\bf \tilde{C}}$ for the semi-Markov processes $_k\tilde{\eta}_\e(t)$.

Let us  also define sets,  $_k\tilde{\YY}_{\e,  i}  = \{ j \in \XX: \, _k\tilde{p}_{\e,  ij} > 0 \}$,   for $i \in \, _k\overline{\DD}$ 
and  $\e \in [0, 1]$. Condition ${\bf B}$ implies that conditions ${\bf \tilde{B}}$, $_k{\bf B}$ and  $_k{\bf \tilde{B}}$  holds and, 
thus,  $_k\tilde{\YY}_{\e, i}  =  \, _k\tilde{\YY}_{1, i}$,  for $i \in \, _k\overline{\DD}$ and $\e \in (0, 1]$.
 Also, conditions ${\bf B}$ and  ${\bf C}$ imply that conditions  ${\bf \tilde{B}}$, ${\bf \tilde{C}}$,  $_k{\bf B}$, 
 $_k{\bf C}$, $_k{\bf \tilde{B}}$ and $_k{\bf \tilde{C}}$ hold, and, thus, $_k\tilde{\YY}_{0, i} \subseteq \, _k\tilde{\YY}_{1, i}$
for $i \in \, _k\overline{\DD}$. The sets $_k\YY_{1, i}$ and $_k\tilde{\YY}_{1, i}$ are connected by  relation, 
$_k\tilde{\YY}_{1, i}   = \, _k\YY_{1, i} \setminus \{ i \}$, for $i \in \, _k\overline{\DD}$.

Note also that, under condition ${\bf \hat{B}}$, sets $_k\tilde{\YY}_{\e, i} = \DD, \e \in [0, 1]$,  for $i \in \, _k\DD$. 
\vspace{1mm} 

{\bf 5.4 Conditions ${\bf D}$ and ${\bf E}$}. Let us assume that conditions ${\bf A}$ -- ${\bf C}$, ${\bf D}$ (or, equivalently ${\bf D}'$), ${\bf E}$ and  ${\bf \tilde{C}}$ (and, thus, ${\bf \tilde{C}}'$), $_k{\bf \tilde{C}}$, $_k{\bf \hat{C}}$, $_k{\bf F}$  hold. In this case, conditions ${\bf \tilde{A}}$, ${\bf \tilde{B}}$, ${\bf \tilde{D}}$ (or, equivalently ${\bf \tilde{D}}'$), ${\bf \tilde{E}}$ hold, and, moreover, conditions $_k{\bf A}$ -- $_k{\bf C}$, $_k{\bf D}$ (or, equivalently $_k{\bf D}'$), $_k{\bf E}$
and $_k{\bf \tilde{A}}$ -- $_k{\bf \tilde{C}}$, $_k{\bf \tilde{C}}'$ also hold.

Therefore, conditions ${\bf D}$ (or, equivalently ${\bf D}'$) and  ${\bf E}$ also hold for  the semi-Markov processes  $_k\tilde{\eta}_\e(t)$. The corresponding proof repeats the proof given in Subsections 4.9 and 4.10, where the semi-Markov processes $\eta_{\e}(t)$ and  $\tilde{\eta}_\e(t)$ should be replaced, respectively,  by the semi-Markov processes $_k\eta_{\e}(t)$ and $_k\tilde{\eta}_{\e}(t)$.
\vspace{1mm}

The local normalising functions $\tilde{v}_{\e, i} = (1 - p_{\e, ii})^{-1} v_{\e, i}, i \in \overline{\DD}$ for the semi-Markov processes  $\tilde{\eta}_{\e}(t)$, 
are defined by relation (\ref{compreg}). The analogue of this relation or the semi-Markov processes  $_k\tilde{\eta}_{\e}(t)$ defining the corresponding local normalising functions takes the following form, for $i \in \, _k\overline{\DD}$, 
\begin{equation}\label{compregase}
 _k\tilde{v}_{\e,  i} = (1 - \, _kp_{\e, ii})^{-1} \, _kv_{\e,  i} = (1 - \, _kp_{\e, ii})^{-1}(1 - p_{\e, ii})^{-1} v_{\e, i}. 
 \end{equation} 
 
 Let $_k\ZZ_{0} = \{i \in \, _k\XX: \, _kp_{0, ii} = 1 \}$ be the set of asymptotically absorbing  states $i \in \, _k\overline{\DD}$ for the semi-Markov processes $_k\eta_{\e}(t)$.

Analogues of relations (\ref{transakai}) and (\ref{trawet}) take the following form, for 
$j \in  \, _k\tilde{\YY}_{1, i} , i \in \, _k\overline{\DD}$,
\begin{align}\label{transakaino}
_{k}\tilde{F}_{\e,   ij}(t) =  \frac{1}{_{k}\tilde{p}_{\e,  ij}} \sum_{n = 0}^\infty \, _{k}F^{(*n)}_{\e, ii}(t) * \, _kF_{\e, ij}(t) 
\, _kp_{\e, ii}^n \, _kp_{\e, ij},  \  t \geq 0,
\end{align}
and 
\begin{align}\label{trawetno}
_k\tilde{\phi}_{\e, ij}(s) &  = \frac{1}{_k\tilde{p}_{\e,  ij}} \cdot   \frac{_k\phi_{\e, ij}(s)  \, _kp_{\e, ij}}{1 - \, _k\phi_{\e, ii}(s)\, _kp_{\e, ii}} 
=  \frac{_k\phi_{\e, ij}(s)(1 - \, _kp_{\e, ii})}{1 - \, _k\phi_{\e, ii}(s) \, _kp_{\e, ii}} \vspace{2mm}\nonumber \\
& =    \frac{_k\phi_{\e, ij}(s)}{1 + \, _kp_{\e, ii} (1 - \, _kp_{\e, ii})^{-1}(1 - \, _k\phi_{\e, ii}(s))}, \  s \geq 0.  
\end{align}

Also,  for $j \notin \, _k\tilde{\YY}_{1,  i} , i \in  \, _k\overline{\DD}$, 
\begin{equation}\label{trawetnomo}
_k\tilde{\phi}_{\e,  ij}(s)   = \, _k\tilde{\phi}_{\e,  i}(s), s \geq 0,  
\end{equation}
where
\begin{align}\label{taswetopmo}
_k\tilde{\phi}_{\e,  i}(s)  = \sum_{j \in \, _k\tilde{\YY}_{1,   i}} \, _k\tilde{\phi}_{\e, ij}(s) \, _k\tilde{p}_{\e,  ij}, \ s \geq 0.
\end{align}

Analogues of asymptotic relations (\ref{trawetba}) and (\ref{trawetbamnop}), which play the role of asymptotic relations appearing in condition ${\bf D}'$ for the semi-Markov processes $\tilde{\eta}_\e(t)$),   take the following forms for the semi-Markov processes $_k\tilde{\eta}_{\e}(t)$.

If $j \in \, _k\tilde{\YY}_{1,  i}, i \in \, _k\overline{\DD} \setminus \, _k\ZZ_{0}$, then,
\begin{equation*}
_k\tilde{\phi}_{\e,  ij}(s /  \, _k\tilde{v}_{\e, i} ) \makebox[90mm]{}
\end{equation*} 
\begin{equation*}
= \left\{
\begin{array}{lll}
\frac{_k\phi_{\e, ij}( (1 - \, _kp_{\e, ii}) s / \, _kv_{\e,  i})(1 - \, _kp_{\e, ii})}{1 - \, _kp_{\e, ii} \,  
_k\phi_{\e, ii}( (1 - \, _kp_{\e, ii})  s / \, _kv_{\e, i})} & \  \text{if} \ i \in \, _k\YY_{1, i}  
\vspace{2mm} \\
_k\phi_{\e, ij}(s / \, _kv_{\e,  i}) & \  \text{if} \ i \in \, _k\XX \setminus \, _k\YY_{1, i}
\end{array}
\right.
\end{equation*}
\begin{equation}\label{trawetbanop}
\to  \, _k\tilde{\phi}_{0,  ij}(s) \ {\rm as} \ \e \to 0, \ {\rm for} \ s \geq 0, \makebox[36mm]{}
\end{equation}
where
\begin{equation}\label{adiinasd}
_k\tilde{\phi}_{0,  ij}(s) = \left\{
\begin{array}{lll}
\frac{_k\phi_{0, ij}((1 - \, _kp_{0, ii})s)(1 - \, _kp_{0, ii})}{1 - \, _kp_{0, ii} \, 
_k\phi_{0, ii}((1 - \, _kp_{0, ii}) s )} & \ \text{if} \ i \in \, _k\YY_{1, i},  
\vspace{2mm} \\
_k\phi_{0, ij}(s) & \ \text{if} \ i \notin \,  _k\XX \setminus \, _k\YY_{1, i}. 
\end{array}
\right.
\end{equation}

If $j \in \, _k\tilde{\YY}_{1,  i}, i \in _k\overline{\DD} \cap  \, _k\ZZ_{0}$,   then,
\begin{align}\label{trawetbamnop}
& \quad _k\tilde{\phi}_{\e, ij}(s /  \, _k\tilde{v}_{\e, i} ) & \vspace{2mm} \nonumber \\ 
& \quad \quad \quad =   \frac{_k\phi_{\e, ij}(  s / (1 - \, _kp_{\e, ii})^{-1} \, _kv_{\e,  i})}{1 
+ \, _kp_{\e, ii} (1 - \, _kp_{\e, ii})^{-1} (1 - \, _k\phi_{\e, ii}(  s / (1 - \, _kp_{\e, ii})^{-1}  \, _kv_{\e, i}))} \vspace{2mm} \nonumber \\
& \quad \quad  \quad \to   \frac{1}{1 +  \, _ke_{0, ii} s} = \, _k\tilde{\phi}_{0,  ij}(s) \ {\rm as} \ \e \to 0, \ {\rm for} \ s \geq 0.   
\end{align}

Analogues of asymptotic relation (\ref{tatop}) takes the following form, for $j \in \, _k\XX \setminus _k\tilde{\YY}_{1,  i}, i \in \, _k\overline{\DD}$,
\begin{align}\label{tatopjut}
_k\tilde{\phi}_{\e,  i}(s)  & = \sum_{j \in \, _k\tilde{\YY}_{1,   i}} \, 
_k\tilde{\phi}_{\e, ij}(s) \, _k\tilde{p}_{\e,  ij} \vspace{2mm} \nonumber \\ 
& \to \sum_{j \in \, _k\tilde{\YY}_{1,   i}} \, _k\tilde{\phi}_{0, ij}(s) \, _k\tilde{p}_{0,  ij} 
= \sum_{j \in \, _k\tilde{\YY}_{0,   i}} \, _k\tilde{\phi}_{0, ij}(s) \, _k\tilde{p}_{0,  ij} 
\vspace{2mm} \nonumber \\ 
& = \,_k\tilde{\phi}_{0,  i}(s) \ {\rm as} \ \e \to 0, {\rm for} \ s \geq 0.
\end{align}   

The corresponding proofs repeat the proofs given in Subsections 3.5, where the semi-Markov processes $\eta_{\e}(t)$ and  $\tilde{\eta}_\e(t)$ should be replaced, respectively,  by the semi-Markov processes $_k\eta_{\e}(t)$ and $_k\tilde{\eta}_{\e}(t)$.
The role of asymptotic relations  appearing in condition ${\bf D}'$, for the semi-Markov processes $\eta_\e(t)$, is played for the semi-Markov processes  $_k\eta_{\e}(t)$, by the asymptotic relations (\ref{trawetbanop}) and (\ref{trawetbamnop}). 

Thus, relations (\ref{trawetbanop}) and  (\ref{trawetbamnop}) imply that condition ${\bf D}'$  and, thus, also condition ${\bf D}$,   holds for  semi-Markov 
processes $_k\tilde{\eta}_{\e}(t)$, with the Laplace transforms of the corresponding limiting distribution functions given in the above relations, and the normalisation functions  $_k\tilde{v}_{\e, i}, i \in \, _k\overline{\DD}$. 
\vspace{1mm}

An analogue of relation (\ref{expofava}) take the following form, for $j \in \, _k\tilde{\YY}_{1,  i}, i \in \, _k\overline{\DD}$, 
\begin{equation*}
\frac{_k\tilde{e}_{\e, ij}}{_k\tilde{v}_{\e,  i}}  = \left\{
\begin{array}{lll}
(1 - \, _kp_{\e, ii}) \frac{_ke_{\e, ij}}{_kv_{\e,  i}} 
+ \, _kp_{\e, ii}  \frac{_ke_{\e, ii}}{_kv_{\e,  i}} & \ \text{if} \  i \in \, _k\YY_{1, i} 
\makebox[10mm]{} \vspace{2mm} \\
\frac{_ke_{\e, ij}}{_kv_{\e,  i}} & \ \text{if} \  i \in \XX_k \setminus \,  _k\YY_{1, i} 
\end{array}
\right.
\end{equation*}
\begin{equation}\label{expofavabeq}
\to  \, _k\tilde{e}_{0,   ij} = \int_0^\infty t \, _k\tilde{F}_{0, ij}(dt) \ {\rm as} \ \e \to 0, \makebox[22mm]{}
\end{equation}
where
\begin{equation}\label{expofavabeqa}
_k\tilde{e}_{0,   ij}  = \left\{
\begin{array}{lll}
(1 -  \, _kp_{0, ii}) \, _ke_{0, ij} + \, _kp_{0, ii} \, _k e_{0, ii} & \ \text{if} \  i \in \, _k\YY_{1, i}, \\
_ke_{0, ij}& \ \text{if} \  i \in \XX_k \setminus \, _k\YY_{1, i}. 
\end{array}
\right.
\end{equation}

Also, the following relation takes place, for $j \in \, _k\XX \setminus \, _k\tilde{\YY}_{1,  i}, i \in \, _k\overline{\DD}$, 
\begin{align}\label{exfavanu}
\frac{_k\tilde{e}_{\e, i}}{_k\tilde{v}_{\e, i}}  &  = \sum_{j \in \, _k\tilde{\YY}_{1,  i}}  
\frac{_k\tilde{e}_{\e, ij}}{_kv_{\e,  i}} \, _k\tilde{p}_{\e, ij}
 \vspace{2mm} \nonumber \\
&  \to \sum_{j \in \, _k\tilde{\YY}_{1,  i}} \, _k\tilde{e}_{0, ij} \, _k\tilde{p}_{0, ij}  
\vspace{2mm} \nonumber \\
& = \sum_{j \in \, _k\tilde{\YY}_{0,  i}} \, _k\tilde{e}_{0, ij} \, _k\tilde{p}_{0, ij}
 \vspace{2mm} \nonumber \\
& = \, _k\tilde{e}_{0, i} = \int_0^\infty t \, _k\tilde{F}_{0, i}(dt) \ {\rm as} \ \e \to 0.
\end{align}

Relation  (\ref{expofavabeq}) plays the role of asymptotic relations appearing in condition ${\bf E}$ for the semi-Markov processes $_k\tilde{\eta}_{\e}(t)$. The corresponding proof repeats the proof given in Subsections 3.6, where the semi-Markov processes $\eta_{\e}(t)$ and  
$\tilde{\eta}_\e(t)$ should be replaced, respectively,  by the semi-Markov processes $_k\eta_{\e}(t)$ and $_k\tilde{\eta}_{\e}(t)$.
The role of asymptotic relations  appearing in condition ${\bf E}$, for the semi-Markov processes $\eta_\e(t)$, is played for the semi-Markov processes  $_k\eta_{\e}(t)$ by the asymptotic relation (\ref{trawetkasnop}). 

Thus, relation (\ref{expofavabeq}) implies that condition ${\bf E}$    holds for  semi-Markov 
processes $_k\tilde{\eta}_{\e}(t)$, with the corresponding limiting expectations given in the above relation, and the normalisation functions  $_k\tilde{v}_{\e, i}, i \in \, _k\overline{\DD}$. 

In what follows, we can also denote by $_k{\bf \tilde{D}}$, $_k{\bf \tilde{D}}'$,  and  $_k{\bf \tilde{E}}$, respectively, conditions ${\bf D}$,  ${\bf D}'$, and  ${\bf E}$  for the semi-Markov processes $_k\eta_{\e}(t)$ 
(expressed, for condition $_k{\bf \tilde{D}}'$,  in the form of relations  (\ref{trawetbanop}) and  (\ref{trawetbamnop}), and, for condition $_k{\bf \tilde{E}}$, in the form of relation (\ref{expofavabeq}). 
 \vspace{1mm}

{\bf 5.5 Summary}. The following lemma summarises the above remarks.  \vspace{1mm}

{\bf Lemma 17}. {\em Let conditions ${\bf A}$ -- ${\bf E}$, and ${\bf \tilde{C}}$, $_k{\bf \tilde{C}}$, $_k{\bf \hat{C}}$, $_k{\bf \hat{F}}$  hold. Then, 
conditions ${\bf A}$ -- ${\bf E}$ and ${\bf \tilde{C}}$ hold for  the  semi-Markov processes 
$_k\tilde{\eta}_{\e}(t)$, respectively,  in the form of conditions  $_k{\bf \tilde{A}}$ -- $_k{\bf \tilde{E}}$ and $_k{\bf \tilde{C}}'$.}
\vspace{1mm}

{\bf Remark 9}. Conditions $_k{\bf \tilde{C}}$ and $_k{\bf \tilde{C}}'$, which is equivalent to condition $_k{\bf \tilde{C}}$,   play, respectively,  the roles of conditions
${\bf C}$ and $_k{\bf \tilde{C}}$  for the semi-Markov processes  $_k\tilde{\eta}_{\e}(t)$.
\vspace{1mm}

{\bf Remark 10}. Condition ${\bf C}_{1}$ is sufficient for holding onditions ${\bf C}$, ${\bf \tilde{C}}$, ${\bf \tilde{C}}'$, $_k{\bf \tilde{C}}$,  $_k{\bf \tilde{C}}'$,  and $_k{\bf \hat{C}}$. \vspace{1mm}

{\bf 5.6 Hitting times for reduced semi-Markov processes with removed virtual transitions}. 
Let us introduce hitting times for semi-Markov 
processes $_k\tilde{\eta}_{\e}(t)$, 
\begin{equation}\label{opopil}
_k\tilde{\tau}_{\e, \DD} = \sum_{n = 1}^{_k\tilde{\nu}_{\e, \DD}} \, _k\kappa_{\e,  n}, \ {\rm where} \  
_k\tilde{\nu}_{\e,  \DD} = \min(n \geq 1: \, _k\tilde{\eta}_{\e, n} \in \DD).
\end{equation}

The following three lemmas are corollaries of Lemma 4 -- 6 and 14 --16. Their proofs can be 
obtained, first, by application of  Lemmas 4 -- 6 to the semi-Markov processes $\eta_{\e}(t)$ and 
$\tilde{\eta}_{\e}(t)$, second, by application of  Lemmas 14 -- 16 to  the semi-Markov processes 
$\tilde{\eta}_{\e}(t)$ and $_k\eta_{\e}(t)$, and, third, by application of  Lemmas  4 -- 6 to the semi-Markov processes  $_k\eta_{\e}(t)$ and $_k\tilde{\eta}_{\e}(t)$.
\vspace{1mm}

{\bf Lemma 18}. {\em Let condition ${\bf A}$ and ${\bf B}$ hold, and, in sequel, conditions ${\bf \tilde{A}}$, ${\bf \tilde{B}}$, $_k{\bf A}$, $_k{\bf B}$, $_k{\bf \tilde{A}}$, $_k{\bf \tilde{B}}$  hold.
Then, the following relation takes place, for  $\e \in (0, 1]$,  
\begin{align*}
& \PP_i \{ \tau_{\e, \DD}  = \,  _k \tau_{\e, \DD} =  \, _k\tilde{\tau}_{\e, \DD}, \makebox[50mm]{}
 %\vspace{2mm} \nonumber \\
 \end{align*}
 \begin{align}\label{ideremopui}
& \quad \  \eta_{\e}(\tau_{\e, \DD}) = \, _k\eta_{\e}(_k\tau_{\e, \DD})  =  
\, _k\tilde{\eta}_{\e}(_k\tilde{\tau}_{\e, \DD})  \} = 1, \ i \in \, _k\XX.   
\end{align}}

Let us also introduce distributions,  for $i \in \, _k\XX$ and $\e \in (0, 1]$,  
\begin{equation}\label{trewqopl}
_k\tilde{G}_{\e, \DD, ij}(t) = \PP_i \{ _k \tilde{\tau}_{\e,  \DD} \leq t, \, _k\tilde{\eta}_{\e}(_k\tilde{\tau}_{\e, \DD}) = j) \}, \, t \geq 0, j \in \DD.
\end{equation}

{\bf Lemma 19}. {\em Let condition ${\bf A}$ and ${\bf B}$ hold, and, in sequel, conditions ${\bf \tilde{A}}$, ${\bf \tilde{B}}$, $_k{\bf A}$, $_k{\bf B}$,  $_k{\bf \tilde{A}}$, $_k{\bf \tilde{B}}$  hold.
  Then, the following relation takes place, for   $\e \in (0, 1]$,  
\begin{align}\label{iderenasbail}
& G_{\e, \DD, ij}(t)   = \, _kG_{\e,  \DD, ij}(t)  =  \, _k\tilde{G}_{\e,  \DD, ij}(t), \  t \geq 0,  j \in \DD, i \in \, _k\XX. 
\end{align}}
\makebox[3mm]{} Let us also introduce Laplace transforms, for $i \in \XX$ and $\e \in (0, 1]$,
\begin{equation}\label{trewqyn}
_k\tilde{\Psi}_{\e, \DD, ij}(s) = \EE_i \exp\{- s \, _k\tilde{\tau}_{\e, \DD}\} {\rm I}(_k\tilde{\eta}_{\e}(_k\tilde{\tau}_{\e, \DD}) = j) \}, \, s \geq 0, j \in \DD.
\end{equation}

{\bf Lemma 20}. {\em Let condition ${\bf A}$ and ${\bf B}$ hold, and, in sequel, conditions ${\bf \tilde{A}}$, ${\bf \tilde{B}}$, $_k{\bf A}$, $_k{\bf B}$, $_k{\bf \tilde{A}}$, $_k{\bf \tilde{B}}$  hold.
 Then, the following relation takes place, for  $\e \in (0, 1]$,
\begin{align}\label{iderenasnopio}
& \Psi_{\e, \DD, ij}(s)   =  \, _k\Psi_{\e, \DD, ij}(s) = \, _k\tilde{\Psi}_{\e, \DD, ij}(s),  \ s \geq 0,  j \in \DD, i \in \, _k\XX. 
\end{align}}
\makebox[3mm]{}
Lemmas 4 -- 6, 14 -- 16, and 18 -- 20  let one reduce study of asymptotics for distributions of hitting times  for  semi-Markov processes $\eta_\e(t)$, $\tilde{\eta}_{\e}(t)$, and $_k\eta_{\e}(t)$  to the case of more simple semi-Markov processes $_k\tilde{\eta}_{\e}(t)$.  

The distributions $_k\tilde{G}_{\e,  \DD, ij}(t), t \geq 0, j \in \DD, i \in \overline{\DD}$ are determined by transition probabilities $_k\tilde{Q}_{\e, kr}(t), t \geq 0,   r \in \, _k\XX, k \in \, _k\overline{\DD}$, which themselves
are determined by transition probabilities $_kQ_{\e, kr}(t), t \geq 0,   r \in \XX, k \in \, _k\overline{\DD}$. The latter transition probabilities  are determined by transition probabilities $\tilde{Q}_{\e, \DD, lr}(t), t \geq 0,   r \in \XX, l \in \, 
\overline{\DD}$, which themselves are determined by transition probabilities $Q_{\e, lr}(t), t \geq 0,   r \in \XX, l \in \, \overline{\DD}$.

This makes it possible to essentially simplify the model in the case, where we are interested to investigate asymptotics of hitting times for some fixed domain $\DD$ and only for initial states $i \in \overline{\DD}$. This case in considered in Sections 2 -- 7. For example, we can assume in what follows in these sections that the transition probabilities  of semi-Markov processes $\eta_\e(t)$ and, in sequel, the transition probabilities of semi-Markov processes  $\tilde{\eta}_\e(t)$, $_k\eta_\e(t)$, $_k\tilde{\eta}_\e(t)$ and satisfy condition 
${\bf \hat{B}}$. \vspace{1mm}

{\bf 5.7 Weak asymptotics of hitting times for the case of two-states domain $\overline{\DD}$}.
The above remarks,  Lemmas 5, 6, 15, 16, 18, and 19 and Theorem 1 let us describe asymptotics for distributions of hitting times for another simplest case, where domain $\overline{\DD} = \{ i, k \}$  is a two-states set.  

We  assume that condition $_k{\bf \hat{F}}$ holds,  i.e., state $k$ is less or equally absorbing with state $i$. According Lemma 18, the following relation takes place, 
\begin{equation}\label{weasd}
G_{\e, \DD, ij}(t) = \, _kG_{\e, \DD, ij}(t) = \, _k\tilde{G}_{\e, ij}(t) ,  \ {\rm for} \ t \geq 0, j \in \DD. 
\end{equation}

It is obvious that, in this case,  domain  $_k\overline{\DD} = \{ i \}$ is a 
one-state set and, thus, the above relation  let us apply Theorem 1 to the semi-Markov processes 
$_k\tilde{\eta}_\e(t)$ instead  of the semi-Markov processes $\eta_\e(t)$.

Relation (\ref{cotyre}) takes the following form, 
\begin{equation}\label{weasderas}
_k\tilde{G}_{\e, \DD, ij}(t) = \, _k\tilde{F}_{\e, ij}(\cdot) \, _k\tilde{p}_{\e, ij},  \ {\rm for} \ t \geq 0, j \in \DD. 
\end{equation} 

Henceforth, we always use notation $_k\check{v}_{\e, i}$ for the final normalising functions.
In the above case, 
\begin{equation}\label{giklo}
_k\check{v}_{\e, i} = \, _k\tilde{v}_{\e, i}.
\end{equation} 
 
The  following theorem takes place. \vspace{1mm}

{\bf Theorem 2}. {\em Let domain $\overline{\DD} = \{ i, k \}$  be a two-states set and conditions ${\bf A}$  -- ${\bf E}$, and ${\bf \tilde{C}}$, $_k{\bf \tilde{C}}$, $_k{\bf \hat{C}}$, $_k{\bf \hat{F}}$   hold.  
Then,  the following asymptotic relation takes place, for $j \in \DD$, 
\begin{equation}\label{cotyrevanio}
G_{\e, \DD, ij}(\cdot \, _k\check{v}_{\e, i}) \Rightarrow  \, G_{0, \DD, ij}(\cdot) = \, _k\tilde{F}_{0, ij}(\cdot) \, _k\tilde{p}_{0, ij} \ {\rm as} \ \e \to 0.  
\end{equation}}
\makebox[4mm]{}{\bf Remark 11}. The distribution functions $_k\tilde{F}_{0, ij}(\cdot), j \in \DD$ are not concentrated in $0$, i.e., $_k\tilde{F}_{0, ij}(0) < 1$, for $j \in \DD$. \vspace{1mm}

The question about asymptotics of distributions $G_{\e, \DD, kj}(t)$ is more complex.

We can use, in this case, relation (\ref{iderenasfat}) given in Lemma 16. This relation takes, in this case, the following form, for $s \geq 0, j \in \DD$, 
\begin{equation}\label{idersfat}
\Psi_{\e, \DD, kj}(s) =  \tilde{\phi}_{\e, kj}(s)\tilde{p}_{\e, kj} + \Psi_{\e, \DD, ij}(s)  \tilde{\phi}_{\e, ki}(s) \tilde{p}_{\e,  ki}. 
\end{equation}

If the limiting probability $\tilde{p}_{0, ki}  > 0$, relation (\ref{idersfat}) and relation (\ref{cotyreva}) given in Theorem 2 imply that the following relation 
takes place,  for $s \geq 0, j \in \DD$, 
\begin{align*}
\Psi_{\e, \DD, kj}(s / \, _k\tilde{v}_{\e, i}) & =  \tilde{\phi}_{\e, kj}(((1 - \, _kp_{\e, ii})  
\frac{\tilde{v}_{\e, k}}{\tilde{v}_{\e, i}}s/\tilde{v}_{\e, k} )\tilde{p}_{\e, kj}  \makebox[30mm]{}
%\vspace{2mm} \nonumber \\
\end{align*}
\begin{align}\label{idersfatnut}
& \quad \quad \quad + \Psi_{\e, \DD, ij}(s / \, _k\tilde{v}_{\e, i}) 
\tilde{\phi}_{\e, ki}((1 - \, _kp_{\e, ii})  \frac{\tilde{v}_{\e, k}}{\tilde{v}_{\e, i}}s/\tilde{v}_{\e, k} )  \tilde{p}_{\e,  ki} 
\vspace{2mm} \nonumber \\
& \quad \quad \to \tilde{\phi}_{0, kj}((1 - \, _kp_{0, ii})  \tilde{w}_{0, ki}s))\tilde{p}_{0, kj} 
\vspace{2mm} \nonumber \\
& \quad \quad \quad + \Psi_{0, \DD, ij}(s)  \tilde{\phi}_{0, ki}((1 - \, _kp_{0, ii})  \tilde{w}_{0, ki}s)  \tilde{p}_{0,  ki}
\vspace{2mm} \nonumber \\ 
& \quad \quad =  \Psi_{0, \DD, kj}(s) \ {\rm as} \ \e \to 0.
\end{align}

If the limiting probability $\tilde{p}_{0, ki}  = 0$,  relation (\ref{idersfat}) and relation (\ref{cotyreva}) given in Theorem 2 imply that the following relation 
takes place,  for $s \geq 0, j \in \DD$, 
\begin{align}\label{idersfatas}
\Psi_{\e, \DD, kj}(s / \tilde{v}_{\e, k}) & =  \tilde{\phi}_{\e, kj}(s/ \tilde{v}_{\e, k})\tilde{p}_{\e, kj} 
\vspace{2mm} \nonumber \\
& \quad + \tilde{\phi}_{\e, ki}( s/\tilde{v}_{\e, k} ) \Psi_{\e, \DD, ij}(s / \tilde{v}_{\e, k}) \tilde{p}_{\e,  ki} \vspace{2mm} \nonumber \\
& \to \Psi_{0, \DD, kj}(s) = \tilde{\phi}_{0, kj}(s)\tilde{p}_{0, kj}  \ {\rm as} \ \e \to 0.
\end{align}

Thus, the corresponding  limiting distribution $G_{0, \DD, kj}(\cdot)$ has the Laplace transform 
$\Psi_{0, \DD, kj}(s) = \int_0^\infty e^{-st} G_{0, \DD, kj}(dt), s \geq 0$ given by relation (\ref{idersfatnut}) and the normalising function is  
$_k\check{v}_{\e, i} = \, _k\tilde{v}_{\e, i}$, if the limiting probability $\tilde{p}_{0, ki}  > 0$. 

However, the corresponding  limiting distribution $G_{0, \DD, kj}(\cdot)$ has the Laplace transform $\Psi_{0, \DD, kj}(\cdot)$ given by 
relation (\ref{idersfatas}) and the normalising function is  $_k\check{v}_{\e, k} = \tilde{v}_{\e, k}$, if the limiting probability $\tilde{p}_{0, ki}  = 0$. 
 
The  following theorem takes place. \vspace{1mm}

{\bf Theorem 3}. {\em Let domain $\overline{\DD} = \{ i, k \}$  be a two-states set and conditions ${\bf A}$ -- ${\bf E}$ and ${\bf \tilde{C}}$, $_k{\bf \tilde{C}}$, $_k{\bf \hat{C}}$, $_k{\bf \hat{F}}$.  
Then,  the following asymptotic relation takes place, for $j \in \DD$, 
\begin{equation}\label{cotyrevaba}
G_{\e, \DD, kj}(\cdot \, _k\check{v}_{\e, k}) \Rightarrow  \, G_{0, \DD, kj}(\cdot) \ {\rm as} \ \e \to 0.
\end{equation}}
\makebox[3mm]{}It is worth to note that in the case, where $_kp_{0, ii} = 1$, the normalising function 
$\tilde{v}_{\e, i} = o(_k\tilde{v}_{\e, i})$ as $\e \to 0$ and, thus,  due to condition $_k{\bf F}$, also,  the normalising function $\tilde{v}_{\e, k} = o(_k\tilde{v}_{\e, i})$ as $\e \to 0$. Thus, in the case, where $_kp_{0, ii} = 1, \tilde{p}_{0, ki}  = 0$,  essentially different normalising functions should be used for distributions $G_{\e, \DD, ij}(\cdot)$ and $G_{\e, \DD, kj}(\cdot)$ (namely, $_k\tilde{v}_{\e, i}$ and  $\tilde{v}_{\e, k}$) in the corresponding weak convergence relations given, respectively, in Theorems 2 and 3. Note that it is possible that the pre-limiting probability $\tilde{p}_{\e, ik}, \tilde{p}_{\e, ki}  > 0$, i.e.,  domain $\overline{\DD}$ is a communicative class of states.

\vspace{1mm}

{\bf Remark 12}. Conditions  ${\bf C}$, ${\bf \tilde{C}}$,   $_k{\bf \tilde{C}}$, and $_k{\bf \hat{C}}$  can be replaced by condition  ${\bf C}_{1}$ in Theorems 2 and 3. \\

{\bf 6.  Algorithm of Multi-Step Reduction of Phase Space  for \\ \makebox[11mm]{}  
Perturbed Semi-Markov Processes} \\ 

This section plays the key role in the paper. Here, we describe  an asymptotic  multi-step phase space reduction procedure for perturbed semi-Markov processes.  What is important that hitting times are asymptotically invariant with respect to this procedure. We also formulate conditions which guarantee that basic perturbation conditions imposed on the initial semi-Markov processes also holds for the semi-Markov processes with reduced phase space resulted by the above multi-step phase space reduction procedure.  Also, we describe a recurrent algorithm for re-calculating  normalisation functions, limiting distributions and expectations  in the corresponding perturbation conditions for  the semi-Markov processes with reduced phase space. The above multi-step phase space reduction procedure let us give a detailed description of weak convergence asymptotics for hitting times for perturbed semi-Markov processes.  \vspace{1mm}

{\bf 6.1 Sequential excluding of  states from the phase space $\XX$}.  An important element of the phase space reduction algorithm presented in the paper is the procedure of sequential excluding some sequence of states $\bar{k}_{\bar{m}-1} = \langle k_1, k_2, \ldots, k_{\bar{m} -1} \rangle$  from domain $\overline{\DD}$ resulting reduction of this domain to a one-state set. We denote by $_{\bar{k}_n}\eta_\e(t)$ the  semi-Markov  process resulted by sequential exclusion from the phase space $\XX$ of semi-Markov processes  $\eta_\e(t)$ the subsequence of states $\bar{k}_n = \langle k_1, \ldots, k_n \rangle$, 
and denote by  $_{\bar{k}_n}\tilde{\eta}_\e(t)$ the  semi-Markov  process resulted by exclusion virtual transitions from trajectories of the 
semi-Markov process $_{\bar{k}_n}\eta_\e(t)$. These semi-Markov processes have  the phase space $_{\bar{k}_n}\XX = \XX \setminus \{ k_1, k_2, \ldots, k_n \}$, while set $_k\overline{\DD} = \overline{\DD} \setminus \{ k_1, k_2, \ldots, k_n \}$ plays  the role of domain $\overline{\DD}$ for these processes. 

In particular, we change notation from $_{k_1}\eta_\e(t)$ and $_{k_1}\tilde{\eta}_\e(t)$, respectively, to $_{\bar{k}_1}\eta_\e(t)$
and $_{\bar{k}_1}\tilde{\eta}_\e(t)$, for the semi-Markov process $_{k_1}\eta_\e(t)$, resulted by the procedure of exclusion from the initial phase space $\XX$  state $k = k_1$ and  the semi-Markov process $_{k_1}\tilde{\eta}_\e(t)$ resulted by removing virtual transitions from trajectories of the semi-Markov process $_{k_1}\eta_\e(t)$.  

In the similar way, we index random variables, in particular hitting times, phase spaces and other related sets, conditions,  probabilities, expectations and other quantities and objects related to the semi-Markov processes $_{\bar{k}_n}\eta_\e(t)$ and $_{\bar{k}_n}\tilde{\eta}_\e(t)$ by the left  lower index $_{\bar{k}_n}$ (in the way such as $_{\bar{k}_n}f$), in order to distingue these random variables, phase spaces and other sets, probabilities, expectations and other quantities and objects,   for $n = 1, \ldots, \bar{m} -1$. 

In order to get the similar forms for  recurrent relations connecting  the above random variables, probabilities and expectations and other quantities and objects,  resulted by  exclusion from the phase space of the subsequences of states $\bar{k}_{n -1}$ and $\bar{k}_{n}$  we use also the left lower index $_{\bar{k}_0}$ (referring to the  ``empty'' subsequence  $\bar{k}_0 = \langle  \rangle$), for $n = 1$.  Thus, we change notations for the initial semi-Markov process $\eta_\e(t)$ and the semi-Markov process 
$\tilde{\eta}_\e(t)$, respectively, to   $_{\bar{k}_0}\eta_\e(t)$ and $_{\bar{k}_0}\tilde{\eta}_\e(t)$ and make the similar indexation for the phase space $_{\bar{k}_0}\XX = \XX$ domain $_{\bar{k}_0}\overline{\DD} = \overline{\DD} $ and other related sets, conditions,  probabilities, expectations and other quantities and objects related to the semi-Markov processes $_{\bar{k}_0}\eta_\e(t) = \eta_\e(t)$ and 
$_{\bar{k}_0}\tilde{\eta}_\e(t) = \tilde{\eta}_\e(t)$.

Let describe in more details the above recurrent algorithm. \vspace{1mm}

{\bf 6.2 Step ${\bf 0}$}.  At the initial step ${\bf 0}$, the semi-Markov process $_{\bar{k}_0}\tilde{\eta}_\e(t) = \tilde{\eta}_\e(t)$ is constructed with the use of the procedure of removing virtual transitions from trajectories of the semi-Markov process $_{\bar{k}_0}\eta_\e(t) = \eta_\e(t)$. This procedure is described in Section 3. 

Conditions $_{\bar{k}_0}{\bf A} = {\bf A}$ -- $_{\bar{k}_0}{\bf E} = {\bf E}$, and  $_{\bar{k}_0}{\bf \tilde{C}} = {\bf \tilde{C}}$ should be assumed to hold for the semi-Markov processes $_{\bar{k}_0}\eta_\e(t) = \eta_\e(t)$. 

According the above remarks about notations, probabilities, $_{\bar{k}_0}p_{\e, ij} = p_{\e, ij},  \, _{\bar{k}_0}\tilde{p}_{\e, ij} = \tilde{p}_{\e, ij}$, distribution functions $_{\bar{k}_0}F_{\e, ij}(\cdot) = F_{\e, ij}(\cdot), \, _{\bar{k}_0}\tilde{F}_{\e, ij}(\cdot) = \tilde{F}_{\e, ij}(\cdot)$, 
Laplace transforms $_{\bar{k}_0}\phi_{\e, ij}(\cdot) = 
\phi_{\e, ij}(\cdot), \, _{\bar{k}_0}\tilde{\phi}_{\e, ij}(\cdot) = 
\tilde{\phi}_{\e, ij}(\cdot)$, and expectations $_{\bar{k}_0}e_{\e, ij} = e_{\e, ij}, \, 
_{\bar{k}_0}\tilde{e}_{\e, ij} = \tilde{e}_{\e, ij}$, for $i, j \in \, _{\bar{k}_0}\XX$. The roles of relations which express the 
above quantities for the semi-Markov processes $_{\bar{k}_0}\tilde{\eta}_\e(t)$ via the corresponding quantities for the 
semi-Markov processes  $_{\bar{k}_0}\eta_\e(t)$ are played by relations (\ref{gopet}), (\ref{transakai}), (\ref{trawet}), and (\ref{twet}).

According Lemma 3,  conditions $_{\bar{k}_0}{\bf \tilde{A}} = {\bf \tilde{A}}$ -- $_{\bar{k}_0}{\bf \tilde{E}} = {\bf \tilde{E}}$ and $_{\bar{k}_0}{\bf \tilde{C}}' = {\bf \tilde{C}}'$ (which is equivalent to conditions $_{\bar{k}_0}{\bf \tilde{C}}$ and plays the role of condition ${\bf \tilde{C}}$
for the  the semi-Markov processes  $_{\bar{k}_0}\tilde{\eta}_\e(t)$) hold.  

Functions $_{\bar{k}_0}v_{\e, i} = v_{\e, i}, i \in \, _{\bar{k}_0}\overline{\DD} = \overline{\DD}$ play the roles of the local normalisation functions in the above conditions $_{\bar{k}_0}{\bf D}$ and  $_{\bar{k}_0}{\bf E}$. 

Functions $_{\bar{k}_0}\tilde{v}_{\e, i} = 
\tilde{v}_{\e, i} = (1  - \, _{\bar{k}_0}p_{\e, ii}) ^{-1}v_{\e, i} = (1  - p_{\e, ii}) ^{-1}v_{\e, i}, i \in \, _{\bar{k}_0}\overline{\DD}$ play the roles of the local normalising functions in the above conditions 
$_{\bar{k}_0}{\bf \tilde{D}}$ and  $_{\bar{k}_0}{\bf \tilde{E}}$. 

The roles of the corresponding asymptotic relations appearing  in 
conditions $_{\bar{k}_0}{\bf \tilde{D}}$ and  $_{\bar{k}_0}{\bf \tilde{E}}$  are played by asymptotic relations (\ref{trawetba}), (\ref{trawetbama}), and (\ref{expofava}).   

Also, Lemmas 4 -- 6 can be reformulated in terms of random variables $_{\bar{k}_0}\tau_{\e, \DD}
 = \tau_{\e, \DD}$, $_{\bar{k}_0}\eta_{\e}(_{\bar{k}_0}\tau_{\e, \DD}) = \eta_{\e}(\tau_{\e, \DD})$, 
 $_{\bar{k}_0}\tilde{\tau}_{\e, \DD} = \tilde{\tau}_{\e, \DD}$, 
 $ _{\bar{k}_0}\tilde{\eta}_{\e}(_{\bar{k}_0}\tilde{\tau}_{\e, \DD}) = \tilde{\eta}_{\e}(\tilde{\tau}_{\e, \DD})$, distributions $_{\bar{k}_0}G_{\e, \DD, ij}(\cdot) = G_{\e, \DD, ij}(\cdot)$, $_{\bar{k}_0}\tilde{G}_{\e, \DD, ij}(\cdot)$ $= \tilde{G}_{\e, \DD, ij}(\cdot)$, and Laplace
transforms $_{\bar{k}_0}\Psi_{\e, \DD, ij}(\cdot) =  \Psi_{\e, \DD, ij}(\cdot)$, \, $_{\bar{k}_0}\tilde{\Psi}_{\e, \DD, ij}(\cdot) = \tilde{\Psi}_{\e, \DD, ij}(\cdot)$,  for $ i \in \XX, j \in \DD$.

Also, it is useful to note that, according Lemmas 1 and 2, condition $_{\bar{k}_0}{\bf C}_1 = {\bf C}_1$ implies that conditions  $_{\bar{k}_0}{\bf C}$, $_{\bar{k}_0}{\bf \tilde{C}}$, and $_{\bar{k}_0}{\bf \tilde{C}}'$ hold. \vspace{1mm}

{\bf 6.3 Step ${\bf 1}$}.   At step ${\bf 1}$, the semi-Markov process $_{\bar{k}_1}\eta_\e(t) = \, _{k_1}\eta_\e(t)$  is constructed, with the use of the procedure of the use of the procedure of exclusion of a specially chosen state $k_1 \in \, _{\bar{k}_0}\overline{\DD} = \overline{\DD}$ from the phase space $_{\bar{k}_0}\XX = \XX$ of the semi-Markov process  $_{\bar{k}_0}\tilde{\eta}_\e(t)$. Then,  
the semi-Markov process $_{\bar{k}_1}\tilde{\eta}_\e(t) = \, _{k_1}\tilde{\eta}_\e(t)$ is constructed,
with the use of  the procedure of removing virtual transitions from trajectories of the semi-Markov process 
$_{\bar{k}_1}\eta_\e(t) = \, _{k_1}\eta_\e(t)$. The above procedures are described, respectively, in  Sections 4 and 5. The semi-Markov processes  $_{\bar{k}_1}\eta_\e(t)$ and $_{\bar{k}_1}\tilde{\eta}_\e(t)$ have the phase space $_{\bar{k}_1}\XX = \, _{\bar{k}_0}\XX \ \setminus \{ k_1 \}$, while domain  $_{\bar{k}_1}\overline{\DD} = \, _{\bar{k}_0}\overline{\DD} \ \setminus \{ k_1 \}$ replaces domain  $_{\bar{k}_0}\overline{\DD}$.

At step ${\bf 1}$, conditions $_{\bar{k}_0}{\bf A}$ -- $_{\bar{k}_0}{\bf E}$, and $_{\bar{k}_0}{\bf \tilde{C}}$ are assumed to hold. Also, condition 
$_{\bar{k}_0}{\bf \tilde{C}}'$ (equivalent to condition $_{\bar{k}_0}{\bf \tilde{C}}$) holds. 

Also, condition  $_{\bar{k}_0}{\bf \tilde{F}} = {\bf \tilde{F}}$ should be additionally assumed to hold.  

Condition $_{\bar{k}_0}{\bf \tilde{F}}$ let us  choose one of the least absorbing states,    $k_1 \in \, _{\bar{k}_0}\overline{\DD}^* = \overline{\DD}^*$, with the  use of the algorithm described in Lemma 11. In fact, the choice of state $k_1 \in  \, _{\bar{k}_0}\overline{\DD}^*$ means that it is assumed  that  condition $_{\bar{k}_1}{\bf \hat{F}} = \, _{k_1}{\bf \hat{F}}$ holds. 

At  the following sub-steps ${\bf 1.1}$ and ${\bf 1.2}$, conditions $_{\bar{k}_1}{\bf \tilde{C}} = \, _{k_1}{\bf \tilde{C}}$,  and $_{\bar{k}_1}{\bf \hat{C}} = \, _{k_1}{\bf \hat{C}}$ should also be assumed to hold. Also, condition $_{\bar{k}_1}{\bf C} = \, _{k_1}{\bf C}$ holds, since it is implied by condition $_{\bar{k}_1}{\bf \tilde{C}}$.

At sub-step ${\bf 1.1}$, the semi-Markov process $_{\bar{k}_1}\eta_\e(t) = \, _{k_1}\eta_\e(t)$ is constructed with the use of the procedure of exclusion state $k_1 \in \, _{\bar{k}_0}\overline{\DD} = \overline{\DD}$ from the phase space $_{\bar{k}_0}\XX = \XX$ of the 
semi-Markov process  $_{\bar{k}_0}\tilde{\eta}_\e(t)$. This procedure is described in Section 4. 

At sub-step ${\bf 1.2}$,  the semi-Markov process $_{\bar{k}_1}\tilde{\eta}_\e(t) = \, _{k_1}\tilde{\eta}_\e(t)$ is constructed  with the use of the procedure of removing virtual transitions from trajectories of the semi-Markov process $_{\bar{k}_1}\eta_\e(t)$. This procedure is described in Section 5. 

According the above remarks about notations, probabilities, $_{\bar{k}_1}p_{\e, ij} = \, _{k_1}p_{\e, ij},  
\, _{\bar{k}_1}\tilde{p}_{\e, ij} = \, _{k_1}\tilde{p}_{\e, ij}$, distribution functions $_{\bar{k}_1}F_{\e, ij}(\cdot) =  \, _{k_1}F_{\e, ij}(\cdot), \, _{\bar{k}_1}\tilde{F}_{\e, ij}(\cdot) = \, _{k_1}\tilde{F}_{\e, ij}(\cdot)$, 
Laplace transforms $_{\bar{k}_1}\phi_{\e, ij}(\cdot) = 
\, _{k_1}\phi_{\e, ij}(\cdot), \, _{\bar{k}_1}\tilde{\phi}_{\e, ij}(\cdot) = 
 \, _{k_1}\tilde{\phi}_{\e, ij}(\cdot)$, and expectations $_{\bar{k}_1}e_{\e, ij} = \, _{k_1}e_{\e, ij}, \, 
_{\bar{k}_1}\tilde{e}_{\e, ij} = \, _{k_1}\tilde{e}_{\e, ij}$, for $i, j \in \, _{\bar{k}_1}\XX$. 
The roles of relations which express the 
above quantities for the semi-Markov processes $_{\bar{k}_1}\eta_\e(t)$ and $_{\bar{k}_1}\tilde{\eta}_\e(t)$ via the 
corresponding quantities for the semi-Markov processes  $_{\bar{k}_0}\eta_\e(t)$ and $_{\bar{k}_0}\tilde{\eta}_\e(t)$ are played by relations (\ref{gopetk}), (\ref{traalk}), (\ref{trawetk}), and (\ref{diotrask}), for the semi-Markov processes $_{\bar{k}_1}\eta_\e(t)$,  and by relations  (\ref{gopetasd}), (\ref{transakaino}), (\ref{trawetno}), and (\ref{trawetkasnop}), for the semi-Markov 
processes  $_{\bar{k}_1}\tilde{\eta}_\e(t)$. 

 According Lemmas 13 and 17, 
conditions  $_{\bar{k}_1}{\bf A}$ -- $_{\bar{k}_1}{\bf E}$  and $_{\bar{k}_1}{\bf \tilde{A}}$ -- $_{\bar{k}_1}{\bf \tilde{E}}$ and   
$_{\bar{k}_1}{\bf \tilde{C}}' = \, _{k_1}{\bf \tilde{C}}'$ (which is equivalent to condition $_{\bar{k}_1}{\bf \tilde{C}}$ and  plays the role of condition ${\bf \tilde{C}}$ for the  the semi-Markov processes  $_{\bar{k}_1}\tilde{\eta}_\e(t)$) hold, respectively,  for the semi-Markov processes  $_{\bar{k}_1}\eta_\e(t)$ and $_{\bar{k}_1}\tilde{\eta}_\e(t)$.

Functions $_{\bar{k}_1}v_{\e, i} = \, _{\bar{k}_0}\tilde{v}_{\e, i} = (1  - \, _{\bar{k}_0}p_{\e, ii}) ^{-1}v_{\e, i} = (1  - p_{\e, ii}) ^{-1}v_{\e, i},  i \in \, _{\bar{k}_1}\overline{\DD} = \, _{k_1}\overline{\DD}$ play the role of the local normalising functions 
in conditions $_{\bar{k}_1}{\bf D}$ and  $_{\bar{k}_1}{\bf E}$.

The roles of the corresponding asymptotic relations appearing  in 
conditions $_{\bar{k}_1}{\bf D}$ and  $_{\bar{k}_1}{\bf E}$   are played by asymptotic relations (\ref{trawetkasopa}), (\ref{trawetkasf}), and (\ref{expofavabeq}).

Functions $_{\bar{k}_1}\tilde{v}_{\e, i} = (1- \, _{\bar{k}_1}p_{\e, ii})^{-1} \, _{\bar{k}_1}v_{\e, i} =  (1- \, _{\bar{k}_1}p_{\e, ii})^{-1}(1 - \, _{\bar{k}_0}p_{\e, ii})^{-1} v_{\e, i} = (1- \, _{k_1}p_{\e, ii})^{-1}(1 - p_{\e, ii})^{-1} v_{\e, i},  i \in \, _{\bar{k}_1}\overline{\DD}$ play the role of the local normalising functions in conditions $_{\bar{k}_1}{\bf \tilde{D}}$ and  
$_{\bar{k}_1}{\bf \tilde{E}}$. 

The roles of the corresponding asymptotic relations appearing  in 
conditions $_{\bar{k}_1}{\bf \tilde{D}}$ and  $_{\bar{k}_1}{\bf \tilde{E}}$
are played by asymptotic relations  (\ref{trawetbanop}), (\ref{trawetbamnop}), and (\ref{expofavabeq}). 

Lemmas \, 14 -- 16 \, can be reformulated  in terms of random variables 
$_{\bar{k}_0}\tau_{\e, \DD}$, $_{\bar{k}_0}\eta_{\e}(_{\bar{k}_0}\tau_{\e, \DD})$, $_{\bar{k}_0}\tilde{\tau}_{\e, \DD}$, 
 $ _{\bar{k}_0}\tilde{\eta}_{\e}(_{\bar{k}_0}\tilde{\tau}_{\e, \DD})$, 
$_{\bar{k}_1}\tau_{\e, \DD} = \, _{k_1}\tau_{\e, \DD}$, 
$_{\bar{k}_1}\eta_{\e}(_{\bar{k}_1}\tau_{\e, \DD}) = \, _{k_1}\eta_{\e}(_{k_1}\tau_{\e, \DD})$, distributions 
$_{\bar{k}_0}G_{\e, \DD, ij}(\cdot)$, \, $_{\bar{k}_0}\tilde{G}_{\e, \DD, ij}(\cdot)$, \,  $_{\bar{k}_1}G_{\e, \DD, ij}(\cdot)$ 
$= \, _{k_1}G_{\e, \DD, ij}(\cdot)$ and Laplace
transforms $_{\bar{k}_0}\Psi_{\e, \DD, ij}(\cdot)$, \, $_{\bar{k}_0}\tilde{\Psi}_{\e, \DD, ij}(\cdot)$,  
$_{\bar{k}_1}\Psi_{\e, \DD, ij}(\cdot)$ $= \, _{k_1}\Psi_{\e, \DD, ij}(\cdot)$, for $ i \in \, _{\bar{k}_1}\XX
= \, _{k_1}\XX, j \in \DD$.

Lemmas \, 18 -- 20 can be reformulated  in  terms of random variables
$_{\bar{k}_0}\tau_{\e, \DD}$, $_{\bar{k}_0}\eta_{\e}(_{\bar{k}_0}\tau_{\e, \DD})$,
$_{\bar{k}_1}\tau_{\e, \DD}$,   $ _{\bar{k}_1}\eta_{\e}(_{\bar{k}_1}\tau_{\e, \DD})$,  \, 
$_{\bar{k}_1}\tilde{\tau}_{\e, \DD} = \, _{k_1}\tilde{\tau}_{\e, \DD}$, $_{\bar{k}_1}\tilde{\eta}_{\e}(_{\bar{k}_1}\tilde{\tau}_{\e, \DD}) = \, _{k_1}\tilde{\eta}_{\e}(_{\bar{k}_1}\tilde{\tau}_{\e, \DD})$, distributions 
$_{\bar{k}_0}G_{\e, \DD, ij}(\cdot)$, $_{\bar{k}_1}G_{\e, \DD, ij}(\cdot)$, \,  $_{\bar{k}_1}\tilde{G}_{\e, \DD, ij}(\cdot) \, = \, _{k_1}\tilde{G}_{\e, \DD, ij}(\cdot)$,  and Laplace transforms $_{\bar{k}_0}\Psi_{\e, \DD, ij}(\cdot)$, \, $_{\bar{k}_1}\Psi_{\e, \DD, ij}(\cdot)$,  \, 
$_{\bar{k}_1}\tilde{\Psi}_{\e, \DD, ij}(\cdot) = \, _{k_1}\tilde{\Psi}_{\e, \DD, ij}(\cdot)$, for $ i \in \, _{\bar{k}_1}\XX 
= \, _{k_1}\XX , j \in \DD$. 

Also, it is useful to note that, according Lemmas 7 -- 10,  condition $_{\bar{k}_0}{\bf C}_1 = {\bf C}_1$ implies that conditions  $_{\bar{k}_1}{\bf C} 
= \, _{k_1}{\bf C}$, $_{\bar{k}_1}{\bf \tilde{C}}$,  $_{\bar{k}_1}{\bf \tilde{C}}'$, $_{\bar{k}_1}{\bf \hat{C}} = \, _{k_1}{\bf \hat{C}}$ hold. Moreover, condition ${\bf C}_1$  is stronger than condition ${\bf C}_0$. That is why, condition ${\bf C}_1$ is also sufficient for holding of conditions $_{\bar{k}_0}{\bf C}$, $_{\bar{k}_0}{\bf \tilde{C}}$, $_{\bar{k}_0}{\bf \tilde{C}}'$.

Also, according Lemma 12, condition  $_{\bar{k}_0}{\bf F}_0 = {\bf F}_0$ implies that condition $_{\bar{k}_0}{\bf \tilde{F}} = {\bf \tilde{F}}$ holds.
\vspace{1mm}

{\bf 6.4 Step ${\bf 2}$}. The above recurrent procedure could be repeated. In order to clarify better the recurrent structure of  the algorithm, let us also describe its next step ${\bf 2}$. 

At step ${\bf 2}$, the semi-Markov process $_{\bar{k}_2}\eta_\e(t)$  is constructed, with the use of  procedure of  exclusion of a specially chosen state $k_2 \in \, _{\bar{k}_1}\overline{\DD}$ from the phase space 
$_{\bar{k}_1}\XX$ of the semi-Markov process  $_{\bar{k}_1}\tilde{\eta}_\e(t)$. Then,  
the semi-Markov process $_{\bar{k}_2}\tilde{\eta}_\e(t)$ is constructed,
with the use of  procedure of removing virtual transitions from trajectories of the semi-Markov process 
$_{\bar{k}_2}\eta_\e(t)$. The above procedures are described, respectively, in  Sections 4 and 5. Now, these procedures should be applied to the
semi-Markov processes $_{\bar{k}_1}\eta_\e(t)$ and $_{\bar{k}_1}\tilde{\eta}_\e(t)$,  instead of the semi-Markov processes
$_{\bar{k}_0}\eta_\e(t)$ and $_{\bar{k}_0}\tilde{\eta}_\e(t)$ used at step ${\bf 1}$.  
The semi-Markov processes  $_{\bar{k}_2}\eta_\e(t)$ and $_{\bar{k}_2}\tilde{\eta}_\e(t)$ have the phase space 
$_{\bar{k}_2}\XX =  \XX \ \setminus \{ k_1, k_2 \}$, while domain  $_{\bar{k}_2}\overline{\DD} = \overline{\DD} \ \setminus \{ k_1, k_2 \}$ 
replaces domain  $_{\bar{k}_1}\overline{\DD}$.

As at step ${\bf 1}$, conditions $_{\bar{k}_0}{\bf A}$ -- $_{\bar{k}_0}{\bf E}$ and $_{\bar{k}_0}{\bf \tilde{C}}$ (and, thus, condition $_{\bar{k}_0}{\bf \tilde{C}}'$ equivalent to condition $_{\bar{k}_0}{\bf \tilde{C}}$)  are assumed to hold. Also, it is assumed that condition 
$_{\bar{k}_0}{\bf \tilde{F}} = {\bf \tilde{F}}$ holds, state $k_1 \in \, _{\bar{k}_0}\overline{\DD}^* = \overline{\DD}^*$ is chosen, and, thus, condition $_{\bar{k}_1}{\bf \hat{F}} = \, _{k_1}{\bf \hat{F}}$ holds. Finally, 
conditions  $_{\bar{k}_1}{\bf \tilde{C}} = \, _{k_1}{\bf \tilde{C}}$ (and, thus, condition $_{\bar{k}_1}{\bf \tilde{C}}'$ equivalent to condition 
$_{\bar{k}_1}{\bf \tilde{C}}$),  and $_{\bar{k}_1}{\bf \hat{C}} = \, _{k_1}{\bf \hat{C}}$ should also be assumed to hold. Also, condition $_{\bar{k}_1}{\bf C} = \, _{k_1}{\bf C}$ holds, since it is implied by condition $_{\bar{k}_0}{\bf \tilde{C}}$. 

Also, as it was pointed out at step ${\bf 1}$, condition $_{\bar{k}_0}{\bf C}_1 = {\bf C}_1$ implies that conditions $_{\bar{k}_0}{\bf C}$, $_{\bar{k}_0}{\bf \tilde{C}}$,  $_{\bar{k}_0}{\bf \tilde{C}}'$ and $_{\bar{k}_1}{\bf C}$, $_{\bar{k}_1}{\bf \tilde{C}}$,  
$_{\bar{k}_1}{\bf \tilde{C}}'$, $_{\bar{k}_1}{\bf \hat{C}} = \, _{k_1}{\bf \hat{C}}$
hold, and condition  $_{\bar{k}_0}{\bf F}_0$ implies that condition $_{\bar{k}_0}{\bf \tilde{F}} = {\bf \tilde{F}}$ holds.

Let us denote by $_{\bar{k}_1}{\bf \tilde{F}}$ and $_{\bar{k}_2}{\bf \hat{F}}$, respectively, conditions
$_{\bar{k}_0}{\bf \tilde{F}}$ and  $_{\bar{k}_1}{\bf \hat{F}}$, in which the characteristics of 
semi-Markov processes $_{\bar{k}_0}\eta_\e(t)$ and $_{\bar{k}_0}\tilde{\eta}_\e(t)$
are replaced by the corresponding characteristics of semi-Markov processes $_{\bar{k}_1}\eta_\e(t)$   
and $_{\bar{k}_1}\tilde{\eta}_\e(t)$.

At step ${\bf 2}$, condition  $_{\bar{k}_1}{\bf \tilde{F}}$ should be additionally  assumed to hold. 

Condition $_{\bar{k}_1}{\bf \tilde{F}}$ let us  choose one of the least absorbing state states   
$k_2 \in \, _{\bar{k}_1}\overline{\DD}^*$
for the semi-Markov processes $_{\bar{k}_1}\eta_\e(t)$, with the  use of the algorithm described in Lemma 11. In fact, the choice of state $k_2 \in  \, _{\bar{k}_1}\overline{\DD}^*$ means that it is assumed  that  condition $_{\bar{k}_2}{\bf \hat{F}}$ holds for the semi-Markov processes  
$_{\bar{k}_1}\eta_\e(t)$.

Let us denote by $_{\bar{k}_2}{\bf C}$, $_{\bar{k}_2}{\bf \tilde{C}}$, $_{\bar{k}_2}{\bf \tilde{C}}'$, $_{\bar{k}_2}{\bf \hat{C}}$,  respectively, conditions $_{\bar{k}_1}{\bf C}$, $_{\bar{k}_1}{\bf \tilde{C}}$, $_{\bar{k}_1}{\bf \tilde{C}}'$, $_{\bar{k}_1}{\bf \hat{C}}$,  in which the characteristics of semi-Markov processes $_{\bar{k}_0}\eta_\e(t)$ and $_{\bar{k}_0}\tilde{\eta}_\e(t)$
are replaced by the corresponding characteristics of semi-Markov processes $_{\bar{k}_1}\eta_\e(t)$   
and $_{\bar{k}_1}\tilde{\eta}_\e(t)$. 

At  the following sub-steps ${\bf 2.1}$ and ${\bf 2.2}$, conditions $_{\bar{k}_2}{\bf \tilde{C}}$  and 
$_{\bar{k}_2}{\bf \hat{C}}$ should also be assumed to hold. Also, condition $_{\bar{k}_2}{\bf C}$ holds, since 
it is implied by condition $_{\bar{k}_1}{\bf \tilde{C}}$. Also, condition $_{\bar{k}_2}{\bf \tilde{C}}'$ (equivalent to condition $_{\bar{k}_2}{\bf \tilde{C}}$) holds. 

At sub-step ${\bf 2.1}$, the semi-Markov process $_{\bar{k}_2}\eta_\e(t)$ is constructed with the use of the procedure of exclusion 
state $k_2 \in \, _{\bar{k}_1}\overline{\DD}$ from the phase space $_{\bar{k}_1}\XX$ of the semi-Markov process  $_{\bar{k}_1}\tilde{\eta}_\e(t)$. This procedure is described in Section 4.  Now, it should be applied to
semi-Markov processes $_{\bar{k}_1}\tilde{\eta}_\e(t)$,  instead of the semi-Markov processes
$_{\bar{k}_0}\tilde{\eta}_\e(t)$ used at step ${\bf 1}$. 

At sub-step ${\bf 2.2}$, the semi-Markov process $_{\bar{k}_2}\tilde{\eta}_\e(t)$ is constructed  with the use of the procedure of removing virtual transitions from trajectories of the 
semi-Markov process $_{\bar{k}_2}\eta_\e(t)$. 
This procedure is described in Section 5. Now, it should be applied to
semi-Markov processes $_{\bar{k}_2}\eta_\e(t)$,  instead of the semi-Markov processes
$_{\bar{k}_1}\eta_\e(t)$ used at step ${\bf 1}$. 

According the above remarks,  the following notations are used for the corresponding characteristics of 
semi-Markov processes  $_{\bar{k}_2}\eta_\e(t)$ and $_{\bar{k}_2}\tilde{\eta}_\e(t)$, namely,  
probabilities, $_{\bar{k}_2}p_{\e, ij},  \, _{\bar{k}_2}\tilde{p}_{\e, ij}$, distribution functions 
$_{\bar{k}_2}F_{\e, ij}(\cdot), \, _{\bar{k}_2}\tilde{F}_{\e, ij}(\cdot)$, 
Laplace transforms $_{\bar{k}_2}\phi_{\e, ij}(\cdot), \, _{\bar{k}_2}\tilde{\phi}_{\e, ij}(\cdot)$, 
and expectations $_{\bar{k}_2}e_{\e, ij}, \, _{\bar{k}_2}\tilde{e}_{\e, ij}$, for $i, j \in \, _{\bar{k}_2}\XX$. 
The roles of relations, which express the 
above quantities for the semi-Markov processes $_{\bar{k}_2}\eta_\e(t)$ and $_{\bar{k}_2}\tilde{\eta}_\e(t)$ via the corresponding quantities for the semi-Markov processes  
$_{\bar{k}_1}\eta_\e(t)$ and $_{\bar{k}_1}\tilde{\eta}_\e(t)$,  
are played by relations (\ref{gopetk}), (\ref{traalk}), (\ref{trawetk}), and (\ref{diotrask}), for the semi-Markov processes $_{\bar{k}_2}\eta_\e(t)$,  and by relations  (\ref{gopetasd}), (\ref{transakaino}), (\ref{trawetno}), and (\ref{expofavabeq}), for the semi-Markov processes  $_{\bar{k}_2}\tilde{\eta}_\e(t)$. 

 According Lemmas 13 and 17, 
conditions  $_{\bar{k}_2}{\bf A}$ -- $_{\bar{k}_2}{\bf E}$  and $_{\bar{k}_2}{\bf \tilde{A}}$ -- $_{\bar{k}_2}{\bf \tilde{E}}$ and  $_{\bar{k}_2}{\bf \tilde{C}}'$ (which plays the role of condition ${\bf \tilde{C}}$
for the  the semi-Markov processes  $_{\bar{k}_2}\tilde{\eta}_\e(t)$) hold, respectively,  for the semi-Markov processes  $_{\bar{k}_2}\eta_\e(t)$ and $_{\bar{k}_2}\tilde{\eta}_\e(t)$.

Functions $_{\bar{k}_2}v_{\e, i} = \, _{\bar{k}_1}\tilde{v}_{\e, i} =  (1- \, _{\bar{k}_1}p_{\e, ii})^{-1}(1 - \, _{\bar{k}_0}p_{\e, ii})^{-1} v_{\e, i},  i \in \, _{\bar{k}_2}\overline{\DD}$ play the role of the local normalisation functions in the above conditions 
$_{\bar{k}_2}{\bf D}$ and  $_{\bar{k}_2}{\bf E}$.

The roles of the corresponding asymptotic relations appearing  in 
conditions $_{\bar{k}_2}{\bf D}$ and  $_{\bar{k}_2}{\bf E}$  are played by the asymptotic relations (\ref{trawetkasopa}), (\ref{trawetkasf}), and (\ref{trawetkasnop}).

Functions  $_{\bar{k}_2}\tilde{v}_{\e, i} = (1- \, _{\bar{k}_2}p_{\e, ii})^{-1} \, _{\bar{k}_2}v_{\e, i} = (1- \, _{\bar{k}_2}p_{\e, ii})^{-1}  (1- \, _{\bar{k}_1}p_{\e, ii})^{-1}(1 - \, _{\bar{k}_0}p_{\e, ii})^{-1} v_{\e, i},  i \in \, _{\bar{k}_2}\overline{\DD}$ play the role of the local normalisation  functions in conditions $_{\bar{k}_2}{\bf \tilde{D}}$ and  $_{\bar{k}_2}{\bf \tilde{E}}$. 

The roles of the corresponding asymptotic relations appearing  in 
conditions $_{\bar{k}_2}{\bf \tilde{D}}$ and  $_{\bar{k}_2}{\bf \tilde{E}}$ are played by the by asymptotic 
relations  (\ref{trawetbanop}), (\ref{trawetbamnop}), and (\ref{expofavabeq}). 

Lemmas \, 14 --16 \, can be reformulated for semi-Markov processes  $_{\bar{k}_0}\eta_\e(t)$, $_{\bar{k}_1}\tilde{\eta}_\e(t)$, and $_{\bar{k}_2}\eta_\e(t)$ replacing, respectively,  
semi-Markov processes  $_{\bar{k}_0}\eta_\e(t)$, $_{\bar{k}_0}\tilde{\eta}_\e(t)$, and $_{\bar{k}_1}\eta_\e(t)$, i.e., in terms of 
random variables  $_{\bar{k}_0}\tau_{\e, \DD}$, $_{\bar{k}_0}\eta_{\e}(_{\bar{k}_0}\tau_{\e, \DD})$, 
$_{\bar{k}_1}\tilde{\tau}_{\e, \DD}$, $_{\bar{k}_1}\tilde{\eta}_{\e}(_{\bar{k}_1}\tilde{\tau}_{\e, \DD})$, 
$_{\bar{k}_2}\tau_{\e, \DD}$, $_{\bar{k}_2}\eta_{\e}(_{\bar{k}_2}\tau_{\e, \DD})$, distributions 
$_{\bar{k}_0}G_{\e, \DD, ij}(\cdot)$, \, $_{\bar{k}_1}\tilde{G}_{\e, \DD, ij}(\cdot)$, \, $_{\bar{k}_2}G_{\e, \DD, ij}(\cdot)$, and Laplace
transforms $_{\bar{k}_0}\Psi_{\e, \DD, ij}(\cdot)$, \, $_{\bar{k}_1}\tilde{\Psi}_{\e, \DD, ij}(\cdot)$,  
$_{\bar{k}_2}\Psi_{\e, \DD, ij}(\cdot)$, for $ i \in \, _{\bar{k}_2}\XX, j \in \DD$.

Lemmas \, 18 -- 20 can be reformulated  for  semi-Markov processes $_{\bar{k}_0}\eta_\e(t)$,  $_{\bar{k}_2}\eta_\e(t)$, and $_{\bar{k}_2}\tilde{\eta}_\e(t)$ replacing, respectively,  
semi-Markov processes  $_{\bar{k}_0}\eta_\e(t)$,  $_{\bar{k}_1}\eta_\e(t)$, and $_{\bar{k}_2}\tilde{\eta}_\e(t)$, i.e., in terms of random variables 
$_{\bar{k}_0}\tau_{\e, \DD}$, $_{\bar{k}_0}\eta_{\e}(_{\bar{k}_0}\tau_{\e, \DD})$,  
$_{\bar{k}_2}\tau_{\e, \DD}$, $_{\bar{k}_2}\eta_{\e}(_{\bar{k}_2}\tau_{\e, \DD})$, 
$_{\bar{k}_2}\tilde{\tau}_{\e, \DD}$, $_{\bar{k}_2}\tilde{\eta}_{\e}(_{\bar{k}_2}\tilde{\tau}_{\e, \DD})$, distributions 
$_{\bar{k}_0}G_{\e, \DD, ij}(\cdot)$, \, $_{\bar{k}_2}G_{\e, \DD, ij}(\cdot)$, \,  $_{\bar{k}_2}\tilde{G}_{\e, \DD, ij}(\cdot) $,  and Laplace
transforms $_{\bar{k}_0}\Psi_{\e, \DD, ij}(\cdot)$, \, $_{\bar{k}_2}\Psi_{\e, \DD, ij}(\cdot)$,  
$_{\bar{k}_2}\tilde{\Psi}_{\e, \DD, ij}(\cdot)$, for $ i \in \, _{\bar{k}_2}\XX, j \in \DD$. 
\vspace{1mm}

{\bf 6.5 Conditions of asymptotic comparability for transition probabilities and normalisation functions}. 
Let us denote by $_{\bar{k}_1}{\bf C}_1 = \, _{k_1}{\bf C}_1$  condition
$_{\bar{k}_0}{\bf C}_1 = {\bf C}_1$, in which the characteristics of 
semi-Markov processes $_{\bar{k}_0}\eta_\e(t)$
are replaced by the corresponding characteristics of semi-Markov process $_{\bar{k}_1}\eta_\e(t)$.

Condition $_{\bar{k}_1}{\bf C}_1$ has the following form:
\begin{itemize}
\item [$_{\bar{k}_1}{\bf C}_{1}$:]  $_{\bar{k}_1}q_{\e}[i i' i'', j j' j''] = 
\frac{_{\bar{k}_1}p_{\e, i i'} \, _{\bar{k}_1}p_{\e, j j'}}{_{\bar{k}_1}p_{\e, i i''} \, _{\bar{k}_1}p_{\e, j j''}}  \to  
\, _{\bar{k}_1}q_{0}[i i' i'', j j' j'']  \in [0, \infty]$ as $\e \to 0$, for  
$i', j' \in \, _{\bar{k}_1}\XX, i'' \in \, _{\bar{k}_1}\YY_{1, i} = \, _{k_1}\YY_{1, i},  j'' \in \, _{\bar{k}_1}\YY_{1, j} = \, _{k_1}\YY_{1, j}, 
i, j \in \, _{\bar{k}_1}\overline{\DD}$.
%, $i \neq j$. 
\end{itemize}

According Lemmas 7 -- 10 (which should be applied to the semi-Markov processes $_{\bar{k}_1}\eta_\e(t)$   
 instead of the semi-Markov processes $_{\bar{k}_0}\eta_\e(t)$)  
condition $_{\bar{k}_1}{\bf C}_1$ implies that conditions  $_{\bar{k}_2}{\bf C}$, $_{\bar{k}_2}{\bf \tilde{C}}$,  and $_{\bar{k}_2}{\bf \hat{C}}$ hold.

Let us introduce condition:  
\begin{itemize}
\item [${\bf C}_{2}$:]  $q_{\e}[i_1i'_1 i''_1, \ldots, i_4 i'_4 i''_4] = 
\frac{p_{\e, i_1i'_1}p_{\e, i_2 i'_2}p_{\e, i_3 i'_3}p_{\e, i_4 i'_4}}{p_{\e, i_1i''_1}p_{\e, i_2 i''_2}p_{\e, i_3 i''_3}p_{\e, i_4 i''_4}}  
\to  q_{0}[i_1i'_1 i''_1, \ldots, i_4 i'_4 i''_4] \in [0, \infty]$ as $\e \to 0$, for  $i'_l  \in \XX, \, 
i''_l \in \, \YY_{1, i_l}, i_l \in \overline{\DD}, l = 1, \ldots, 4$. 
%, \,  i_1, i_2, i_3 = i_4 \in \overline{\DD}, i_1 \neq i_2 \neq i_3$.
\end{itemize}

{\bf Lemma 21}. {\em Condition ${\bf C}_{2}$ is sufficient for holding of condition 
$_k{\bf C}_1$, for any $k \in \overline{\DD}$.} \vspace{1mm}

{\bf Proof}. Since condition ${\bf B}$ holds for the semi-Markov processes $\eta_\e(t)$, and, thus, probabilities 
$1 - p_{\e, ii}, 1 - p_{\e, jj}, 1 - p_{\e, kk} \in (0, 1], \e \in (0, 1]$, for $i, j \in \, _k\overline{\DD}, i \neq j$. Thus, the following relation takes place,  for $i', i'' \in \, _k\YY_{1, i},  j', j'' \in \, _k\YY_{1, j}, i, j \in \, _k\overline{\DD}$, 
%i \neq j$,
\begin{align*}
_kq_{\e}[ii'i'', jj'j''] & =  \frac{_kp_{\e, ij'}\, _kp_{\e, jj'}}{_kp_{\e, ij''} \, _kp_{\e, jj''}} \makebox[60mm]{}
\vspace{2mm} \nonumber \\
& = \frac{\tilde{p}_{\e, ii'} + \tilde{p}_{\e, ik}\tilde{p}_{\e, ki'}}{\tilde{p}_{\e, ii''} + \tilde{p}_{\e, ik}\tilde{p}_{\e, ki''}}
\frac{\tilde{p}_{\e, jj'} + \tilde{p}_{\e, jk}\tilde{p}_{\e, kj'}}{\tilde{p}_{\e, jj''} + \tilde{p}_{\e, jk}\tilde{p}_{\e, kj''}} 
\vspace{2mm} \nonumber \\
\end{align*}
\begin{align}\label{fourth}
& = \frac{\sum_{l' \neq k}p_{\e, ii'}p_{\e, kl'}  + p_{\e, ik}p_{\e, ki'}}{\sum_{l'' \neq k}p_{\e, ii''}p_{\e, kl''} + p_{\e, ik}p_{\e, ki''}} 
\vspace{2mm} \nonumber \\
& \quad \times \frac{\sum_{m' \neq k}p_{\e, jj'}p_{\e, km'}  + 
p_{\e, rk}p_{\e, kj'}}{\sum_{m'' \neq k}p_{\e, jj''}p_{\e, km''} + p_{\e, jk}p_{\e, kj''}}, 
\end{align}

Since condition ${\bf B}$ holds for the semi-Markov processes $\eta_\e(t)$, $\tilde{\eta}_\e(t)$ and $_k\eta_\e(t)$, probabilities $_kp_{\e, ii'}, \, _kp_{\e, ii''}, \, _kp_{\e, jj'}, \, 
_kp_{\e, jj''} \in (0, 1], \e \in (0, 1]$, for $i', i'' \in \,  _k\YY_{1, i},  j', j'' \in \, _k\YY_{1, j}, 
i, j \in \, _k\overline{\DD}, i \neq j$.   That is why,  $\sum_{l' \neq k}p_{\e, ii'}p_{\e, kl'}  + p_{\e, ik}p_{\e, ki'}$,    
$\sum_{l'' \neq k}p_{\e, ii''}p_{\e, kl''} + p_{\e, ik}p_{\e, ki''}$, \ $\sum_{m' \neq k}p_{\e, jj'}p_{\e, km'}  + 
p_{\e, jk}p_{\e, kj'}$, \ $\sum_{m'' \neq k}p_{\e, rj''}p_{\e, km''}$ $+ p_{\e, jk}p_{\e, kj'} > 0,  \e \in (0, 1]$, for $i', i'' \in \,  _k\YY_{1, i},  j', j'' \in \, _k\YY_{1, j}, i, j \in \, 
_k\overline{\DD}$.
%, i \neq j$.  
Also,  every product  $p_{\e, ii'}p_{\e, kl'}, l' \neq k$ and $p_{\e, ik}p_{\e, ki'}$, as well as, every product $p_{\e, rj'}p_{\e, km'}$, $m' \neq k$ and $p_{\e, jk}p_{\e, kj'}$,
either takes positive value,  for every  $\e \in (0, 1]$, or equals $0$, for every  $\e \in (0, 1]$.
Let us introduce indicators  ${\rm I}_{ii'kl'} =  {\rm I}(p_{1, ii'}p_{1, kl'} > 0), l' \neq k, {\rm I}_{ikki'} =  {\rm I}(p_{1, ik}p_{1, ki'} > 0), {\rm I}_{rj'km'} =  {\rm I}(p_{1, rj'}p_{1, km'} > 0), m' \neq k, {\rm I}_{rkkj'} =  {\rm I}(p_{1, jk}p_{1, kj'} > 0)$.
By the above remarks, at least one of these indicators take value $1$.

Using the above remarks and relation  (\ref{fourth}), we get the following relation,  for $i', i'' \in \,  _k\YY_{1, i},  j', j'' \in \, _k\YY_{1, j}, i, j \in \, _k\overline{\DD}$, 
%i \neq j$,
\begin{align*}
& _kq_{\e}[ii'i'', jj'j'']  \vspace{2mm} \nonumber \\
& \quad = \Big( \sum_{l' \neq k} {\rm I}_{ii'kl'} 
\big( \sum_{l'' \neq k}  \frac{p_{\e, ii''} p_{\e, kl''}}{p_{\e, ii'} p_{\e, kl'}} 
+ \frac{p_{\e, kl''}}{p_{\e, kl'} } \big)^{-1}  \vspace{2mm} \nonumber \\
& \quad \ \  + {\rm I}_{ikki'} 
\big( \sum_{l'' \neq k}  \frac{p_{\e, ii''} p_{\e, kl''}}{p_{\e, ik}p_{\e, ki'}} 
+ \frac{p_{\e, kl''}}{p_{\e, ki'} } \big)^{-1} \Big) \vspace{2mm} \nonumber \\
& \quad \ \ \times \Big( \sum_{m' \neq k} {\rm I}_{rj'km'} 
\big( \sum_{m'' \neq k}  \frac{p_{\e, jj''} p_{\e, km''}}{p_{\e, jj'} p_{\e, km'}} 
+ \frac{p_{\e, km''}}{p_{\e, km'} } \big)^{-1}  \vspace{2mm} \nonumber \\
& \quad \ \  + {\rm I}_{rkkj'} \big( \sum_{m'' \neq k}  \frac{p_{\e, jj''} p_{\e, km''}}{p_{\e, jk}p_{\e, ks'}} 
+ \frac{p_{\e, km''}}{p_{\e, ks'} } \big)^{-1} \Big) \vspace{2mm} \nonumber \\
& \quad = \sum_{l', m' \neq k} {\rm I}_{ii'kl'}  {\rm I}_{jj'km'} \Big( \sum_{l'', m'' \neq k}  
\frac{p_{\e, ij''} p_{\e, kl''} p_{\e, jj''} p_{\e, km''}}{p_{\e, ii'} p_{\e, kl'}p_{\e, js'} p_{\e, km'}} \vspace{2mm} \nonumber \\
& \quad \ \ + \sum_{l'' \neq k}  \frac{p_{\e, ii''} p_{\e, kl''}p_{\e, km''}}{p_{\e, ii'} p_{\e, kl'}p_{\e, km'}} 
+  \sum_{m'' \neq k} \frac{p_{\e, kl''}p_{\e, js''} p_{\e, km''}}{p_{\e, kl'}p_{\e, js'} p_{\e, km'}} 
 +  \frac{p_{\e, kl''}p_{\e, km''}}{p_{\e, kl'} p_{\e, km'}} \Big)^{-1} 
  \vspace{2mm} \nonumber \\
   & \quad \ \ + \sum_{l' \neq k} {\rm I}_{ii'kl'}  {\rm I}_{jkks'} \Big( \sum_{l'', m'' \neq k} 
\frac{p_{\e, ii''} p_{\e, kl''}p_{\e, js''} p_{\e, km''}}{p_{\e, ii'} p_{\e, kl'}p_{\e, jk}p_{\e, kj'}} 
\vspace{2mm} \nonumber \\
  \end{align*}
\begin{align}\label{forthen}
& \quad \ \ + \sum_{l'' \neq k} \frac{p_{\e, ii''} p_{\e, kl''}p_{\e, km''}}{p_{\e, ii'} p_{\e, kl'}p_{\e, kj'} } 
+ \sum_{m'' \neq k}  \frac{p_{\e, kl''}p_{\e, js''} p_{\e, km''}}{p_{\e, kl'}p_{\e, jk}p_{\e, km'}}  
+  \frac{p_{\e, kl''}p_{\e, km''}}{p_{\e, kl'} p_{\e, kj'}} \Big)^{-1} 
\vspace{2mm} \nonumber \\
& \quad \quad + \sum_{m' \neq k} {\rm I}_{ikki'} {\rm I}_{jj'km'} \Big( \sum_{l'', m'' \neq k} 
\frac{p_{\e, ii''} p_{\e, kl''}p_{\e, jj''} p_{\e, km''}}{p_{\e, ik}p_{\e, ki'}p_{\e, jj'} p_{\e, km'}}  \vspace{2mm} \nonumber \\
& \quad \ \  + \sum_{l'' \neq k}  \frac{p_{\e, ii''} p_{\e, kl''}p_{\e, km''}}{p_{\e, ik}p_{\e, ki'}p_{\e, km'}}  
+ \sum_{m'' \neq k}  \frac{p_{\e, kl''}p_{\e, jj''} p_{\e, km''}}{p_{\e, ki'} p_{\e, jj'} p_{\e, km'}} 
+ \frac{p_{\e, kl''}p_{\e, km''}}{p_{\e, ki'}p_{\e, km'} } \Big)^{-1}  \vspace{2mm} \nonumber \\
& \quad = {\rm I}_{ikki'}  {\rm I}_{rkks'} 
\big( \sum_{l'', m'' \neq k}  \frac{p_{\e, ii''} p_{\e, kl''}p_{\e, jj''} p_{\e, km''}}{p_{\e, ik}p_{\e, ki'}p_{\e, jk}p_{\e, kj'}}   
\vspace{2mm} \nonumber \\
& \quad \ \  + \sum_{l'' \neq k} \frac{p_{\e, ii''} p_{\e, kl''}p_{\e, km''}}{p_{\e, ik}p_{\e, ki'}p_{\e, kj'} } 
+ \sum_{m'' \neq k} \frac{p_{\e, kl''}p_{\e, jj''} p_{\e, km''}}{p_{\e, ki'} p_{\e, jk}p_{\e, kj'}}  
+ \frac{p_{\e, kl''}p_{\e, km''}}{p_{\e, ki'} p_{\e, kj'}} \Big)^{-1},
\end{align}
where every product of the form ${\rm I}_{\mathbf{\cdot}} \, {\rm I}_{\mathbf{\cdot}} \, ( \mathbf{\cdot} )^{-1}$ should be counted as $0$ in the above sums,   if the corresponding product of indicators ${\rm I}_{\mathbf{\cdot}} \, {\rm I}_{\mathbf{\cdot}} = 0$.

Relation (\ref{forthen}) and condition ${\bf C}_{2}$ imply that the following relation holds, for $i', i'' \in \,  _k\YY_{1, i},  j', j'' \in \, _k\YY_{1, j}, i, j \in \, _k\overline{\DD}$, 
%i \neq j$,
\begin{align*}
&_kq_{\e}[ii'i'', jj'j'']  \vspace{2mm} \nonumber \\ 
& \quad \to   \sum_{l', m' \neq k} {\rm I}_{ii'kl'}  {\rm I}_{jj'km'} \Big( \sum_{l'', m'' \neq k}  q_{0}[ii''i', kl''l', jj''j', km''m']  
\vspace{2mm} \nonumber \\ 
& \quad \ \ + \sum_{l'' \neq k} q_{0}[ii''j', kl''l', km''m']  + \sum_{m'' \neq k} q_{0}[kl''l', jj''j', km''m']  
\vspace{2mm} \nonumber \\ 
& \quad \ \ + q_{0}[kl''l', km''m'] \Big)^{-1} \vspace{2mm} \nonumber \\
& \quad \ \  + \sum_{l' \neq k} {\rm I}_{ii'kl'}  {\rm I}_{jkkj'} 
\Big( \sum_{l'', m'' \neq k}  q_{0}[ii''j', kl''l', jj''k, km''j']  
\vspace{2mm} \nonumber \\ 
& \quad \ \ + \sum_{l'' \neq k} q_{0}[ii''i', kl''l', km''j']  + \sum_{m'' \neq k} q_{0}[kl''l', jj''k, km''m']  
\vspace{2mm} \nonumber \\
&  \quad  \ \ + q_{0}[kl''l', km''j'] \Big)^{-1} \vspace{2mm} \nonumber \\ 
&  \quad  \ \ + \sum_{m' \neq k} {\rm I}_{ikki'}  {\rm I}_{jj'km'} 
\Big( \sum_{l'', m'' \neq k}  q_{0}[ii''k, kl''i', jj''j', km''m']  
\vspace{2mm} \nonumber \\ 
& \quad \ \ + \sum_{l'' \neq k} q_{0}[ii''k, kl''i', km''m']  + \sum_{m'' \neq k} q_{0}[kl''i', jj''j', km''m']  \vspace{2mm} \nonumber \\ 
\end{align*}
\begin{align}\label{forthenas}
& \quad \ \ + q_{0}[kl''i', km''m'] \Big)^{-1}
\vspace{2mm} \nonumber \\ 
&  \quad  \ \ +  {\rm I}_{ikki'}  {\rm I}_{jkkj'} \Big( \sum_{l'', m'' \neq k}  q_{0}[ii''k, kl''i', jj''k, km''j']  
\vspace{2mm} \nonumber \\   
& \quad \ \ + \sum_{l'' \neq k} q_{0}[ii''k, kl''i', km''j']  + \sum_{m'' \neq k} q_{0}[kl''i', jj''k, km''j']  
\vspace{2mm} \nonumber \\ 
& \quad \ \ + q_{0}[kl''i', km''j'] \Big)^{-1}
 = \, _kq_{0}[ii'i'', jj'j''] \in [0, \infty] \ {\rm as} \ \e \to 0,
\end{align}  
where every product of the form ${\rm I}_{\mathbf{\cdot}} {\rm I}_{\mathbf{\cdot}} ( \mathbf{\cdot} )^{-1}$ should be counted as $0$ in the above sums , if the corresponding product of indicators ${\rm I}_{\mathbf{\cdot}} \, {\rm I}_{\mathbf{\cdot}} = 0$.

Relation  (\ref{forthenas}) implies that condition ${\bf C}_{1}$  is sufficient for holding of condition 
$_k{\bf C}_{1}$. $\Box$

Lemma 21 and the the above remarks imply that condition ${\bf C}_{2}$ is sufficient for holding conditions $_{\bar{k}_2}{\bf C}$, $_{\bar{k}_2}{\bf \tilde{C}}$,  and $_{\bar{k}_2}{\bf \hat{C}}$. Moreover, condition ${\bf C}_{2}$ is stronger than conditions ${\bf C}_1$ and ${\bf C}_0$. That is why, this condition is also sufficient for holding of conditions  $_{\bar{k}_1}{\bf C}$, 
$_{\bar{k}_1}{\bf \tilde{C}}$, $_{\bar{k}_1}{\bf \tilde{C}}'$, $_{\bar{k}_1}{\bf \hat{C}}$ and $_{\bar{k}_0}{\bf C}$, $_{\bar{k}_0}{\bf \tilde{C}}$, $_{\bar{k}_0}{\bf \tilde{C}}'$. 

Let us denote by  $_{\bar{k}_1}{\bf F}_0 = \, _{k_1}{\bf F}_0$ condition
$\, _{k_0}{\bf F}_0 = {\bf F}_0$, in which the characteristics of 
semi-Markov processes $_{\bar{k}_0}\eta_\e(t)$ and $_{\bar{k}_0}\tilde{\eta}_\e(t)$
are replaced by the corresponding characteristics of semi-Markov processes $_{\bar{k}_1}\eta_\e(t)$   
and $_{\bar{k}_1}\tilde{\eta}_\e(t)$.
Condition $_{\bar{k}_1}{\bf F}_0$ has the following form:
\begin{itemize}
\item [$_{\bar{k}_1}{\bf F}_0$:] $_{\bar{k}_1}u_{\e}[i i', j j'']  =   
\frac{_{\bar{k}_1}p_{\e, i i'} \, _{\bar{k}_1}v_{\e, i}^{-1}}{_{\bar{k}_1}p_{\e, j j''} \, _{\bar{k}_1}v_{\e, j}^{-1}} \to \, 
_{\bar{k}_1}u_{0}[i i', j j'']  \in [0, \infty]$  as  $\e \to 0$, for $i' \in \, _{\bar{k}_1}\XX$, 
$j'' \in \, _{\bar{k}_1}\YY_{1, j},  i, j \in \, _{\bar{k}_1}\overline{\DD}$.
%, i \neq j$.
\end{itemize}

According Lemma 12 (which should be applied to the semi-Markov processes $_{\bar{k}_1}\eta_\e(t)$   
and $_{\bar{k}_1}\tilde{\eta}_\e(t)$ instead of the semi-Markov processes $_{\bar{k}_0}\eta_\e(t)$ and $_{\bar{k}_0}\tilde{\eta}_\e(t)$)  
condition  $_{\bar{k}_1}{\bf F}_0$ implies that condition $_{\bar{k}_1}{\bf \tilde{F}}$ holds.

Let us introduce condition:  
\begin{itemize}
\item [${\bf F}_1$:] $u_{\e}[i_1 i_1', i_2 i_2'', i_3 i_3' i_3'']  =   
\frac{p_{\e, i_1 i_1'}p_{\e, i_3 i_3'} v_{\e, i_1}^{-1}}{p_{\e, i_2 i_2'} p_{\e, i_3 i_3''}v_{\e, j}^{-1}} \to u_{0}[i_1 i_1', i_2 i_2'', i_3 i_3' i_3'']  \in [0, \infty]$  as  $\e \to 0$, for $i_1', i_3' \in \XX, i_2'' \in \YY_{1, i_2}, 
i_3'' \in \YY_{1, i_3},  i_1, i_2, i_3 \in \overline{\DD}$.
%, i_1 \neq i_2 \neq i_3$.
\end{itemize}

{\bf Lemma 22}. {\em Condition ${\bf F}_{1} = \, _{\bar{k}_0}{\bf F}_{1}$ is sufficient for holding of condition $_k{\bf F}_0$, for any $k \in \overline{\DD}$.} \vspace{1mm}

{\bf Proof}. Since condition ${\bf B}$ holds for the semi-Markov processes $\eta_\e(t)$, and, thus, probabilities 
$1 - p_{\e, ii}, 1 - p_{\e, jj}, 1 - p_{\e, kk} \in (0, 1], \e \in (0, 1]$, for $i, j \in \, _k\overline{\DD}, i \neq j$. Thus, the following relation takes place,  for $j' \in \, _k\YY_{1, j},  i' \in \, _k\YY_{1, i}, i, j \in \, _k\overline{\DD}, i \neq j$, 
\begin{align}\label{afasdi}
_{k}u_{\e}[ii', jj']  & =   
\frac{_{k}p_{\e, ii'} \, _{k}v_{\e, i}^{-1}}{_{k}p_{\e, jj''} \, _{k}v_{\e, j}^{-1}} \makebox[50mm]{}
\vspace{2mm} \nonumber \\
%\end{align*}
%\begin{align}
& = \frac{(\tilde{p}_{\e, ii'} + \tilde{p}_{\e, ik}\tilde{p}_{\e, ki'}) 
(1 - p_{\e, ii})v_{\e, i}^{-1}}{(\tilde{p}_{\e, jj''} + \tilde{p}_{\e, jk}\tilde{p}_{\e, kj''}) 
(1 - p_{\e, jj})v_{\e, j}^{-1}}  \vspace{2mm} \nonumber \\
& = \frac{(\sum_{l' \neq k}p_{\e, ii'}p_{\e, kl'}  + p_{\e, ik} p_{\e, ki'}) 
v_{\e, i}^{-1}}{(\sum_{m' \neq k}p_{\e, jj''}p_{\e, km'}  + p_{\e, jk} p_{\e, kj''}) v_{\e, j}^{-1}}
\end{align}

Since condition ${\bf B}$ holds for the semi-Markov processes $\eta_\e(t)$, $\tilde{\eta}_\e(t)$ and $_k\eta_\e(t)$, probabilities $_kp_{\e, ii'}, \, _kp_{\e, jj''}, \in (0, 1], \e \in (0, 1]$, for $i'  \in \,  _k\YY_{1, i},  j'',  \in \, _k\YY_{1, j}$, 
$i, j \in \, _k\overline{\DD}, i \neq j$.   That is why,  $\sum_{l' \neq k}p_{\e, ii'}p_{\e, kl'}  + p_{\e, ik}p_{\e, ki'}$, \ 
$\sum_{m' \neq k}p_{\e, jj''}p_{\e, km'}  + p_{\e, jk}p_{\e, kj''}  > 0,  \e \in (0, 1]$, for $i' \in \,  _k\YY_{1, i},  
j'' \in \, _k\YY_{1, j}, i, j \in \, _k\overline{\DD}, i \neq j$. Also,  every product  $p_{\e, ii'}p_{\e, kl'}, l' \neq k$ and $p_{\e, ik}p_{\e, ki'}$, as well as, every product $p_{\e, jj''}p_{\e, km'}$, $m' \neq k$ and $p_{\e, jk}p_{\e, kj'}$,
either takes positive value,  for every  $\e \in (0, 1]$, or equals $0$, for every  $\e \in (0, 1]$.
Let us introduce indicators  ${\rm I}_{ii'kl'} =  {\rm I}(p_{1, ii'}p_{1, kl'} > 0), l' \neq k, {\rm I}_{ikki'} =  {\rm I}(p_{1, ik}p_{1, ki'} > 0)$. By the above remarks, at least one of these indicators take value $1$.

Using the above remarks and relation  (\ref{afasdi}), we get the following relation,  for $i' \in \,  _k\YY_{1, i},  j' \in \, _k\YY_{1, j}, i, j \in \, _k\overline{\DD}$, 
%i \neq j$,
\begin{align}\label{afasdin}
_{k}u_{\e}[ii', jj'']  & =   
 \sum_{l' \neq k} {\rm I}_{ii'kl'} \big(\sum_{m' \neq k}\frac{p_{\e, jj''}p_{\e, km'}v_{\e, j}^{-1}}{p_{\e, ii'}p_{\e, kl'}v_{\e, i}^{-1}}  
 + \frac{p_{\e, jk} p_{\e, kj''}v_{\e, j}^{-1}}{p_{\e, ii'}p_{\e, kl'}v_{\e, i}^{-1}} \big)^{-1} 
\vspace{2mm} \nonumber \\
%\end{align*}
%\begin{align}
& \quad + {\rm I}_{ikki'}\big( \sum_{m' \neq k}\frac{p_{\e, jj''}p_{\e, km'}v_{\e, j}^{-1}}{p_{\e, ik} p_{\e, ki'}v_{\e, i}^{-1}} + \frac{p_{\e, jk} p_{\e, kj'''}v_{\e, j}^{-1}}{p_{\e, ik} p_{\e, ki'}v_{\e, i}^{-1}} \big)^{-1},
\end{align}
where every product of the form ${\rm I}_{\mathbf{\cdot}}  \, ( \mathbf{\cdot} )^{-1}$ should be counted as $0$ in the above sums,   if the corresponding indicator ${\rm I}_{\mathbf{\cdot}}  = 0$.

Relation (\ref{afasdin}) and condition ${\bf F}_{1}$ imply that the following relation holds, for $i' \in \,  _k\YY_{1, i},  j'' \in \, _k\YY_{1, j}, i, j \in \, _k\overline{\DD}$, 
%i \neq j$,
\begin{align*}
_{k}u_{\e}[ii', jj'']  & \to   
 \sum_{l' \neq k} {\rm I}_{ii'kl'} \big(\sum_{m' \neq k} u_{0}[jj'', ii', km'l']  
 + u_{0}[jk, ii', kj''l'] \big)^{-1}  \vspace{2mm} \nonumber \\
 \end{align*}
 \begin{align}\label{afasdinno}
& \quad + {\rm I}_{ikki'}\big( \sum_{m' \neq k} u_{0}[jj'', ik, km' i'] +  u_{0}[jk, ik, kj''i'] \big)^{-1}  \vspace{2mm} \nonumber \\
& = \,  _{k}u_{0}[ii', jj''] \in [0, \infty] \  {\rm as} \ \e \to 0, 
\end{align}
where every product of the form ${\rm I}_{\mathbf{\cdot}}  \, ( \mathbf{\cdot} )^{-1}$ should be counted as $0$ in the above sums,   if the corresponding indicator ${\rm I}_{\mathbf{\cdot}}  = 0$.

Relation  (\ref{afasdinno}) implies that condition ${\bf F}_{1}$  is sufficient for holding of 
condition $_k{\bf F}_{0}$. $\Box$

Lemma 23 and the  above remarks imply that condition ${\bf F}_{1}$ is sufficient for holding condition $_{\bar{k}_1}{\bf \tilde{F}}$. Moreover, condition ${\bf F}_{1}$ is stronger than conditions ${\bf F}_0$. That is why, this condition is also sufficient for holding of condition  
$_{\bar{k}_0}{\bf \tilde{F}}$.

\vspace{1mm}

{\bf 6.6 A general description of recurrent asymptotic algorithm of phase space reduction}. We are now prepared to describe an arbitrary ${\bf n}$-th step of the phase space reduction algorithm, for $1 \leq n \leq \bar{m}- 1$. Therefore, we assume that $n-1$ steps of  this algorithms have been already realised, i.e., that the initial semi-Markov process $_{\bar{k}_{0}}\eta_\e(t)$ have been transformed with the use of procedure of removing virtual transitions from trajectories in this semi-Markov process  in the semi-Markov process $_{\bar{k}_{0}}\tilde{\eta}_\e(t)$, and, then, the above pair of semi-Markov processes $_{\bar{k}_{0}}\eta_\e(t), \, _{\bar{k}_{0}}\tilde{\eta}_\e(t)$ has been sequentially, for $r = 1, \ldots, n-1$, transformed  in the  pairs of reduced semi-Markov processes  $_{\bar{k}_{r}}\eta_\e(t), \, _{\bar{k}_{r}}\tilde{\eta}_\e(t)$,  with the use of procedure of removing virtual transitions from trajectories of the semi-Markov processes $_{\bar{k}_{r}}\eta_\e(t)$,  then the the exclusion of state $k_r$ from the phase space of the semi-Markov process  $_{\bar{k}_{r}}\tilde{\eta}_\e(t)$,  The above procedures are described, respectively, in  Sections 4 and 5.

Thus, it is assumed that initial basic conditions $_{\bar{k}_0}{\bf A}$ -- $_{\bar{k}_0}{\bf E}$,  and $_{\bar{k}_0}{\bf \tilde{C}}$ (and, thus, condition  $_{\bar{k}_0}{\bf \tilde{C}}'$ equivalent to condition $_{\bar{k}_0}{\bf \tilde{C}}$) hold.

Let us denote, for $r = 1, \ldots, n-1$, by $_{\bar{k}_r}{\bf C}$, $_{\bar{k}_r}{\bf \tilde{C}}$, $_{\bar{k}_r}{\bf \tilde{C}}'$, and $_{\bar{k}_r}{\bf \hat{C}}$, respectively, conditions $_{\bar{k}_1}{\bf C}$, $_{\bar{k}_1}{\bf \tilde{C}}$, $_{\bar{k}_1}{\bf \tilde{C}}'$, $_{\bar{k}_1}{\bf \hat{C}}$,    in which the characteristics of semi-Markov processes $_{\bar{k}_0}\eta_\e(t)$ and $_{\bar{k}_0}\tilde{\eta}_\e(t)$
are replaced by the corresponding characteristics of semi-Markov processes $_{\bar{k}_r}\eta_\e(t)$   
and $_{\bar{k}_r}\tilde{\eta}_\e(t)$.  

The realisation of $n -1$ steps for  the phase space reduction algorithm pointed above requires to assume that conditions 
$_{\bar{k}_r}{\bf \tilde{C}}$, $_{\bar{k}_r}{\bf \hat{C}}, r = 1, \ldots, n-1$ are assumed to hold.  Also conditions $_{\bar{k}_r}{\bf C}, r = 1, \ldots, n-1$ hold (since condition $_{\bar{k}_r}{\bf C}$ is implied by condition $_{\bar{k}_{r-1}}{\bf \tilde{C}}$, for every  $r = 1, \ldots, n-1$), and conditions 
$_{\bar{k}_r}{\bf \tilde{C}}',  r = 1, \ldots, n- 1$ hold (since condition $_{\bar{k}_r}{\bf \tilde{C}}'$ is equivalent to condition $_{\bar{k}_r}{\bf \tilde{C}}$
 for every  $r = 1, \ldots, n-1$). 

Let us also denote, for $r = 1, \ldots, n$, by $_{\bar{k}_{r-1}}{\bf \tilde{F}}$ condition  $_{\bar{k}_{0}}{\bf \tilde{F}}$, in which the characteristics of semi-Markov processes $_{\bar{k}_0}\eta_\e(t)$ and $_{\bar{k}_0}\tilde{\eta}_\e(t)$
are replaced by the corresponding characteristics of semi-Markov processes $_{\bar{k}_{r-1}}\eta_\e(t)$   
and $_{\bar{k}_{r-1}}\tilde{\eta}_\e(t)$.  

Also conditions 
$_{\bar{k}_{r-1}}{\bf \tilde{F}}, r = 1, \ldots, n-1$ are assumed to hold, and states $k_r \in \, 
_{\bar{k}_{r-1}}\overline{\DD}^*$ are chosen, for $r = 1, \ldots, n-1$. This implies that conditions
$_{\bar{k}_{r}}{\bf \hat{F}}, r = 1, \ldots, n-1$ hold.

According  Lemmas 13 and 17 (applied sequentially, for $r = 1, \ldots, n-1$,  to  the semi-Markov processes $_{\bar{k}_{r}}\eta_\e(t)$ and $_{\bar{k}_{r}}\tilde{\eta}_\e(t)$), conditions $_{\bar{k}_{r}}{\bf A}$ -- $_{\bar{k}_r}{\bf E}$ and 
$_{\bar{k}_{r}}{\bf \tilde{A}}$ -- $_{\bar{k}_r}{\bf \tilde{E}}$, and $_{\bar{k}_r}{\bf \tilde{C}}'$ hold, for $r = 1, \ldots, n -1$.

At ${\bf n}$-th step, condition  $_{\bar{k}_{n-1}}{\bf \tilde{F}}$ should be additionally  assumed to hold. 

Condition $_{\bar{k}_{n-1}}{\bf \tilde{F}}$ let us  choose one of the least absorbing state states   $k_n \in \, _{\bar{k}_{n-1}}\overline{\DD}^*$
for the semi-Markov processes $_{\bar{k}_{n-1}}\eta_\e(t)$, with the  use of the algorithm described in Lemma 11. In fact, the choice of state 
$k_n \in  \, _{\bar{k}_{n-1}}\overline{\DD}^*$ means that it is assumed  that  condition $_{\bar{k}_n}{\bf \hat{F}}$ holds for the semi-Markov processes  $_{\bar{k}_{n-1}}\eta_\e(t)$.

Let us denote by $_{\bar{k}_n}{\bf C}$, $_{\bar{k}_n}{\bf \tilde{C}}$, $_{\bar{k}_n}{\bf \tilde{C}}'$, $_{\bar{k}_n}{\bf \hat{C}}$,  respectively, conditions $_{\bar{k}_1}{\bf C}$, $_{\bar{k}_1}{\bf \tilde{C}}$, $_{\bar{k}_1}{\bf \tilde{C}}'$, $_{\bar{k}_1}{\bf \hat{C}}$,  in which the characteristics of semi-Markov processes $_{\bar{k}_0}\eta_\e(t)$, $_{\bar{k}_0}\tilde{\eta}_\e(t)$, and $_{\bar{k}_1}\eta_\e(t)$
are replaced by the corresponding characteristics of semi-Markov processes $_{\bar{k}_{n-1}}\eta_\e(t)$,   $_{\bar{k}_{n -1}}\tilde{\eta}_\e(t)$, and  $_{\bar{k}_{n}}\eta_\e(t)$. 

At  the following sub-steps ${\bf n.1}$ and ${\bf n.2}$, conditions $_{\bar{k}_n}{\bf \tilde{C}}$  and 
$_{\bar{k}_n}{\bf \hat{C}}$ should also be assumed to hold. Also, condition $_{\bar{k}_n}{\bf C}$ holds, since 
it is implied by condition $_{\bar{k}_{n-1}}{\bf \tilde{C}}$, and condition $_{\bar{k}_n}{\bf \tilde{C}}'$ (equivalent to condition $_{\bar{k}_n}{\bf \tilde{C}}$) holds.

At sub-step ${\bf n.1}$, the semi-Markov process $_{\bar{k}_n}\eta_\e(t)$ is constructed with the use of the procedure of exclusion 
state $k_n \in \, _{\bar{k}_{n-1}}\overline{\DD}$ from the phase space $_{\bar{k}_{n-1}}\XX = \XX \setminus \{ k_1, \ldots, k_{n-1} \}$ of the semi-Markov process  $_{\bar{k}_{n-1}}\tilde{\eta}_\e(t)$. This procedure is described in Section 4.  Now, it should be applied to
semi-Markov processes $_{\bar{k}_{n-1}}\tilde{\eta}_\e(t)$,  instead of the semi-Markov processes
$_{\bar{k}_0}\tilde{\eta}_\e(t)$ used at sub-step ${\bf 1.1}$. 

At sub-step ${\bf n.2}$,  the semi-Markov process $_{\bar{k}_n}\tilde{\eta}_\e(t)$ is constructed  with the use of the procedure of removing virtual transitions from trajectories of the semi-Markov process 
$_{\bar{k}_n}\eta_\e(t)$. 
This procedure is described in Section 5. Now, it should be applied to
semi-Markov processes $_{\bar{k}_n}\eta_\e(t)$,  instead of the semi-Markov processes
$_{\bar{k}_1}\eta_\e(t)$ used at sub-step ${\bf 1.2}$. 

The semi-Markov processes $_{\bar{k}_n}\eta_\e(t)$ and $_{\bar{k}_n}\tilde{\eta}_\e(t)$ have the phase space 
$_{\bar{k}_{n}}\XX = \XX \setminus \{ k_1, \ldots, k_{n} \}$. The role of domain $\overline{\DD}$ is played by the domain 
$_{\bar{k}_{n}}\overline{\DD} = \overline{\DD} \setminus \{ k_1, \ldots, k_{n} \}$.  

According the above remarks,  the following notations are used for the corresponding characteristics of 
semi-Markov processes  $_{\bar{k}_n}\eta_\e(t)$ and $_{\bar{k}_n}\tilde{\eta}_\e(t)$, namely,  
probabilities, $_{\bar{k}_n}p_{\e, ij},  \, _{\bar{k}_n}\tilde{p}_{\e, ij}$, distribution functions 
$_{\bar{k}_n}F_{\e, ij}(\cdot), \, _{\bar{k}_n}\tilde{F}_{\e, ij}(\cdot)$, 
Laplace transforms $_{\bar{k}_n}\phi_{\e, ij}(\cdot), \, _{\bar{k}_n}\tilde{\phi}_{\e, ij}(\cdot)$, 
and expectations $_{\bar{k}_n}e_{\e, ij}, \, _{\bar{k}_n}\tilde{e}_{\e, ij}$, for $i, j \in \, _{\bar{k}_n}\XX$. 

The roles of relations, which express the 
above quantities for the semi-Markov processes $_{\bar{k}_n}\eta_\e(t)$ and $_{\bar{k}_n}\tilde{\eta}_\e(t)$ via the corresponding quantities for the semi-Markov processes  $_{\bar{k}_{n-1}}\eta_\e(t)$ and 
'$_{\bar{k}_{n-1}}\tilde{\eta}_\e(t)$,  
are played by relations (\ref{gopetk}), (\ref{traalk}), (\ref{trawetk}), and (\ref{diotrask}), for the semi-Markov processes $_{\bar{k}_n}\eta_\e(t)$,  and by relations  (\ref{gopetasd}), (\ref{transakaino}), (\ref{trawetno}), and (\ref{expofavabeq}), for the semi-Markov processes  $_{\bar{k}_n}\tilde{\eta}_\e(t)$. 

 According Lemmas 13 and 17, 
conditions  $_{\bar{k}_n}{\bf A}$ -- $_{\bar{k}_n}{\bf E}$  and $_{\bar{k}_n}{\bf \tilde{A}}$ -- $_{\bar{k}_n}{\bf \tilde{E}}$ and   $_{\bar{k}_n}{\bf \tilde{C}}'$ (which plays the role of condition ${\bf \tilde{C}}$
for the  the semi-Markov processes  $_{\bar{k}_n}\tilde{\eta}_\e(t)$) hold, respectively,  for the semi-Markov processes  $_{\bar{k}_n}\eta_\e(t)$ and $_{\bar{k}_n}\tilde{\eta}_\e(t)$.

Functions $_{\bar{k}_n}v_{\e, i} = \, _{\bar{k}_{n-1}}\tilde{v}_{\e, i} =  \prod_{l = 0}^{n-1} (1- \, _{\bar{k}_l}p_{\e, ii})^{-1} v_{\e, i},  i \in \, _{\bar{k}_n}\overline{\DD}$ play the role of the local normalising functions in the above conditions 
$_{\bar{k}_n}{\bf D}$ and  $_{\bar{k}_n}{\bf E}$.

The roles of the corresponding asymptotic relations appearing  in 
conditions $_{\bar{k}_n}{\bf D}$ and  $_{\bar{k}_n}{\bf E}$  are played by the asymptotic relations (\ref{trawetkasopa}), (\ref{trawetkasf}), and (\ref{trawetkasnop}), where characteristics of the semi-Markov processes 
$_{\bar{k}_{0}}\tilde{\eta}_\e(t)$ and $_{\bar{k}_1}\eta_\e(t)$ should be replaced by the corresponding characteristics of the semi-Markov processes $_{\bar{k}_{n-1}}\tilde{\eta}_\e(t)$ and $_{\bar{k}_n}\eta_\e(t)$. 

Functions  $_{\bar{k}_n}\tilde{v}_{\e, i} = (1- \, _{\bar{k}_n}p_{\e, ii})^{-1} \, _{\bar{k}_n}v_{\e, i} = \prod_{l = 0}^{n}  (1- \, _{\bar{k}_l}p_{\e, ii})^{-1}   v_{\e, i},  i \in \, _{\bar{k}_n}\overline{\DD}$ play the role of the local normalising functions in conditions $_{\bar{k}_n}{\bf \tilde{D}}$ and  
$_{\bar{k}_n}{\bf \tilde{E}}$. 

The roles of the corresponding asymptotic relations appearing  in 
conditions $_{\bar{k}_n}{\bf \tilde{D}}$ and  $_{\bar{k}_n}{\bf \tilde{E}}$ are played by the by asymptotic 
relations  (\ref{trawetbanop}), (\ref{trawetbamnop}), and (\ref{expofavabeq}), where characteristics of the semi-Markov processes $_{\bar{k}_1}\eta_\e(t)$ and $_{\bar{k}_{1}}\tilde{\eta}_\e(t)$ should be replaced by the corresponding characteristics of the semi-Markov processes  $_{\bar{k}_n}\eta_\e(t)$ and $_{\bar{k}_{n}}\tilde{\eta}_\e(t)$. 

Lemmas \, 14 -- 16 \, can be reformulated for semi-Markov processes  $_{\bar{k}_0}\eta_\e(t)$, $_{\bar{k}_{n-1}}\tilde{\eta}_\e(t)$, and $_{\bar{k}_n}\eta_\e(t)$ replacing, respectively,  
semi-Markov processes  $_{\bar{k}_0}\eta_\e(t)$, $_{\bar{k}_0}\tilde{\eta}_\e(t)$, and $_{\bar{k}_1}\eta_\e(t)$, i..e., in terms of random variables  $_{\bar{k}_0}\tau_{\e, \DD}$, $_{\bar{k}_0}\eta_{\e}(_{\bar{k}_0}\tau_{\e, \DD})$, $_{\bar{k}_{n-1}}\tilde{\tau}_{\e, \DD}$, \, $_{\bar{k}_{n-1}}\tilde{\eta}_{\e}(_{\bar{k}_{n-1}}\tilde{\tau}_{\e, \DD})$, \, $_{\bar{k}_n}\tau_{\e, \DD}$, \,  $_{\bar{k}_n}\eta_{\e}(_{\bar{k}_n}\tau_{\e, \DD})$, distributions 
\, $_{\bar{k}_0}G_{\e, \DD, ij}(\cdot)$, \, $_{\bar{k}_{n-1}}\tilde{G}_{\e, \DD, ij}(\cdot)$, \, $_{\bar{k}_n}G_{\e, \DD, ij}(\cdot)$, and Laplace
transforms $_{\bar{k}_0}\Psi_{\e, \DD, ij}(\cdot)$, \, $_{\bar{k}_{n-1}}\tilde{\Psi}_{\e, \DD, ij}(\cdot)$,  
$_{\bar{k}_n}\Psi_{\e, \DD, ij}(\cdot)$, for $ i \in \, _{\bar{k}_n}\XX, j \in \DD$.

Lemmas \, 18 -- 20 can be reformulated  for  semi-Markov processes $_{\bar{k}_0}\eta_\e(t)$,  $_{\bar{k}_n}\eta_\e(t)$, and $_{\bar{k}_n}\tilde{\eta}_\e(t)$ replacing, respectively,  
semi-Markov processes  $_{\bar{k}_0}\eta_\e(t)$,  $_{\bar{k}_1}\eta_\e(t)$, and $_{\bar{k}_2}\tilde{\eta}_\e(t)$,  i.e., in  terms of random variables 
$_{\bar{k}_0}\tau_{\e, \DD}$, $_{\bar{k}_0}\eta_{\e}(_{\bar{k}_0}\tau_{\e, \DD})$,  
$_{\bar{k}_n}\tau_{\e, \DD}$, $_{\bar{k}_n}\eta_{\e}(_{\bar{k}_n}\tau_{\e, \DD})$, 
$_{\bar{k}_n}\tilde{\tau}_{\e, \DD}$, $_{\bar{k}_n}\tilde{\eta}_{\e}(_{\bar{k}_n}\tilde{\tau}_{\e, \DD})$, distributions 
$_{\bar{k}_0}G_{\e, \DD, ij}(\cdot)$, \, $_{\bar{k}_n}G_{\e, \DD, ij}(\cdot)$,  \,  $_{\bar{k}_n}\tilde{G}_{\e, \DD, ij}(\cdot) $,  and Laplace
transforms $_{\bar{k}_0}\Psi_{\e, \DD, ij}(\cdot)$, \, $_{\bar{k}_n}\Psi_{\e, \DD, ij}(\cdot)$,  
$_{\bar{k}_n}\tilde{\Psi}_{\e, \DD, ij}(\cdot)$, for $ i \in \, _{\bar{k}_n}\XX, j \in \DD$. 
\vspace{1mm}

{\bf 6.7 Summary}. The following lemma summarises the above remarks. \vspace{1mm}

{\bf Lemma 23}.  {\em Let conditions ${\bf A}$ -- ${\bf E}$, ${\bf \tilde{C}}$, $_{\bar{k}_{r}}{\bf \tilde{C}}$, $_{\bar{k}_{r}}{\bf \hat{C}}$, $_{\bar{k}_{r-1}}{\bf \tilde{F}}$, $r = 1, \ldots, n$ and hold, where the  states $k_r \in \, _{\bar{k}_{r-1}} \overline{\DD}^*, r = 1, \ldots,  n -1$ are chosen in such way that condition
$_{\bar{k}_{r}}{\bf \hat{F}}$ holds, for  $r = 1, \ldots,  n$. Then, 
conditions ${\bf A}$ -- ${\bf E}$ and ${\bf \tilde{C}}$ hold for  the  semi-Markov processes $_{\bar{k}_n}\eta_{\e}(t)$ and 
$_{\bar{k}_n}\tilde{\eta}_{\e}(t)$, respectively,  in the form of conditions $_{\bar{k}_n}{\bf A}$ -- $_{\bar{k}_n}{\bf E}$, $_{\bar{k}_n}{\bf \tilde{C}}'$ and  
$_{\bar{k}_n}{\bf \tilde{A}}$ -- $_{\bar{k}_n}{\bf \tilde{E}}$, $_{\bar{k}_n}{\bf \tilde{C}}'$.} \vspace{1mm}

{\bf 6.7 Asymptotic comparability for transition probabilities and normalising functions}. 
Let us denote by $_{\bar{k}_{n-1}}{\bf C}_2$  condition
$_{\bar{k}_0}{\bf C}_2$, in which the characteristics of 
semi-Markov processes $_{\bar{k}_0}\eta_\e(t)$
are replaced by the corresponding characteristics of semi-Markov process $_{\bar{k}_{n-1}}\eta_\e(t)$.

Let us denote, for $i \in \, _{\bar{k}_{n-1}}\overline{\DD}$,
\begin{equation}\label{nota}
_{\bar{k}_{n-1}}\YY_{\e, i} = \{r \in \, _{\bar{k}_{n-1}}\overline{\DD}:  \, _{\bar{k}_{n-1}}p_{\e, ir} > 0 \}
\end{equation}

Condition $_{\bar{k}_{n-1}}{\bf B}$ implies that sets $_{\bar{k}_{n-1}}\YY_{\e, i}  =  \, _{\bar{k}_{n-1}}\YY_{1, i}, \e \in (0, 1]$, 
for $i \in \, _{\bar{k}_{n-1}}\overline{\DD}$.

Condition $_{\bar{k}_{n-1}}{\bf C}_1$ has the following form: \vspace{2mm}

\noindent $_{\bar{k}_{n-1}}{\bf C}_{1}$: $_{\bar{k}_{n-1}}q_{\e}[ii'i'',jj'j''] = 
\frac{_{\bar{k}_{n-1}}p_{\e, ii'} \, _{\bar{k}_{n-1}}p_{\e, jj'}}{_{\bar{k}_{n-1}}p_{\e, ii''} \, _{\bar{k}_{n-1}}p_{\e, jj''}} 
 \to  \, _{\bar{k}_{n-1}}q_{0}[ii'i'', jj'j'']  \in [0, \infty]$ \vspace{1mm} \\  
\makebox[13mm]{} as $\e \to 0$, for  $i', j' \in \, _{\bar{k}_{n -1}}\XX, i'' \in \, _{\bar{k}_{n-1}}\YY_{1, i}$,  $j'' \in  \, _{\bar{k}_{n-1}}\YY_{1, j}$, $i, j \in \, _{\bar{k}_{n-1}}\overline{\DD}$.
%\makebox[13mm]{} $i \neq j$.  
\vspace{2mm}

According Lemmas 11--14 (which should be applied to the semi-Markov processes $_{\bar{k}_{n-1}}\eta_\e(t)$   
 instead of the semi-Markov processes $_{\bar{k}_0}\eta_\e(t)$)  condition $_{\bar{k}_{n-1}}{\bf C}_{1}$ implies that  
 conditions $_{\bar{k}_n}{\bf C}$, $_{\bar{k}_n}{\bf \tilde{C}}$, $_{\bar{k}_n}{\bf \tilde{C}}'$, and $_{\bar{k}_n}{\bf \hat{C}}$ hold.

Let us introduce, for $n \geq 3$,  condition:  
\begin{itemize}
\item [${\bf C}_{n}$:]  $q_{\e}[i_1 i'_1, i''_1, \ldots, i_{2^n} i'_{2^n} i''_{2^n}] =  \prod_{l = 0}^{2^n} 
\frac{p_{\e, i_l i'_l}}{p_{\e, i_li''_l}}  
\to  q_{0}[i_1 i'_1, i''_1, \ldots, i_{2^n} i'_{2^n} i''_{2^n}]  \in [0, \infty]$  as $\e \to 0$, for  $i'_l \in \, \XX, 
i''_l \in \, \YY_{1, i_l}, i_l \in \overline{\DD}, l = 1, \ldots, 2^n$. 
%, \ i_1, i_2, i_3 = i_4, i_5 = i_6 = i_7 = i_8, \ldots, i_{2^{n-1} + 1} = \cdots =   i_{2^{n}} \in \overline{\DD}, \, i_1 \neq i_2 \neq i_3 \neq \cdots \neq i_{2^{n-1} + 1}$.
\end{itemize}

{\bf Lemma 24}. {\em Condition ${\bf C}_{n}$ is sufficient for holding of condition 
$_{\bar{k}_{n-1}}{\bf C}_1$, for any $k_1, \ldots, k_{n-1} \in \overline{\DD},k_1 \neq \cdots \neq k_{n-1}$.} \vspace{1mm}

The proof of Lemma 24 is analogous to  the proof of Lemma 21. 

Lemma 24 and the above remarks imply that condition ${\bf C}_{n}$ is sufficient for holding of conditions $_{\bar{k}_n}{\bf C}$, 
$_{\bar{k}_n}{\bf \tilde{C}}$, $_{\bar{k}_n}{\bf \tilde{C}}'$, and $_{\bar{k}_n}{\bf \hat{C}}$. Moreover, since condition ${\bf C}_{n}$ is stronger than condition ${\bf C}_{r-1}$, for $r = 1, \ldots, n$, condition 
${\bf C}_{n}$ is, in fact, sufficient for holding of all conditions $_{\bar{k}_r}{\bf C}$, 
$_{\bar{k}_r}{\bf \tilde{C}}$, $_{\bar{k}_r}{\bf \tilde{C}}'$, $_{\bar{k}_r}{\bf \hat{C}}, r = 0, \ldots, n$.

The comparability condition ${\bf C}_{n}$ is expressed via transition probabilities $p_{\e, ij}$ in much more explicit form than conditions   $_{\bar{k}_r}{\bf C}$, $_{\bar{k}_r}{\bf \tilde{C}}$, $_{\bar{k}_r}{\bf \tilde{C}}'$, $_{\bar{k}_r}{\bf \hat{C}}, r = 0, \ldots, n$. 

A disadvantage  of condition ${\bf C}_{n}$ is that it   requires existence of  a large number of limits. 

However,  some effective sufficient (for holding of condition ${\bf C}_{n}$ conditions based on the notion of a complete family of asymptotically comparable functions can be formulated.

A family of positive functions ${\cal H} = \{ h(\cdot) \}$ defined on interval $(0, 1]$ is a complete family of asymptotically comparable functions if: {\bf (1)} it is closed with respect to operation of summation, multiplication and devision, and 
{\bf (2)} there exist $\lim_{\e \to 0} h(\e) = a[h(\cdot)] \in [0, \infty]$, for any 
function $h(\cdot) \in {\cal H}$.

An example is the family ${\cal H}_1$ of positive functions such that, for any  $h(\cdot)  \in {\cal H}_1$,  there exist constants $a_h > 0$ and $b_h \in (-\infty, \infty)$ such that, 
\begin{equation}\label{byrce}
\frac{h(\e)}{a_h \e^{b_h}} \to 1 \ {\rm  as} \ \e \to 0.
\end{equation}

The detailed discussion related to the above notion of complete family of asymptotically comparable functions and other examples of such families are given in Appendix A.

This condition automatically hold and limits appearing in it  are easily computable, for example,  if condition 
${\bf B}$  holds and transition probabilities $p_{\e, ij}, \e \in (0, 1], j \in \YY_{1, i}, i \in \overline{\DD}$ satisfy the following condition:
\begin{itemize}
\item [${\bf G}$:] Functions $p_{\cdot, ij}, j \in \YY_{1, i}, i \in \overline{\DD}$ belong to some complete family of asymptotically comparable 
functions.    
\end{itemize} 

{\bf Lemma 25}. {\em Condition ${\bf G}$ is sufficient for holding of condition 
${\bf C}_{\bar{m} -1}$ and, thus,  for holding of conditions $_{\bar{k}_r}{\bf C}$, $_{\bar{k}_r}{\bf \tilde{C}}$, $_{\bar{k}_r}{\bf \tilde{C}}'$, $_{\bar{k}_r}{\bf \hat{C}}, r = 0, \ldots, \bar{m} -1$, for any $k_1, \ldots, k_{\bar{m}-1} \in \overline{\DD}, k_1 \neq \cdots \neq k_{\bar{m} -1}$.} \vspace{1mm}
 
The proof of this lemma and the discussion connected  additional restrictions caused  by stochasticity of matrix $\| p_{\e, ij} \|$ are given in Appendix A.

Here, we also would like to mention that only limits appearing in conditions $_{\bar{k}_r}{\bf C}$, $_{\bar{k}_r}{\bf \tilde{C}}$, $_{\bar{k}_r}{\bf \tilde{C}}'$, $_{\bar{k}_r}{\bf \hat{C}}, r = 0, \ldots, n$ are actually  used for computing limiting Laplace transforms in conditions $_{\bar{k}_0}{\bf \tilde{D}}'$, $_{\bar{k}_r}{\bf D}'$, 
$_{\bar{k}_r}{\bf \tilde{D}},  r = 1, \ldots, n$. As far as limits in condition ${\bf C}_{2^{n}}$ are concerned, their existence is sufficient for existence, but these limits itself are not directly used  for computing the  above mentioned limits appearing in conditions $_{\bar{k}_r}{\bf C}$, $_{\bar{k}_r}{\bf \tilde{C}}$, $_{\bar{k}_r}{\bf \tilde{C}}'$, $_{\bar{k}_r}{\bf \hat{C}}, r = 0, \ldots, n$.

Let us denote by  $_{\bar{k}_{n-1}}{\bf F}_0$ condition
$_{\bar{k}_{0}}{\bf F}_0$, in which the characteristics of 
semi-Markov processes $_{\bar{k}_0}\eta_\e(t)$ and $_{\bar{k}_0}\tilde{\eta}_\e(t)$
are replaced by the corresponding characteristics of semi-Markov processes $_{\bar{k}_{n-1}}\eta_\e(t)$   
and $_{\bar{k}_{n-1}}\tilde{\eta}_\e(t)$.
Condition $_{\bar{k}_{n -1}}{\bf F}_0$ has the following form: \vspace{2mm}

\noindent $_{\bar{k}_{n-1}}{\bf F}_{0}$:  $_{\bar{k}_{n-1}}u_{\e}[jj', ii']  =   
\frac{_{\bar{k}_{n-1}}p_{\e, jj'} \, _{\bar{k}_{n-1}}v_{\e, j}^{-1}}{_{\bar{k}_{n-1}}p_{\e, ii'} \, _{\bar{k}_{n-1}}v_{\e, i}^{-1}} \to 
\, _{\bar{k}_{n-1}}u_{0}[jj', ii']  \in [0, \infty]$  as  $\e \to 0$, \\ 
\makebox[13mm]{} for  $j' \in \, _{\bar{k}_{n-1}}\XX, 
i' \in \, _{\bar{k}_{n-1}}\YY_{1, i}$, $j, i \in \, _{\bar{k}_{n-1}}\overline{\DD}$.
%j \neq i, 
\vspace{2mm}

According Lemma 12 (which should be applied to the semi-Markov processes $_{\bar{k}_{n-1}}\eta_\e(t)$   
and $_{\bar{k}_{n-1}}\tilde{\eta}_\e(t)$ instead of the semi-Markov processes $_{\bar{k}_0}\eta_\e(t)$ and $_{\bar{k}_0}\tilde{\eta}_\e(t)$)  condition  $_{\bar{k}_{n-1}}{\bf F}_0$ implies that condition $_{\bar{k}_{n-1}}{\bf \tilde{F}}$ holds.

Let us introduce, for $n  \geq 3$, condition:  \vspace{2mm}

\noindent ${\bf F}_{n-1}$: $u_{\e}[i_1 i_1', i_2 i_2'', i_3 i_3' i_3'', \ldots, i_{2^{n-1} +1} i_{2^{n-1}+1}' i_{2^{n-1}+1}'']$ 
$=  \frac{p_{\e, i_1i_1'}}{p_{\e, i_2 i_2''}} \prod_{l = 3}^{2^{n-1} + 1} \frac{p_{\e, i_l i_l'}}{p_{\e, i_l i_l''}} \frac{v_{\e, i_1}^{-1}}{v_{\e, i_2}^{-1}}$ 

\noindent \makebox[10mm]{} $\to u_{0}[i_1 i_1', i_2 i_2'', i_3 i_3' i_3'', \ldots, i_{2^{n-1} +1} i_{2^{n-1}+1}' i_{2^{n-1}+1}''] \in [0, \infty]$ as $\e \to 0$, for

\noindent \makebox[10mm]{}  $i'_1, i_l' \in \, \XX, i''_2, i_l'' \in \, \YY_{1, i_l}, i_1, i_2,  i_l \in \overline{\DD}, 
l = 3, \ldots, 2^{n-1} + 1$. 
%, \ i_1, i_2, i_3,  i_4 = i_5, i_6 = i_7 = i_8$ 
%\noindent \makebox[10mm]{}  $=  i_9, \ldots, i_{2^{n-2} + 2} = \cdots =   i_{2^{n-1} + 1} \in \overline{\DD}, \, 
%i_1 \neq i_2 \neq i_3 \neq i_4 \neq \cdots \neq i_{2^{n-2} + 2}$. 
\vspace{2mm}

{\bf Lemma 26}. {\em Condition ${\bf F}_{n-1}$ is sufficient for holding of condition 
$_{\bar{k}_{n-1}}{\bf F}_1$, for any $k_1, \ldots, k_{n-1} \in \overline{\DD}, k_1 \neq \cdots \neq k_{n-1}$.} 
\vspace{1mm}

The proof of Lemma 26 is analogous with the proof of Lemma 22. 

Lemma 26 and the above remarks imply that condition ${\bf F}_{n-1}$ is sufficient for holding of conditions 
condition $_{\bar{k}_{n-1}}{\bf \tilde{F}}$. 
Moreover, since condition ${\bf F}_{n-1}$ is stronger than condition ${\bf F}_{r-1}$, for $r = 1, \ldots, n-1$, condition ${\bf F}_{n-1}$ is, in fact, sufficient for holding of all conditions $_{\bar{k}_{r-1}}{\bf \tilde{F}}$
for $r = 1, \ldots, n-1$. 

The comparability condition ${\bf F}_{n-1}$ is expressed via transition probabilities $p_{\e, ij}$ and normalising functions $v_{\e, i}$ in much more explicit form than conditions $_{\bar{k}_{r-1}}{\bf \tilde{F}}$ for $r = 1, \ldots, n$. 

A disadvantage  of condition ${\bf F}_{n-1}$   requires existence of  a large number of limits. 

However,  this condition automatically hold,  if condition 
${\bf B}$ holds, and transition probabilities $p_{\e, ij},  \e \in (0, 1], j \in \YY_{1, i}, i \in \overline{\DD}$  and the normalisation functions $v_{\e, i},  \e \in (0, 1], i \in \overline{\DD}$ satisfy the following condition:  
\begin{itemize}
\item [${\bf H}$:] Functions $p_{\cdot, ij}, j \in \YY_{1, i}, i \in \overline{\DD}$ and $v_{\cdot, i}, i \in i \in \overline{\DD}$ belong to some complete family of asymptotically comparable functions.   
\end{itemize} 

{\bf Lemma 27}. {\em Condition ${\bf H}$ is sufficient for holding 
of condition ${\bf F}_{\bar{m} -2}$ and, thus,  for holding of conditions 
$_{\bar{k}_{r-1}}{\bf \tilde{F}}$, $r = 1, \ldots, \bar{m} -1$, for any $k_1, \ldots, k_{\bar{m}-1} \in \overline{\DD}, k_1 \neq \cdots \neq k_{\bar{m}-1}$.} 
\vspace{1mm} 

The proof of the lemma is given in Appendix A.

Here, we also would like to mention that only limits appearing in conditions $_{\bar{k}_{r-1}}{\bf \tilde{F}}$
for $r = 1, \ldots, n$ are actually  used for computing limiting Laplace transforms and expectations in conditions $_{\bar{k}_0}{\bf \tilde{D}}'$, $_{\bar{k}_r}{\bf D}'$, $_{\bar{k}_r}{\bf \tilde{E}},  r = 1, \ldots, n$. As far as limits in condition 
${\bf F}_{n-1}$ are concerned, their existence is sufficient for existence, but these limits itself are not directly used  for computing the  above mentioned limits appearing in conditions $_{\bar{k}_{r-1}}{\bf \tilde{F}}, r = 1, \ldots, n$.

In conclusion, we would like also to mention that the number of limits, which should be computed in every  above conditions  $_{\bar{k}_0}{\bf C}$, $_{\bar{k}_r}{\bf C}$, $_{\bar{k}_r}{\bf \tilde{C}}$, $_{\bar{k}_r}{\bf \tilde{C}}'$, $_{\bar{k}_r}{\bf \hat{C}}, r = 1, \ldots, n$ does not exceed $(m + \bar{m}) \bar{m}$ (recall that $m$ and $\bar{m}$ are  the numbers of states, respectively, in domain $\DD$ and domain $\overline{\DD}$). Also, the number of limits, which should be computed in every  above conditions $_{\bar{k}_{r-1}}{\bf \tilde{F}}, r = 1, \ldots, n$ doed not exceed $\bar{m}$.

Thus, for $n = \bar{m} -1$, the total number of limits in all conditions listed above  
do not exceed $(m + \bar{m} +1) \bar{m}^2 $. 

In the case, where condition ${\bf H}$ holds, the pre-limiting quantities in conditions $_{\bar{k}_0}{\bf C}$, $_{\bar{k}_r}{\bf C}$, $_{\bar{k}_r}{\bf \tilde{C}}$, $_{\bar{k}_r}{\bf \tilde{C}}'$, $_{\bar{k}_r}{\bf \hat{C}}, r = 1, \ldots, n$ and $_{\bar{k}_{r-1}}{\bf \tilde{F}}, r = 1, \ldots, n$ belong to the complete family of asymptotically comparable functions pointed in condition ${\bf H}$. 

The above algorithms became more effective, in the case where the asymptotic comparability condition ${\cal H}$ is based  on some concrete  complete family of asymptotically comparable functions  admitting effective computation of limits related to summation, multiplication and division operations. For example,  this relates to the above mentioned family ${\cal H}_1$, and families  ${\cal H}_2$ and  ${\cal H}_3$ described in Appendix A.
  
In this case, computations of limits in conditions $_{\bar{k}_0}{\bf C}$, $_{\bar{k}_r}{\bf C}$, $_{\bar{k}_r}{\bf \tilde{C}}$, $_{\bar{k}_r}{\bf \tilde{C}}'$, $_{\bar{k}_r}{\bf \hat{C}}, r = 1, \ldots, n$ and $_{\bar{k}_{r-1}}{\bf \tilde{F}}, r = 1, \ldots, n$ do require to perform the finite number of arithmetic operations. Moreover, due to the recurrent character of conditions $_{\bar{k}_r}{\bf C}$, $_{\bar{k}_r}{\bf \tilde{C}}$, $_{\bar{k}_r}{\bf \tilde{C}}'$, $_{\bar{k}_r}{\bf \hat{C}}$ and 
$_{\bar{k}_{r-1}}{\bf \tilde{F}}$, the numbers of operations  is the same for each of the above conditions, for every $r = 1, \ldots, n$. This implies that, in the case $n = \bar{m} -1$, the total number of operations required for computing all limits in conditions   $_{\bar{k}_0}{\bf C}$, $_{\bar{k}_r}{\bf C}$, $_{\bar{k}_r}{\bf \tilde{C}}$, $_{\bar{k}_r}{\bf \tilde{C}}'$, $_{\bar{k}_r}{\bf \hat{C}}, r = 1, \ldots, \bar{m} -1$ and $_{\bar{k}_{r-1}}{\bf \tilde{F}}, r = 1, \ldots, \bar{m} -1$  is  of the order $O(\bar{m}^3)$.   \\

{\bf 7 Weak Asymptotics for Distributions of Hitting Times} \\

In this section, we formulate and prove theorems about weak convergence of hitting times as well as  
describe the recurrent algorithm for computing the corresponding limiting Laplace transforms and 
normalisation functions. \vspace{1mm} 

{\bf 7.1 Recurrent relations for hitting times, their distributions and Laplace transforms}. Let $\bar{k}_{\bar{m}} = \langle k_1, k_2, \ldots, k_{\bar{m}} \rangle$ be an arbitrary sequence of different states from domain $\overline{\DD}$.

The following Lemmas 28 -- 33 are analogs of Lemmas  14 -- 16 and 18 -- 20. They  play the key role in the phase space reduction algorithm and getting recurrent weak convergence relations for hitting times.  

The proofs of Lemmas 28 -- 33  can be obtained by recurrent realisation of the following steps. At step ${\bf 0}$,  Lemmas 4 -- 6 should be applied to the semi-Markov processes $_{\bar{k}_0}\eta_{\e}(t)$ and $_{\bar{k}_0}\tilde{\eta}_{\e}(t)$. At step ${\bf 1}$, Lemmas 14 -- 16 should be applied to  the semi-Markov processes 
$_{\bar{k}_0}\tilde{\eta}_{\e}(t)$ and $_{\bar{k}_1}\eta_{\e}(t)$ and, then, Lemmas  18 -- 20 should be applied to the semi-Markov processes  $_{\bar{k}_1}\eta_{\e}(t)$ and $_{\bar{k}_1}\tilde{\eta}_{\e}(t)$. At step ${\bf 2}$, Lemmas 14 -- 16 should be applied to  the semi-Markov processes 
$_{\bar{k}_1}\tilde{\eta}_{\e}(t)$ and $_{\bar{k}_2}\eta_{\e}(t)$ and, then, Lemmas  18 -- 20 should be applied to the semi-Markov processes  $_{\bar{k}_2}\eta_{\e}(t)$ and $_{\bar{k}_2}\tilde{\eta}_{\e}(t)$, etc. Finally, 
at step ${\bf n}$, Lemmas 14 -- 16 should be applied to  the semi-Markov processes 
$_{\bar{k}_{n - 1}}\tilde{\eta}_{\e}(t)$ and $_{\bar{k}_n}\eta_{\e}(t)$ and, then, Lemmas  18 -- 20 should be applied to the semi-Markov processes  $_{\bar{k}_n}\eta_{\e}(t)$ and $_{\bar{k}_n}\tilde{\eta}_{\e}(t)$.
\vspace{1mm}

{\bf Lemma 28}. {\em Let conditions $_{\bar{k}_0}{\bf A}$, $_{\bar{k}_0}{\bf B}$ hold. Then, for every $n = 0, \ldots, \bar{m} -1$, conditions   $_{\bar{k}_n}{\bf A}$,  $_{\bar{k}_n}{\bf B}$, $_{\bar{k}_n}{\bf \tilde{A}}$, 
$_{\bar{k}_{n}}{\bf \tilde{B}}$ hold, and the following relations take place, for $\e \in (0, 1]$, 
\begin{align}\label{iderebaikokom}
& \PP_i \{ \tau_{\e, \DD} =  \, _{\bar{k}_n}\tau_{\e,  \DD}  =  \, _{\bar{k}_{n}}\tilde{\tau}_{\e, \DD} , 
\makebox[22mm]{} \vspace{2mm} \nonumber \\
& \quad  \     \eta_{\e}(\tau_{\e \DD}) =  \, _{\bar{k}_n}\eta_{\e}(_{\bar{k}_n}\tau_{\e, \DD}) =  \, _{\bar{k}_{n}}\tilde{\eta}_{\e}(_{\bar{k}_{n}}\tilde{\tau}_{\e \DD}) \}
= 1, i \in \, _{\bar{k}_n}\XX.    
\end{align}

{\bf Lemma 29}. {\em Let conditions $_{\bar{k}_0}{\bf A}$, $_{\bar{k}_0}{\bf B}$ hold. Then, for every $n = 0, \ldots, \bar{m} -1$, conditions   $_{\bar{k}_n}{\bf A}$,  $_{\bar{k}_n}{\bf B}$, $_{\bar{k}_n}{\bf \tilde{A}}$, 
$_{\bar{k}_n}{\bf \tilde{B}}$ hold, and the following relations take place, for $\e \in (0, 1]$,  
\begin{equation}\label{iderebanomaokom}
G_{\e, \DD,  ij}(t)  =  \, _{\bar{k}_n}G_{\e, \DD, ij}(t)  = \, _{\bar{k}_{n}}\tilde{G}_{\e, \DD,  ij}(t), 
\, t \geq 0,  j \in \DD, i \in \, _{\bar{k}_n}\XX.       
\end{equation}}
%\vspace{-2mm}

{\bf Lemma 30}. {\em Let conditions $_{\bar{k}_0}{\bf A}$, $_{\bar{k}_0}{\bf B}$ hold. Then, for every $n = 0, \ldots, \bar{m} -1$, conditions   $_{\bar{k}_n}{\bf A}$,  $_{\bar{k}_n}{\bf B}$, $_{\bar{k}_{n}}{\bf \tilde{A}}$, $_{\bar{k}_{n}}{\bf \tilde{B}}$ hold, and the following relations take place, for $\e \in (0, 1]$,  
 \begin{equation}\label{mokubadaokuom}
\Psi_{\e, \DD, ij}(s)  =   \,  _{\bar{k}_{n}}\Psi_{\e, \DD, ij}(s) = 
\, _{\bar{k}_{n}}\tilde{\Psi}_{\e, \DD, ij}(s), \, s \geq 0,  j \in \DD, i \in \, _{\bar{k}_n}\XX. 
\end{equation}
\vspace{1mm} 

{\bf Lemma 31}. {\em Let conditions $_{\bar{k}_0}{\bf A}$, $_{\bar{k}_0}{\bf B}$ hold. Then, for every $n = 1, \ldots, \bar{m} -1$, conditions 
$_{\bar{k}_{n-1}}{\bf \tilde{A}}$, $_{\bar{k}_{n-1}}{\bf \tilde{B}}$, $_{\bar{k}_n}{\bf A}$,  $_{\bar{k}_n}{\bf B}$ hold, and the following relations take place,  for  $\e \in (0, 1]$,
\begin{align}\label{iderebaikom}
& \PP_i \{ \tau_{\e, \DD} =  \, _{\bar{k}_{n-1}}\tilde{\tau}_{\e, \DD} =  \, _{\bar{k}_n}\tau_{\e,  \DD},  \makebox[18mm]{} \vspace{2mm} \nonumber \\
& \quad \   \eta_{\e}(\tau_{\e \DD}) = \, _{\bar{k}_{n-1}}\tilde{\eta}_{\e}(_{\bar{k}_{n-1}}\tilde{\tau}_{\e \DD}) = \, _{\bar{k}_n}\eta_{\e}(_{\bar{k}_n}\tau_{\e, \DD})  \} = 1, \, j \in \DD,  i \in \, _{\bar{k}_n}\XX,  
\end{align}
and
\begin{align}\label{iderebasikom}
& \quad \PP_{k_n} \{ \tau_{\e, \DD} 
 =  \,  _{\bar{k}_{n-1}}\tilde{\tau}_{\e, \DD} = \, _{\bar{k}_{n-1}}\tilde{\kappa}_{\e,  1} {\rm I}(_{\bar{k}_{n-1}}\tilde{\eta}_{\e, 1} \in \DD) \makebox[50mm]{}
 \vspace{2mm} \nonumber \\
&  \quad \quad   \quad \quad
+  (_{\bar{k}_{n-1}}\tilde{\kappa}_{\e, 1}   
+ \, _{\bar{k}_n}\tau_{\e, \DD})  {\rm I}(_{\bar{k}_{n-1}}\tilde{\eta}_{\e, 1} \in \, _{\bar{k}_n}\overline{\DD}), 
\vspace{2mm} \nonumber \\
&  \quad \quad \quad 
\eta_{\e}(\tau_{\e, \DD}) = \, _{\bar{k}_{n-1}}\tilde{\eta}_{\e}(_{\bar{k}_{n-1}}\tilde{\tau}_{\e, \DD}) = 
\, _{\bar{k}_{n-1}}\tilde{\eta}_{\e, 1}{\rm I}(_{\bar{k}_{n-1}}\tilde{\eta}_{\e, 1} \in \DD) \vspace{2mm} \nonumber \\
&  \quad \quad   \quad \quad
+  \,  _{\bar{k}_n}\eta_{\e}( _{\bar{k}_n}\tau_{\e,  \DD})
{\rm I}(_{\bar{k}_{n-1}}\tilde{\eta}_{\e, 1} \in \, _{\bar{k}_n}\overline{\DD}) \} = 1, \, j \in \DD.  
\end{align}}

{\bf Lemma 32}. {\em Let conditions $_{\bar{k}_0}{\bf A}$, $_{\bar{k}_0}{\bf B}$ hold. Then, for every $n = 1, \ldots, \bar{m} -1$, conditions $_{\bar{k}_{n-1}}{\bf \tilde{A}}$, $_{\bar{k}_{n-1}}{\bf \tilde{B}}$, $_{\bar{k}_n}{\bf A}$,  $_{\bar{k}_n}{\bf B}$ hold, and the following relations take place, for  $\e \in (0, 1]$, 
\begin{align}\label{iderebanomaom}
G_{\e, \DD,  ij}(t)  & = \, _{\bar{k}_{n-1}}\tilde{G}_{\e, \DD,  ij}(t) \vspace{1mm} \nonumber \\
& =  \, _{\bar{k}_n}G_{\e, \DD, ij}(t), \, t \geq 0,  j \in \DD, i \in \, _{\bar{k}_n}\XX,   
\end{align}
and 
\begin{align}\label{iderebasnomom}
& \quad G_{\e, \DD, k_n j}(t)  = \, _{\bar{k}_{n-1}}\tilde{G}_{\e, \DD, k_n j}(t) = \, _{\bar{k}_{n-1}}\tilde{F}_{\e, k_n j}(t) \, 
_{\bar{k}_{n-1}}\tilde{p}_{\e,  k_n j}  \makebox[31mm]{} \vspace{3mm} \nonumber \\
&  \quad \quad \quad   + \sum_{r \in \, _{\bar{k}_n}\overline{\DD}} \, _{\bar{k}_{n-1}}\tilde{F}_{\e, k_n r}(t) * \,  _{\bar{k}_n}G_{\e, \DD, rj}(t) \, _{\bar{k}_{n-1}}\tilde{p}_{\e, k_n r} \vspace{2mm} \nonumber \\
&  \quad \quad  = \, _{\bar{k}_{n-1}}\tilde{F}_{\e, k_n j}(t) \, _{\bar{k}_{n-1}}\tilde{p}_{\e,  k_n j}  \vspace{2mm} \nonumber \\
& \quad \quad \quad    +  \sum_{r \in \, _{\bar{k}_n}\overline{\DD}} \, _{\bar{k}_{n-1}}\tilde{F}_{\e, k_n r}(t) * G_{\e, \DD, rj}(t) \, _{\bar{k}_{n-1}}\tilde{p}_{\e, k_n r}, \, t \geq 0,  j \in \DD. 
\end{align}}

{\bf Lemma 33}. {\em Let conditions $_{\bar{k}_0}{\bf A}$, $_{\bar{k}_0}{\bf B}$ hold. Then, for every $n = 1, \ldots, \bar{m} -1$, conditions 
${\bf _{\bar{k}_{n-1}}\tilde{A}}$, $_{\bar{k}_{n-1}}{\bf \tilde{B}}$, $_{\bar{k}_n}{\bf A}$,  $_{\bar{k}_n}{\bf B}$ hold,  and the following relations take place,  for $\e \in (0, 1]$,
 \begin{align}\label{mokubadaom}
\Psi_{\e, \DD, ij}(s) & =  \, _{\bar{k}_{n-1}}\tilde{\Psi}_{\e, \DD, ij}(s) \vspace{1mm} \nonumber \\
& = \,  _{\bar{k}_{n}}\Psi_{\e, \DD, ij}(s), \, s \geq 0,  j \in \DD, i \in \, _{\bar{k}_n}\XX, 
\end{align} 
and
\begin{align}\label{iderenasfatom}
& \ \ \Psi_{\e, \DD, k_n j}(s)  = \, _{\bar{k}_{n-1}}\tilde{\Psi}_{\e, \DD,  k_n j}(s)  =  \, _{\bar{k}_{n-1}}\tilde{\phi}_{\e, k_n j}(s) \, _{\bar{k}_{n-1}}\tilde{p}_{\e, k_n j} 
 \makebox[27mm]{} \vspace{3mm} \nonumber \\ 
&   \quad \quad \quad \quad  + \sum_{r \in \, _{\bar{k}_{n}}\overline{\DD}} \, 
_{\bar{k}_{n}}\Psi_{\e, \DD, rj}(s) \, _{\bar{k}_{n-1}}\tilde{\phi}_{\e, k_n r}(s) \,  _{\bar{k}_{n-1}}\tilde{p}_{\e,  k_n r} 
\vspace{2mm} \nonumber \\ 
&  \quad \quad \quad   =  \, _{\bar{k}_{n-1}}\tilde{\phi}_{\e, k_n j}(s) \, _{\bar{k}_{n-1}}\tilde{p}_{\e, k_n j}  \vspace{2mm} \nonumber \\ 
&  \quad \quad \quad \quad  
+ \sum_{r \in \, _{\bar{k}_{n}}\overline{\DD}} \Psi_{\e, \DD, rj}(s)  \, _{\bar{k}_{n-1}}\tilde{\phi}_{\e, k_n r}(s) \, _{\bar{k}_{n-1}}\tilde{p}_{\e,  k_n r}, \, s \geq 0,  j \in \DD. 
\end{align}} 

In what follows, we assume that the asymptotic recurrent algorithm of phase space reduction described in Section 6 is realised for $n = \bar{m} -1$. 

In the following Subsections 7.2 -- 7.8, we assume that conditions ${\bf A}$ -- ${\bf E}$, ${\bf \tilde{C}}$, $_{\bar{k}_{r}}{\bf \tilde{C}}$, $_{\bar{k}_{r}}{\bf \hat{C}}$, $_{\bar{k}_{r-1}}{\bf \tilde{F}}$, $r = 1, \ldots, \bar{m} -1$ hold, where the   
states $k_r \in \, _{\bar{k}_{r-1}} \overline{\DD}^*, r = 1, \ldots,  \bar{m} -1$ are chosen in such way that condition
$_{\bar{k}_{r}}{\bf \hat{F}}$ holds, for  $r = 1, \ldots,  \bar{m} -1$. 

In this case, domain $_{\bar{k}_{\bar{m}-1}}\overline{\DD} = 
\overline{\DD} \setminus \{ k_1, \ldots, k_{\bar{m} -1} \} = \{ k_{\bar{m}} \}$ is a one-state set, which can be referred as one of the most absorbing states in domain $\overline{\DD}$. 

In what follows, we denote as $\bar{k}_m = \langle k_1, \ldots, k_m \rangle$ the sequence of states  constructed with the use of the above algorithm,  and use for the corresponding final normalisation function for hitting times 
$\tau_{\e, \DD}$ used for the case, where an initial state is $k_i \in  \overline{\DD}$, the notation 
$_{\bar{k}_m}\check{v}_{\e, k_i}$.  \vspace{1mm}

{\bf 7.2 Weak asymptotic for distributions of hitting times for the case with the most absorbing initial state}.  Let us now consider the case, where the initial state  is $k_{\bar{m}}$. 

The reduced semi-Markov processes  $_{\bar{k}_{\bar{m}-1}}\eta_\e(t)$ and $_{\bar{k}_{\bar{m}
-1}}\tilde{\eta}_\e(t)$ have the phase space $_{\bar{k}_{\bar{m}-1}}\XX = \XX \setminus \{k_1, \ldots, k_{\bar{m} -1} \} = \DD \cup \{ k_{\bar{m}} \}$.

Obviously, the hitting time $_{\bar{k}_{\bar{m}-1}}\tilde{\tau}_{\e, \DD} = \, _{\bar{k}_{\bar{m}-1}}\tilde{\kappa}_{\e, 1}$ and 
$_{\bar{k}_{\bar{m}-1}}\tilde{\eta}_{\e}(_{\bar{k}_{\bar{m}-1}}\tilde{\tau}_{\e, \DD}) = \, _{\bar{k}_{\bar{m}-1}}\tilde{\eta}_{\e, 1}$, 
if $\eta_\e(0) = k_{\bar{m}}$.

Lemma 30 implies the following relation takes place,
\begin{equation}\label{cotyren}
G_{\e, \DD, k_{\bar{m}} j}(t)  = \, _{\bar{k}_{\bar{m} -1}}\tilde{F}_{\e, k_{\bar{m}} j}(t) \, _{\bar{k}_{\bar{m}-1}}\tilde{p}_{\e, k_{\bar{m}}j}, t \geq 0, j \in \DD.
\end{equation}

The corresponding final normalisation function takes the following form,
\begin{equation}\label{norm}
 _{\bar{k}_{\bar{m}}}\check{v}_{\e, k_{\bar{m}}} =  
\, _{\bar{k}_{\bar{m}-1}}\tilde{v}_{\e, k_{\bar{m}}} = \prod_{l = 0}^{\bar{m} -1}(1 - \, _{\bar{k}_l}p_{\e, k_{\bar{m}} k_{\bar{m}}})^{-1} v_{\e, k_{\bar{m}}}.
\end{equation} 

The  following theorem, which is the direct  corollary of Lemma 23 and relation (\ref{cotyren}),  takes place. \vspace{1mm}

{\bf Theorem 4}. {\em Let conditions ${\bf A}$ -- ${\bf E}$, ${\bf \tilde{C}}$, $_{\bar{k}_{r}}{\bf \tilde{C}}$, $_{\bar{k}_{r}}{\bf \hat{C}}$, $_{\bar{k}_{r-1}}{\bf \tilde{F}}$, $r = 1, \ldots, \bar{m} -1$ hold, where the   
states $k_r \in \, _{\bar{k}_{r-1}} \overline{\DD}^*, r = 1, \ldots,  \bar{m} -1$ are chosen in such way that condition
$_{\bar{k}_{r}}{\bf \hat{F}}$ holds, for  $r = 1, \ldots,  \bar{m} -1$. Then,  the following asymptotic 
relation takes place, for $j \in \DD$, 
\begin{align}\label{cotyrevaner}
G_{\e, \DD, k_{\bar{m}} j}(\cdot \, _{\bar{k}_{\bar{m}}} \check{v}_{\e, k_{\bar{m}} }) & \Rightarrow  G_{0, \DD, k_{\bar{m}} j}(\cdot) \vspace{1mm} \nonumber \\ 
& = \, _{\bar{k}_{\bar{m}-1}}\tilde{F}_{0, k_{\bar{m}}j}(\cdot) \, _{\bar{k}_{\bar{m}-1}}\tilde{p}_{0, k_{\bar{m}}j} \ {\rm as} \ \e \to 0. 
\end{align}} 
\makebox[3mm]{} Note that, in this case,  the limiting  distribution $G_{0, \DD, k_{\bar{m}} j}(\cdot)$ has, for every $j \in \DD$,  the Laplace transform,
\begin{align}\label{gasd} 
\Psi_{0, \DD, k_{\bar{m}} j}(s) & =  \int_0^\infty e^{-st} G_{0, \DD, k_{\bar{m}} j}(dt) \vspace{2mm} \nonumber \\ 
& = \tilde{\phi}_{0, k_{\bar{m}} j}(s)\tilde{p}_{0, k_{\bar{m}}j}, s \geq 0.
\end{align}

{\bf 7.3 Weak asymptotic for distributions of hitting times for the case with the second most absorbing initial state}. 
Let us now consider the case, where the initial 
state  is $k_{\bar{m} -1}$. 

The reduced semi-Markov processes  $_{\bar{k}_{\bar{m}-2}}\eta_\e(t)$ and $_{\bar{k}_{\bar{m}
-2}}\tilde{\eta}_\e(t)$ have the phase space $_{\bar{k}_{\bar{m}-2}}\XX = \XX \setminus \{k_1, \ldots, k_{\bar{m} -2} \} = 
\DD \cup \{ k_{\bar{m}-1}, k_{\bar{m}} \}$, while domain  $_{\bar{k}_{\bar{m}-2}}\bar{\DD} = 
\overline{\DD} \setminus \{k_1, \ldots, k_{\bar{m} -2} \} = \{ k_{\bar{m} -1},  k_{\bar{m}} \}$ is a two-states set. 

We can use, in this case, relation (\ref{iderenasfatom}) given in Lemma 33. This relation takes, in this case, the following form, for $s \geq 0, j \in \DD$, 
\begin{align}\label{sfatomnega}
\Psi_{\e, \DD, k_{\bar{m} -1} j}(s) & =   \, _{\bar{k}_{\bar{m} -2}}\tilde{\phi}_{\e, k_{\bar{m} - 1} j}(s) \, 
_{\bar{k}_{\bar{m} - 2}}\tilde{p}_{\e, k_{\bar{m} - 1} j}  \vspace{2mm} \nonumber \\ 
&  \quad   + \Psi_{\e, \DD, k_{\bar{m}} j}(s)   \, _{\bar{k}_{\bar{m} -2}}\tilde{\phi}_{\e, k_{\bar{m} - 1} 
k_{\bar{m}}}(s) \, _{\bar{k}_{\bar{m} - 2}}\tilde{p}_{\e,  k_{\bar{m} - 1} k_{\bar{m}}}. 
\end{align} 

If the limiting probability $\, _{\bar{k}_{\bar{m} - 2}}\tilde{p}_{0,  k_{\bar{m} - 1} k_{\bar{m}}} > 0$, the corresponding final normalisation function take the form,
\begin{equation}\label{norma}
 _{\bar{k}_{\bar{m}}}\check{v}_{\e, k_{\bar{m} -1}} =  
\, _{\bar{k}_{\bar{m}-1}}\tilde{v}_{\e, k_{\bar{m} }} = \prod_{l = 0}^{\bar{m} -1}(1 - \, _{\bar{k}_l}p_{\e, k_{\bar{m}} k_{\bar{m}}})^{-1} v_{\e, k_{\bar{m}}}.
\end{equation} 

Relation (\ref{cotyrevaner}) given in Theorem 4 and relation  (\ref{sfatomnega})  imply that the following relation 
takes place,  for $s \geq 0, j \in \DD$, 
\begin{align}\label{idersfatnen}
& \Psi_{\e, \DD, k_{\bar{m} -1} j}(s / \, _{\bar{k}_{\bar{m}-1}}\tilde{v}_{\e, k_{\bar{m}}})  
\vspace{2mm} \nonumber \\
& \quad =  \, _{\bar{k}_{\bar{m} -2}}\tilde{\phi}_{\e, k_{\bar{m} - 1} j}((1 - \, _{\bar{k}_{\bar{m} - 1}}p_{\e, k_{\bar{m}} k_{\bar{m}}})
\frac{_{\bar{k}_{\bar{m}-2}}\tilde{v}_{\e, k_{\bar{m} -1}}}{_{\bar{k}_{\bar{m}-2}}\tilde{v}_{\e, k_{\bar{m}}}} s/ \, _{\bar{k}_{\bar{m}-2}}\tilde{v}_{\e, k_{\bar{m} -1}})  \vspace{2mm} \nonumber \\ 
& \quad \quad  \times \, _{\bar{k}_{\bar{m}-2}}\tilde{p}_{\e, k_{\bar{m} - 1} j} 
\vspace{2mm} \nonumber \\
& \quad \quad + \Psi_{\e, \DD, \, _{\bar{k}_{\bar{m}}} j}(s / \, _{\bar{k}_{\bar{m}-1}}\tilde{v}_{\e, k_{\bar{m}}}) \vspace{2mm} \nonumber \\
%\end{align*}
%\begin{align}
& \quad \quad \times \, _{\bar{k}_{\bar{m} -2}}\tilde{\phi}_{\e,  k_{\bar{m} -1} k_{\bar{m}}}((1 - \, _{\bar{k}_{\bar{m} - 1}}p_{\e, k_{\bar{m}} k_{\bar{m}}})
\frac{_{\bar{k}_{\bar{m}-2}}\tilde{v}_{\e, k_{\bar{m} -1}}}{_{\bar{k}_{\bar{m}-2}}\tilde{v}_{\e, k_{\bar{m}}}} s/ \, _{\bar{k}_{\bar{m}-2}}\tilde{v}_{\e, k_{\bar{m} -1}} ) \vspace{2mm} \nonumber \\
& \quad \quad \times 
\, _{\bar{k}_{\bar{m} - 2}}\tilde{p}_{\e,  k_{\bar{m} - 1} k_{\bar{m}}} 
\vspace{2mm} \nonumber \\
&  \quad \to \, _{\bar{k}_{\bar{m} -2}}\tilde{\phi}_{0,   k_{\bar{m} -1} j}((1 - \, _{\bar{k}_{\bar{m} - 1}}p_{0, k_{\bar{m}} k_{\bar{m}}}) \, _{k_{\bar{m} - 2}}\tilde{w}_{0, k_{\bar{m} - 1} k_{\bar{m}}}s)  \, _{\bar{k}_{\bar{m}-2}}\tilde{p}_{0, k_{\bar{m} - 1} j} 
\vspace{2mm} \nonumber \\
& \quad \quad + \Psi_{0, \DD, \, _{\bar{k}_{\bar{m}}} j}(s) \vspace{2mm} \nonumber \\
& \quad \quad  \times \, _{\bar{k}_{\bar{m} -2}}\tilde{\phi}_{0,  k_{\bar{m} -1} k_{\bar{m}}}((1 - \, _{\bar{k}_{\bar{m} - 1}}p_{0, k_{\bar{m}} k_{\bar{m}}}) \, _{k_{\bar{m} - 2}}\tilde{w}_{0, k_{\bar{m} - 1} k_{\bar{m}}}s) \vspace{2mm} \nonumber \\
& \quad \quad \times \, _{\bar{k}_{\bar{m} - 2}}\tilde{p}_{0,  k_{\bar{m} - 1} k_{\bar{m}}} \vspace{2mm} \nonumber \\
& \quad = \Psi_{0, \DD,  k_{\bar{m} -1} j}(s) \ {\rm as} \ \e \to 0.
\end{align}

In this case, the corresponding  limiting distribution $G_{0, \DD, k_{\bar{m} -1} j}(\cdot)$ has, for every $j \in \DD$,  the Laplace transform $\Psi_{0, \DD, k_{\bar{m} -1} j}(s) = 
\int_0^\infty e^{-st} G_{0, \DD, k_{\bar{m} -1} j}(dt), s \geq 0$ given by relation (\ref{idersfatnen}). 

If $_{\bar{k}_{\bar{m} - 2}}\tilde{p}_{0,  k_{\bar{m} - 1} k_{\bar{m}}} = 0$, the corresponding final normalisation function take the form,
\begin{equation}\label{normal}
 _{\bar{k}_{\bar{m}}}\check{v}_{\e, k_{\bar{m} -1}} =  
\, _{\bar{k}_{\bar{m}-2}}\tilde{v}_{\e, k_{\bar{m}- 1}} = \prod_{l = 0}^{\bar{m} -2}(1 - \, 
_{\bar{k}_l}p_{\e, k_{\bar{m}-1} k_{\bar{m} -1}})^{-1} v_{\e, k_{\bar{m}-1}}.
\end{equation} 

Since, probabilities $_{\bar{k}_{\bar{m} - 2}}\tilde{p}_{\e,  k_{\bar{m} - 1} k_{\bar{m}}} \to \, _{\bar{k}_{\bar{m} - 2}}\tilde{p}_{0,  k_{\bar{m} - 1} k_{\bar{m}}}   = 0$ as $\e \to 0$, relation (\ref{cotyrevaner}) given in Theorem 4 and relation (\ref{sfatomnega}) imply that the following relation 
takes place,  for $s \geq 0, j \in \DD$, 
\begin{align*}
& \Psi_{\e, \DD, k_{\bar{m} -1} j}(s / \, _{\bar{k}_{\bar{m}-2}}\tilde{v}_{\e, k_{\bar{m}-1}})  
\vspace{2mm} \nonumber \\
& \quad \quad =  \, _{\bar{k}_{\bar{m} -2}}\tilde{\phi}_{\e, k_{\bar{m} - 1} j}( s/ \, _{\bar{k}_{\bar{m}-2}}\tilde{v}_{\e, k_{\bar{m} -1}})  \, _{\bar{k}_{\bar{m}-2}}\tilde{p}_{\e, k_{\bar{m} - 1} j} 
\vspace{2mm} \nonumber \\
& \quad \quad \quad + \Psi_{\e, \DD, \, _{\bar{k}_{\bar{m}}} j}(s / \, _{\bar{k}_{\bar{m}-2}}\tilde{v}_{\e, k_{\bar{m}-1}}) \vspace{2mm} \nonumber \\
\end{align*}
\begin{align}\label{idersfa}
& \quad \quad \quad \quad \times \, _{\bar{k}_{\bar{m} -2}}\tilde{\phi}_{\e,  k_{\bar{m} -1} k_{\bar{m}}}( s/ \, _{\bar{k}_{\bar{m}-2}}\tilde{v}_{\e, k_{\bar{m} -1}} ) 
\, _{\bar{k}_{\bar{m} - 2}}\tilde{p}_{\e,  k_{\bar{m} - 1} k_{\bar{m}}} 
\vspace{2mm} \nonumber \\
&\quad  \quad \to \, _{\bar{k}_{\bar{m} -2}}\tilde{\phi}_{0, k_{\bar{m} - 1} j}( s)  \, _{\bar{k}_{\bar{m}-2}}\tilde{p}_{0, k_{\bar{m} - 1} j} \vspace{2mm} \nonumber \\
& \quad \quad = \Psi_{0, \DD,  k_{\bar{m} -1} j}(s) \ {\rm as} \ \e \to 0.
\end{align}

In this case, the corresponding  limiting distribution $G_{0, \DD, k_{\bar{m} -1} j}(\cdot)$ has, for every $j \in \DD$,  the Laplace 
transform $\Psi_{0, \DD, k_{\bar{m} -1} j}(s) = \int_0^\infty e^{-st} G_{0, \DD, k_{\bar{m} -1} j}(dt), s \geq 0$ given by 
relation (\ref{idersfa}).  
 
The  following theorem takes place. \vspace{1mm}

{\bf Theorem 5}. {\em Let conditions ${\bf A}$ -- ${\bf E}$, ${\bf \tilde{C}}$, $_{\bar{k}_{r}}{\bf \tilde{C}}$, $_{\bar{k}_{r}}{\bf \hat{C}}$, $_{\bar{k}_{r-1}}{\bf \tilde{F}}$, $r = 1, \ldots, \bar{m} -1$ hold, where the   
states $k_r \in \, _{\bar{k}_{r-1}} \overline{\DD}^*, r = 1, \ldots,  \bar{m} -1$ are chosen in such way that condition
$_{\bar{k}_{r}}{\bf \hat{F}}$ holds, for  $r = 1, \ldots,  \bar{m} -1$. Then,  the following asymptotic 
relation takes place, for $j \in \DD$, 
\begin{equation}\label{cotyrevanert}
G_{\e, \DD, k_{\bar{m}-1} j}(\cdot \, _{\bar{k}_{\bar{m}}} \check{v}_{\e, k_{\bar{m}-1} }) 
\Rightarrow  G_{0, \DD, k_{\bar{m}- 1}j}(\cdot) \ {\rm as} \ \e \to 0. 
\end{equation}} 
\makebox[3mm]{} According Theorems 4 and 5,  the normalisation function for hitting times $\tau_{\e, \DD}$ is the same  
for both cases, where initial state  is $k_{\bar{m}}$ or $k_{\bar{m}-1}$, if probability $_{\bar{k}_{\bar{m} - 2}}\tilde{p}_{0,  k_{\bar{m} - 1} k_{\bar{m}}} > 0$. However, the normalisation functions for hitting times $\tau_{\e, \DD}$ can differ 
for  cases, where initial state  is $k_{\bar{m}}$ or $k_{\bar{m}-1}$, if probability $_{\bar{k}_{\bar{m} - 2}}\tilde{p}_{0,  k_{\bar{m} - 1} k_{\bar{m}}} = 0$. \vspace{1mm}

{\bf 7.4 Weak asymptotics for distributions of hitting times for the case with an arbitrary initial state from sequence 
$\bar{k}_{\bar{m}}$}. Asymptotic relations (\ref{cotyrevaner}) and (\ref{cotyrevanert}) given, respectively, in Theorems 4 and 5 can be considered as the results of  first and second steps in some backward asymptotic recurrent algorithm for computing weak limits for distributions of hitting times $G_{\e, \DD, k_n  j}(\cdot \, _{\bar{k}_{\bar{m}}}\check{v}_{\e, k_{n} })$, for $n = \bar{m}, \ldots, 1$. 

In what follows, we can assume that $\bar{m} \geq 3$. 

Let us assume that, we have already realised $\bar{m} - n$ steps in this backward algorithm resulted by the following asymptotic relations, for $j \in \DD, l = \bar{m}, \ldots, n +1$,
\begin{equation}\label{resulna}
G_{\e, \DD, k_l j}(\cdot \, _{\bar{k}_{\bar{m}}}\check{v}_{\e, k_l})
\Rightarrow  G_{0, \DD, k_l j}(\cdot) \ {\rm as} \ \e \to 0. 
\end{equation}
with some normalisation functions $_{\bar{k}_{\bar{m}}}\check{v}_{\e, k_{l}}$, (given by relations (\ref{norm}), for $l = \bar{m}$, or
by relations (\ref{norma}) and (\ref{normal}), for $l = \bar{m} -1$) and limiting distributions $G_{0, \DD, k_{l} j}(t), t \geq 0,  j \in \DD$, which Laplace transforms  
$\Psi_{0, \DD, k_{l} j}(s)$ $= \int_0^\infty e^{-st} G_{0, \DD, k_{l} j}(dt), s \geq 0,  j \in \DD$ (given by relation (\ref{gasd}), for $l = \bar{m}$ or relations (\ref{idersfatnen}) and (\ref{idersfa}), for $l = \bar{m} -1$).

The assumptions  given in the form of relations (\ref{resulna}) can be also expressed in the equivalent form of the following relations expressed in terms of the corresponding Laplace transforms and assumed to hold, for  $j \in \DD, l = \bar{m}, \ldots, n +1$,
 \begin{equation}\label{resul}
\Psi_{\e, \DD, k_{l} j}(s / \, _{\bar{k}_{\bar{m}}}\check{v}_{\e, k_{l} })
\Rightarrow  \Psi_{0, \DD, k_{l} j}(s) \ {\rm as} \ \e \to 0, \ {\rm for} \ s \geq 0. 
\end{equation}

In what follows, we  consider the case, where the initial state  is $k_{n}$, for some  $n \leq  \bar{m} - 2$.

The reduced semi-Markov processes  $_{\bar{k}_{n -1}}\eta_\e(t)$ and 
$_{\bar{k}_{n -1}}\tilde{\eta}_\e(t)$ have the phase space $_{\bar{k}_{n -1}}\XX = \XX \setminus \{k_1, \ldots, k_{n -1} \} = 
\DD \cup \{ k_{n}, \ldots, k_{\bar{m}} \}$, while domain  $_{\bar{k}_{n -1}}\bar{\DD} = 
\overline{\DD} \setminus \{k_1, \ldots, k_{n - 1} \} = \{ k_{n},  \ldots, k_{\bar{m}} \}$. 

We can use, in this case, relation (\ref{iderenasfatom}) given in Lemma 28. This relation takes, in this case, 
can be written in the following form, for $s \geq 0, j \in \DD$, 
\begin{align}\label{sfatomne}
\Psi_{\e, \DD, k_{n} j}(s) & =   \, _{\bar{k}_{n -1 }}\tilde{\phi}_{\e, k_{n} j}(s) \, 
_{\bar{k}_{n -1}}\tilde{p}_{\e, k_{n} j}  + \sum_{r \in \, _{\bar{k}_{n}}\overline{\DD} }\Psi_{\e, \DD, r j}(s)  \vspace{2mm} \nonumber \\ 
&  \quad \quad \times   \, _{\bar{k}_{n -1}}
\tilde{\phi}_{\e, k_{n} r}(s) \, _{\bar{k}_{n -1}}\tilde{p}_{\e,  k_{n} r} 
\vspace{2mm} \nonumber \\
%\end{align*} 
%\begin{align}
& =   \, _{\bar{k}_{n -1}}\tilde{\phi}_{\e, k_{n} j}(s) \, _{\bar{k}_{n -1}}\tilde{p}_{\e, k_{n} j}  
+ \sum_{n +1 \leq l \leq \bar{m}}\Psi_{\e, \DD, k_{l} j}(s)  
\vspace{2mm} \nonumber \\ 
&  \quad \quad \times   \, _{\bar{k}_{n -1}}
\tilde{\phi}_{\e, k_{n} k_{l}}(s) \, _{\bar{k}_{n -1}}\tilde{p}_{\e,  k_{n} k_{l}}.  
\end{align} 

Note that Lemma 23 and the above assumptions imply that probabilities  
$_{\bar{k}_{n -1}}\tilde{p}_{\e,  k_{n} k_{l}}  \to 0$ as $\e \to 0$,  if 
$_{\bar{k}_{n -1}}\tilde{p}_{0,  k_{n} k_{l}}  = 0$.

The corresponding normalisation functions take the following forms,  for $n =  
\bar{m}, \ldots, 1$,
\begin{align}\label{normaba}
_{\bar{k}_{\bar{m}}}\check{v}_{\e, k_{n}}  & =  
\, _{\bar{k}_{\bar{n}(\bar{k}_{\bar{m}}, n) -1}}\tilde{v}_{\e, k_{\bar{n}(\bar{k}_{\bar{m}}, n)}} \vspace{2mm} \nonumber \\ 
&  = \prod_{q = 0}^{\bar{n}(\bar{k}_{\bar{m}}, n) -1}(1 - \, 
_{\bar{k}_q}p_{\e, k_{\bar{n}(\bar{k}_{\bar{m}}, n)} k_{\bar{n}(\bar{k}_{\bar{m}}, n)}})^{-1}  v_{\e, k_{\bar{n}(\bar{k}_{\bar{m}}, n)}},
\end{align} 
where
\begin{equation}\label{erty}
\bar{n}(\bar{k}_{\bar{m}}, n) = \left \{
\begin{array}{ll}
\bar{n} & \ \text{if there exists} \ n +1 \leq \bar{n} \leq \bar{m} \ \text{such that} \\
&  \ _{\bar{k}_{n -1}}\tilde{p}_{0,  k_{n} k_{l}}  = 0,  \bar{n} + 1 \leq l \leq \bar{m}, \vspace{2mm} \\ 
&  \  \text{and} \ _{\bar{k}_{n -1}}\tilde{p}_{0,  k_{n} k_{\bar{n}}}   > 0,   \vspace{2mm} \\
n  & \ \text{if} \ _{\bar{k}_{n -1}}\tilde{p}_{0,  k_{n} k_{l}}  = 0, \  n +1 \leq l \leq \bar{m}.
\end{array}
\right.
\end{equation}

First, let us assume that,
\begin{equation}\label{assump}
\bar{n}(\bar{k}_{\bar{m}}, n)  = \bar{n}, \ {\rm where} \  n +1 \leq \bar{n} \leq \bar{m}.
 \end{equation}

An analogue of relation (\ref{idersfatnen})
takes the following form, for $s \geq 0, j \in \DD$, 
\begin{align}\label{sfatomnema}
& \Psi_{\e, \DD, k_{n} j}(s /  _{\bar{k}_{\bar{m}}}\check{v}_{\e, k_{n}} ) \makebox[80mm]{}
\vspace{2mm} \nonumber \\ 
& \quad =   \, _{\bar{k}_{n -1}}\tilde{\phi}_{\e, k_{n} j}(\prod_{q = n}^{\bar{n} -1}(1 - \, 
_{\bar{k}_q}p_{\e, k_{\bar{n}} k_{\bar{n}}}) 
\frac{_{\bar{k}_{n -1}}\tilde{v}_{\e, k_{n}}}{_{\bar{k}_{n - 1}}\tilde{v}_{\e, k_{\bar{n}}}}s / 
_{\bar{k}_{n -1}}\tilde{v}_{\e, k_{n}}) \, _{\bar{k}_{n -1}}\tilde{p}_{\e, k_{n} j}  
\vspace{2mm} \nonumber \\ 
&  \quad  \quad + \sum_{\bar{n} + 1 \leq l \leq \bar{m}}
\Psi_{\e, \DD, k_{l} j}(s / \, _{\bar{k}_{\bar{m}}}\check{v}_{\e, k_{n}})  
 \, _{\bar{k}_{n -1}}
\tilde{\phi}_{\e, k_{n} k_{l}}(s /  \, _{\bar{k}_{\bar{m}}}\check{v}_{\e, k_{n}}) \, 
_{\bar{k}_{n -1}}\tilde{p}_{\e,  k_{n} k_{l}}.  
\vspace{2mm} \nonumber \\
&  \quad  \quad  + \Psi_{\e, \DD, k_{\bar{n}} j}(s/ _{\bar{k}_{\bar{n} -1}}
\tilde{v}_{\e, k_{\bar{n}}})  
\vspace{2mm} \nonumber \\ 
&  \quad \quad \times   \, _{\bar{k}_{n -1}}\tilde{\phi}_{\e, k_{n} k_{\bar{n}}}
(\prod_{q = n}^{\bar{n} -1}(1 - \, _{\bar{k}_q}p_{\e, k_{\bar{n}} k_{\bar{n}}}) 
 \frac{_{\bar{k}_{n - 1}}\tilde{v}_{\e, k_{n}}}{_{\bar{k}_{n -1}}\tilde{v}_{\e, k_{\bar{n}}}}s / _{\bar{k}_{n -1}}\tilde{v}_{\e, k_{n}}) \, 
_{\bar{k}_{n -1}}\tilde{p}_{\e,  k_{n} k_{\bar{n}}}. 
 \vspace{2mm} \nonumber \\
 &  \quad  \quad  + \sum_{ n  + 1 \leq l < \bar{n}}\Psi_{\e, \DD, k_{l} j}
(\prod_{q = l}^{\bar{n} -1}(1 - \, 
_{\bar{k}_q}p_{\e, k_{\bar{n}} k_{\bar{n}}}) 
 \frac{_{\bar{k}_{l -1}}\tilde{v}_{\e, k_{l}}}{_{\bar{k}_{l -1}}
\tilde{v}_{\e, k_{\bar{n}}}}s / _{\bar{k}_{\bar{m} - l}}\tilde{v}_{\e, k_{l}})
\vspace{2mm} \nonumber \\ 
&  \quad \quad \times   \, _{\bar{k}_{n -1}}
\tilde{\phi}_{\e, k_{n} k_{l}}(\prod_{q = n}^{\bar{n} -1}(1 - \, _{\bar{k}_q}p_{\e, k_{\bar{n}} k_{\bar{n}}}) 
\frac{_{\bar{k}_{n -1}}\tilde{v}_{\e, k_{n}}}{_{\bar{k}_{n -1}}
\tilde{v}_{\e, k_{\bar{n}}}}s / _{\bar{k}_{n -1}}\tilde{v}_{\e, k_{n}}) \, 
 _{\bar{k}_{n -1}}\tilde{p}_{\e,  k_{n} k_{l}} \vspace{2mm} \nonumber \\ 
% \end{align*}
%\begin{align}
 & \quad \to   \, _{\bar{k}_{n -1}}\tilde{\phi}_{0, k_{n} j}(\prod_{q = n}^{\bar{n} -1}(1 - \, 
_{\bar{k}_q}p_{0, k_{\bar{n}} k_{\bar{n}}}) 
\, _{\bar{k}_{n -1}}\tilde{w}_{k_{n}, k_{\bar{n}}}
s) \, _{\bar{k}_{n -1}}\tilde{p}_{0, k_{n} j} 
\vspace{2mm} \nonumber \\
&  \quad \quad   + \Psi_{0, \DD, k_{\bar{n}} j}(s)  \vspace{2mm} \nonumber \\
&  \quad \quad \times \, _{\bar{k}_{n -1}}\tilde{\phi}_{0, k_{n} k_{\bar{n}}}(\prod_{q = n}^{\bar{n} -1}
(1 - \, _{\bar{k}_q}p_{0, k_{\bar{n}} k_{\bar{n}}}) 
\, _{\bar{k}_{n -1}}\tilde{w}_{k_{n}, k_{\bar{n}}}s) \, _{\bar{k}_{n -1}}\tilde{p}_{0,  k_{n} k_{\bar{n}}} 
\vspace{2mm} \nonumber \\ 
&  \quad  \quad  + \sum_{n + 1 \leq l < \bar{n}}\Psi_{0, \DD, k_{l} j}(\prod_{q = l}^{\bar{n} - 1}
(1 - \, _{\bar{k}_q}p_{0, k_{\bar{n}} k_{\bar{n}}}) 
\, _{\bar{k}_{l -1}}\tilde{w}_{0, k_{l}, k_{\bar{n}}}s)
\vspace{2mm} \nonumber \\ 
&  \quad \quad \times   \, _{\bar{k}_{n - 1}}
\tilde{\phi}_{0, k_{n} k_{l}}(\prod_{q = n}^{\bar{n} - 1}(1 - \, _{\bar{k}_q}p_{0, k_{\bar{n}} k_{\bar{n}}}) 
 \, _{\bar{k}_{n - 1}}\tilde{w}_{k_{n}, k_{\bar{n}}} s) \,  
_{\bar{k}_{n -1}}\tilde{p}_{0,  k_{n} k_{l}}   
\vspace{2mm} \nonumber \\ 
& \quad =  \Psi_{0, \DD, k_{n} j}(s) \ {\rm as} \ \e \to 0.
\end{align}

In this case, the limiting distribution $G_{0, \DD, k_{n} j}(\cdot)$ has, for every $j \in \DD, n = 1, \ldots, \bar{m}$,  the Laplace transform $\Psi_{0, \DD, k_{n} j}(s) = 
\int_0^\infty e^{-st} G_{0, \DD, k_{n} j}(dt)$, $s \geq 0$ given by relation (\ref{sfatomnema}). 

Second, let us assume that 
\begin{equation} \label{sanaset}
\bar{n}(\bar{k}_{\bar{m}}, n)  = n.
\end{equation}

In this case, an analogue of relation (\ref{idersfa})
takes the following form, for $s \geq 0, j \in \DD$, 
\begin{align}\label{sfatomnekan}
& \Psi_{\e, \DD, k_{n } j}(s / _{\bar{k}_{\bar{m}}}\check{v}_{\e, k_{n}}) \vspace{2mm} \nonumber \\ 
& \quad =   \, _{\bar{k}_{n -1}}\tilde{\phi}_{\e, k_{n} j}(s / _{\bar{k}_{n -1}}\tilde{v}_{\e, k_{n}} ) \, 
_{\bar{k}_{n - 1}}\tilde{p}_{\e, k_{n} j}  
\vspace{2mm} \nonumber \\ 
&  \quad \quad   + \sum_{n + 1 \leq l \leq \bar{m}}\Psi_{\e, \DD, k_{l} j}(s / _{\bar{k}_{\bar{m}}}\check{v}_{\e, k_{n}})  \vspace{2mm} \nonumber \\ 
&  \quad \quad \quad \times   
\, _{\bar{k}_{n -1}}\tilde{\phi}_{\e, k_{n} k_{l}}(s / _{\bar{k}_{\bar{m}}}\check{v}_{\e, k_{n}}) \, 
_{\bar{k}_{n - 1}}\tilde{p}_{\e,  k_{n} k_{l}}  
\vspace{2mm} \nonumber \\ 
& \quad \to  \, _{\bar{k}_{n - 1}}\tilde{\phi}_{0, k_{n} j}(s) \, _{\bar{k}_{n - 1}}\tilde{p}_{0, k_{n} j}  \vspace{2mm} \nonumber \\ 
& \quad = \Psi_{0, \DD, k_{n} j}(s)  \ {\rm as} \ \e \to 0.
\end{align} 

In this case, the corresponding  limiting distribution $G_{0, \DD, k_{n} j}(\cdot)$ has, for every $j \in \DD$,  the Laplace transform 
$\Psi_{0, \DD, k_{n} j}(s) = \int_0^\infty e^{-st} G_{0, \DD, k_{n} j}(dt)$, $s \geq 0$ given by relation (\ref{sfatomnekan}). 

As it was mentioned above, the above asymptotic relations (\ref{sfatomnema}) or (\ref{sfatomnekan}) hold, due to Theorems 4 and 5, 
for $n = \bar{m}$ and $n = \bar{m} -1$. The recurrent application of asymptotic relations (\ref{sfatomnema}) and (\ref{sfatomnekan})  let one sequentially get them for $n = \bar{m} - 3, \ldots, 1$. 

The formulas for the limiting Laplace transforms $\Psi_{0, \DD, k_{n} j}(s), n = 1, \ldots, \bar{m}$ given by relations 
(\ref{sfatomnema}) and (\ref{sfatomnekan}) have the recurrent character. 
In the case, $j \in \DD, n = \bar{m}$, the Laplace transform  $\Psi_{0, \DD, k_{\bar{m}} j}(s)$ given by relation (\ref{sfatomnekan})  coincides with the Laplace transform given, according Theorem 4, by relation (\ref{gasd}). 
In the case, $j \in \DD, n = \bar{m} -1$, the Laplace transform  $\Psi_{0, \DD, k_{\bar{m} -1} j}(s)$ given by relations  
(\ref{sfatomnema}) and (\ref{sfatomnekan}) coincide with the Laplace transform given, according Theorem 5,
by relations (\ref{idersfatnen}) and (\ref{idersfa}).
The sequential application of formulas given by relations (\ref{sfatomnema}) and (\ref{sfatomnekan}), for $j \in \DD, n = \bar{m} -2, \ldots$
let one compute the limiting Laplace transforms $\Psi_{0, \DD, k_{n} j}(s)$, for $j \in \DD, n = \bar{m} - 2, \ldots, 1$.  

The  above remarks imply that the following theorem takes place. \vspace{1mm}

{\bf Theorem 6}. {\em Let conditions ${\bf A}$ -- ${\bf E}$, ${\bf \tilde{C}}$, $_{\bar{k}_{r}}{\bf \tilde{C}}$, $_{\bar{k}_{r}}{\bf \hat{C}}$, $_{\bar{k}_{r-1}}{\bf \tilde{F}}$, $r = 1, \ldots, \bar{m} -1$ hold, where the   
states $k_r \in \, _{\bar{k}_{r-1}} \overline{\DD}^*, r = 1, \ldots,  \bar{m} -1$ are chosen in such way that condition
$_{\bar{k}_{r}}{\bf \hat{F}}$ holds, for  $r = 1, \ldots,  \bar{m} -1$. Then,  the following asymptotic 
relation takes place, for $j \in \DD, n = 1, \ldots, \bar{m}$,  
\begin{equation}\label{cotyrevanertew}
G_{\e, \DD, k_{n} j}(\cdot  \, _{\bar{k}_{\bar{m}}} \check{v}_{\e, k_{n} }) 
\Rightarrow  G_{0, \DD, k_{n} j}(\cdot) \ {\rm as} \ \e \to 0. 
\end{equation}}
\vspace{-1mm}

{\bf 7.5 Weak asymptotics for distributions  of hitting times with an arbitrary initial state in domain $\overline{\DD}$}. It is useful to note that, for any $i \in \overline{\DD}$,  there exist the unique $1 \leq \, N_{\bar{k}_m, i} \leq \bar{m}$ such  that,
\begin{equation}\label{hopertino}
k_{N_{\bar{k}_m, i}} = i.
\end{equation}

Note that,  $N_{\bar{k}_m, k_n} = n$, and, thus,  $k_{N_{\bar{k}_m, k_n}} = k_n$, for $n = 1, \ldots, \bar{m}$. 

Respectively, the corresponding normalising functions take 
the form, for $i \in \overline{\DD}$ and $\e \in (0, 1]$, 
\begin{equation}\label{sopitr}
_{\bar{k}_{\bar{m}}}\check{v}_{\e, i} = \,  _{\bar{k}_{\bar{m}}}\check{v}_{\e, k_{N_{\bar{k}_m, i}}}
= \, _{\bar{k}_{\bar{n}(\bar{k}_{\bar{m}}, N_{\bar{k}_m, i}) -1}}\tilde{v}_{\e, k_{\bar{n}(\bar{k}_{\bar{m}}, N_{\bar{k}_m, i})}} 
\end{equation}
where the functions $_{\bar{k}_{\bar{m}}}\check{v}_{\e, k_{n}}, n = 1, \ldots, \bar{m}$  are given by relations 
(\ref{normaba}) and (\ref{erty}).

It is useful to note that the normalisation functions $_{\bar{k}_{\bar{m}}}\check{v}_{\e, i}, i \in \overline{\DD}$ are asymptotically comparable, in the sense that the following relation holds, for every $j \in \overline{\DD}^*_i,  i \in \overline{\DD}$, 
\begin{align}\label{compas}
\frac{_{\bar{k}_{\bar{m}}}\check{v}_{\e, j}}{_{\bar{k}_{\bar{m}}}\check{v}_{\e, i}} &
= \frac{_{\bar{k}_{\bar{n}(\bar{k}_{\bar{m}}, N_{\bar{k}_m, j}) -1}}
\tilde{v}_{\e, k_{\bar{n}(\bar{k}_{\bar{m}}, N_{\bar{k}_m, j})}}}{_{\bar{k}_{\bar{n}(\bar{k}_{\bar{m}}, N_{\bar{k}_m, i}) -1}}
\tilde{v}_{\e, k_{\bar{n}(\bar{k}_{\bar{m}}, N_{\bar{k}_m, i})}}} \vspace{2mm} \nonumber \\
& =  \prod_{q = \bar{n}(\bar{k}_{\bar{m}}, N_{\bar{k}_m, j})}^{\bar{n}(\bar{k}_{\bar{m}}, N_{\bar{k}_m, i}) - 1}
(1 - \, _{\bar{k}_q}p_{\e, k_{\bar{n}(\bar{k}_{\bar{m}}, N_{\bar{k}_m, i})} k_{\bar{n}(\bar{k}_{\bar{m}}, N_{\bar{k}_m, i})}}) 
\vspace{2mm} \nonumber \\
%\end{align*} 
%\begin{align}
&  \quad \times \frac{_{\bar{k}_{\bar{n}(\bar{k}_{\bar{m}}, N_{\bar{k}_m, j}) -1}}\tilde{v}_{\e, k_{\bar{n}(\bar{k}_{\bar{m}}, N_{\bar{k}_m, j})}}}{_{\bar{k}_{\bar{n}(\bar{k}_{\bar{m}}, N_{\bar{k}_m, j}) -1}}\tilde{v}_{\e, k_{\bar{n}(\bar{k}_{\bar{m}}, N_{\bar{k}_m, i})}}} 
\vspace{2mm} \nonumber \\ 
& \to  \prod_{q =  \bar{n}(\bar{k}_{\bar{m}}, N_{\bar{k}_m, j})}^{\bar{n}(\bar{k}_{\bar{m}}, N_{\bar{k}_m, i}) -1}(1 - \, _{\bar{k}_q}p_{0, k_{\bar{n}(\bar{k}_{\bar{m}}, N_{\bar{k}_m, i})} k_{\bar{n}(\bar{k}_{\bar{m}}, N_{\bar{k}_m, i})}}) 
\vspace{2mm} \nonumber  \\
&  \quad \times \, _{\bar{k}_{\bar{n}(\bar{k}_{\bar{m}}, N_{\bar{k}_m, j}) -1}}
\tilde{w}_{0, k_{\bar{n}(\bar{k}_{\bar{m}}, 
N_{\bar{k}_m, j})} k_{\bar{n}(\bar{k}_{\bar{m}}, N_{\bar{k}_m, i})}} 
\vspace{2mm} \nonumber \\
& = \, _{\bar{k}_{\bar{m}}}\check{w}_{0, ji} \in [0, \infty) \ {\rm as} \ \e \to 0,   
\end{align}
where, for $i \in \overline{\DD}$, 
\begin{equation}\label{jokl}
\overline{\DD}^*_i = \{ j \in \overline{\DD}: \bar{n}(\bar{k}_{\bar{m}}, N_{\bar{k}_m, j}) \leq \bar{n}(\bar{k}_{\bar{m}}, N_{\bar{k}_m, i}) \}.
\end{equation}

The above remarks, let us reformulate the weak asymptotic relation (\ref{cotyrevanert}) given in Theorem 6 in the following 
simpler  form,  for $i \in \overline{\DD}, j \in \DD$,  
\begin{equation}\label{cotyrevanertop}
G_{\e, \DD, i j}(\cdot  \, _{\bar{k}_{\bar{m}}}\check{v}_{\e, i}) 
\Rightarrow  G_{0, \DD, i j}(\cdot) \ {\rm as} \ \e \to 0,  
\end{equation}
where the limiting distributions and their Laplace transforms are given by the following relations, for $i \in \overline{\DD}, j \in \DD$, 
\begin{equation}\label{fity}
G_{0, \DD, i j}(\cdot) =  G_{0, \DD, k_{N_{\bar{k}_m, i}}, j}(\cdot),     
\end{equation} 
and
\begin{equation}\label{fityh}
 \Psi_{0, \DD, i j}(\cdot) =  \Psi_{0, \DD, k_{N_{\bar{k}_m, i}}, j}(\cdot).    
\end{equation} 
 
{\bf 7.6 Convergence of hitting probabilities}. 
Let us introduce hitting probabilities, for $i \in \XX, j \in \DD$ and $\e \in [0, 1]$, 
\begin{equation}\label{hittata}
P_{\e, \DD, i j}  = \PP_i \{ \eta_\e(\tau_{\e, \DD}) = j \}  =  G_{\e, \DD, i j}(\infty) = \Psi_{\e, \DD, i j}(0).  
\end{equation}

Note that condition ${\bf B}$ and Theorem 6 imply that, for any $i \in \overline{\DD}$ and $\e \in [0, 1]$,
\begin{equation}\label{suma}
P_{\e, \DD, i j} \geq 0, j \in \DD, \ \sum_{j \in \DD} P_{\e, \DD, i j} = 1.
\end{equation}

By taken $s = 0$ in relations (\ref{sfatomne}), we get, for every $j \in \DD$ and $\e \in (0, 1]$,  the following backward 
recurrent relation for hitting probabilities, 
for $n = \bar{m}, \ldots, 1$, 
\begin{equation}\label{somne}
P_{\e, \DD, k_{n} j}  =  \, _{\bar{k}_{n - 1}}\tilde{p}_{\e, k_{n} j}  
+ \sum_{n + 1 \leq l \leq \bar{m}} P_{\e, \DD, k_{l} j} \, _{\bar{k}_{n -1}}\tilde{p}_{\e,  k_{n} k_{l}}.  
\end{equation} 

Let us introduce the following sets, for $j \in \XX$ and $\e \in [0, 1]$,
\begin{equation}\label{setrabok}
\YY_{\e, \DD, i} = \{ j \in \DD: P_{\e, \DD, i j}  > 0 \}. 
\end{equation}

{\bf Lemma 34}. {\em Let conditions of Theorem 6 hold. Then, sets, $\YY_{\e, \DD, i}, i  \in \overline{\DD}$ do not depend on $\e \in (0, 1]$, i.e., for $i \in \overline{\DD}$,
\begin{equation}\label{seter}
\YY_{\e, \DD, i} = \YY_{1, \DD, i}, \e \in (0, 1].
\end{equation}}
\makebox[3mm]{} {\bf Proof}. According condition ${\bf B}$ and  Lemma 23, either  
 $_{\bar{k}_{n -1}}\tilde{p}_{\e, k_{n} j} = 0$, $\e \in (0, 1]$, or $_{\bar{k}_{n - 1}}\tilde{p}_{\e, k_{n} j} > 0, \e \in (0, 1]$, 
 for every $n = 1, \ldots, \bar{m}$ and,
either   $_{\bar{k}_{n - 1}}\tilde{p}_{\e,  k_{n} k_{l}} = 0, \e \in (0, 1]$,  or
$_{\bar{k}_{n - 1}}\tilde{p}_{\e,  k_{n} k_{l}} > 0, \e \in (0, 1]$, for 
 every  $l =  n + 1, \ldots, \bar{m}, n = 1, \ldots, \bar{m}$.  

The above remarks and relation (\ref{somne}) taken for  $n = \bar{m}$, imply that sets $\YY_{\e, \DD, k_{\bar{m}}} = \YY_{1, \DD, k_{\bar{m}}}, \e \in (0, 1]$. Then, the above remarks and relation (\ref{somne}) taken for  $n = \bar{m} - 1$, imply that sets  $\YY_{\e, \DD, k_{\bar{m} -1}} = \YY_{1, \DD, k_{\bar{m} - 1}}, \e \in (0, 1]$. By continuing in this way, we get that, $\YY_{\e, \DD, k_{n}} = 
\YY_{1, \DD, k_{n}}$, for any $n = \bar{m}, \ldots, 1$ and $\e \in (0, 1]$. 

It remains to note that, for $i \in \overline{\DD}, j \in \DD$ and  $\e \in (0, 1]$,
\begin{equation}\label{koplk}
P_{\e, \DD, i j} =  P_{\e, \DD, k_{N_{\bar{k}_m, i}}, j}, 
\end{equation}
and, 
\begin{equation}\label{mukl}
 \YY_{\e, \DD, i} =  \YY_{\e, \DD, k_{N_{\bar{k}_m, i}}}. 
\end{equation}

Relations (\ref{cotyrevanertop}), (\ref{hittata}) and (\ref{koplk}) imply that the 
the following lemma takes place. \vspace{1mm}

{\bf Lemma 35}. {\em Let conditions of Theorem 6 holds. Then,  the following relation takes place, for $i \in \overline{\DD}, j \in \DD$,
\begin{equation}\label{cotyrevanertopf}
P_{\e, \DD, i j} \to  P_{0, \DD, i j}  \ {\rm as} \ \e \to 0. 
\end{equation}}
\makebox[3mm]{} Note that relations (\ref{seter}) and  (\ref{cotyrevanertopf}) imply that sets $\YY_{0, \DD, i} \subseteq \YY_{1, \DD, i}$, 
for $i \in \overline{\DD}$. 

It is also useful to note that the sequential application of formulas given by relations (\ref{sfatomnema}) and (\ref{sfatomnekan}), where one should choose $s = 0$, let one compute the limiting probabilities  $P_{0, \DD, k_{n} j} = \Psi_{0, \DD, k_{n} j}(0)$, for 
$j \in \DD, n =  \bar{m}, \ldots, 1$. Then, the following relation let compute the limiting hitting probabilities $P_{0, \DD, i j}, j \in \DD, i \in \overline{\DD}$,
\begin{equation}\label{gopi}
P_{0, \DD, i j} =  P_{0, \DD, k_{N_{\bar{k}_m, i}}, j}.  
\end{equation}

{\bf 7.7 Weak asymptotics for conditional distributions of hitting times}.  Let us introduce conditional distribution functions, for $i \in \XX, j \in \DD$ and $\e \in [0, 1]$,   
\begin{equation}\label{hittatanas}
F_{\e, \DD, i j}(t)  = \PP_i  \{ \tau_{\e, \DD} \leq t / \eta_\e(\tau_{\e, \DD}) = j \}, \ t \geq 0.  
\end{equation}

If $P_{\e, \DD, i j} > 0$, then, obviously, 
\begin{equation}\label{okure}
F_{\e, \DD, i j}(t) =  P_{\e, \DD, i j}^{-1} G_{\e, \DD, i j}(t), \ t \geq 0.
\end{equation}

If $P_{\e, \DD, i j} = 0$, one can define the distribution function $F_{\e, \DD, i j}(t)$ in the standard way,  as,
\begin{align}\label{kasd}
F_{\e, \DD, i j}(t)  & = F_{\e, \DD, i}(t)  = \sum_{j \in \YY_{\e, i}} G_{\e, \DD, i j}(t) \vspace{2mm} \nonumber \\
&  = \sum_{j \in \YY_{\e, \DD, i}} F_{\e, \DD, i j}(t) P_{\e, \DD, i j}, \ t \geq 0.
\end{align}  

Also, relations (\ref{cotyrevanertop}), (\ref{cotyrevanertopf}) and (\ref{okure}) imply that the following asymptotic relation holds, for $j \in \YY_{1, \DD, i}, i \in \overline{\DD}$, 
\begin{equation}\label{cotyrebad}
F_{\e, \DD, i j}(\cdot \, _{\bar{k}_{\bar{m}}}\check{v}_{\e, i}) 
\Rightarrow  F_{0, \DD, i j}(\cdot) \ {\rm as} \ \e \to 0. 
\end{equation} 

Relations (\ref{cotyrevanertopf}) and (\ref{cotyrebad}) imply that, for $j \notin \YY_{1, \DD, i}, i \in \overline{\DD}$,
\begin{align}\label{cotyrebe}
F_{\e, \DD, i j}(\cdot \, _{\bar{k}_{\bar{m}}}\check{v}_{\e, i})  & = \sum_{j \in \YY_{1, \DD,  i}} F_{\e, \DD, i j}(\cdot /  
\, _{\bar{k}_{\bar{m}}}\check{v}_{\e, i}) P_{\e, \DD, i j}  \vspace{2mm} \nonumber \\
& \Rightarrow  \sum_{j \in \YY_{1, \DD, i}} F_{0, \DD, i j}(\cdot) P_{0, \DD, i j} 
= \sum_{j \in \YY_{0, \DD, i}} F_{0, \DD, i j}(\cdot) P_{0, \DD, i j} \vspace{2mm} \nonumber \\
&  =  F_{0, \DD, i}(\cdot)  =  F_{0, \DD, i j}(\cdot) \ {\rm as} \ \e \to 0. 
\end{align} 

Let us also introduce Laplace transforms, for $i \in \XX, j \in \DD$ and $\e \in [0, 1]$,   
\begin{equation}\label{hiatanas}
\Phi_{\e, \DD, i j}(s)  = \int_0^\infty e^{-st} F_{\e, \DD, i j}(dt), \ s \geq 0.  
\end{equation}

If $P_{\e, \DD, i j} > 0$, then, obviously, 
\begin{equation}\label{nouret}
\Phi_{\e, \DD, i j}(s) =  P_{\e, \DD, i j}^{-1} \Psi_{\e, \DD, i j}(s), \ s \geq 0.
\end{equation}

If $P_{\e, \DD, i j} = 0$, then,
\begin{align}\label{cotyrebena}
\phi_{\e, \DD, i j}(s) & = \phi_{\e, \DD, i}(s)  \vspace{2mm} \nonumber \\
& = \sum_{j \in \YY_{\e, \DD, i}} \Phi_{\e, \DD, i j}(s) P_{\e, \DD, i j}, \ s \geq 0.  
\end{align} 

Relations (\ref{cotyrebad}) and (\ref{cotyrebe}) can be expressed in the equivalent form in terms of the above Laplace transforms. It takes the 
following form, for $j \in \YY_{1, \DD, i}, i \in \overline{\DD}$, 
\begin{equation}\label{cotyreba}
\Phi_{\e, \DD, i j}(s / \, _{\bar{k}_{\bar{m}}}\check{v}_{\e, i}) 
\to  \Phi_{0, \DD, i j}(s) \ {\rm as} \ \e \to 0, \ {\rm for} \ s \geq 0.  
\end{equation} 
and, for $j \notin \YY_{1, \DD, i}, i \in \overline{\DD}$,
\begin{align}\label{cotyrebebaf}
\Phi_{\e, \DD, i j}(s / \,    _{\bar{k}_{\bar{m}}}\check{v}_{\e, i})  & = \sum_{j \in \YY_{1, i}} \Phi_{\e, \DD, i j}(s /  
\, _{\bar{k}_{\bar{m}}}\check{v}_{\e, i}) P_{\e, \DD, i j}  \vspace{2mm} \nonumber \\
& \Rightarrow   \sum_{j \in \YY_{1, i}} \Phi_{0, \DD, i j}(s) P_{0, \DD, i j} 
= \sum_{j \in \YY_{0, i}} \Phi_{0, \DD, i j}(s) P_{0, \DD, i j} \vspace{2mm} \nonumber \\
&  =  \Phi_{0, \DD, i}(s)  =  \Phi_{0, \DD, i j}(s) \ {\rm as} \ \e \to 0, \ {\rm for} \ s \geq 0. 
\end{align} 

Let also point out that the distribution functions $F_{0, \DD, i j}(\cdot), i \in \overline{\DD}, j \in \DD$ are not concentrated in $0$, i.e., 
for  $i \in \overline{\DD}, j \in \DD$, 
\begin{equation}\label{oput}
F_{0, \DD, i j}(0) < 1,
\end{equation}
or,   equivalently, for $s > 0$ and $i \in \overline{\DD}, j \in \DD$, 
\begin{equation}\label{oputew}
\Phi_{0, \DD, i j}(s) < 1. 
\end{equation}

Indeed, let us assume that assumption (\ref{sanaset}) holds. 

In this case, 
relation (\ref{sfatomnekan}) implies that $F_{0, \DD, k_{n} j}(\cdot) = \,  _{\bar{k}_{n -1}}\tilde{F}_{0, k_{n} j}(\cdot)$, 
for $j \in \DD, n = \bar{m}, \ldots, 1$, and, thus, by Lemma 23,  $F_{0, \DD, k_{n} j}(0) < 1$, for $j \in \DD, n =  \bar{m}, \ldots, 1$.

Let us now assume  that,  assumption (\ref{assump}) holds,  for some $1 \leq n \leq \bar{m}$, and 
$F_{0, \DD, \bar{k}_{l} j}(0) < 1, \, j \in \DD,  n + 1 \leq l \leq \bar{m}$.  

In this, case,  $F_{0, \DD, \bar{k}_{n} j}(0) < 1, j \in \DD$.
Indeed, relations (\ref{sfatomnema})  and (\ref{somne}) imply that, for every  the following inequality holds, for every $s > 0, j \in \DD$,
\begin{align}\label{inequal}
\Psi_{0, \DD, k_{n} j}(s) & \leq  \,  _{\bar{k}_{n - 1}}\tilde{p}_{0, k_{n} j} 
+ \Psi_{0, \DD, k_{\bar{n}} j}(s) \,  _{\bar{k}_{n - 1}}\tilde{p}_{0,  k_{n} k_{\bar{n}}} 
\vspace{2mm} \nonumber \\ 
&  \quad   +  \sum_{n  +1 \leq l <  \bar{n}} P_{0, \DD, k_{l} j} \, _{\bar{k}_{n - 1}}\tilde{p}_{0,  k_{n} k_{l}} \vspace{2mm} \nonumber \\ 
& <   \,  _{\bar{k}_{\bar{m} - n}}\tilde{p}_{0, k_{\bar{m} - n +1} j} 
+ \,   _{\bar{k}_{n - 1}}\tilde{p}_{0,  k_{n} k_{\bar{n}}} \vspace{2mm} \nonumber \\ 
&  \quad   +  \sum_{n  + 1 \leq l <  \bar{n}} P_{0, \DD, k_{l} j} \, _{\bar{k}_{n -1}}\tilde{p}_{0,  k_{n} k_{l}}  = P_{0, \DD, k_{n} j}. 
\end{align} 

Relation (\ref{inequal}) obviously implies that, for $s > 0$ and $j \in \YY_{0, \DD, k_{n}}, n = \bar{m}, \ldots, 1$.
\begin{equation}\label{noplk}
\Phi_{0, \DD, k_{n} j}(s) =  
P_{0, \DD, k_{n} j}^{-1}\Psi_{0, \DD, k_{n} j}(s) < 1.
\end{equation}

In sequel, for $s > 0, j \notin \YY_{0, \DD, k_{n}}, n = \bar{m}, \ldots, 1$,  
\begin{align}\label{noplkasty}
\Phi_{0, \DD, k_{n} j}(s) & = \Phi_{0, \DD, k_{n}}(s)  = \sum_{r \in \YY_{0, \DD, k_{n}}} \Phi_{0, \DD, k_{n} r}(s) 
P_{0, \DD, k_{n} r} < 1.
\end{align}

The above remarks, prove relations (\ref{oput}). \vspace{1mm}

{\bf 7.8 Summary of weak asymptotics for distributions of hitting times}. The remarks made in Subsections 7.5 -- 7.7 let us reformulate Theorem 6 in the following equivalent but more informative final form. \vspace{1mm}

{\bf Theorem 7}. {\em Let conditions ${\bf A}$ -- ${\bf E}$, ${\bf \tilde{C}}$, $_{\bar{k}_{r}}{\bf \tilde{C}}$, $_{\bar{k}_{r}}{\bf \hat{C}}$, $_{\bar{k}_{r-1}}{\bf \tilde{F}}$, $r = 1, \ldots, \bar{m} -1$ hold, where the   
states $k_r \in \, _{\bar{k}_{r-1}} \overline{\DD}^*, r = 1, \ldots,  \bar{m} -1$ are chosen in such way that condition
$_{\bar{k}_{r}}{\bf \hat{F}}$ holds, for  $r = 1, \ldots,  \bar{m} -1$. Then,  the following asymptotic 
relation takes place, for $i \in \overline{\DD}, j \in \DD$,  
\begin{align}\label{cotyr}
G_{\e, \DD, i j}(\cdot \, _{\bar{k}_{\bar{m}}}\check{v}_{\e, i}) &
\Rightarrow  G_{0, \DD, i j}(\cdot) \vspace{2mm} \nonumber \\
& = F_{0, \DD, i j}(\cdot) P_{0, \DD, i j}\ {\rm as} \ \e \to 0, 
\end{align}
where the limiting distribution functions $F_{0, \DD, i j}(\cdot)$ are not concentrated in zero, i.e.,  
$F_{0, \DD, i j}(0) < 1$, for $i \in \overline{\DD}, j \in \DD$.} \vspace{1mm}

{\bf Remark 13}. The normalisation functions $_{\bar{k}_{\bar{m}}}\check{v}_{\e, i}$ can be found by applications of recurrent relations
(\ref{normaba}) and (\ref{erty}), and, then, relations (\ref{hopertino}), and (\ref{sopitr}). The limiting Laplace transforms $\Psi_{0, \DD, i j}(s) = \int_0^\infty e^{-st}G_{0, \DD, i j}(dt), s \geq 0$, probabilities $P_{0, \DD, i j}$ and Laplace transforms $\Phi_{0, \DD, i j}(s) =  
\int_0^\infty e^{-st}F_{0, \DD, i j}(dt)$, $s \geq 0$ can be found by the recurrent application of relations  (\ref{sfatomnema}) and (\ref{sfatomnekan}), and, then relations (\ref{hopertino}), (\ref{fity}), (\ref{nouret}), and (\ref{cotyrebena}).  \vspace{1mm}

 {\bf Remark 14}. Condition ${\bf C}_{\bar{m}-1}$ is sufficient for holding of  conditions ${\bf \tilde{C}}$, $_{\bar{k}_{r}}{\bf \tilde{C}}$, $_{\bar{k}_{r}}{\bf \hat{C}}$, $r = 1, \ldots, \bar{m} -1$, and, thus, can  replace them in Theorems 4 -- 7. \vspace{1mm}

{\bf Remark 15}. Condition ${\bf G}$  is sufficient for holding of condition ${\bf C}_{\bar{m}-1}$. 

\vspace{1mm}

{\bf Remark 16}. Condition ${\bf F}_{\bar{m}-2}$  is sufficient for holding of  conditions $_{\bar{k}_{r-1}}{\bf \tilde{F}}$, $r = 1, \ldots, \bar{m} -1$, and, thus, can replace them in Theorems 4 -- 7. \vspace{1mm}

{\bf Remark 17}. Condition ${\bf H}$  are sufficient for holding of conditions ${\bf C}_{\bar{m}-1}$ and ${\bf F}_{\bar{m}-2}$. \\
 
 {\bf 8  Weak Convergence of Hitting Times, for  the Case  Where \\ \makebox[9mm]{} an Initial State Belongs  
 to Domain $\DD$} \\
 
In this section, we extend the results presented in Section 7 to the model of perturbed semi-Markov processes  with an initial state that belongs to domain $\DD$. \vspace{1mm}

 {\bf 8.1  Weak asymptotics of distributions of  hitting times for the case with an initial state in domain $\DD$}. In this subsection,  we present 
 conditions of weak convergence for distributions $G_{\e, \DD. ij}(\cdot)$ for the case, where the  initial state $i \in \DD$. 
 
 The key role is played by the following relation, which takes place for $t \geq 0, r, j \in \DD$ 
 and $\e \in (0,1]$, 
 \begin{align}\label{inita}
 G_{\e, \DD, rj}(t)  & = F_{\e, rj}(t) p_{\e, rj} 
 + \sum_{i \in \overline{\DD}} G_{\e, \DD, ij}(t)*F_{\e, ri}(t)  p_{\e, ri} 
 \vspace{3mm} \nonumber \\
 & = F_{\e, rj}(t) p_{\e, rj} + \sum_{i \in \overline{\DD}} F_{\e, \DD, ij}(t)*F_{\e, ri}(t)  P_{\e, \DD, ij} p_{\e, ri}.  
 \end{align}

The above relation can be re-written in terms of the corresponding Laplace transforms, for $s \geq 0, i, j \in \DD$ and $\e \in (0,1]$, 
 \begin{align}\label{initamok}
  \Psi_{\e, \DD, rj}(s)  & = \phi_{\e, rj}(s) p_{\e, rj} 
 + \sum_{i \in \overline{\DD}}  \Psi_{\e, \DD, ij}(s)  \phi_{\e, ri}(s) p_{\e, ri} 
 \vspace{3mm} \nonumber \\
% \end{align*}
%\begin{align}
 & = \phi_{\e, rj}(s) p_{\e, rj} 
 + \sum_{i \in \overline{\DD}}  \Phi_{\e, \DD, ij}(s)  \phi_{\e, ri}(s) P_{\e, \DD, ij}  p_{\e, ri}.  
 \end{align}

These relations hint us that, in this case, analogous of basic conditions ${\bf B}$ {\bf (a)},  ${\bf C}$ -- ${\bf E}$ (for the case where domain $\overline{\DD}$ is replaced by domain $\DD$) should be assumed to hold: 
\begin{itemize}
\item[${\bf \dot{B}}$:]  {\bf (a)} $p_{\e, ij} > 0, \e \in (0, 1]$, for $j \in \YY_{1, i}, i \in \DD$ and 
 $p_{\e, ij} = 0, \e \in (0, 1]$, for every $j \in \overline{\YY}_{1, i}, i \in \DD$, 
\end{itemize}
where sets $\YY_{1, i}$ have been defined in relation (\ref {hore}).

\begin{itemize}
\item[${\bf  \dot{C}}$:]  $p_{\e, ij} \to p_{0, ij}$ as $\e \to 0$, for $i \in \DD, j \in \XX$. \vspace{2mm}
\end{itemize}
\begin{itemize}
\vspace{1mm}

\item[${\bf  \dot{D}}$:]   {\bf (a)} $F_{\e, ij}(\cdot \, v_{\e,  i}) = \PP \{ \kappa_{\e, 1}/ v_{\e,  i}   \leq \cdot /  \eta_{\e, 0} = i, \eta_{\e, 1} = j  \} \Rightarrow F_{0, ij}(\cdot)$ as $\e \to 0$,  for $j \in \YY_{1, i}, i \in \DD$, {\bf (b)} $F_{0, ij}(\cdot), j \in \YY_{1, i}, i \in \DD$ are proper distribution functions such that $F_{0, ij}(0) < 1, j \in \YY_{1, i}, i \in \DD$,   
{\bf (c)} $1 \leq v_{\e, i} \to v_{0, i} \in [1, \infty]$ as $\e \to 0$, for $i \in \DD$. 
\end{itemize}
and
\begin{itemize}
\item[${\bf  \dot{E}}$:]  $e_{\e, ij}/v_{\e,  i} = \int_0^{\infty} t F_{\e, ij}(dt)/ v_{\e,  i}   \to 
e_{0, ij} = \int_0^{\infty} t F_{0, ij}(dt)$ as $\e \to 0$,  for $j \in \YY_{1, i}, i \in \DD$. 
\end{itemize}

Condition ${\bf \dot{B}}$, ${\bf \dot{C}}$ and ${\bf \dot{D}}$ imply that the following relation of weak convergence holds, for $i \in \DD$,
\begin{align}\label{alsoknew}
F_{\e, i}( \cdot v_{\e,  i}) & = \sum_{j  \in \YY_{1, i}}  F_{\e, ij}(\cdot v_{\e,  i} ) p_{\e, ij} \vspace{2mm} \nonumber \\
& \Rightarrow \sum_{j  \in \YY_{1, i}}  F_{0, ij}( \cdot) p_{0, ij} \vspace{2mm} \nonumber \\
& = \sum_{j  \in \YY_{0, i}}  F_{0, ij}( \cdot) p_{0, ij} =  F_{0, i}(\cdot) \ {\rm as} \ \e \to 0.
\end{align}

Obviously, $F_{0, i}(t)$ is a proper distribution function such that  $F_{0, i}(0) < 1$, for $i \in  \overline{\DD}$.

Conditions ${\bf \dot{B}}$, ${\bf \dot{C}}$, ${\bf \dot{D}}$ and ${\bf \dot{E}}$ also imply the following asymptotic relation  holds, for $i \in \XX$,
\begin{align}\label{alsokalan}
\frac{e_{\e, i}}{v_{\e,  i}}  & = \sum_{j  \in \YY_{1, i}}  \frac{e_{\e, ij}}{v_{\e,  i}} p_{\e, ij}  
 \to \sum_{j  \in \YY_{1, i}}  e_{0, ij} p_{0, ij} \vspace{2mm} \nonumber \\
%\end{align*}
%\begin{align}
& =   \sum_{j  \in \YY_{0, i}}  e_{0, ij} p_{0, ij} = e_{0, i} \ {\rm as} \ \e \to 0.
\end{align}

It is useful to note that expectation $e_{0, i} \in (0, \infty)$, for $i \in \overline{\DD}$.

It is useful to note that condition ${\bf \check{D}}$ can be formulated in the following equivalent form:
\begin{itemize}
\item[${\bf  \dot{D}}'$:]   {\bf (a)} $\phi_{\e, ij}(\cdot / v_{\e,  i}) = \EE \{ e^{-s \kappa_{\e, 1}/ v_{\e,  i} }  /  \eta_{\e, 0} = i, \eta_{\e, 1} = j  \} \Rightarrow \phi_{0, ij}(s)$ as $\e \to 0$,  for $j \in \YY_{1, i}, i \in \DD$, {\bf (b)} $\phi_{0, ij}(s)= \int_0^\infty e^{-st}F_{0, ij}(dt), s \geq 0$, for $j \in \YY_{1, i}, i \in \DD$, where  $F_{0, ij}(\cdot)$ are proper distribution functions such that $F_{0, ij}(0) < 1, j \in \YY_{1, i}, i \in \DD$,   
{\bf (c)} $1 \leq v_{\e, i} \to v_{0, i} \in [1, \infty]$ as $\e \to 0$, for $i \in \DD$. 
\end{itemize}

In Sections 8.2 and 8.3,  we assume that conditions of Theorem 7 hold. \vspace{1mm} 

{\bf 8.2 Convergence of hitting probabilities}.  Relation (\ref{initamok}) implies that the following relation holds, for $i, j \in \DD$ and $\e \in (0, 1]$,  
\begin{align}\label{hopkufa}
 P_{\e, \DD, ij} & = \PP_i \{ \eta_\e(\tau_{\e, \DD}) = j \} = \Psi_{\e, \DD, ij}(0)    
 \vspace{2mm} \nonumber \\
 &  = p_{\e, ij}  + \sum_{l \in \overline{\DD}}  P_{\e, \DD, lj} p_{\e, il}. 
\end{align}

Note that condition ${\bf \dot{B}}$, Lemmas 34 and 35  and relation (\ref{hopkufa})  imply that, for any $i \in \DD$ and $\e \in [0, 1]$,
\begin{equation}\label{sumat}
P_{\e, \DD, i j} \geq 0, j \in \DD, \ \sum_{j \in \DD} P_{\e, \DD, i j} = 1.
\end{equation}

Condition ${\bf \dot{B}}$, Lemmas 34 and 35 and  relation (\ref{hopkufa}) also   imply that sets $\YY_{\e, \DD, i} = \{ j \in \DD: P_{\e, \DD, i j}  > 0 \}, i \in \DD$ do not depend on $\e \in (0, 1]$, i.e., for $i \in \DD$,
\begin{equation}\label{setera}
\YY_{\e, \DD, i} = \YY_{1, \DD, i}, \, \e \in (0, 1].
\end{equation}

Also, condition ${\bf \dot{C}}$, Lemma 35 and relation (\ref{hopkufa}) imply that the following relation holds, for $i, j \in \DD$, 
\begin{align}\label{hopkufno}
 P_{\e, \DD, ij} & = p_{\e, rj}  + \sum_{l \in \overline{\DD}}  P_{\e, \DD, lj} p_{\e, il} 
\vspace{2mm} \nonumber \\
& \to p_{0, ij}  + \sum_{l \in \overline{\DD}}  P_{0, \DD, lj} p_{0, il}  
=  P_{0, \DD, ij} \ {\rm as} \ \e \to 0.
\end{align}

{\bf 8.3 Weak asymptotics for distributions of hitting times}.  
Let us also assume that the following comparability condition  holds:
\vspace{2mm}

\noindent$_{\bar{k}_{\bar{m}}}{\bf \dot{F}}$:  $\frac{_{\bar{k}_{\bar{m}}}\check{v}_{\e, l}}{v_{\e, i}} \to \, _{\bar{k}_{\bar{m}}}\dot{w}_{0, li} \in [0, \infty]$ as $\e \to 0$, for $l \in \overline{\DD}, i \in \DD$.  

\vspace{2mm} 

It is useful to note that condition $_{\bar{k}_{\bar{m}}}{\bf \dot{F}}$ holds if condition ${\bf H}$ holds and,  additionally, the following condition holds:
\begin{itemize}
\item [${\bf \dot{H}}$:] Functions $p_{\cdot, ij}, j \in \YY_{1, j}, i \in \DD$ and $v_{\cdot, i}, 
i \in \DD$ beong to some complete family of asymptotically comparable functions.   
\end{itemize} 

It turns out that, in this case,  a more complex normalising functions should be used and, moreover, these functions should depend on an entrance  state  to domain $\DD$. Let us define sets,
\begin{equation}\label{normadsr}
\overline{\DD}_{ij} = \{ l \in \overline{\DD}: P_{0, \DD, lj} \, p_{0, il} > 0 \}. 
\end{equation}

Relation (\ref{hopkufno}) implies that, in the case, where $P_{0, \DD, ij} > 0$, it is impossible that probability $p_{0, ij} = 0$ and set $\overline{\DD}_{ij} = \emptyset$,  simultaneously. 

Let us define the normalisation functions, for $i, j \in \DD$ such that $P_{0, \DD, ij} > 0$, 
\begin{equation}\label{nortypir}
_{\bar{k}_{\bar{m}}}\dot{v}_{\e, i j} = 
v_{\e, i}{\rm I}(p_{0, ij} > 0) + \sum_{l \in \overline{\DD}_{ij}}  \, _{\bar{k}_{\bar{m}}}\check{v}_{\e, l}  
\in [1, \infty), \ \e \in (0, 1].
\end{equation}

Condition $_{\bar{k}_{\bar{m}}}{\bf \check{F}}$ imply that the following relation holds, for 
$r, j \in \DD$, 
 \begin{align}\label{opluty}
\frac{v_{\e, i}}{_{\bar{k}_{\bar{m}}}\dot{v}_{\e, i j}} & = \big( {\rm I}(p_{0, ij} > 0) 
+ \sum_{l \in \overline{\DD}_{ij}} \,  \frac{_{\bar{k}_{\bar{m}}}\check{v}_{\e, l}}{v_{\e, i}} \big)^{-1} 
\vspace{2mm} \nonumber \\
& \to \big( {\rm I}(p_{0, ij} > 0)  + \sum_{l \in \overline{\DD}_{ij}} \, _{\bar{k}_{\bar{m}}}\dot{w}_{0, li} \big)^{-1} \vspace{2mm} \nonumber \\
& = \, _{\bar{k}_{\bar{m}}}\dot{u}_{0}[iij] \in [0, 1] \ {\rm as} \ \e \to 0. 
 \end{align}
 
Also, relations (\ref{compas}), (\ref{jokl}) and condition $_{\bar{k}_{\bar{m}}}{\bf \dot{F}}$ imply that the following relations holds, for $l \in \overline{\DD}_{ij}, i, j \in \DD$,
\begin{align}\label{oplutyk}
\frac{_{\bar{k}_{\bar{m}}}\check{v}_{\e, l}}{_{\bar{k}_{\bar{m}}}\dot{v}_{\e, ij}} & = 
\big( \frac{v_{\e, i}}{_{\bar{k}_{\bar{m}}}\check{v}_{\e, l}} {\rm I}(p_{0, ij} > 0) 
+ \sum_{l' \in \overline{\DD}_{ij}} \,  
\frac{_{\bar{k}_{\bar{m}}}\check{v}_{\e, l'}}{_{\bar{k}_{\bar{m}}}\check{v}_{\e, l}} \big)^{-1} 
\vspace{2mm} \nonumber \\
& \to  \big( \, _{\bar{k}_{\bar{m}}}\dot{w}_{0, li}^{-1} {\rm I}(p_{0, ij} > 0) + 
\sum_{l' \in \overline{\DD}_{ij} \cap \overline{\DD}^*_l} \,  _{\bar{k}_{\bar{m}}}\check{w}_{0, l'l}  + \sum_{l' \in \overline{\DD}_{ij} \setminus \overline{\DD}^*_l}  \, _{\bar{k}_{\bar{m}}}\check{w}_{0, l l'}^{-1}\big)^{-1} \vspace{2mm} \nonumber \\
& = \, _{\bar{k}_{\bar{m}}}\dot{u}_{0}[l i j] \in [0, 1] \ {\rm as} \ \e \to 0, 
 \end{align}
 where one should count  $_{\bar{k}_{\bar{m}}}\dot{w}_{0, l i}^{-1}{\rm I}(p_{0, ij} > 0) = 0$, if ${\rm I}(p_{0, ij} > 0) = 0$. 
 
 Remind that sets $\overline{\DD}^*_l, l \in \overline{\DD}$ have been defined in relation 
 (\ref{jokl}).
 
Obviously, the following relation holds, for $i, j \in \DD$,
\begin{equation}\label{simplesr}
_{\bar{k}_{\bar{m}}}\dot{u}_{0, ii j}  + \sum_{l \in \overline{\DD}_{ij}} \, _{\bar{k}_{\bar{m}}}\dot{u}_{0, l i j} = 1. 
\end{equation}

In this case, the following asymptotic relation takes place, for $s \geq 0, i, j \in \DD$, 
\begin{align*}
 \Psi_{\e, \DD, ij}(s / \, _{\bar{k}_{\bar{m}}}\dot{v}_{\e, ij}) & 
 = \phi_{\e, ij}(s  / \, _{\bar{k}_{\bar{m}}}\dot{v}_{\e, ij}) p_{\e, ij} \vspace{2mm} \nonumber \\
 &  \quad + \sum_{l \in \overline{\DD}_{ij}}  \Psi_{\e, \DD, lj}(s  / \, _{\bar{k}_{\bar{m}}}\dot{v}_{\e, ij})  
 \phi_{\e, il}(s / \, _{\bar{k}_{\bar{m}}}\dot{v}_{\e, ij}) p_{\e, il} \vspace{2mm} \nonumber \\
 &  \quad + \sum_{l \in \overline{\DD} \setminus \overline{\DD}_{ij}}  \Psi_{\e, \DD, lj}(s  / \, _{\bar{k}_{\bar{m}}}\dot{v}_{\e, ij})  
 \phi_{\e, il}(s / \, _{\bar{k}_{\bar{m}}}\dot{v}_{\e, ij}) p_{\e, il} \vspace{2mm} \nonumber \\
  & = \phi_{\e, ij}(s  \frac{v_{\e, i}}{_{\bar{k}_{\bar{m}}}\dot{v}_{\e, ij}} / v_{\e, i}) p_{\e, ij} \vspace{2mm} \nonumber \\
 \end{align*}
\begin{align}\label{hopku}
 & \quad + \sum_{l \in \overline{\DD}_{ij}}  
 \Psi_{\e, \DD, lj}(s  \frac{_{\bar{k}_{\bar{m}}} \check{v}_{\e, l}}{_{\bar{k}_{\bar{m}}}\dot{v}_{\e, ij}}
 / \, _{\bar{k}_{\bar{m}}} \check{v}_{\e, l}) \phi_{\e, il}(s  \frac{v_{\e, i}}{_{\bar{k}_{\bar{m}}}\dot{v}_{\e, ij}} / 
 v_{\e, i})p_{\e, il} \makebox[3mm]{} \vspace{2mm} \nonumber \\
 &  \quad + \sum_{l \in \overline{\DD} \setminus \overline{\DD}_{ij}}  \Psi_{\e, \DD, lj}(s  / \, _{\bar{k}_{\bar{m}}}\dot{v}_{\e, ij})  
 \phi_{\e, il}(s / \, _{\bar{k}_{\bar{m}}}\dot{v}_{\e, ij}) p_{\e, il} \vspace{2mm} \nonumber \\
& \to \phi_{0, ij}(_{\bar{k}_{\bar{m}}}\dot{u}_{0}[ii j] s)p_{0, ij} \vspace{2mm} \nonumber \\
& \quad + \sum_{l \in \overline{\DD}_{ij}} \Psi_{0, lj}(_{\bar{k}_{\bar{m}}}\dot{u}_{0}[lij] s) 
\phi_{0, il}(_{\bar{k}_{\bar{m}}}\dot{u}_{0}[ii j]s)p_{0, il} \vspace{2mm} \nonumber \\
& =  \Psi_{0, \DD, ij}(s) \ {\rm as} \ \e \to 0.
 \end{align}
 
 In this case, the corresponding  limiting distribution $G_{0, \DD, i j}(\cdot)$ has, for every 
 $i, j \in \DD$,  the Laplace transform $\Psi_{0, \DD, ij}(s) = 
\int_0^\infty e^{-st} G_{0, \DD, ij}(dt)$, $s \geq 0$ given by relation (\ref{hopku}). 

Relation (\ref{hopku}) implies that, for $s \geq 0$ and $i, j \in \DD$ such that  $P_{0, \DD, ij} > 0$, 
\begin{align}\label{noplinio}
\Phi_{0, \DD, ij}(s) & =  P_{0, \DD, ij}^{-1} \Psi_{0, \DD, ij}(s) = P_{0, \DD, ij}^{-1} 
\big( \phi_{0, ij}(_{\bar{k}_{\bar{m}}}\dot{u}_{0}[ii j]s) \, p_{0, ij}  \vspace{2mm} \nonumber \\
%\end{align*}
%\begin{align}
& \quad   + \sum_{l \in \overline{\DD}_{ij}} \Phi_{0, lj}(_{\bar{k}_{\bar{m}}}\dot{u}_{0}[lij] s) 
\phi_{0, il}(_{\bar{k}_{\bar{m}}}\dot{u}_{0}[ii j]s)P_{0, \DD, lj} \, p_{0, il} \big). 
\end{align}

Note that $\phi_{0, ij}(s), \Phi_{0, lj}( s), \phi_{0, il}(s) < 1$, for $s > 0, i, j \in \DD, l \in \overline{\DD}$. Also, the assumption, $P_{0, \DD, ij} > 0$, implies that, either (a) $p_{0, ij} > 0$ and 
$\overline{\DD}_{ij} = \emptyset$, or
(b) $p_{0, ij} = 0$ and $\overline{\DD}_{ij} \neq \emptyset$, or (c) $p_{0, ij} > 0$ and $\overline{\DD}_{ij} \neq \emptyset$.
In the case (a), $_{\bar{k}_{\bar{m}}}\dot{u}_{0}[ii j] = 1$. In the case (b), 
$_{\bar{k}_{\bar{m}}}\dot{u}_{0}[ii j] = 0$ and there exists $l \in \overline{\DD}_{ij}$ such that 
$_{\bar{k}_{\bar{m}}}\dot{u}_{0}[li j]  > 0$. In the case (c), either $_{\bar{k}_{\bar{m}}}\dot{u}_{0}[ii j] > 0$ or there exists $l \in \overline{\DD}_{ij}$ such that $_{\bar{k}_{\bar{m}}}\dot{u}_{0}[li j]  > 0$. 

 These remarks and relation (\ref{noplinio}) imply that, 
for $s > 0$ and $i, j \in \DD$ such that  $P_{0, \DD, ij} > 0$, 
\begin{equation}\label{noplikuneba}
\Phi_{0, \DD, ij}(s)  <  P_{0, \DD, ij}^{-1} \big( p_{0, ij}  + 
\sum_{i \in \overline{\DD}} P_{0, \DD, lj} p_{0, il} \big) = 1. 
\end{equation}

The above remarks imply that the following theorem takes place.  \vspace{1mm}

{\bf Theorem 8}. {\em Let conditions ${\bf A}$ -- ${\bf E}$, ${\bf \tilde{C}}$, $_{\bar{k}_{r}}{\bf \tilde{C}}$, $_{\bar{k}_{r}}{\bf \hat{C}}$, $_{\bar{k}_{r-1}}{\bf \tilde{F}}$, $r = 1, \ldots, \bar{m} -1$ hold, where the  states $k_r \in \, _{\bar{k}_{r-1}} \overline{\DD}^*, r = 1, \ldots,  \bar{m} -1$ are chosen in such way that condition $_{\bar{k}_{r}}{\bf \hat{F}}$ holds, for  $r = 1, \ldots,  \bar{m} -1$. Let also conditions ${\bf \dot{B}}$, 
${\bf \dot{C}}$, ${\bf \dot{D}}$ and $_{\bar{k}_{\bar{m}}}{\bf \dot{F}}$
hold. Then,  the following asymptotic relation takes place, 
for $i j \in \DD$ such that  $P_{0, \DD, ij} > 0$,   
\begin{equation}\label{cotyrevanertopi}
G_{\e, \DD, i j}(\cdot \, _{\bar{k}_{\bar{m}}}\dot{v}_{\e, ij}) 
\Rightarrow  G_{0, \DD, i j}(\cdot)  = F_{0, \DD, i j}(\cdot) P_{0, \DD, i j} \ {\rm as} \ \e \to 0, 
\end{equation}
where the limiting distribution function $F_{0, \DD, i j}(\cdot)$ is not concentrated in zero, i.e.,  
$F_{0, \DD, i j}(0) < 1$ and has the Laplace transform $\Phi_{0, \DD, i j}(s), 
s \geq 0$, for $i, j \in \DD$.} 

\vspace{1mm}

{\bf Remark 18}.  Conditions ${\bf H}$ and ${\bf \dot{H}}$  are sufficient for holding of conditions  
${\bf C}$, ${\bf \tilde{C}}$, $_{\bar{k}_{r}}{\bf \tilde{C}}$, $_{\bar{k}_{r}}{\bf \hat{C}}$, $_{\bar{k}_{r-1}}{\bf \tilde{F}}$, $r = 1, \ldots, \bar{m} -1$,  ${\bf \dot{C}}$, and $_{\bar{k}_{\bar{m}}}{\bf \dot{F}}$. \\

{\bf 9. Asymptotics for Expectations of Hitting Times} \\

In this section, we present results concerning convergence of expectations of hitting times for regularly and singularly perturbed semi-Markov processes. These results supplement results concerning weak asymptotics of distributions of  hitting times for regularly and singularly perturbed semi-Markov processes. \vspace{1mm}

{ \bf 9.1  Recurrent relations for expectations for hitting times}. Let us introduce expectations, for $i \in \XX, j \in \DD$ and $\e \in [0, 1]$, 
\begin{equation}\label{expart}
E_{\e, \DD, ij} = \int_0^\infty t G_{\e, \DD, ij} (dt),
\end{equation}
and, for $i \in \overline{\DD}, j \in \DD, n = 0, \dots, \bar{m}-1$ and $\e \in [0, 1]$, 
\begin{equation}\label{expartoki}
_{\bar{k}_{n}}E_{\e, \DD, ij} =  \int_0^\infty t \, _{\bar{k}_{n -1}}G_{\e, \DD, ij} (dt). 
\end{equation}
and
\begin{equation}\label{exparture}
_{\bar{k}_{n}}\tilde{E}_{\e, \DD, ij} =  \int_0^\infty t \, _{\bar{k}_{n -1}}\tilde{G}_{\e, \DD, ij} (dt). 
\end{equation}

The following lemmas supplement Lemmas  28 -- 33. As in the above lemmas, let $\bar{k}_{\bar{m}} = \langle k_1, k_2, \ldots, k_{\bar{m}} \rangle$ be an arbitrary sequence of different states from domain $\overline{\DD}$. \vspace{1mm}

{\bf Lemma 36}. {\em Let conditions $_{\bar{k}_0}{\bf A}$, $_{\bar{k}_0}{\bf B}$, and $_{\bar{k}_0}{\bf E}$ {\bf (a)}   hold. Then, for every $n = 0, \ldots, \bar{m} -1$, conditions   $_{\bar{k}_n}{\bf A}$,  $_{\bar{k}_n}{\bf B}$, $_{\bar{k}_n}{\bf \tilde{A}}$, 
$_{\bar{k}_n}{\bf \tilde{B}}$, $_{\bar{k}_{n}}{\bf E}$ {\bf (a)}, $_{\bar{k}_n}{\bf \tilde{E}}$ {\bf (a)} hold, and the following relations take place, for  $\e \in (0, 1]${\rm :} 
\begin{equation}\label{aokom}
E_{\e, \DD,  ij} =  \, _{\bar{k}_n}E_{\e, \DD, ij}  = \, _{\bar{k}_{n}}\tilde{E}_{\e, \DD,  ij} < \infty,  j \in \DD, i \in \, _{\bar{k}_n}\overline{\DD}.         
\end{equation}}
 {\bf Lemma 37}. {\em Let conditions $_{\bar{k}_0}{\bf A}$, $_{\bar{k}_0}{\bf B}$, $_{\bar{k}_0}{\bf E}$ {\bf (a)}  hold. Then, for every $n = 1, \ldots, \bar{m} -1$, conditions $_{\bar{k}_{n-1}}{\bf \tilde{A}}$, $_{\bar{k}_{n-1}}{\bf \tilde{B}}$, $_{\bar{k}_n}{\bf A}$,  $_{\bar{k}_n}{\bf B}$, 
$_{\bar{k}_{n-1}}{\bf \tilde{E}}$ {\bf (a)}, $_{\bar{k}_n}{\bf E}$ {\bf (a)} hold, and the following relations take place, for  $\e \in (0, 1]$,  
\begin{equation}\label{anomaom}
E_{\e, \DD,  ij}  = \, _{\bar{k}_{n-1}}\tilde{E}_{\e, \DD,  ij} = 
\, _{\bar{k}_n}E_{\e, \DD, ij}  < \infty,  j \in \DD, i \in \, _{\bar{k}_{n}}\overline{\DD}, 
\end{equation}
and
\begin{align}\label{nomome}
& \quad E_{\e, \DD, k_n j}  = \, _{\bar{k}_{n-1}}\tilde{E}_{\e, \DD, k_n j} = \, _{\bar{k}_{n-1}}\tilde{e}_{\e, k_n j} \, 
_{\bar{k}_{n-1}}\tilde{p}_{\e,  k_n j}  \makebox[22mm]{} \vspace{3mm} \nonumber \\
&  \quad \quad    + \sum_{r \in \, _{\bar{k}_n}\overline{\DD}}  (_{\bar{k}_{n-1}}\tilde{e}_{\e, k_n r} +  \, _{\bar{k}_n}E_{\e, \DD, rj}) \, _{\bar{k}_{n-1}}\tilde{p}_{\e, k_n r} \vspace{2mm} \nonumber \\
%\end{align}
%\begin{align}
&  \quad \quad   = \, _{\bar{k}_{n-1}}\tilde{e}_{\e, k_n j} \, _{\bar{k}_{n-1}}\tilde{p}_{\e,  k_n j}  \vspace{2mm} \nonumber \\
& \quad \quad     +  \sum_{r \in \, _{\bar{k}_n}\overline{\DD}} ( _{\bar{k}_{n-1}}\tilde{e}_{\e, k_n r} + E_{\e, \DD, rj}) \, 
_{\bar{k}_{n-1}}\tilde{p}_{\e, k_n r}  < \infty,  j \in \DD, i \in \, _{\bar{k}_{n}}\overline{\DD}. 
\end{align}}

\makebox[3mm]{} {\bf Proof}. Condition $_{\bar{k}_0}{\bf B}$ implies that probability $_{\bar{k}_0}p_{\e, ii} < 1$, for $i \in \overline{\DD}$ and $\e \in (0, 1]$. Thus, conditions $_{\bar{k}_0}{\bf A}$, $_{\bar{k}_0}{\bf B}$, $_{\bar{k}_0}{\bf E}$ {\bf (a)} and relations (\ref{twet}), (\ref{twetok}) imply that conditional expectation $_{\bar{k}_0}\tilde{e}_{\e, ij} < \infty$, for $j \in \XX, i \in \overline{\DD}$ and $\e \in (0, 1]$.

Also, conditions $_{\bar{k}_0}{\bf A}$, $_{\bar{k}_0}{\bf B}$,  $_{\bar{k}_0}{\bf E}$ {\bf (a)} and relations (\ref{diotrask}), (\ref{againab})  imply that conditional expectation $_{\bar{k}_1}e_{\e, ij} < \infty$, for  $j \in \XX, i \in \, _{\bar{k}_1}\overline{\DD}, k_1 \in \overline{\DD}$ and $\e \in (0, 1]$.

By repeating recurrent application procedures (described in Sections 3 -- 5) of removing virtual transitions and exclusion states belonging to domain $\overline{\DD}$ from the phase space $\XX$, we get that conditional expectations    
$_{\bar{k}_{n-1}}\tilde{e}_{\e, ij}   < \infty$, for $j \in  \, _{\bar{k}_{n-1}}\XX, i \in \, _{\bar{k}_{n-1}}\DD, n = 1, \ldots, \bar{m} -1$
and  $\e \in (0, 1]$.

Relations given in Lemmas 36 and 37 follow from relations given in Lemmas 29 and 32, 
for both cases, where the 
corresponding expectations $E_{\e, \DD, k_{\bar{m}} j}$ are finite 
nite.
However, relation (\ref{cotyre}) and Lemma 24 imply  that, for $j \in \DD$ and $\e \in (0, 1]$, 
\begin{equation}\label{gopko}
E_{\e, \DD, k_{\bar{m}} j} = \, _{\bar{k}_{\bar{m}-1}}\tilde{e}_{\e, k_{\bar{m}} j} 
\ _{\bar{k}_{\bar{m}-1}}\tilde{p}_{\e,  k_{\bar{m}} j}  < \infty.  
\end{equation}

Relation (\ref{nomome}), taken for  $n = \bar{m} - 1$,  implies that, for $j \in \DD$ and $\e \in (0, 1]$, 
\begin{align}\label{gopkom}
 E_{\e, \DD, k_{\bar{m} -1} j} & = \, _{\bar{k}_{\bar{m}-2}}\tilde{e}_{\e, k_{\bar{m}-1} j} \, 
 _{\bar{k}_{\bar{m}-2}}\tilde{p}_{\e,  k_{\bar{m} -1} j} \vspace{2mm} \nonumber \\
 & \quad +   (_{\bar{k}_{\bar{m} -2}}\tilde{e}_{\e, k_{\bar{m} -1} k_{\bar{m}}} + E_{\e, \DD, k_{\bar{m}}j}) \, 
_{\bar{k}_{\bar{m}-2}}\tilde{p}_{\e, k_{\bar{m} - 1} k_{\bar{m}}}  < \infty.
\end{align}

By continuing the recurrent use of relation (\ref{gopko}), for $n =  \bar{m} - 2, \dots$, we prove that all expectations 
$E_{\e, \DD, k_n j} < \infty, j \in \DD, n = \bar{m}, \ldots, 1$ and $\e \in (0, 1]$.  
$\Box$
\vspace{1mm}

{ \bf 9.2  Asymptotics for expectations of hitting times in the case with the most absorbing initial state}. Let us consider the case, where the initial state is the most absorbing state $k_m$. 

In this case, we can use the normalisation function, given in relation (\ref{norm}), i.e., the same  as in 
the corresponding weak convergence relation for hitting times given in Theorem 4,
\begin{equation}\label{hopun}
_{\bar{k}_{\bar{m}}} \bar{v}_{\e, k_{\bar{m}}} = \, _{\bar{k}_{\bar{m}}} \check{v}_{\e, k_{\bar{m}}} = \, _{\bar{k}_{\bar{m}-1}}\tilde{v}_{\e, k_{\bar{m}}}.
\end{equation} 

The following theorem follows from relation (\ref{gopko}) and Lemma 23. \vspace{1mm}

{\bf Theorem 9}. {\em Let conditions ${\bf A}$ -- ${\bf E}$, ${\bf \tilde{C}}$, $_{\bar{k}_{r}}{\bf \tilde{C}}$, $_{\bar{k}_{r}}{\bf \hat{C}}$, $_{\bar{k}_{r-1}}{\bf \tilde{F}}$, $r = 1, \ldots, \bar{m} -1$ hold, where the   
states $k_r \in \, _{\bar{k}_{r-1}} \overline{\DD}^*, r = 1, \ldots,  \bar{m} -1$ are chosen in such way that condition
$_{\bar{k}_{r}}{\bf \hat{F}}$ holds, for  $r = 1, \ldots,  \bar{m} -1$. Then,  the following asymptotic 
relation takes place, for $j \in \DD$, 
\begin{align}\label{cotynerva}
E_{\e, \DD, k_{\bar{m}} j} / \, _{\bar{k}_{\bar{m}}} \bar{v}_{\e, k_{\bar{m}} } & \to  \bar{E}_{0, \DD, k_{\bar{m}} j} = E_{0, \DD, k_{\bar{m}} j} 
\vspace{1mm} \nonumber \\ 
& = \, _{\bar{k}_{\bar{m}-1}}\tilde{e}_{0, k_{\bar{m}}j}
\, _{\bar{k}_{\bar{m}-1}}\tilde{p}_{0, k_{\bar{m}}j} < \infty \ {\rm as} \ \e \to 0. 
\end{align}} 
\makebox[3mm]{} Note that the limit $E_{0, \DD, k_{\bar{m}} j}$ in the
asymptotic relation (\ref{cotynerva}) is the first moment for the limiting distribution  $G_{\e, \DD, k_{\bar{m}} j}(\cdot)$ given in Theorem 4, for $j \in \DD$. Note also that
$E_{0, \DD, k_{\bar{m}} j} \in (0, \infty)$, for $j \in \DD$ such that  
$_{\bar{k}_{\bar{m}-1}}\tilde{p}_{0, k_{\bar{m}}j} > 0$. 
\vspace{1mm}

{ \bf 9.3  Asymptotics for expectations of hitting times in the case of the second most absorbing initial state}. In this section, we assume that conditions of Theorem 9 hold.

Let us  consider the case, where the initial state is the second most absorbing state $k_{m -1}$.  

Let us define the following normalisation function,
\begin{align}\label{vopuer}
_{\bar{k}_{\bar{m}}}\bar{v}_{\e, k_{\bar{m} -1}} 
& = \, _{\bar{k}_{\bar{m}-2}}\tilde{v}_{\e, k_{\bar{m}- 1}} 
+ \, _{\bar{k}_{\bar{m}}}\bar{v}_{\e, k_{\bar{m} }} \ _{\bar{k}_{\bar{m} - 2}}\tilde{p}_{\e,  k_{\bar{m} - 1} k_{\bar{m}}} \vspace{2mm} \nonumber \\
& = \, _{\bar{k}_{\bar{m} -2}}\tilde{v}_{\e, k_{\bar{m}- 1}} 
+ \, _{\bar{k}_{\bar{m}-1}}\tilde{v}_{\e, k_{\bar{m}}} \ _{\bar{k}_{\bar{m} - 2}}
\tilde{p}_{\e,  k_{\bar{m} - 1} k_{\bar{m}}} 
\end{align}

Obviously, $_{\bar{k}_{\bar{m}}}\bar{v}_{\e, k_{\bar{m} -1}}  \in [1, \infty), \, \e \in (0, 1]$. 

Let us assume that the following asymptotic comparability condition holds:
\begin{equation}\label{vopuerta}
\frac{_{\bar{k}_{\bar{m}}}\bar{v}_{\e, k_{\bar{m} }} \ 
_{\bar{k}_{\bar{m} - 2}}\tilde{p}_{\e,  k_{\bar{m} - 1} k_{\bar{m}}}}{_{\bar{k}_{\bar{m}}}\bar{v}_{\e, k_{\bar{m}- 1}}} 
 \to \, _{\bar{k}_{\bar{m}}}\bar{u}_0[k_{\bar{m} - 1} k_{\bar{m}}] \in [0, 1] \ {\rm as} \ \e \to 0.
\end{equation}

It is useful to note that condition ${\bf H}$ is sufficient for holding of condition (\ref{vopuerta}). 

Relation (\ref{gopkom}), (\ref{cotynerva}),  and condition (\ref{vopuerta})  imply that the following relation 
takes place,  for $j \in \DD$, 
\begin{align*}
\frac{E_{\e, \DD, k_{\bar{m} -1} j}}{_{\bar{k}_{\bar{m}}}\bar{v}_{\e, k_{\bar{m} -1}}}  
&  =  \frac{_{\bar{k}_{\bar{m}-2}}\tilde{v}_{\e, k_{\bar{m} -1}}}{_{\bar{k}_{\bar{m}}}\bar{v}_{\e, k_{\bar{m} -1}}} \frac{_{\bar{k}_{\bar{m} -2}}\tilde{e}_{\e, k_{\bar{m} - 1} j}}{_{\bar{k}_{\bar{m}-2}}\tilde{v}_{\e, k_{\bar{m} -1}}}  \,  
_{\bar{k}_{\bar{m}-2}}\tilde{p}_{\e, k_{\bar{m} - 1} j} 
\makebox[20mm]{}\vspace{2mm} \nonumber \\
\end{align*}
\begin{align}\label{idernent}
& \quad +  
\frac{_{\bar{k}_{\bar{m}-2}}\tilde{v}_{\e, k_{\bar{m} -1}}}{_{\bar{k}_{\bar{m}}}\bar{v}_{\e, k_{\bar{m} -1}}} \, 
\frac{_{\bar{k}_{\bar{m} -2}}\tilde{e}_{\e,  k_{\bar{m} -1} k_{\bar{m}}}}{_{\bar{k}_{\bar{m}-2}}\tilde{v}_{\e, k_{\bar{m} -1}}} \, 
_{\bar{k}_{\bar{m}-2}}\tilde{p}_{\e,  k_{\bar{m} - 1} k_{\bar{m}}}
\vspace{2mm} \nonumber \\
& \quad   + \frac{E_{\e, \DD, k_{\bar{m}} j}}
{_{\bar{k}_{\bar{m}}}\bar{v}_{\e, k_{\bar{m}}}} \frac{_{\bar{k}_{\bar{m}}}\bar{v}_{\e, k_{\bar{m}}} \, _{\bar{k}_{\bar{m} - 2}}\tilde{p}_{\e,  k_{\bar{m} - 1} k_{\bar{m}}}}{_{\bar{k}_{\bar{m}}}\bar{v}_{\e, k_{\bar{m} -1}} } 
\vspace{2mm} \nonumber \\
&  \to (1 - \,  _{\bar{k}_{\bar{m}}}\bar{u}_0[k_{\bar{m} - 1} k_{\bar{m}}] ) \,
_{\bar{k}_{\bar{m} -2}}\tilde{e}_{0,   k_{\bar{m} -1} j} \, 
_{\bar{k}_{\bar{m}-2}}\tilde{p}_{0, k_{\bar{m} - 1} j} 
\vspace{2mm} \nonumber \\
& \quad   +  (1 - \,  _{\bar{k}_{\bar{m}}}\bar{u}_0[k_{\bar{m} - 1} k_{\bar{m}}] ) \,
_{\bar{k}_{\bar{m}-2}}\tilde{e}_{0,  k_{\bar{m} -1} k_{\bar{m}}} \,
_{\bar{k}_{\bar{m}-2}}\tilde{p}_{0,  k_{\bar{m} - 1} k_{\bar{m}}}
\vspace{2mm} \nonumber \\
& \quad   + \,  _{\bar{k}_{\bar{m}}}\bar{u}_0[k_{\bar{m} - 1} k_{\bar{m}}] \bar{E}_{0, \DD, k_{\bar{m}} j}  
\vspace{2mm} \nonumber \\
&   = \bar{E}_{0, \DD,  k_{\bar{m} -1} j}   < \infty \ {\rm as} \ \e \to 0.
\end{align}

Note that $\bar{E}_{0, \DD,  k_{\bar{m} -1} j}  \in  (0, \infty)$ if $\bar{E}_{0, \DD,  k_{\bar{m}} j}  \in  (0, \infty)$.

In the following theorem, limit $\bar{E}_{0, \DD, k_{\bar{m} -1} j}$ is given by relation (\ref{idernent}). \vspace{1mm} 

{\bf Theorem 10}. {\em Let conditions ${\bf A}$ -- ${\bf E}$, ${\bf \tilde{C}}$, $_{\bar{k}_{r}}{\bf \tilde{C}}$, $_{\bar{k}_{r}}{\bf \hat{C}}$, $_{\bar{k}_{r-1}}{\bf \tilde{F}}$, $r = 1, \ldots, \bar{m} -1$ hold, where the   
states $k_r \in \, _{\bar{k}_{r-1}} \overline{\DD}^*, r = 1, \ldots,  \bar{m} -1$ are chosen in such way that condition
$_{\bar{k}_{r}}{\bf \hat{F}}$ holds, for  $r = 1, \ldots,  \bar{m} -1$, and the additional condition {\rm (\ref{vopuerta})}  holds. Then,  the following asymptotic 
relation takes place, for $j \in \DD$, 
\begin{equation}\label{cotynervamol}
E_{\e, \DD, k_{\bar{m} -1} j} / \, _{\bar{k}_{\bar{m}}}\bar{v}_{\e, k_{\bar{m} -1} }  \to  \bar{E}_{0, \DD, k_{\bar{m} -1} j} < \infty \ {\rm as} \ \e \to 0. 
\end{equation}}
\vspace{-1mm}  

{ \bf 9.4  Conditions of simultaneous weak convergence of distributions and expectations for hitting times}.   
The question arises about conditions, under which expectations $E_{\e, \DD, k_{\bar{m} -1} j}$, normalised  by  the function $_{\bar{k}_{\bar{m}}} \check{v}_{\e, k_{\bar{m} -1}}$ used as the normalisation function in the weak convergence relation for distributions $G_{\e, \DD, k_{\bar{m} -1} j}(\cdot)$  given in Theorem 5, converge to the first moment of the corresponding limiting distribution $G_{0, \DD, k_{\bar{m} -1} j}(\cdot)$.

First, let us consider the case, where $k_{\bar{m}} \, \in \ _{\bar{k}_{\bar{m} - 2}}\tilde{\YY}_{1, k_{\bar{m} - 1}}$ and $_{\bar{k}_{\bar{m} - 2}}\tilde{p}_{0,  k_{\bar{m} - 1} k_{\bar{m}}}$ $> 0$, i.e., the following condition holds:
\begin{equation}\label{cond}
_{\bar{k}_{\bar{m} - 2}}\tilde{p}_{\e,  k_{\bar{m} - 1} k_{\bar{m}}} > 0, \  {\rm for} \ \e \in [0, 1].
\end{equation}

In this case, condition (\ref{vopuerta}) holds. Indeed,
\begin{align*}
& \frac{_{\bar{k}_{\bar{m}}}\bar{v}_{\e, k_{\bar{m} }} \ 
_{\bar{k}_{\bar{m} - 2}}\tilde{p}_{\e,  k_{\bar{m} - 1} 
k_{\bar{m}}}}{_{\bar{k}_{\bar{m}}}\bar{v}_{\e, k_{\bar{m}- 1}}} 
=   \big(1 + \frac{_{\bar{k}_{\bar{m}-2}}\tilde{v}_{\e, k_{\bar{m}- 1}}}{_{\bar{k}_{\bar{m}-1}}\tilde{v}_{\e, k_{\bar{m} }} \ _{\bar{k}_{\bar{m} - 2}}\tilde{p}_{\e,  k_{\bar{m} - 1} k_{\bar{m}}}} \big)^{-1}
\vspace{2mm} \nonumber \\
& \quad \quad =  \big(1 + \frac{(1 - \, _{\bar{k}_{\bar{m} - 1}}p_{\e, k_{\bar{m}} k_{\bar{m}}}) \, _{\bar{k}_{\bar{m}-2}}\tilde{v}_{\e, k_{\bar{m}- 1}}}{_{\bar{k}_{\bar{m}-2}}\tilde{v}_{\e, k_{\bar{m} }} \ _{\bar{k}_{\bar{m} - 2}}\tilde{p}_{\e,  k_{\bar{m} - 1} k_{\bar{m}}}} \big)^{-1}
\vspace{2mm} \nonumber \\
\end{align*}
\begin{align}\label{vopuertok}
& \quad \quad \to  \big( 1 + \frac{(1 - \, _{\bar{k}_{\bar{m} - 1}}p_{0, k_{\bar{m}} k_{\bar{m}}}) \, 
_{\bar{k}_{\bar{m}-2}}\tilde{w}_{0, k_{\bar{m} - 1} k_{\bar{m}} }}{_{\bar{k}_{\bar{m} - 2}}\tilde{p}_{0,  k_{\bar{m} - 1} k_{\bar{m}} }} \big)^{-1}
\vspace{2mm} \nonumber \\
& \quad \quad =  \, _{\bar{k}_{\bar{m}}}\bar{u}_0[k_{\bar{m} - 1} k_{\bar{m}}] \in (0, 1] \ {\rm as} \ \e \to 0.
\end{align}

In the case where assumption (\ref{cond}) holds, the normalisation function 
$_{\bar{k}_{\bar{m}}} \check{v}_{\e, k_{\bar{m} -1}} = \, _{\bar{k}_{\bar{m}-1}}\tilde{v}_{\e, k_{\bar{m}}}$ is used in the weak convergence relation given in Theorem 5. Obviously,
\begin{align}\label{bopl}
\frac{_{\bar{k}_{\bar{m}}}\bar{v}_{\e, k_{\bar{m} -1}}}{_{\bar{k}_{\bar{m}}} \check{v}_{\e, k_{\bar{m} -1}}}
& = \frac{_{\bar{k}_{\bar{m}}}\bar{v}_{\e, k_{\bar{m} -1}}}{_{\bar{k}_{\bar{m}-1}}\tilde{v}_{\e, k_{\bar{m}}}}
\vspace{2mm} \nonumber \\ 
& = \frac{_{\bar{k}_{\bar{m}}}\bar{v}_{\e, k_{\bar{m} -1}}
\, _{\bar{k}_{\bar{m}-2}}\tilde{p}_{\e, k_{\bar{m} - 1} k_{\bar{m}}} }{_{\bar{k}_{\bar{m}-1}}\tilde{v}_{\e, k_{\bar{m}}}  
\, _{\bar{k}_{\bar{m}-2}}\tilde{p}_{\e, k_{\bar{m} - 1} k_{\bar{m}}}}
\vspace{2mm} \nonumber \\
& \to \frac{_{\bar{k}_{\bar{m}-2}}\tilde{p}_{0, k_{\bar{m} - 1} k_{\bar{m}}} }{_{\bar{k}_{\bar{m}}}\bar{u}_0[k_{\bar{m} - 1} k_{\bar{m}}]} \ {\rm as} \ \e \to 0.
\end{align}

Also, relation (\ref{vopuertok}) implies that,
\begin{equation}\label{fureq}
\frac{1 - \, _{\bar{k}_{\bar{m}}}\bar{u}_0[k_{\bar{m} - 1} k_{\bar{m}}]}{_{\bar{k}_{\bar{m}}}\bar{u}_0[k_{\bar{m} - 1} k_{\bar{m}}]} =
\frac{(1 - \, _{\bar{k}_{\bar{m} - 1}}p_{0, k_{\bar{m}} k_{\bar{m}}}) \, 
_{\bar{k}_{\bar{m}-2}}\tilde{w}_{0, k_{\bar{m} - 1} k_{\bar{m}} }}{_{\bar{k}_{\bar{m} - 2}}\tilde{p}_{0,  k_{\bar{m} - 1} k_{\bar{m}} }}.
\end{equation}

Finally, relations (\ref{cotynerva}), (\ref{idernent}), (\ref{bopl}), and (\ref{fureq}) imply that, in the case where assumption (\ref{cond}) holds, the following relation takes place, for $j \in \DD$,
\begin{align}\label{idernentolp}
& \frac{E_{\e, \DD, k_{\bar{m} -1} j}}{_{\bar{k}_{\bar{m}}} \check{v}_{\e, k_{\bar{m} -1}}} 
= \frac{E_{\e, \DD, k_{\bar{m} -1} j}}{_{\bar{k}_{\bar{m}}}\bar{v}_{\e, k_{\bar{m} -1}}}
\frac{_{\bar{k}_{\bar{m}}}\bar{v}_{\e, k_{\bar{m} -1}}}{_{\bar{k}_{\bar{m}}} \check{v}_{\e, k_{\bar{m} -1}}} \vspace{2mm} \nonumber \\
& \quad \to \bar{E}_{0, \DD,  k_{\bar{m} -1} j} \frac{_{\bar{k}_{\bar{m}-2}}\tilde{p}_{0, k_{\bar{m} - 1} k_{\bar{m}}} }{_{\bar{k}_{\bar{m}}}\bar{u}_0[k_{\bar{m} - 1} k_{\bar{m}}]} \makebox[50mm]{} \vspace{2mm} \nonumber \\
%\end{align*}
%\begin{align}
& \quad = (1 - \, _{\bar{k}_{\bar{m} - 1}}p_{0, k_{\bar{m}} k_{\bar{m}}}) \, 
_{\bar{k}_{\bar{m}-2}}\tilde{w}_{0, k_{\bar{m} - 1} k_{\bar{m}} } \, 
_{\bar{k}_{\bar{m} -2}}\tilde{e}_{0,   k_{\bar{m} -1} j} \, _{\bar{k}_{\bar{m}-2}}\tilde{p}_{0, k_{\bar{m} - 1} j} \vspace{2mm} \nonumber \\
&  \quad \quad + (1 - \, _{\bar{k}_{\bar{m} - 1}}p_{0, k_{\bar{m}} k_{\bar{m}}})  \,
_{\bar{k}_{\bar{m}-2}}\tilde{w}_{0, k_{\bar{m} - 1} k_{\bar{m}} } \,
_{\bar{k}_{\bar{m} -2}}\tilde{e}_{0,  k_{\bar{m} -1} k_{\bar{m}}} \, 
\, _{\bar{k}_{\bar{m}-2}}\tilde{p}_{0,  k_{\bar{m} - 1} k_{\bar{m}}} \vspace{2mm} \nonumber \\
& \quad \quad + E_{0, \DD, k_{\bar{m}} j} \, _{\bar{k}_{\bar{m}-2}}\tilde{p}_{0,  k_{\bar{m} - 1} k_{\bar{m}}} \vspace{2mm} \nonumber \\
& \quad  = E_{0, \DD, k_{\bar{m} -1} j} \ {\rm as} \ \e \to 0.
\end{align}

Note that $E_{0, \DD,  k_{\bar{m} -1} j}  \in (0, \infty)$ if $E_{0, \DD,  k_{\bar{m}} j} \in (0, \infty)$, i.e., if $_{\bar{k}_{\bar{m}-1}}\tilde{p}_{0, k_{\bar{m}} j} > 0$.

Relation (\ref{idernentolp}) means that, in the case where assumption (\ref{cond}) holds, the expectations $E_{\e, \DD, k_{\bar{m} -1} j}$ normalised  by  function $_{\bar{k}_{\bar{m}}} \check{v}_{\e, k_{\bar{m} -1}} =  \, _{\bar{k}_{\bar{m}-1}}\tilde{v}_{\e, k_{\bar{m}}}$ (used as the normalisation function in the weak convergence relation for distributions $G_{\e, \DD, k_{\bar{m} -1} j}(\cdot)$  given in Theorem 5, in the case where probability $_{\bar{k}_{\bar{m} - 2}}\tilde{p}_{0,  k_{\bar{m} - 1} k_{\bar{m}}} > 0$), converge to the first moment of distribution $G_{0, \DD, k_{\bar{m} -1} j}(\cdot)$.

Second, let us assume that state $k_{\bar{m}} \notin \, _{\bar{k}_{\bar{m} - 2}}\tilde{\YY}_{1, k_{\bar{m} - 1}}$, i.e., the following condition holds:
\begin{equation}\label{conda} 
_{\bar{k}_{\bar{m} - 2}}\tilde{p}_{\e,  k_{\bar{m} - 1} k_{\bar{m}}} = 0, \ \e \in [0, 1]. 
\end{equation}

In this case, the normalisation function,
\begin{equation}\label{hopunal}  
_{\bar{k}_{\bar{m}}}\bar{v}_{\e, k_{\bar{m} -1}} = \,  
_{\bar{k}_{\bar{m}}}\check{v}_{\e, k_{\bar{m} -1}} =  \, _{\bar{k}_{\bar{m}-2}}\tilde{v}_{\e, k_{\bar{m}- 1}}.
\end{equation}

Condition (\ref{vopuerta}) obviously holds, with the limit,
\begin{equation}\label{notyre} 
_{\bar{k}_{\bar{m}}}\bar{u}_0[k_{\bar{m} - 1} k_{\bar{m}}] = 0.
\end{equation}

The normalisation function $_{\bar{k}_{\bar{m}}}\bar{v}_{\e, k_{\bar{m} -1}}$ coincides with the normalisation function $_{\bar{k}_{\bar{m}}}\check{v}_{\e, k_{\bar{m} -1}}$ used in the corresponding weak convergence relation for distributions of hitting times given in Theorem 5, in the case where probability
$_{\bar{k}_{\bar{m} - 2}}\tilde{p}_{0,  k_{\bar{m} - 1} k_{\bar{m}}} = 0$. This equality is included in the 
condition (\ref{conda}).

Relation (\ref{idernent})  takes, in this case,  the following form,  for $j \in \DD$, 
\begin{align}\label{idersfase}
\frac{E_{\e, \DD, k_{\bar{m} -1} j}}{_{\bar{k}_{\bar{m}}}\check{v}_{\e, k_{\bar{m}-1}}} 
& =  \frac{_{\bar{k}_{\bar{m} -2}}\tilde{e}_{\e, k_{\bar{m} - 1} j}}{_{\bar{k}_{\bar{m}-2}}\tilde{v}_{\e, k_{\bar{m} -1}}}  \, _{\bar{k}_{\bar{m}-2}}\tilde{p}_{\e, k_{\bar{m} - 1} j} 
\vspace{2mm} \nonumber \\
&  \to \, _{\bar{k}_{\bar{m} -2}}\tilde{e}_{0, k_{\bar{m} - 1} j}  \ _{\bar{k}_{\bar{m}-2}}\tilde{p}_{0, k_{\bar{m} - 1} j}  \vspace{2mm} \nonumber \\
& = \bar{E}_{0, \DD,  k_{\bar{m} -1} j}  = E_{0, \DD,  k_{\bar{m} -1} j} < \infty \ {\rm as} \ \e \to 0.
\end{align}

Note that $\bar{E}_{0, \DD,  k_{\bar{m} -1} j}  \in (0, \infty)$, for $j \in \DD$ such that 
$_{\bar{k}_{\bar{m}-2}}\tilde{p}_{0, k_{\bar{m} - 1} j}  > 0$.

Relation (\ref{idersfase}) means that, in the case where condition (\ref{conda}) holds, the expectations $E_{\e, \DD, k_{\bar{m} -1} j}$ normalised  by  function $_{\bar{k}_{\bar{m}}} \check{v}_{\e, k_{\bar{m} -1}} =  \, _{\bar{k}_{\bar{m}-2}}\tilde{v}_{\e, k_{\bar{m}- 1}}$ (used as the normalisation function in the weak convergence relation for distributions $G_{\e, \DD, k_{\bar{m} -1} j}(\cdot)$  given in Theorem 5, in the case where probability $_{\bar{k}_{\bar{m} - 2}}\tilde{p}_{0,  k_{\bar{m} - 1} k_{\bar{m}}} = 0$), converge to the first moment of  distribution $G_{0, \DD, k_{\bar{m} -1} j}(\cdot)$.

Third, let us assume that state \ $k_{\bar{m}} \, \in \ _{\bar{k}_{\bar{m} - 2}}\tilde{\YY}_{1, k_{\bar{m} - 1}}$ and $_{\bar{k}_{\bar{m} - 2}}\tilde{p}_{0,  k_{\bar{m} - 1} k_{\bar{m}}} = 0$, i.e., the following condition holds,
\begin{equation}\label{condas}
_{\bar{k}_{\bar{m} - 2}}\tilde{p}_{\e,  k_{\bar{m} - 1} k_{\bar{m}}} > 0, \ \e \in (0, 1], \ {\rm while} \ _{\bar{k}_{\bar{m} - 2}}\tilde{p}_{0,  k_{\bar{m} - 1} k_{\bar{m}}} = 0. 
\end{equation}

Note that,  in this case, 
\begin{equation}\label{condasra}
_{\bar{k}_{\bar{m} - 2}}\tilde{p}_{\e,  k_{\bar{m} - 1} k_{\bar{m}}}  \to 0 \ {\rm as} \ \e \to 0. 
\end{equation}

In this case, condition (\ref{vopuerta}) should be assumed to hold. Three cases should be considered, when
the limit $_{\bar{k}_{\bar{m}}}\bar{u}_0[k_{\bar{m} - 1} k_{\bar{m}}]$  equals $0$, or takes value in interval $(0, 1)$, or equals $1$. 

If limit $_{\bar{k}_{\bar{m}}}\bar{u}_0[k_{\bar{m} - 1} k_{\bar{m}}] = 0$, then the asymptotic relation
(\ref{idernent}) takes the following form, for $j \in \DD$,
\begin{align}\label{idernentyw}
\frac{E_{\e, \DD, k_{\bar{m} -1} j}}{_{\bar{k}_{\bar{m}}}\bar{v}_{\e, k_{\bar{m} -1}}}  
&  \to  \, _{\bar{k}_{\bar{m} -2}}\tilde{e}_{0, k_{\bar{m} - 1} j}  \ _{\bar{k}_{\bar{m}-2}}\tilde{p}_{0, k_{\bar{m} - 1} j}  \vspace{2mm} \nonumber \\
& = \bar{E}_{0, \DD,  k_{\bar{m} -1} j}  = E_{0, \DD,  k_{\bar{m} -1} j} < \infty \ {\rm as} \ \e \to 0.
\end{align}

In this case, 
\begin{equation}\label{futyw}
\frac{_{\bar{k}_{\bar{m}}}\bar{v}_{\e, k_{\bar{m}-1}}}{_{\bar{k}_{\bar{m}-2}}\tilde{v}_{\e, k_{\bar{m} -1}}} \to 1 \ {\rm  as} \  \e \to 0. 
\end{equation}

That is why relation (\ref{idernentyw}) implies that, for $j \in \DD$,
\begin{align}\label{idernentywo}
\frac{E_{\e, \DD, k_{\bar{m} -1} j}}{_{\bar{k}_{\bar{m}-2}}\tilde{v}_{\e, k_{\bar{m} -1}}}  
&  \to  \, _{\bar{k}_{\bar{m} -2}}\tilde{e}_{0, k_{\bar{m} - 1} j}  \ _{\bar{k}_{\bar{m}-2}}\tilde{p}_{0, k_{\bar{m} - 1} j}  \vspace{2mm} \nonumber \\
& = \bar{E}_{0, \DD,  k_{\bar{m} -1} j}  = E_{0, \DD,  k_{\bar{m} -1} j} < \infty \ {\rm as} \ \e \to 0.
\end{align}

Relation (\ref{idernentywo}) means that, in the case where condition (\ref{vopuerta})  holds, with the limit \, $_{\bar{k}_{\bar{m}}}\bar{u}_0[k_{\bar{m} - 1} k_{\bar{m}}] = 0$, and condition (\ref{condas}) holds,  expectations $E_{\e, \DD, k_{\bar{m} -1} j}$ normalised  by  function $_{\bar{k}_{\bar{m}}} \check{v}_{\e, k_{\bar{m} -1}} =  \, _{\bar{k}_{\bar{m}-2}}\tilde{v}_{\e, k_{\bar{m}- 1}}$ (used as the normalisation function in the weak convergence relation for distributions $G_{\e, \DD, k_{\bar{m} -1} j}(\cdot)$  given in Theorem 5, in the case where probability $_{\bar{k}_{\bar{m} - 2}}\tilde{p}_{0,  k_{\bar{m} - 1} k_{\bar{m}}} = 0$), converge to the first moment of distribution $G_{0, \DD, k_{\bar{m} -1} j}(\cdot)$.

If limit $_{\bar{k}_{\bar{m}}}\bar{u}_0[k_{\bar{m} - 1} k_{\bar{m}}] \in (0, 1)$, then the asymptotic relation
(\ref{idernent}) takes the following form, for $j \in \DD$,
\begin{align}\label{idernentbyr}
\frac{E_{\e, \DD, k_{\bar{m} -1} j}}{_{\bar{k}_{\bar{m}}}\bar{v}_{\e, k_{\bar{m} -1}}}  
&  \to (1 - \,  _{\bar{k}_{\bar{m}}}\bar{u}_0[k_{\bar{m} - 1} k_{\bar{m}}] ) \,
_{\bar{k}_{\bar{m} -2}}\tilde{e}_{0,   k_{\bar{m} -1} j} \, 
_{\bar{k}_{\bar{m}-2}}\tilde{p}_{0, k_{\bar{m} - 1} j} 
\vspace{2mm} \nonumber \\
& \quad   + \,  _{\bar{k}_{\bar{m}}}\bar{u}_0[k_{\bar{m} - 1} k_{\bar{m}}] \bar{E}_{0, \DD, k_{\bar{m}} j}  
\vspace{2mm} \nonumber \\
&   = \bar{E}_{0, \DD,  k_{\bar{m} -1} j}   < \infty \ {\rm as} \ \e \to 0.
\end{align}

In this case, 
\begin{equation}\label{futywa}
\frac{_{\bar{k}_{\bar{m}}}\bar{v}_{\e, k_{\bar{m}-1}}}{_{\bar{k}_{\bar{m}-2}}\tilde{v}_{\e, k_{\bar{m} -1}}} \to (1 - \, _{\bar{k}_{\bar{m}}}\bar{u}_0[k_{\bar{m} - 1} k_{\bar{m}}])^{-1} \ {\rm as} \ \e \to 0. 
\end{equation}

That is why relation (\ref{idernentbyr}) implies that, for $j \in \DD$,
\begin{align}\label{idernentbyrbo}
\frac{E_{\e, \DD, k_{\bar{m} -1} j}}{_{\bar{k}_{\bar{m}-2}}\tilde{v}_{\e, k_{\bar{m} -1}} }  
&  \to  \,
_{\bar{k}_{\bar{m} -2}}\tilde{e}_{0,   k_{\bar{m} -1} j} \, 
_{\bar{k}_{\bar{m}-2}}\tilde{p}_{0, k_{\bar{m} - 1} j} 
\vspace{2mm} \nonumber \\
& \quad   + \,  \frac{_{\bar{k}_{\bar{m}}}\bar{u}_0[k_{\bar{m} - 1} k_{\bar{m}}]}{1 - \,  _{\bar{k}_{\bar{m}}}\bar{u}_0[k_{\bar{m} - 1} k_{\bar{m}}]} E_{0, \DD, k_{\bar{m}} j}  
 \ {\rm as} \ \e \to 0.
\end{align}

If $E_{0, \DD,  k_{\bar{m}} j} \in (0, \infty)$, i.e., if $_{\bar{k}_{\bar{m}-1}}\tilde{p}_{0, k_{\bar{m}} j} > 0$, then
\begin{align}\label{byrbo}
& _{\bar{k}_{\bar{m} -2}}\tilde{e}_{0,   k_{\bar{m} -1} j} \, 
_{\bar{k}_{\bar{m}-2}}\tilde{p}_{0, k_{\bar{m} - 1} j} 
  + \,  \frac{_{\bar{k}_{\bar{m}}}\bar{u}_0[k_{\bar{m} - 1} k_{\bar{m}}]}{1 - \,  _{\bar{k}_{\bar{m}}}\bar{u}_0[k_{\bar{m} - 1} k_{\bar{m}}]} E_{0, \DD, k_{\bar{m}} j}  \vspace{2mm} \nonumber \\
& \quad \quad  > \, _{\bar{k}_{\bar{m} -2}}\tilde{e}_{0,   k_{\bar{m} -1} j} \, 
_{\bar{k}_{\bar{m}-2}}\tilde{p}_{0, k_{\bar{m} - 1} j}  = E_{0, \DD,  k_{\bar{m} -1} j}.  
\end{align}

Relation (\ref{byrbo}) means that, in the case where condition (\ref{vopuerta})  holds, with the limit \, $_{\bar{k}_{\bar{m}}}\bar{u}_0[k_{\bar{m} - 1} k_{\bar{m}}]$ $\in \,  (0, 1)$, condition (\ref{condas}) holds, and $E_{0, \DD,  k_{\bar{m}} j}$ $\in \, (0, \infty)$,  expectations $E_{\e, \DD, k_{\bar{m} -1} j}$ normalised  by  function $_{\bar{k}_{\bar{m}}} \check{v}_{\e, k_{\bar{m} -1}} =  \, _{\bar{k}_{\bar{m}-2}}\tilde{v}_{\e, k_{\bar{m}- 1}}$ (used as the normalisation function in the weak convergence relation for distributions $G_{\e, \DD, k_{\bar{m} -1} j}(\cdot)$  given in Theorem 5, in the case where probability $_{\bar{k}_{\bar{m} - 2}}\tilde{p}_{0,  k_{\bar{m} - 1} k_{\bar{m}}} = 0$), converge to the limit, which differs of the first moment of distribution $G_{0, \DD, k_{\bar{m} -1} j}(\cdot)$.

If limit $_{\bar{k}_{\bar{m}}}\bar{u}_0[k_{\bar{m} - 1} k_{\bar{m}}] = 1$, then the asymptotic relation
(\ref{idernent}) takes the following form, for $j \in \DD$,
\begin{align}\label{byrnop}
\frac{E_{\e, \DD, k_{\bar{m} -1} j}}{_{\bar{k}_{\bar{m}}}\bar{v}_{\e, k_{\bar{m} -1}}}  
&  \to \bar{E}_{0, \DD, k_{\bar{m}} j}  
\vspace{2mm} \nonumber \\
&   = \bar{E}_{0, \DD,  k_{\bar{m} -1} j}   < \infty \ {\rm as} \ \e \to 0.
\end{align}

In this case, 
\begin{equation}\label{futywan}
\frac{_{\bar{k}_{\bar{m}}}\bar{v}_{\e, k_{\bar{m}-1}}}{_{\bar{k}_{\bar{m}-2}}\tilde{v}_{\e, k_{\bar{m} -1}}} \to \infty
 \ {\rm as} \ \e \to 0. 
\end{equation}

That is why,  relation (\ref{byrnop}) implies that, for $j \in \DD$ such that  
$\bar{E}_{0, \DD,  k_{\bar{m} -1} j}  \in (0, \infty)$,
\begin{equation}\label{byrnopas}
\frac{E_{\e, \DD, k_{\bar{m} -1} j}}{_{\bar{k}_{\bar{m}-2}}\tilde{v}_{\e, k_{\bar{m} -1}} }  
 \to \infty \ {\rm as} \ \e \to 0.
\end{equation}

Relation (\ref{byrnopas}) means that, in the case where condition (\ref{vopuerta})  holds, with the limit $_{\bar{k}_{\bar{m}}}\bar{u}_0[k_{\bar{m} - 1} k_{\bar{m}}] = 1$, condition (\ref{condas}) holds, and $E_{0, \DD,  k_{\bar{m}} j} \in (0, \infty)$, the expectations $E_{\e, \DD, k_{\bar{m} -1} j}$ normalised  by  function $_{\bar{k}_{\bar{m}}} \check{v}_{\e, k_{\bar{m} -1}} =  \, _{\bar{k}_{\bar{m}-2}}\tilde{v}_{\e, k_{\bar{m}- 1}}$ (used as the normalisation function in the weak convergence relation for distributions $G_{\e, \DD, k_{\bar{m} -1} j}(\cdot)$  given in Theorem 5, in the case where probability $_{\bar{k}_{\bar{m} - 2}}\tilde{p}_{0,  k_{\bar{m} - 1} k_{\bar{m}}} = 0$), converge to $\infty$ as $\e \to 0$.
 \vspace{1mm}

{ \bf 9.5  Asymptotics for expectations of hitting times in the case 
 with an arbitrary initial state from domain $\overline{\DD} \setminus \{k_m \}$}.
Let now, consider the general case, with the initial state is $k_{n}$, for some 
$n = 1, \ldots, \bar{m} -1$. 

We can use, in this case, relation (\ref{nomome}) given in Lemma 37. This relation takes, in this case, 
can be written in the following form, for $j \in \DD$ and $n = 1, \ldots, \bar{m} -1$,  
\begin{align}\label{sfatomnem}
E_{\e, \DD, k_{n} j}  & =   \, _{\bar{k}_{n -1}}\tilde{e}_{\e, k_{n} j} \, 
_{\bar{k}_{n - 1}}\tilde{p}_{\e, k_{n} j}  
\vspace{2mm} \nonumber \\ 
&  \quad  + \sum_{n +1 \leq l \leq \bar{m}}     
\ ( \, _{\bar{k}_{n - 1}}\tilde{e}_{\e, k_{n} k_{l}} \vspace{2mm} \nonumber \\ 
&  \quad  \quad + E_{\e, \DD, k_{l} j} ) 
\ _{\bar{k}_{n - 1}}\tilde{p}_{\e,  k_{n} k_{l}}.  
\end{align}

Let us define recurrently  the following normalisation function, for $n = 1, \ldots, \bar{m}$,
\begin{equation}\label{vopueraso}
_{\bar{k}_{\bar{m}}}\bar{v}_{\e, k_{n}}   = \, _{\bar{k}_{\bar{n}- 1}}\tilde{v}_{\e, k_{n}}  
+ \sum_{n + 1 \leq l \leq \bar{m}} \, _{\bar{k}_{\bar{m}}}\bar{v}_{\e, k_{l}} \
_{\bar{k}_{n - 1}}\tilde{p}_{\e,  k_{n} k_{l}}.
\end{equation}

Obviously, $_{\bar{k}_{\bar{m}}}\bar{v}_{\e, k_{n}}  \in [1, \infty), \, \e \in (0, 1]$, for  $n = 1, \ldots, \bar{m} -1$.

Let us assume that the following condition holds: \\
\vspace{2mm}

\noindent $_{\bar{k}_{\bar{m}}}{\bf  \bar{F}}$: $\frac{_{\bar{k}_{\bar{m}}}\bar{v}_{\e, k_{l}} \, 
_{\bar{k}_{n - 1}}\tilde{p}_{\e,  k_{n} k_{l}}}{_{\bar{k}_{\bar{m}}}\bar{v}_{\e, k_{n}}} 
\to \, _{\bar{k}_{\bar{m}}}\bar{u}_{0} [k_n \, k_l] \in [0, 1]$ 
as $\e \to 0$, for $l = n + 1$, $\ldots, \bar{m}$,  \makebox[10mm]{} $n = 1, \ldots, \bar{m} -1$.  
\vspace{2mm}

It is useful to note that condition ${\bf H}$ is sufficient for holding of 
condition $_{\bar{k}_{\bar{m}}}{\bf  \bar{F}}$. 

Condition $_{\bar{k}_{\bar{m}}}{\bf  \bar{F}}$ obviously implies that, for $n = 1, \ldots, \bar{m} -1$,
\begin{align}\label{unmo}
\frac{_{\bar{k}_{\bar{n}- 1}}\tilde{v}_{\e, k_{n}}}{_{\bar{k}_{\bar{m}}}\bar{v}_{\e, k_{n}}} & \to \, _{\bar{k}_{\bar{m}}}\bar{u}_{0} [k_n \, k_n] \vspace{2mm} \nonumber \\
& = 1 - \sum_{l = n +1}^{\bar{m}} \, _{\bar{k}_{\bar{m}}}\bar{u}_{0} [k_n \, k_l] \in  [0, 1] \ {\rm as} \ \e \to 0.
\end{align}

The asymptotic relation (\ref{cotynerva}) and (\ref{cotynervamol})  given, respectively,  in Theorems 9  and 10 can be considered as the results of  first  step in some backward asymptotic recurrent algorithm for computing  limits for expectations of hitting times, 
$E_{\e, \DD, k_{n}, j} / \, _{\bar{k}_{\bar{m}}}\bar{v}_{\e, k_{n} }$, for $n = \bar{m}, \ldots, 1$, where the corresponding normalisation functions are given by relations  (\ref{hopun}) and (\ref{vopuer}). 

Let us assume that, we have already realised $\bar{m} -  n$ steps in this backward algorithm resulted by the following asymptotic 
relations, for $j \in \DD, l = \bar{m}, \ldots, n + 1$, for some $1 \leq n \leq \bar{m} - 1$, 
\begin{equation}\label{resulnol}
E_{\e, \DD, k_{l}, j} / \, _{\bar{k}_{\bar{m}}}\bar{v}_{\e, k_{l} }
\to  \bar{E}_{0, \DD, k_{n} j} < \infty \ {\rm as} \ \e \to 0. 
\end{equation}
with some limits $\bar{E}_{0, \DD, k_{l} j} < \infty,  j \in \DD, l = \bar{m}, \ldots, n + 1$. 

Then, using relation (\ref{sfatomnem}), we get the following relations, for $j \in \DD$, 
\begin{align}\label{sfatomnemas}
 E_{\e, \DD, k_{n} j} / _{\bar{k}_{\bar{m}}}\bar{v}_{\e, k_{n}} & =   
 \frac{_{\bar{k}_{\bar{n}- 1}}\tilde{v}_{\e, k_{n}}}{_{\bar{k}_{\bar{m}}}\bar{v}_{\e, k_{n}}}
\frac{_{\bar{k}_{n -1}}\tilde{e}_{\e, k_{n} j}}{_{\bar{k}_{\bar{n}- 1}}\tilde{v}_{\e, k_{n}}}
\, _{\bar{k}_{n - 1}}\tilde{p}_{\e, k_{n} j}  
\vspace{2mm} \nonumber \\ 
& \quad  + \sum_{n +1 \leq l \leq \bar{m}}
 \frac{_{\bar{k}_{\bar{n}- 1}}\tilde{v}_{\e, k_{n}}}{_{\bar{k}_{\bar{m}}}\bar{v}_{\e, k_{n}}}
\frac{_{\bar{k}_{n -1}}\tilde{e}_{\e, k_{n} l}}{_{\bar{k}_{\bar{n}- 1}}\tilde{v}_{\e, k_{n}}}      
\ _{\bar{k}_{n - 1}}\tilde{p}_{\e,  k_{n} k_{l}} 
\vspace{2mm} \nonumber \\ 
& \quad  + \sum_{n +1 \leq l \leq \bar{m}}     
\ \frac{_{\bar{k}_{\bar{m}}}\bar{v}_{\e, k_{l}} \ _{\bar{k}_{n - 1}}\tilde{p}_{\e,  k_{n} k_{l}}}{_{\bar{k}_{\bar{m}}}\bar{v}_{\e, k_{n}}} 
\frac{ E_{\e, \DD, k_{l} j}}{_{\bar{k}_{\bar{m}}}\bar{v}_{\e, k_{l}}} 
\vspace{2mm} \nonumber \\ 
& = (1 - \, _{\bar{k}_{\bar{m}}}\bar{u}_{0} [k_n \, k_{n}] ) \, 
( \, _{\bar{k}_{n -1}}\tilde{e}_{0, k_{n} j} \, _{\bar{k}_{n - 1}}\tilde{p}_{0, k_{n} j} 
\vspace{2mm} \nonumber \\ 
& \quad  + \sum_{n +1 \leq l \leq \bar{m}}  \, _{\bar{k}_{n -1}}\tilde{e}_{0, k_{n} l} \, _{\bar{k}_{n - 1}}\tilde{p}_{0, k_{n} l} )
\vspace{2mm} \nonumber \\ 
& \quad  + \sum_{n +1 \leq l \leq \bar{m}} \,    _{\bar{k}_{\bar{m}}}\bar{u}_{0} [k_n \, k_{l}] \bar{E}_{0, \DD, k_{l} j} \vspace{2mm} \nonumber \\ 
& = \bar{E}_{0, \DD, k_{n} j}  < \infty \ {\rm as} \ \e \to 0. 
\end{align} 

By induction, relation (\ref{sfatomnemas}) holds for any $j \in \DD, n = \bar{m} -1, \ldots, 1$. 
Moreover, (\ref{sfatomnemas}) gives the explicit recurrent formulas for computing limits 
 $\bar{E}_{0, \DD, k_{n} j} < \infty,  j \in \DD, n = \bar{m} -1, \ldots, 1$. 
 
The following theorem takes place. \vspace{1mm}

{\bf Theorem 11}. {\em Let conditions ${\bf A}$ -- ${\bf E}$, ${\bf \tilde{C}}$, $_{\bar{k}_{r}}{\bf \tilde{C}}$, $_{\bar{k}_{r}}{\bf \hat{C}}$, $_{\bar{k}_{r-1}}{\bf \tilde{F}}$, $r = 1, \ldots, \bar{m} -1$ hold, where the   
states $k_r \in \, _{\bar{k}_{r-1}} \overline{\DD}^*, r = 1, \ldots,  \bar{m} -1$ are chosen in such way that condition 
$_{\bar{k}_{r}}{\bf \hat{F}}$ holds, for  $r = 1, \ldots,  \bar{m} -1$ hold, and additionally
condition $_{\bar{k}_{\bar{m}}}{\bf  \bar{F}}$ holds. Then,  the following asymptotic 
relation takes place, for $j \in \DD, n = 1, \ldots, \bar{m} - 1$, 
\begin{equation}\label{cotynervamoty}
E_{\e, \DD, k_{n} j} / \, _{\bar{k}_{\bar{m}}}\bar{v}_{\e, k_{n} }  \to  \bar{E}_{0, \DD, k_{n} j} 
< \infty \ {\rm as} \ \e \to 0. 
\end{equation}} 
\makebox[3mm]{}{\bf Remark 19}.  Conditions ${\bf H}$ is sufficient for holding of conditions  
${\bf C}$, ${\bf \tilde{C}}$, $_{\bar{k}_{r}}{\bf \tilde{C}}$, $_{\bar{k}_{r}}{\bf \hat{C}}$, $_{\bar{k}_{r-1}}{\bf \tilde{F}}$, $r = 1, \ldots, \bar{m} -1$, and $_{\bar{k}_{\bar{m}}}{\bf \bar{F}}$.
\vspace{1mm}

The normalisation function $_{\bar{k}_{\bar{m}}}\bar{v}_{\e, k_{n}}$ can be expressed in terms of 
normalisation functions  $_{\bar{k}_{l- 1}}\tilde{v}_{\e, k_{l}}, l = n,  \ldots, \bar{m}$, for $n = 1, \ldots, \bar{m}$.

In particular, according relations (\ref{hopun}) and  (\ref{vopuer}), the normalisation functions
$_{\bar{k}_{\bar{m}}}\bar{v}_{\e, k_{\bar{m}}} = \, _{\bar{k}_{\bar{m} -1}}\tilde{v}_{\e, k_{\bar{m}}}$ and  
$_{\bar{k}_{\bar{m}}}\bar{v}_{\e, k_{\bar{m} -1}} =  \, _{\bar{k}_{\bar{m} -2}}\tilde{v}_{\e, k_{\bar{m} -1}}  + 
\, _{\bar{k}_{\bar{m} -1}}\tilde{v}_{\e, k_{\bar{m}}} \,  _{\bar{k}_{\bar{m} - 2}}\tilde{p}_{\e,  k_{\bar{m} -1} k_{\bar{m}}}$.
By continuing the recurrent substitution in relation (\ref{vopueraso}), we get, 
\begin{align}\label{nolky}
_{\bar{k}_{\bar{m}}}\bar{v}_{\e, k_{\bar{m} -2}} & = \, _{\bar{k}_{\bar{m} - 3}}\tilde{v}_{\e, k_{\bar{m} -2}} 
\vspace{2mm} \nonumber \\
& + (_{\bar{k}_{\bar{m} -2}}\tilde{v}_{\e, k_{\bar{m} -1}}  + 
\, _{\bar{k}_{\bar{m} -1}}\tilde{v}_{\e, k_{\bar{m}}} \,  _{\bar{k}_{\bar{m} - 2}}\tilde{p}_{\e,  k_{\bar{m} -1} k_{\bar{m}}})
\,  _{\bar{k}_{\bar{m} - 3}}\tilde{p}_{\e,  k_{\bar{m} -2} k_{\bar{m}-1}} 
\vspace{2mm} \nonumber \\
& + \, _{\bar{k}_{\bar{m} -1}}\tilde{v}_{\e, k_{\bar{m}}} \,  _{\bar{k}_{\bar{m} - 3}}\tilde{p}_{\e,  k_{\bar{m} -2} k_{\bar{m}}} 
\vspace{2mm} \nonumber \\
& =    \, _{\bar{k}_{\bar{m} - 3}}\tilde{v}_{\e, k_{\bar{m} -2}} \vspace{2mm} \nonumber \\
& + \, _{\bar{k}_{\bar{m} -2}}\tilde{v}_{\e, k_{\bar{m} -1}} \,  _{\bar{k}_{\bar{m} - 2}}\tilde{p}_{\e,  k_{\bar{m} -1} k_{\bar{m}}} 
\vspace{2mm} \nonumber \\
& + \, _{\bar{k}_{\bar{m} -1}}\tilde{v}_{\e, k_{\bar{m}}} ( _{\bar{k}_{\bar{m} - 3}}\tilde{p}_{\e,  k_{\bar{m} -2} k_{\bar{m}-1}} 
\,  _{\bar{k}_{\bar{m} - 2}}\tilde{p}_{\e,  k_{\bar{m} -1} k_{\bar{m}}} + \,  _{\bar{k}_{\bar{m} - 3}}\tilde{p}_{\e,  k_{\bar{m} -2} k_{\bar{m}}} ).
\end{align}

Analogous relations can be obtained for $n = \bar{m} - 3, \ldots, 1$.

The following analysis  (concerned conditions, under which expectations $E_{\e, \DD, k_{\bar{m} -1} j}$, normalised  by  the function $_{\bar{k}_{\bar{m}}} \check{v}_{\e, k_{\bar{m} -1}}$ used as the normalisation function in the weak convergence relation for distributions $G_{\e, \DD, k_{\bar{m} -1} j}(\cdot)$  given in Theorem 6, converge to the first moment of the corresponding limiting distribution 
$G_{0, \DD, k_{\bar{m} -1} j}(\cdot)$)  can be performed in the way analogous to those used in Subsection 9.4. 
We, leave this for the future publications. \vspace{1mm} 

{\bf 9.6  Asymptotics for expectations of  hitting times for the case with an initial state in domain $\DD$}. In this subsection we present conditions of convergence for normalised expectations of hitting times  for the case, where the  initial state $i \in \DD$. 

We assume that conditions of Theorems 9 and  10 holds and additionally assume that conditions ${\bf \dot{B}}$, ${\bf \dot{C}}$, ${\bf \dot{D}}$, ${\bf \dot{E}}$, and $_{\bar{k}_{\bar{m}}}{\bf  \dot{F}}$ hold.

The key role is played by the following relation, which takes place, for $i, j \in \DD$ and $\e \in (0,1]$, 
 \begin{align}\label{initabert}
 E_{\e, \DD, ij}  & = e_{\e, ij} p_{\e, ij} 
 + \sum_{k_n \in \overline{\DD}} (E_{\e, \DD, k_n j} + e_{\e, i k_n} ) p_{\e, i k_n} 
 \vspace{2mm} \nonumber \\
 & = e_{\e, ij} p_{\e, ij}  + \sum_{k_n \in \overline{\DD}} e_{\e, i k_n}  p_{\e, i k_n} 
 +  \sum_{k_n \in \overline{\DD}} E_{\e, \DD, k_n j} p_{\e, i k_n}. 
 \end{align}

Note, first of all, that the above assumptions and relation (\ref{initabert}) imply that $E_{\e, \DD, ij} < \infty$, for 
$i, j \in \DD$ and $\e \in (0,1]$. 

Let us introduce the following normalising functions, for $i \in \DD$,  
\begin{equation}\label{nortypirnas}
_{\bar{k}_{\bar{m}}}\ddot{v}_{\e, i} = v_{\e, i} + \sum_{k_n \in \overline{\DD}}  \, 
_{\bar{k}_{\bar{m}}}\bar{v}_{\e, k_n} 
\end{equation}

Let us also assume that the following comparability condition  holds:
\vspace{2mm}

\noindent$_{\bar{k}_{\bar{m}}}{\bf \ddot{F}}$:  
$\frac{_{\bar{k}_{\bar{m}}}\bar{v}_{\e, k_n}}{_{\bar{k}_{\bar{m}}}\ddot{v}_{\e, i}} \to 
\, _{\bar{k}_{\bar{m}}}\ddot{u}_{0}[k_n \,  i] \in [0, \infty]$ as $\e \to 0$, for $k_n \in \overline{\DD},  i \in \DD$.  

\vspace{2mm} 

Condition $_{\bar{k}_{\bar{m}}}{\bf \ddot{F}}$ imply that the following relations holds, for 
$r \in \DD$, 
 \begin{align}\label{oplutynastyk}
\frac{v_{\e, i}}{_{\bar{k}_{\bar{m}}}\ddot{v}_{\e, i}} & =  1 
- \sum_{k_n \in \overline{\DD}} \,  
\frac{_{\bar{k}_{\bar{m}}}\bar{v}_{\e, k_n} }{_{\bar{k}_{\bar{m}}}\ddot{v}_{\e, i}}  
\vspace{2mm} \nonumber \\
& \to 1  - \sum_{k_n \in \overline{\DD}} \, _{\bar{k}_{\bar{m}}}\ddot{u}_{0}[k_n \, i]  = \, _{\bar{k}_{\bar{m}}}\ddot{u}_{0}[ii] \in [0, 1] \ {\rm as} \ \e \to 0. 
 \end{align}
 
Obviously, the following relation, for $i \in \DD$,
\begin{equation}\label{simplesrty}
_{\bar{k}_{\bar{m}}}\ddot{u}_{0}[ii]  + \sum_{k_n \in \overline{\DD}} \, 
_{\bar{k}_{\bar{m}}}\ddot{u}_{0}[k_n \, i] = 1. 
\end{equation}

Theorems 9  and 10, condition $_{\bar{k}_{\bar{m}}}{\bf \ddot{F}}$ and relations (\ref{initabert}) and (\ref{oplutynastyk}) imply that the following relation hold, for $i, j \in \DD$, 
\begin{align}\label{initabertulk}
 E_{\e, \DD, ij} / \, _{\bar{k}_{\bar{m}}}\ddot{v}_{\e, i} & =  \frac{v_{\e, i}}{_{\bar{k}_{\bar{m}}}\ddot{v}_{\e, i}} \frac{e_{\e, ij}}{v_{\e, i}} p_{\e, ij}  
  + \sum_{k_n \in \overline{\DD}}  \frac{v_{\e, i}}{_{\bar{k}_{\bar{m}}}\ddot{v}_{\e, i}} 
\frac{e_{\e, i k_n}}{v_{\e, i}}   p_{\e, i k_n}  \makebox[20mm]{}
\vspace{2mm} \nonumber \\
% \end{align*}
% \begin{align}
& \quad +  \sum_{k_n \in \overline{\DD}} 
\frac{_{\bar{k}_{\bar{m}}}\bar{v}_{\e, k_n}}{_{\bar{k}_{\bar{m}}}\ddot{v}_{\e, i}} 
\frac{E_{\e, \DD, k_n j}}{_{\bar{k}_{\bar{m}}}\bar{v}_{\e, k_n}} p_{\e, i k_n} 
\vspace{2mm} \nonumber \\ 
& \to \, _{\bar{k}_{\bar{m}}}\ddot{u}_{0}[ii] e_{0, ij} 
+ \sum_{k_n \in \overline{\DD}} \, _{\bar{k}_{\bar{m}}}\ddot{u}_{0}[ii] e_{0, r k_n} 
\vspace{2mm} \nonumber \\ 
& \quad + \sum_{k_n \in \overline{\DD}}  \, _{\bar{k}_{\bar{m}}}\ddot{u}_{0} [k_n i] \bar{E}_{0, \DD, k_n j}
 = \ddot{E}_{0, \DD, i j} \ {\rm as} \ \e \to 0. 
 \end{align}
 
 The following theorem takes place. \vspace{1mm}

{\bf Theorem 12}. {\em Let conditions ${\bf A}$ -- ${\bf E}$, ${\bf \tilde{C}}$, $_{\bar{k}_{r}}{\bf \tilde{C}}$, $_{\bar{k}_{r}}{\bf \hat{C}}$, $_{\bar{k}_{r-1}}{\bf \tilde{F}}$, $r = 1, \ldots, \bar{m} -1$ hold, where the   
states $k_r \in \, _{\bar{k}_{r-1}} \overline{\DD}^*, r = 1, \ldots,  \bar{m} -1$ are chosen in such way that condition $_{\bar{k}_{r}}{\bf \hat{F}}$ holds, for  $r = 1, \ldots,  \bar{m} -1$ hold, and, also,   
condition $_{\bar{k}_{\bar{m}}}{\bf  \bar{F}}$, ${\bf \dot{B}}$ -- ${\bf \dot{E}}$, 
$_{\bar{k}_{\bar{m}}}{\bf  \dot{F}}$, and $_{\bar{k}_{\bar{m}}}{\bf  \ddot{F}}$ hold. Then,  the following asymptotic relation takes place, for $i, j \in \DD$, 
\begin{equation}\label{cotynervamolo}
E_{\e, \DD, ij} / \, _{\bar{k}_{\bar{m}}}\ddot{v}_{\e, i }  \to  \ddot{E}_{0, \DD, i j} 
< \infty \ {\rm as} \ \e \to 0. 
\end{equation}} 
\makebox[3mm]{}{\bf Remark 20}.  Conditions ${\bf H}$ and ${\bf \dot{H}}$ are sufficient for holding of conditions  
${\bf C}$, ${\bf \tilde{C}}$, $_{\bar{k}_{r}}{\bf \tilde{C}}$, $_{\bar{k}_{r}}{\bf \hat{C}}$, $_{\bar{k}_{r-1}}{\bf \tilde{F}}$, $r = 1, \ldots, \bar{m} -1$, and $_{\bar{k}_{\bar{m}}}{\bf \bar{F}}$, ${\bf \dot{C}}$, 
$_{\bar{k}_{\bar{m}}}{\bf  \dot{F}}$, and $_{\bar{k}_{\bar{m}}}{\bf  \ddot{F}}$.  \\

{\bf 10. Comments,  Generalisations and Examples} \\

In this section, we comment some natural generalisations of asymptotic results presented in Sections 2 -- 9. In particular, we discuss possibilities of reward interpretation for hitting type functionals as well as consideration of vector and real-valued rewards. In the last subsection, we also present numerical examples illustrating the main theoretical results of the paper. \vspace{1mm} 

{ \bf 10.1  Hitting rewards  and regularity conditions}. Let us begin from the remark that the regularity condition ${\bf A}$ can be, in fact, omitted. In this case, condition ${\bf B}$ still guarantee that relation (\ref{hitta}) holds, i.e. the hitting times $\nu_{\e, \DD}$ and $\tau_{\e, \DD}$ are proper random variables.

In this case, it is possible that the random variable $\nu_{\e}(t)$ can be improper random variables, i.e. take value $\infty$ with positive probabilities, and, thus, the semi-Markov process $\eta_\e (t)$ is not well defined.

Despite of this, one can interpret the random variable $\tau_{\e, \DD}$ as a Markov reward accumulated at  a trajectory of the Markov chain $\eta_{\e, n}$ up to the first hitting of domain $\DD$ by this Markov chain.

The recurrent algorithms presented in Sections 1 -- 9 and the asymptotic results formulated in Theorems 1 -- 12 remain to be valid. Condition ${\bf A}$ can be omitted in the corresponding algorithms and theorems.
\vspace{1mm} 

{ \bf 10.2   Atoms at zero for limiting distributions}.  These  generalisations are concerned condition ${\bf D}$ {\bf (b)}, which requires that the corresponding limiting distributions for transition times are not concentrated at zero.

This condition, together with other basic conditions, guarantees that the corresponding limiting conditional distributions of hitting times $F_{0, \DD, ij}(\cdot)$ also are not concentrated at zero. 

%One of the possible generalisaion is connected with replacing condition ${\bf D}$ {\bf (b)} by stronger %condition, which would require that limiting distribution functions of hitting times $F_{0,  ij}(\cdot)$  have %not atoms at zero, i.e., $F_{0, ij}(0) = 0, j \in \YY_{1, i}, i \in \overline{\DD}$. 

%It is readily seen that the asymptotic procedure of removing virtual transitions, described in Subsection  %3, as well as the asymptotic procedure of  one-state reduction of phase space, described in Section 4, %yield, the new distributions of transition times,  $\tilde{F}_{0,  ij}(\cdot)$ and $_kF_{0,  ij}(\cdot)$, which, %also,  have not atoms at $0$. This property remains to hold for the  limiting distribution functions of %transition times $_{\bar{k}_n}\tilde{F}_{0,  ij}(\cdot)$ and $_{\bar{k}_n}F_{0,  ij}(\cdot)$ resulted by the %corresponding recurrent application of the above mentioned  procedures of asymptotic removing virtual %transitions and asymptotic  reduction of phase space.

%Due to this, the resulting limiting conditional distributions of hitting times $F_{0, \DD, k_m j}(\cdot)$ Ifor %the most absorbing state in domain $\overline{\DD}  also have no atoms at $0$. 

One possible  generalisation is connected with omitting condition ${\bf D}$ {\bf (b)}. The asymptotic recurrent algorithms presented in the paper  can be realised without changes. The only difference is that, in this case, one can not guarantee  that the corresponding limiting distributions $_{\bar{k}_n}\tilde{F}_{0,  ij}(\cdot)$ and $_{\bar{k}_n}F_{0,  ij}(\cdot)$ are not concentrated at zero, and, in sequel, that
limiting conditional distributions of hitting times $F_{0, \DD, ij}(\cdot)$ are not concentrated at zero.  In this case, one should compute the limiting Laplace transforms $\Psi_{0, \DD, ij}(\cdot)$  and the limits
of expectations $\hat{E}_{0, \DD, ij}$  using the same asymptotic recurrent formulas given in  Theorems 1 -- 12, and then  check that $\Psi_{0, \DD, ij}(s_0) < 1$, for some $s_0 > 0$ or, respectively, that 
$\hat{E}_{0, \DD, ij} > 0$.   

An important example is connected with the functional, which can be defined for any domain $\CC \subseteq \overline{\DD}$,  
\begin{equation}\label{timesd}
\tau_{\e, \CC, \DD} = \sum_{n = 1}^{\nu_{\e, \DD}} \kappa_{\e, n} {\rm I}(\eta_{\e, n -1} \in \CC). 
\end{equation}

The random variable $\tau_{\e, \CC, \DD}$ can be interpreted  the reward accumulated in states from domain $\CC$ at a trajectory of the Markov chain $\eta_{\e, n}$ up to the first hitting of domain $\DD$ by this Markov chain. 

In this case, one can consider the new Markov renewal process $(\eta_{\e,  \CC, n}$,  
$\kappa_{\e, \CC, n}$, where $\eta_{\e,  \CC, n} = \eta_{\e, n}$ and  $\kappa_{\e, \CC, n} = \kappa_{\e, n} {\rm I}(\eta_{\e, n -1} \in \CC)$. Obviously, the basic 
conditions  ${\bf B}$ -- ${\bf E}$ and ${\bf \tilde{C}}$ remain to hold for the new Markov renewal processes, wth the only change that the distribution functions $F_{\e, ij}(t) = {\rm I}(t \geq 0), t \geq 0$ and expectations $e_{\e, ij} = 0$, for $i \in \CC, j \in \XX, \e \in (0, 1]$. One can choose 
the initial normalisation functions $v_{\e, i} \equiv 1$, for $i \in \CC$ . In this case, the limiting distributions  for transition times $F_{0, ij}(t) = {\rm I}(t \geq 0), t \geq 0$ and expectations $e_{0, ij} = 0$, for $i \in \CC, j \in \XX$. 

\vspace{1mm}

{ \bf 10.3   Normalisation functions}.   It also should be mention that, in fact,   the normalisation functions $\check{v}_{\e, i} = \, _{\bar{k}_{\bar{m}}}\check{v}_{\e, i}, i \in \XX$ depend on the choice of  sequence  $\bar{k}_{\bar{m}} = \langle k_1, \ldots k_{\bar{m}} \rangle$ chosen according the corresponding algorithm of phase space  reduction.  However,  the property  of non-concentration at zero for the corresponding limiting distributions $G_{0 \DD, ij}(\cdot)$ implies that the above normalisation functions $_{\bar{k}'_{\bar{m}}}\check{v}_{\e, i}$ and, $_{\bar{k}''_{\bar{m}}}\check{v}_{\e, i}$ for any two admissible sequences $_{\bar{k}'_{\bar{m}}}$ and $_{\bar{k}''_{\bar{m}}}$ are asymptotically comparable (their quotient $_{\bar{k}'_{\bar{m}}}\check{v}_{\e, i}/\, _{\bar{k}''_{\bar{m}}}\check{v}_{\e, i}$ should converge to some constant $w_{0, i}[\bar{k}'_{\bar{m}}, \bar{k}''_{\bar{m}}] \in (0, \infty)$, as $\e \to 0$).  Condition ${\bf H}$ (or ${\bf H}$ and ${\bf \dot{H}}$) make it possible  recurrent computing of these constants.
The corresponding limiting distributions  for hitting times differ only by the scaling factor $w_{0, i}[\bar{k}'_{\bar{m}}, \bar{k}''_{\bar{m}}]$.

It also should be mentioned that the assumption that the initial normalisation functions, $1 \leq v_{\e, i} \to v_{0, i} \in [1, \infty]$ as $\e \to 0$, can be weaken and replaced by the assumption that 
$0 < v_{\e, i} \to v_{0, i} \in (0, \infty]$ as $\e \to 0$.   All algorithms described in Sections 1 -- 9 and results formulated in Theorems 1 -- 12 remain to be valid. \vspace{1mm}

{\bf 10.4 Limiting distributions for hitting times}.  The limiting conditional distributions 
$F_{0, \DD, i j}(\cdot)$, appearing in Theorem 7 and other weak convergence theorems for hitting times for perturbed semi-Markov processes, have semi-Markov phase-type, i.e., themselves are conditional distributions of hitting types for semi-Markov processes, possibly, with a smaller phase space than the phase space of the initial semi-Markov processes. 
 
In particular, it is possible that the limiting conditional distributions  for hitting times  belong to the well known class of phase-type distributions (conditional distributions of hitting times for continuous time Markov chains). We refer to woks [8 - 10], where a detailed description of phase-type distributions can be found. 
\vspace{1mm}

{\bf 10.5 Weak asymptotics for distributions and expectations  of  hitting times for an arbitrary domain 
$\DD$}.  It is readily seen that, in order to provide holding of conditions of Theorems 1 -- 12, for any domain $\DD \subseteq \XX$, one should, first require holding of the basic conditions ${\bf B}$ -- ${\bf E}$ and ${\bf \tilde{C}}$ in the following ``maximal''   form:

\vspace{2mm}

\noindent ${\bf B}_{\max}$:  {\bf (a)} $p_{\e, ij} > 0, \e \in (0, 1]$ or $p_{\e, ij} = 0, \e \in (0, 1]$, for every $i \in \DD, j \in \XX$, 
\makebox[11mm]{} {\bf (b)} there exists, for every $i, j  \in \XX$, a chain of states $i = j_0, j_1, \ldots, j_{n_{ij}- 1}$,  
\makebox[11mm]{}  $j_{n_{ij}} = j$ such that $\prod_{1 \leq l \leq n_{ij}} p_{1, j_{l-1} j_l} > 0$.
\vspace{2mm}

\noindent ${\bf C}_{\max}$: $p_{\e, ij} \to p_{0, ij}$ as $\e \to 0$, for $i, j \in \XX$. 
\vspace{2mm}

\noindent ${\bf D}_{\max}$:   {\bf (a)} $F_{\e, ij}(\cdot \, v_{\e,  i}) = \PP \{ \kappa_{\e, 1}/ v_{\e,  i}   \leq \cdot /  \eta_{\e, 0} = i, \eta_{\e, 1} = j  \} \Rightarrow F_{0, ij}(\cdot)$ as 
$\e \to 0$,  \makebox[11mm]{}   for $j \in \YY_{1, i}, \, i \in \XX$, {\bf (b)} $F_{0, ij}(\cdot), j \in \YY_{1, i}, i \in \XX$ are proper 
 distribution \\ \makebox[11mm]{}  functions  such that $F_{0, ij}(0) < 1, \, j \in \YY_{1, i}, i \in \XX$,   
{\bf (c)} $1 \leq v_{\e, i} \to v_{0, i} \in$ \\  \makebox[12mm]{}$[1, \infty]$ as $\e \to 0$, for $i \in \XX$. 
\vspace{2mm}

\noindent ${\bf E}_{\max}$:  $e_{\e, ij}/v_{\e,  i} \, = \int_0^{\infty} t F_{\e, ij}(dt)/ v_{\e,  i}   \to 
e_{0, ij} \, = \int_0^{\infty} t F_{0, ij}(dt)$ as $\e \to 0$,  for $j \in$ \\ \makebox[11mm]{}  $\YY_{1, i}, i \in \XX$. 
\vspace{2mm}

\noindent and \vspace{2mm}

\noindent ${\bf \tilde{C}}_{\max}$: $\tilde{p}_{\e, ij}  = {\rm I}(j \neq i)\frac{p_{\e, ij}}{1 - p_{\e, ii}}  \to \tilde{p}_{0, ij} \in [0, 1]$ as $\e \to 0$, for $i, j \in \XX$. 
\vspace{2mm}

It is useful to note that, under condition ${\bf B}_{\max}$  {\bf (a)},   condition ${\bf B}_{\max}$  {\bf (b)}  is equivalent to the assumption that the phase space $\XX$ is one class of communicative states for the Markov chain $\eta_{\e, n}$, for every $\e \in (0, 1]$. 

In this case, probabilities $p_{\e, ii} < 1, i \in \XX$,  for $\e \in [0, 1]$.

However, the phase space $\XX$ can possess an arbitrary communicative structure for the Markov chain $\eta_{0, n}$, i.e.  split in one or several closed classes of communicative states plus possibly a class of transient states.

Condition ${\bf H}$ should also be replaced  by its maximal variant: \vspace{2mm}

\noindent ${\bf H}_{\max}$: Functions $p_{\cdot, ij}, i, j \in \XX$ and  $v_{\cdot, i},  i \in \XX$ belong to some complete family \makebox[11mm]{}  of asymptotically comparable functions.
\vspace{2mm}

Condition ${\bf H}_{\max}$ is sufficient for holding of conditions ${\bf C}_{n}$ and ${\bf F}_{n}$, for any nonempty domain $\DD \subseteq \XX$ and $n \geq 0$. In sequel,  condition ${\bf H}_{\max}$ is sufficient for holding of  conditions $_{\bar{k}_{r}}{\bf \tilde{C}}$, $_{\bar{k}_{r}}{\bf \hat{C}}$, $_{\bar{k}_{r-1}}{\bf \tilde{F}}$, $r = 1, \ldots, \bar{m} -1$ and $_{\bar{k}_{\bar{m}}}{\bf \dot{F}}$ (where the  states $k_r \in \, _{\bar{k}_{r-1}} \overline{\DD}^*, r = 1, \ldots,  \bar{m} -1$  are chosen in such way that condition $_{\bar{k}_{r}}{\bf \hat{F}}$ holds, for  $r = 1, \ldots,  \bar{m} -1$),  
for any nonempty domain $\DD \subseteq \XX$. \vspace{1mm}

{ \bf 10.6  Vector hitting rewards}. In this case, one consider the Markov renewal process $(\eta_{\e, n},  \bar{\kappa}_{\e, n}) = 
(\eta_{\e, n},  (\kappa_{\e, 1, n}, \ldots, \kappa_{\e, L, n})) , n = 0, 1, \ldots$
with a phase space $\XX \times \RR^+_l$, where $\XX = \{1, \ldots, M \}$ is a finite set and $\RR^+_l = [0, \infty) \times \cdots \times  [0, \infty)$ is 
the $l$-product of the interval $[0, \infty)$.

The corresponding vector hitting reward functional is defined in the following way, for $\DD \subseteq \XX$, 
\begin{equation}\label{timesderd}
\bar{\tau}_{\e, \DD} = \sum_{n = 1}^{\nu_{\e, \DD}} \bar{\kappa}_{\e, n} = (\tau_{\e, 1, \DD}, \ldots,  \tau_{\e, l, \DD}),   \ \, {\rm where} \ \, 
\tau_{\e, l, \DD} =  \sum_{n = 1}^{\nu_{\e, \DD}} \kappa_{\e, l, n}.  
\end{equation}

In this case, one can reduce the problem to scalar case by using the method known as the Wold-Cram\'{e}r device.  Let us introduce the following scalar hitting reward, for $\bar{z} = (z_1, \ldots, z_l) \in  \RR^+_l$,
\begin{equation}\label{timesded}
\tau^{(\bar{z})}_{\e, \DD} = \sum_{n = 1}^{\nu_{\e, \DD}}  \kappa^{(\bar{z})}_{\e, n},  \ \, {\rm where} \ \,   \kappa^{(\bar{z})}_{\e, n} =  \sum_{r = 1}^l  z_r \kappa_{\e, r, n}, n \geq 1.
\end{equation}

Condition  ${\bf D}$ {\bf (a)} takes the form of weak convergence relations for $l$-dimen\-sional  distribution functions,   $F_{\e, ij}(\cdot \, v_{\e,  i}) = \PP \{ v_{\e,  i}^{-1}\bar{\kappa}_{\e, 1}   \leq \cdot /  \eta_{\e, 0} = i, \eta_{\e, 1} = j  \}$. It is equivalent to the assumption that the  weak convergence relations holds for the corresponding one-dimensional distribution functions,  $F^{(\bar{z})}_{\e, ij}(\cdot \, v_{\e,  i}) = \PP \{ v_{\e,  i}^{-1}\kappa^{(\bar{z})}_{\e, 1}  \leq \cdot /  \eta_{\e, 0} = i, \eta_{\e, 1} = j  \}$, for $ \bar{z} \in \RR^+_l$. 

%Analogously, condition  ${\bf E}$ {\bf (b)} takes the form of asymptotic relations for expectations   
%$v_{\e,  i}^{-1}e_{\e, r, ij} = \EE \{  v_{\e,  i}^{-1}\kappa_{\e, r, 1} /  \eta_{\e, 0} = i, \eta_{\e, 1} = j  \}$, for %$r = 1, \ldots, l$  that is equivalent to the corresponding asymptotic relations for expectations  $v_{\e,  i}%^{-1}e^{(\bar{z})}_{\e, ij} = \EE \{ v_{\e,  i}^{-1}\kappa^{(\bar{z})}_{\e, 1} /  \eta_{\e, 0} = i, \eta_{\e, 1} = j  %\}$, for $\bar{z} \in \RR_l$. 

The weak convergence relations for the corresponding distributions of  ``scalar''  hitting reward functionals, $G^{(\bar{z})}_{\e, \DD, ij} (\cdot \, \check{v}_{\e, i})$, for $\bar{z} \in \RR_l$ are equivalent to the  weak convergence relations for $l$-dimensional  distributions, $G_{\e, \DD, ij} (\cdot \, \check{v}_{\e, i}) = \PP_i \{ v_{\e, i}^{-1} \bar{\tau}_{\e, \DD}  \leq \cdot, \eta_{\e, \nu_{\e, \DD}} = j \}$. \vspace{1mm}

{ \bf 10.7  Real-valued hitting rewards}. In this case, one consider the Markov renewal process $(\eta_{\e, n},  \kappa_{\e, n}), n = 0, 1, \ldots$ with a phase space $\XX \times \RR$, where $\XX = \{1, \ldots, M \}$ is a finite set, and $ \RR = (-\infty, \infty)$ is a real line. 

The corresponding real-valued hitting reward functional can be defined in the standard way and represented as the difference of two non-negative hitting reward functionals, for $\DD \subseteq \XX$, 
\begin{equation}\label{timesduki}
\tau_{\e, \DD} = \sum_{n = 1}^{\nu_{\e, \DD}} \kappa_{\e, n} = \tau^+_{\e,\DD} - 
\tau^-_{\e, \DD},   
\end{equation}
where 
\begin{equation}\label{timesddew}
\tau^{\pm}_{\e, \DD} =  \sum_{n = 1}^{\nu_{\e, \DD}} \kappa^{\pm}_{\e, n} \  {\rm and} \   
\kappa^{\pm}_{\e, n} =  \pm \kappa_{\e, n} {\rm I}( \pm \kappa_{\e, n} \geq 0), \, n \geq 1.  
\end{equation}

One can consider the vector hitting reward functional $\bar{\tau}_{\e, \DD} = (\tau^{+}_{\e, \DD}, \tau^{-}_{\e, \DD})$ and to apply to this vector hitting reward the Wold-Cram\'{e}r device as it is described in Subsection 10.6,  for getting weak convergence relations of the type, $G_{\e, \DD, ij} (\cdot \, _{\bar{k}_{\bar{m}}}v_{\e, i}) = \PP_i \{ \, _{\bar{k}_{\bar{m}}}v_{\e, i}^{-1} \bar{\tau}_{\e, \DD}  \leq \cdot, 
\eta_{\e, \nu_{\e, \DD}} = j \}$.  Such relations in an obvious way imply the corresponding weak convergence relations real-valued hitting reward functionals $\tau_{\e, \DD} = \tau^+_{\e,\DD} 
- \tau^-_{\e, \DD}$. 

In order to escape some possible side effects, which can be caused by discontinuity of distribution functions of transition times at point $0$, one can also use a more general splitting procedure for hitting times, which is  based on random variables 
$\kappa^{\pm}_{\e, c,  n} =  \pm (\kappa_{\e, n}  - c) {\rm I}( \pm \kappa_{\e, n} \geq \pm c), n \geq 1$, for some $c \geq 0$. 

\vspace{1mm}

{ \bf 10.8  Examples}. Let us consider an  example, which let us illustrate results presented in the paper. 

In what follows, denote  $Exp_\lambda(t) = 1 - e^{-\lambda t}, t \geq 0$ the exponential distribution function with parameter $\lambda > 0$,  $Geo_{p}(\cdot)$ a geometrical distribution function with parameter $p \in (0, 1)$ (i.e. the distribution function of a random variable taking value $n$ with probability $p(1- p)^{n-1}$, for $n = 1, 2, \ldots$),  and 
$Con_a(t) = {\rm I}(t \geq a)$ the distribution function with unit jump at a point $a \geq 0$ (the distribution function of a random variable taking value $a$ with probability $1$).

In this example, the phase space $\XX = \{1, 2, 3 \}$, and  the matrix $\mathbf{P}_\e  = 
\| p_{\e, ij} \|$ of transition probabilities for the embedded Markov chain $\eta_{\e, n}$ has the following form, for $\e \in (0, 1]$,  
\begin{equation}\label{pitoma}
 \mathbf{P}_\e = \left\|
\begin{array}{ccc}
1 - \frac{1}{2} \e^{\alpha} - \frac{1}{2}\e^{\beta} & \frac{1}{2}\e^{\alpha}  & \frac{1}{2}\e^{\beta}  \vspace{1mm} \\
\frac{1}{2} \e & 1 - \e & \frac{1}{2} \e  \vspace{1mm} \\
0 & 0 &  1 
\end{array}
\right\|.
\end{equation}

Here,  parameters $\alpha, \beta \in [0, \infty)$.

Also,  we assume that distributions of transition times take the following forms,  for 
 $\e \in (0, 1]$, 
\begin{equation} \label{distra}
F_{\e, ij}(t) = \left\{
\begin{array}{ll}
Con_{\e^{-\gamma}}(t), t \geq 0 & \ \text{for} \ i = 1, \,  j = 1, 2, 3, \\
Con_1(t), t \geq 0 & \ \text{for} \ i = 2, 3, \,  j = 1, 2, 3, 
\end{array}
\right.
\end{equation}
where parameter $\gamma \in [0, \infty)$.

This means that the transition times from state 1 take value $ \e^{-\gamma}$, while transition times from states 2 and 3 take value 1.

Note that in the case $\gamma = 0$,  the semi-Markov process $\eta_{\e}(t)$ is a discrete time Markov chain embedded 
in continuous time.

Let us consider the case, where domain $\DD = \{ 3 \}$ and, thus, domain $\overline{\DD} = \{ 1, 2 \}$, i.e., it  is a two-states set.  

We shall describe weak convergence asymptotics for distributions  of hitting times 
$G_{\e, i, 3}(\cdot \check{v}_{\e, i}), i = 1, 2$ that is find proper normalisation 
functions $\check{v}_{\e, i}, i = 1, 2$ and  the corresponding limiting distributions 
$G_{0, i, 3}(\cdot), i = 1, 2$, which are not concentrated at $0$, and, also, describe 
the corresponding asymptotics for the expectations of hitting times.

In this case, conditions ${\bf A}$ and ${\bf B}$ obviously hold.

Moreover, in this case condition ${\bf \hat{B}}$ also holds.

Since, $\DD$ is a one-state set, the indicator ${\rm I}(\eta_\e(\tau_{\e, \DD} = 3) = 1$, and, thus, for $i = 1, 2$, 
\begin{equation}\label{vutyr}
G_{\e, \DD, i3}(t) = \PP \{\tau_{\e, \DD} \leq t \}, \ t \geq 0. 
\end{equation}

Condition ${\bf C}$ also holds. The corresponding limiting matrix of transition probabilities for the embedded Markov chain $\eta_{0, n}$ takes different forms, for cases: {\bf (1)} $\alpha = \beta = 0$, {\bf (2)} $\alpha > 0, \beta = 0$, {\bf (3)} $\alpha = 0, \beta > 0$, 
{\bf (4)} $\alpha = \beta > 0$, {\bf (5)} $\alpha > \beta > 0$, and {\bf (6)}  $\beta > \alpha > 0$. We denote matrix $\mathbf{P}_0$ for these three cases, respectively,  as $\mathbf{P}_0^{(h)}$, for $h = 1, \ldots, 6$.   Obviously, 
\begin{equation}\label{pitomak}
 \mathbf{P}^{(1)}_0 = \left\|
\begin{array}{ccc}
0 & \frac{1}{2} & \frac{1}{2}    \vspace{1mm} \\
0  & 1 &  0  \vspace{1mm} \\
0 & 0 &  1 
\end{array}
\right\|, \, \mathbf{P}^{(2)}_0 = \left\|
\begin{array}{ccc}
\frac{1}{2}   & 0 & \frac{1}{2}   \vspace{1mm} \\
0  & 1 &  0  \vspace{1mm} \\
0 & 0 &  1 
\end{array}
\right\|, \, \mathbf{P}^{(3)}_0 = \left\|
\begin{array}{ccc}
\frac{1}{2}   & \frac{1}{2}  & 0   \vspace{1mm} \\
0  & 1 &  0  \vspace{1mm} \\
0 & 0 &  1 
\end{array}
\right\|,
\end{equation}
and
\begin{equation}\label{pitomakio}
 \mathbf{P}^{(h)}_0  = \left\|
\begin{array}{ccc}
1 & 0 & 0    \vspace{1mm} \\
0  & 1 &  0  \vspace{1mm} \\
0 & 0 &  1 
\end{array}
\right\|, \ h = 4, 5, 6.
\end{equation} 

In this model,  state $3$ is an absorbing state, state $2$ is asymptotically absorbing 
state, while state $1$ is asymptotically non-absorbing state, for cases $h = 1, 2, 3$ or  
asymptotically absorbing state, for cases $h = 4, 5, 6$.

In this case, it is natural to choose the initial normalisation functions,
\begin{equation}\label{nouty}
v_{\e, 1} = \e^{-\gamma}, \ v_{\e, 2} = 1, \ \e \in (0, 1].
\end{equation}

Note, first of all that, in this case, the condition of asymptotic comparability ${\bf H}$ holds for transition probabilities $p_{\e, ij}, i \in \overline{\DD},  j \in \XX$ and the normalisation functions $v_{\e, i}, i \in \overline{\DD}$, which belong to the complete family of asymptotically comparable functions ${\cal H}_1$ defined by 
relation (\ref{byrce}).

Conditions ${\bf D}$ obviously holds, with the corresponding limiting distribution functions $F_{0, ij}(t)$ and  the corresponding Laplace transforms $\phi_{0, ij}(s)$, given by the following relations, for  $i = 1, 2, j = 1, 2, 3$,
\begin{equation}\label{ferta}
F_{0, ij}(t) = Con_1(t), \ t \geq 0,
\end{equation}
and
\begin{equation}\label{fertava}
\phi_{0, ij}(s) = e^{-s}, \ s \geq 0. 
\end{equation} 

Conditions ${\bf E}$ also holds.  The corresponding limiting expectations,
take the following forms, for $ i = 1, 2, j = 1, 2, 3$.
\begin{equation}\label{fertabva}
e_{0, ij} = 1. 
\end{equation} 

Since domain $\overline{\DD}$ is a two-states set, the only steps ${\bf 0}$ and ${\bf 1}$ of the asymptotic recurrent algorithm described in Sections 3 -- 9 should be realised.

At the step ${\bf 0}$, the procedure of asymptotic removing virtual transitions described in Section 3 should be applied to the semi-Markov processes $\eta_\e(t)$. 
 
The Markov chain $\tilde{\eta}_{\e, n}$ has the transition probabilities  $\tilde{p}_{\e, ij}  = {\rm I}(j \neq i)\frac{p_{\e, ij}}{1 - p_{\e, ii}}$, $i = 1, 2, j = 1, 2, 3$, which given by  relation (\ref{gopet}),  and the matrix of its transition probabilities takes the 
following form, for $\e \in (0, 1]$,
\begin{equation}\label{pitomake}
 \tilde{\mathbf{P}}_\e = \left\|
\begin{array}{ccc}
0  & \frac{\e^{\alpha}}{\e^{\alpha} + e^{\beta}}   & \frac{\e^{\beta}}{\e^{\alpha} + e^{\beta}}    
\vspace{1mm} \\
\frac{1}{2}  & 0 & \frac{1}{2}   \vspace{1mm} \\
0 & 0 &  1 
\end{array}
\right\|.
\end{equation}

According Lemma 3, conditions ${\bf \tilde{A}}$ and ${\bf \tilde{B}}$ hold.
 
Condition  ${\bf \tilde{C}}$  also holds. The corresponding limiting matrix of 
transition probabilities $\tilde{\mathbf{P}}_0$ for the Markov chain $\tilde{\eta}_{0, n}$ takes different forms for cases {\bf (1)} -- {\bf (6)}  listed above. 
We denote matrix $\tilde{\mathbf{P}}_0$ for these  cases, respectively,  as $\tilde{\mathbf{P}}_0^{(h)}$, for $h = 1, \ldots, 6$.   Obviously, 
\begin{equation}\label{pitomaker}
 \tilde{\mathbf{P}}^{(1)}_0 = \left\|
\begin{array}{ccc}
0  & \frac{1}{2}    & \frac{1}{2}    \vspace{1mm} \\
\frac{1}{2}  & 0 & \frac{1}{2}   \vspace{1mm} \\
0 & 0 &  1 
\end{array}
\right\|, \ \tilde{\mathbf{P}}^{(2)}_0 = \left\|
\begin{array}{ccc}
0  & 0   & 1     \vspace{1mm} \\
\frac{1}{2}  & 0 & \frac{1}{2}   \vspace{1mm} \\
0 & 0 &  1 
\end{array}
\right\|, \ \tilde{\mathbf{P}}^{(3)}_0 = \left\|
\begin{array}{ccc}
0  & 1   & 0   \vspace{1mm} \\
\frac{1}{2}  & 0 & \frac{1}{2}   \vspace{1mm} \\
0 & 0 &  1 
\end{array}
\right\|, 
\end{equation}
and 
\begin{equation}\label{pitomakerop}
 \tilde{\mathbf{P}}^{(4)}_0 = \left\|
\begin{array}{ccc}
0  & \frac{1}{2}    & \frac{1}{2}    \vspace{1mm} \\
\frac{1}{2}  & 0 & \frac{1}{2}   \vspace{1mm} \\
0 & 0 &  1 
\end{array}
\right\|, \ \tilde{\mathbf{P}}^{(5)}_0 = \left\|
\begin{array}{ccc}
0  & 0   & 1     \vspace{1mm} \\
\frac{1}{2}  & 0 & \frac{1}{2}   \vspace{1mm} \\
0 & 0 &  1 
\end{array}
\right\|, \ \tilde{\mathbf{P}}^{(6)}_0 = \left\|
\begin{array}{ccc}
0  & 1   & 0   \vspace{1mm} \\
\frac{1}{2}  & 0 & \frac{1}{2}   \vspace{1mm} \\
0 & 0 &  1 
\end{array}
\right\|.
\end{equation}

In this case, the normalisation functions, $\tilde{v}_{\e, i} = (1 - p_{\e, ii})^{-1} v_{\e, i}, i = 1, 2$,  take, according relations (\ref{compreg}) the following forms, 
\begin{equation}\label{nortena}
\tilde{v}_{\e, 1} =  2\e^{-\gamma}(\e^{\alpha} + \e^{\beta})^{-1}, \  
\tilde{v}_{\e, 2} = \e^{-1},  \  \e \in (0, 1]. 
\end{equation}

According Lemma 3, conditions ${\bf \tilde{D}}$ holds.

Let us define parameter,
\begin{equation}\label{nortenas}
\delta = \min(\alpha, \beta) + \gamma.
\end{equation}

If $\delta = 0$, then: $\tilde{v}_{\e, 1} \equiv 1$, in the case {\bf (1)}; $\tilde{v}_{\e, 1} \to \frac{1}{2}$ as $\e \to 0$, in the cases {\bf (2)} and {\bf (3)}. 

If $\delta > 0$ then $\tilde{v}_{\e, 1} \to \infty$ as $\e \to 0$. 

Also, $\tilde{v}_{\e, 2} \to \infty$ as $\e \to 0$. 

Probability $p_{0, 11}  < 1$ in the cases {\bf (1)} -- {\bf (3)}, i.e., if $\alpha \wedge \beta = 0$, while $p_{0, 11}  = 1$ in the cases {\bf (4)} -- {\bf (6)}, i.e., if $\alpha \wedge \beta >  0$. 

The Laplace transform $\tilde{\phi}_{0, 1j}(\cdot), j = 1, 2, 3$ for corresponding limiting distribution function $\tilde{F}_{0, 1j}(\cdot)$ appearing in condition ${\bf \tilde{D}}$  should be computed with the use of  relation (\ref{trawetba}) and (\ref{tatop}), if  $p_{0, 11}  < 1$, or (\ref{trawetbama}) and (\ref{tatop}), if $p_{0, 11}  =1$.

Thus, relations (\ref{pitomak})  and (\ref{fertava}) imply that, for $s \geq 0$ and  $j = 1, 2, 3$,
\begin{equation*}
\tilde{\phi}_{0,  1j}(s)   =  \frac{\phi_{0, 1j}((1 - p_{0, 11})s)(1 - p_{0, 11})}{1 - p_{0, 11} \phi_{0, 1j}((1 - p_{0, 11})s ))} \makebox[57mm]{}
\end{equation*}
\begin{equation}\label{traetba}
\makebox[12mm]{} = \left\{
\begin{array}{ll}
e^{-s} & \ \text{if} \ \alpha = \beta = 0, \vspace{2mm} \\
\frac{1}{2}e^{- \frac{s}{2}} / (1 - \frac{1}{2}e^{- \frac{s}{2}})  &  \ \text{if} \ \alpha > 0, \beta = 0 \ \text{or} \ \alpha = 0, \beta > 0. 
\end{array}
\right.
\end{equation}

Note that the  distribution functions corresponding to the Laplace transforms appearing in relation (\ref{traetba}) are, respectively, 
$Con_1(t)$, for the case $\alpha = \beta = 0$,  and $Geo_{\frac{1}{2}}(2t)$, for the cases $\alpha > 0, \beta = 0$ and $\alpha = 0, \beta > 0$. 

In the case,  where $\alpha \wedge \beta > 0$, relations (\ref{pitomak})  and (\ref{fertabva}) imply that, for  $s \geq 0$ and $j = 1, 2, 3$,
\begin{equation}\label{traetbak}
\tilde{\phi}_{0,  1j}(s)   = (1 + e_{0, 11}s)^{-1} = (1 + s)^{-1}.
\end{equation}

Probability $p_{0, 22}  = 1$. 

Thus, the Laplace transform $\tilde{\phi}_{0, 2j}(s), s \geq 0$ for corresponding limiting distribution function $\tilde{F}_{0, 2j}(\cdot), j = 1, 2, 3$ appearing in condition ${\bf \tilde{D}}$  should be computed with the use of  relation  (\ref{trawetbama}) and (\ref{tatop}).

In the case, relations (\ref{pitomak})  and (\ref{fertabva}) imply that, for $s \geq 0$ and $j = 1, 2, 3$,
\begin{equation}\label{traetbakn}
\tilde{\phi}_{0,  2j}(s)   = (1 + e_{0, 22}s)^{-1} = (1 + s)^{-1}.
\end{equation}

Note that the distribution function corresponding the Laplace transform appearing in relations (\ref{traetbak}) and (\ref{traetbakn}) 
is $Exp_1(\cdot)$. 

According Lemma 3, conditions  ${\bf \tilde{E}}$ also holds.

The corresponding limiting expectations take the following forms, for $ i = 1, 2, j = 1, 2, 3$.
\begin{equation}\label{fertabvan}
\tilde{e}_{0, ij} = 1. 
\end{equation} 

Note also that the distribution function $\tilde{F}_{\e, 3j}(t) = Con_1(t), t \geq 0$, for $j  = 1, 2, 3$ and the corresponding expectation $\tilde{e}_{\e, 3j} = 1$, for $j  = 1, 2, 3$ and $\e \in [0, 1]$. 

At the  sub-step ${\bf 1.1}$,  the procedure of asymptotic one-state reduction of phase space, described in Section 4 should be applied to the semi-Markov processes 
$\tilde{\eta}_\e(t)$.

Relations (\ref{nortena}) and  (\ref{nortenas}) imply that, 
\begin{equation}\label{simaba}
\tilde{v}_{\e, 1} = 2\e^{-\gamma}(\e^{\alpha} + \e^{\beta})^{-1}
 \sim \frac{2}{1 + {\rm I}(\alpha =  \beta)} 
\e^{- \delta} \ {\rm as} \ \e \to 0, 
\end{equation}

Three cases should be considered, where {\bf (i)} $\delta \in [0, 1)$, {\bf (ii)} $\delta = 1$, and {\bf (iii)} $\delta \in (1, \infty)$. Obviously, 
 \begin{equation}\label{simany}
\tilde{w}_{\e, 12} = \frac{\tilde{v}_{\e, 1}}{\tilde{v}_{\e, 2}} \to \tilde{w}_{0, 12} \
\ {\rm as} \ \e \to 0,
 \end{equation}
 where 
 \begin{equation}\label{vopire}
 \tilde{w}_{0, 12}  = \left\{
 \begin{array}{ll}
0 & \ \text{if} \ \delta \in [0, 1), \\
\frac{2}{1 + {\rm I}(\alpha =  \beta)}  & \ \text{if} \ \delta = 1, \\
 \infty & \ \text{if} \ \delta \in (1, \infty). 
 \end{array}
 \right. 
 \end{equation}
 
Let us consider the case {\bf (i)}, where $\delta \in [0, 1)$.

In this case, conditions ${\bf \tilde{F}}$ and $_1{\bf \hat{F}}$ holds, i.e., state $1$ is asymptotically less absorbing than state $2$. 

That is why,  state $1$ should be chosen for exclusion it from the phase space $\XX$, according the procedure described in Section 4. In this case, the reduced phase space is $_1\XX = \{ 2, 3 \}$ and domain $_1\overline{\DD} = \{ 2 \}$ is a one-state set.

The Markov chain $_1\eta_{\e, n}$ has the transition probabilities  
$_1p_{\e, ij}  =  \tilde{p}_{\e, ij} + \tilde{p}_{\e, i1} \tilde{p}_{\e, 1j}$, $i, j = 2, 3$, which given by  relation (\ref{gopetk}),  and the matrix of its transition probabilities takes the following form, for $\e \in (0, 1]$,
\begin{equation}\label{pitomakbat}
_1\mathbf{P}_\e =  \left\|
\begin{array}{cc}
\frac{\e^{\alpha}}{2(\e^{\alpha} + \e^{\beta})}  & \frac{\e^{\alpha} + 2\e^{\beta}}{2(\e^{\alpha} + \e^{\beta})}      \vspace{1mm} \\
0 & 1  
\end{array}
\right\|.
\end{equation}

According Lemma 13,  conditions $_1{\bf A}$ and $_1{\bf B}$  hold.
 
Condition  $_1{\bf C}$  also holds. The corresponding limiting matrix of 
transition probabilities $_1\mathbf{P}_0$ for the Markov chain $_1\eta_{0, n}$ takes three different forms  in the cases: {\bf (1)}, {\bf (4)}, i.e., if $\alpha = \beta$; {\bf (2)}, 
{\bf (5)}, i.e., if   $\alpha > \beta$; and {\bf (3)}, {\bf (6)}, i.e., if  $\beta > \alpha$.  We denote matrix $_1\mathbf{P}_0$ for these  cases, respectively,  as 
$_1\mathbf{P}^{(h)},  h = 1, \ldots, 6$. Obviously, 
\begin{equation*}
 _1\mathbf{P}^{(1)}_0 = \left\|
\begin{array}{cc}
\frac{1}{4}  & \frac{3}{4}      \vspace{1mm} \\
0  & 1    
\end{array}
\right\|, \ _1\mathbf{P}^{(2)}_0 = \left\|
\begin{array}{cc}
 0   & 1     \vspace{1mm} \\
0  & 1 
\end{array}
\right\|, \ _1\mathbf{P}^{(3)}_0 = \left\|
\begin{array}{cc}
 \frac{1}{2}    &  \frac{1}{2}    \vspace{1mm} \\
 0 &  1 
\end{array}
\right\|.
\end{equation*}
\begin{equation}\label{pitomakerno}
\makebox[8mm]{} _1\mathbf{P}^{(4)}_0 = \left\|
\begin{array}{cc}
\frac{1}{4}  & \frac{3}{4}      \vspace{1mm} \\
0  & 1    
\end{array}
\right\|, \ _1\mathbf{P}^{(5)}_0 = \left\|
\begin{array}{cc}
 0   & 1     \vspace{1mm} \\
0  & 1 
\end{array}
\right\|, \ _1\mathbf{P}^{(6)}_0 = \left\|
\begin{array}{cc}
 \frac{1}{2}    &  \frac{1}{2}    \vspace{1mm} \\
 0 &  1 
\end{array}
\right\|.
\end{equation} 

In this case, the normalisation function, $_1v_{\e, 2}$ takes, according relation (\ref{normali}),  the following form, 
\begin{equation}\label{norte}
_1v_{\e, 2} = \tilde{v}_{\e, 2} = \e^{-1},  \ \e \in (0, 1]. 
\end{equation}

By Lemma 9, condition $_1{\bf \tilde{C}}$ holds and, for $j = 2, 3$,
\begin{equation}\label{honpl}
_1\tilde{p}_{0, 2 j}  = {\rm I}(j \neq 2)\frac{_1p_{0, 2 j}}{1 -  \, _1p_{0, 2 2}} = {\rm I}(j \neq 2). 
\end{equation}

Also, by Lemma 10, condition $_1{\bf \hat{C}}$  holds, and, 
\begin{equation}\label{nesaw}
_1\hat{q}_{0}[2j] =  {\rm I}(j \neq 2) \frac{\tilde{p}_{0, 2 j}}{_1p_{0, 2 j}}  = \left\{
\begin{array}{ll}
0 & \text{if} \ j = 2, \vspace{2mm} \\
 \frac{2}{3} &  \text{if} \ j = 3, \, \alpha = \beta, \vspace{2mm} \\
\frac{1}{2}  &  \text{if} \ j = 3, \, \alpha > \beta, \vspace{2mm} \\
 1 &  \text{if} \ j = 3, \, \beta > \alpha.  
\end{array}
\right. 
\end{equation}

According Lemma 13,  conditions $_1{\bf D}$ and $_1{\bf E}$ hold.

The Laplace transform $_1\phi_{0, 2j}(s), s \geq 0$ for corresponding limiting distribution function $_1F_{0, 2j}(\cdot)$ appearing in condition $_1{\bf D}$  should be computed with the use of  relations (\ref{trawetkasopa}) and (\ref{trawetkasf}), for $j = 2, 3$. 

Using relations (\ref{trawetkasopa}), (\ref{trawetkasf}), (\ref{traetbakn}), (\ref{simany}), (\ref{vopire}), and  (\ref{nesaw}) we get, taking into account that $\tilde{w}_{0,  12} = 0$, that, for $s \geq 0$ and $j = 3$,
\begin{align}\label{nopur}
_1\phi_{0, 23}(s) & =  \tilde{\phi}_{0,  23}(s) \, _1q_{0, 23} +  
\tilde{\phi}_{0, 21}(s) \tilde{\phi}_{0, 13}(\tilde{w}_{0,  12} s) (1 - \, _1q_{0, 23}) 
\vspace{2mm} \nonumber \\
& = (1 + s)^{-1} \, _1q_{0, 23} + (1 + s)^{-1} (1 - \, _1q_{0, 23})  = (1 + s)^{-1}, 
\end{align}
and, for $s \geq 0$ and $j = 2$,
\begin{equation}\label{nopura}
_1\phi_{0, 22}(s)  =  \tilde{\phi}_{0,  21}(s) \tilde{\phi}_{0,  12}(\tilde{w}_{0,  12} s) = 
(1 + s)^{-1}. 
\end{equation}

Thus, the distribution functions $_1F_{0, 2j}(\cdot) = Exp_1(\cdot)$, for $j = 2, 3$.

Also, the limiting expectation $_1e_{0, 2j}$ appearing in condition $_1{\bf E}$ is the first moment of the distribution function 
$_1F_{0, 2j}(\cdot)$ and, thus, for $j = 2, 3$,
\begin{equation}\label{noplur}
_1e_{0, 2j} = 1.
\end{equation}

Note also that the distribution function $_1F_{\e, 3j}(t) \equiv Con_1(t \geq 1), t \geq 0$, for $j  = 2, 3$ and the corresponding expectation $_1e_{\e, 3j} \equiv 1$, for $j  = 2, 3$ and $\e \in [0, 1]$. 

At step ${\bf 1.2}$, the procedure of asymptotic removing virtual transitions described in Section 5 should be applied to the semi-Markov processes $_1\eta_\e(t)$. 

The Markov chain $_1\tilde{\eta}_{\e, n}$ has the transition probabilities  
$_1p_{\e, ij}, i, j = 2, 3$, which given by  relation (\ref{gopetasd}),  and the matrix of its transition probabilities takes the following form, for $\e \in (0, 1]$,
\begin{equation}\label{pitomakbatni}
_1\tilde{\mathbf{P}}_\e =  \left\|
\begin{array}{cc}
0  & 1      \vspace{1mm} \\
0 & 1  
\end{array}
\right\|.
\end{equation}

In this case, the normalisation function, $_1\tilde{v}_{\e, 2}$ takes, according 
relation (\ref{compregase}), (\ref{nortena}), (\ref{pitomakbat}), and  (\ref{pitomakerno}) the following form, 
\begin{align}\label{nortever}
_1\tilde{v}_{\e, 2} & = (1 - \, _1p_{\e, 22})^{-1}(1 - p_{\e, 22})^{-1} v_{\e, 2} 
\vspace{2mm} \nonumber \\
& = (1 - \frac{\e^{\alpha}}{2(\e^{\alpha} +  e^{\beta})})^{-1} \e^{-1} = 
\frac{2(\e^\alpha + \e^{\beta})}{\e(\e^{\alpha} + 2 \e^{\beta})} 
\vspace{2mm} \nonumber \\ 
& \sim \, _1p_{0, 23}^{-1} \, \e^{-1} \ {\rm as} \ \e \to 0,
\end{align}
where the coefficient $_1p_{0, 23}^{-1}$ takes, according relation (\ref{pitomakerno}),   the following value, 
\begin{equation}\label{coef}
_1p_{0, 23}^{-1} = \frac{4}{3} \cdot {\rm I}(\alpha = \beta) + 1 \cdot {\rm I}(\alpha > \beta) + 
2 \cdot {\rm I}(\beta > \alpha). 
\end{equation} 

According Lemma 17,  conditions $_1{\bf \tilde{A}}$ and $_1{\bf \tilde{B}}$  hold. 

As it was mentioned above condition $_1{\bf \tilde{C}}$ (which plays the role of condition
${\bf C}$ for the semi-Markov processes $_1\tilde{\eta}_\e(t)$) holds. Also, as was mentioned 
in Subsection 5.6, condition $_1{\bf \tilde{C}}'$ (which is equivalent to condition
$_1{\bf \tilde{C}}$ and plays the role of condition ${\bf \tilde{C}}$ for the semi-Markov processes $_1\tilde{\eta}_\e(t)$) holds.

According Lemma 17,  conditions $_1{\bf \tilde{D}}$ and $_1{\bf \tilde{E}}$ hold.

The Laplace transform $_1\tilde{\phi}_{0, 23}(s), s \geq 0$ for corresponding limiting distribution function $_1\tilde{F}_{0, 23}(\cdot)$ appearing in condition $_1{\bf \tilde{D}}$  should be computed with the use of  relation (\ref{trawetbanop}), if probability $_1p_{0, 22} < 1$, or relation
(\ref{trawetbamnop}), if probability $_1p_{0, 22} = 1$.  Relation (\ref{pitomakerno}) shows that 
the first variant takes place. Using relations  (\ref{trawetbamnop}), (\ref{pitomakerno}), (\ref{nopur}), and (\ref{nopura}), we get, for $s \geq 0$,
\begin{align}\label{tbanop}
_1\tilde{\phi}_{0,  23}(s) & =  \frac{_1\phi_{0, 23}((1 - \, _1p_{0, 22})s)(1 - 
\, _1p_{0, 22})}{1 - \, _1p_{0, 22} \, _1\phi_{0, 23}((1 - \, _1p_{0, 22})s ))} 
\vspace{2mm} \nonumber \\
& =  \frac{1 - \, _1p_{0, 22}}{1 + (1 - \, _1p_{0, 22})s} \, \big(1 - 
\frac{_1p_{0, 22}}{1 + (1 - \, _1p_{0, 22})s} \big)^{-1} 
= \frac{1}{1 + s}.
\end{align}

Also, according relation (\ref{pitomakerno}) and (\ref{noplur}),
\begin{equation}\label{expofan}
_1\tilde{e}_{0,   23}  =   (1 -  \, _1p_{0, 22}) \, _1e_{0, 23} + \, _1p_{0, 22}  \, _1e_{0, 22} =  1.
\end{equation}

According, Theorem 2, $_1\tilde{v}_{\e, 2}$, given by relation (\ref{norte}) should be taken as the  normalisation function and the limiting Laplace transform $\Psi_{0, \DD, 23}(s) = \, _1\tilde{\phi}_{0,  23}(s)$, $s \geq 0$. Thus, Theorem 2, yields, 
according relation (\ref{tbanop}), the following asymptotic relation,
\begin{equation}\label{hopceb}
G_{\e, \DD, 23}(\cdot \, _1\tilde{v}_{\e, 2}) \Rightarrow G_{0, \DD, 23}(\cdot) =  Exp_1(\cdot) \ {\rm as} \ \e \to 0.
\end{equation}

Now the backward recurrent algorithm described in Theorems 2 and 3, which are given in Section 5,
should be applied.

According relations (\ref{pitomaker}) and (\ref{pitomakerop}), probability $\tilde{p}_{0, 12}$ is takes the following value, 
\begin{equation}\label{gerd}
\tilde{p}_{0, 12} =  \frac{1}{2} \cdot {\rm I}(\alpha = \beta) + 0 \cdot {\rm I}(\alpha > \beta) + 
1 \cdot {\rm I}(\beta > \alpha). 
\end{equation} 

In the cases $\alpha = \beta$ and  $\beta > \alpha$, probability $\tilde{p}_{0, 12}$ 
takes positive values. 

In this case, according Theorem 3, $_1\tilde{v}_{\e, 2}$ should be taken as the  normalisation function. Also, 
relations (\ref{idersfat}), (\ref{simany}), (\ref{vopire}), and (\ref{gerd}) yield, in this case, that the limiting Laplace transform $\Psi_{0, \DD, 13}(s)$ takes the following form, for $s \geq 0$,
\begin{align*}
\Psi_{0, \DD, 13}(s) & = \tilde{\phi}_{0, 13}((1 - \, _1p_{0, 22})  \tilde{w}_{0, 12}s))\tilde{p}_{0, 13} 
\vspace{2mm} \nonumber \\
& \quad + \Psi_{0, \DD, 23}(s)  \tilde{\phi}_{0, 12}((1 - \, _1p_{0, 22})  \tilde{w}_{0, 12}s)  \tilde{p}_{0,  12} 
\makebox[30mm]{} 
\end{align*} 
\begin{equation}\label{idersfatolta}
 =  \tilde{p}_{0, 13}  + \frac{1}{1 + s} \tilde{p}_{0,  12} = \left\{
\begin{array}{ll}
\frac{1}{2} +  \frac{1}{2}\frac{1}{1 + s}   & \text{if} \  \alpha = \beta, \vspace{2mm} \\
\frac{1}{1 + s}  & \text{if} \  \beta > \alpha. \makebox[8mm]{}
\end{array}
\right.
\end{equation}

Thus, Theorem 2, yields, 
according relation (\ref{idersfatolta}), the following asymptotic relation,
\begin{equation}\label{hopcena}
G_{\e, \DD, 13}(\cdot \, _1\tilde{v}_{\e, 2}) \Rightarrow G_{0, \DD, 13}(\cdot)  \ {\rm as} \ \e \to 0,
\end{equation}
where
\begin{equation}\label{hopcenaba}
G_{0, \DD, 13}(t)  = \left\{
\begin{array}{ll}
\frac{1}{2} Con_0(t) +  \frac{1}{2}Exp_1(t)   & \text{if} \  \alpha = \beta, 
\vspace{2mm} \\
Exp_1(t) & \text{if} \ \beta > \alpha.
\end{array}
\right.
\end{equation}

In the case $\alpha > \beta$, probability $\tilde{p}_{0, 12}  = 0$. 

According to  Theorem 3, 
function $\tilde{v}_{\e, 1}$, given by relation (\ref{simaba}),  should be taken as the  normalisation function. Also, relations (\ref{idersfatas}) yields, in this case, that the limiting Laplace transform $\Psi_{0, \DD, 13}(s)$ takes the following form, for $s \geq 0$,
\begin{equation*}
\Psi_{0, \DD, 13}(s) = \tilde{\phi}_{0, 13}(s)\tilde{p}_{0, 13}  \makebox[61mm]{} 
\end{equation*}
\begin{equation}\label{idelta}
  = \left\{
\begin{array}{ll}
\frac{1}{2}e^{- \frac{s}{2}} / (1 - \frac{1}{2}e^{- \frac{s}{2}})   & \text{if} \  \alpha > \beta = 0, \vspace{2mm} \\
\frac{1}{1 + s}  & \text{if} \  \alpha > \beta > 0. 
\end{array}
\right.
\end{equation}

Thus, Theorem 3, yields, 
according relation (\ref{idelta}), the following asymptotic relation,
\begin{equation}\label{hopcenanu}
G_{\e, \DD, 13}(\cdot \, \tilde{v}_{\e, 1}) \Rightarrow G_{0, \DD, 13}(\cdot)  \ {\rm as} \ \e \to 0,
\end{equation}
where
\begin{equation}\label{hopcenabam}
G_{0, \DD, 13}(\cdot)  = \left\{
\begin{array}{ll}
Geo_{\frac{1}{2}}(2t)   & \text{if} \  \alpha > \beta = 0, 
\vspace{2mm} \\
Exp_1(t) & \text{if} \ \alpha > \beta > 0.
\end{array}
\right.
\end{equation}

Remind that we consider the case {\bf (i)}, where $\delta  = \min(\alpha, \beta) + 
\gamma \in [0, 1)$. In this case,  the normalisation function $_1\tilde{v}_{\e, 2}$, used in the weak convergence relations 
(\ref{hopceb}),  (\ref{hopcena}), and the normalisation function $\tilde{v}_{\e, 1}$, used in the weak convergence relation (\ref{hopcenanu}),  are connected by relation,
 \begin{equation}\label{osmal}
 \tilde{v}_{\e, 1} = o(_1\tilde{v}_{\e, 2}) \ {\rm as}  \ \e \to 0. 
\end{equation}

Finally, according to Theorem 9,   the following asymptotic relation takes place for expectations of hitting times,
\begin{equation}\label{hopcebnip}
E_{\e, \DD, 23} / \,  _1\tilde{v}_{\e, 2} \to  \bar{E}_{0, \DD, 23} = E_{0, \DD, 23} = 1 \ {\rm as} \ \e \to 0.
\end{equation}

The situation is more complex in the case, where the initial state is $1$.

In this case, the normalisation function proposed in relation (\ref{vopuer}) takes the following form,
\begin{equation}\label{hurewq}
\bar{v}_{\e, 1} =  \tilde{v}_{\e, 1}  + \, _1\tilde{v}_{\e, 2} \, \tilde{p}_{\e, 12} = 
\frac{2}{\e^{\gamma}(\e^{\alpha} + \e^{\beta})} + 
\frac{2(\e^\alpha + \e^{\beta})}{\e(\e^{\alpha} + 2 \e^{\beta})}  \frac{\e^{\alpha}}{\e^{\alpha} + e^{\beta}}   .   
\end{equation}

The asymptotic comparability condition represented by relation (\ref{vopuerta}) holds and this relation takes the following form,
\begin{equation}\label{asymbolt}
\frac{_1\tilde{v}_{\e, 2} \, \tilde{p}_{\e, 12}}{\bar{v}_{\e, 1}}  = 
\big(1 + \frac{\e(\e^{\alpha} + 2 \e^{\beta})}{\e^{\alpha + \gamma}(\e^\alpha + \e^\beta)} \big)^{-1} 
 \to \, \bar{u}_{\alpha, \beta, \gamma} \ {\rm as} \ \e \to 0, 
\end{equation}
where
\begin{equation}\label{asymbult}
\bar{u}_{\alpha, \beta, \gamma}  = \left\{
\begin{array}{ll}
1 & \text{if} \ \alpha + \gamma < 1,  \vspace{2mm}  \\
\frac{1}{3}{\rm I}(\alpha > \beta) + \frac{2}{5}{\rm I}(\alpha = \beta) +  \frac{1}{2}{\rm I}(\alpha < \beta) & \text{if} \ \alpha +  \gamma = 1,  \vspace{2mm}  \\
0 & \text{if} \ \alpha + \gamma > 1.
\end{array}
\right.
\end{equation}

Remind that we consider the case, where parameter $\delta = \min(\alpha, \beta) + \gamma \in [0, 1)$. 

If $\alpha \leq \beta$, then  $\delta = \alpha + \gamma < 1$, and, thus, $\bar{u}_{\alpha, \beta, \gamma}  = 1$. In this case,  
probability $\tilde{p}_{0, 12} = \frac{1}{2}{\rm I}(\alpha = \beta) + 1{\rm I}(\alpha < \beta) > 0$ and $ \bar{v}_{\e, 1}/ \, _1\tilde{v}_{\e, 2}  \to 
\tilde{p}_{0, 12}$ as $\e \to 0$.  Therefore, according Theorem 10 and the remarks made in Subsection 9.4, the following equivalent relations take place,
\begin{equation}\label{hopcebnipabo}
E_{\e, \DD, 13} / \,  \bar{v}_{\e, 1} \to   \bar{E}_{0, \DD, 13} =  \bar{E}_{0, \DD, 23} = 1 \ {\rm as} \ \e \to 0, 
\end{equation}
and  
\begin{equation}\label{hopcebnipa}
E_{\e, \DD, 13} / \,  _1\tilde{v}_{\e, 2} \to   E_{0, \DD, 13} \ {\rm as} \ \e \to 0.
\end{equation}
where
\begin{equation}\label{cenaba}
E_{0, \DD, 13}  = \left\{
\begin{array}{ll}
\frac{1}{2}     & \text{if} \ \alpha = \beta, \vspace{2mm} \\
1  & \text{if} \  \alpha < \beta.
\end{array}
\right.
\end{equation}

In the above relations (\ref{hopcebnip}) and (\ref{hopcebnipa}), the expectations of hitting times $E_{\e, \DD, 13}$ normalised by  function $_1\tilde{v}_{\e, 2}$ (used in the corresponding weak convergence relations for distributions of hitting times, for the case where 
$\tilde{p}_{0, 12} > 0$) converge to the first moment $E_{0, \DD, 13}$ for the corresponding limiting distribution $G_{0, \DD, 13}(\cdot)$.

If $\alpha > \beta$, then  $\delta = \beta + \gamma < 1$, and, thus, $\bar{u}_{\alpha, \beta, \gamma}  = 1 {\rm I}(\alpha + \gamma < 1) + 
\frac{1}{3} {\rm I}(\alpha + \gamma = 1) + 0 {\rm I}(\alpha + \gamma > 1)$. In this case, 
probability $\tilde{p}_{0, 12} = 0$, and  $ \bar{v}_{\e, 1}/ \tilde{v}_{\e, 1}  \to (1 - \bar{u}_{\alpha, \beta, \gamma})^{-1}$ as $\e \to 0$. Therefore, according Theorem 10 and the remarks made in Subsection 9.4, the following 
relations take place, 
\begin{align}\label{honipa}
E_{\e, \DD, 13} / \bar{v}_{\e, 1} &  \to \bar{E}_{0, \DD, 13} \vspace{2mm} \nonumber \\
& =  (1 - \bar{u}_{\alpha, \beta, \gamma})
 \tilde{e}_{0, 13} \tilde{p}_{0, 13} + \bar{u}_{\alpha, \beta, \gamma} \bar{E}_{0, \DD, 23}
 \vspace{2mm} \nonumber \\
 & = (1 - \bar{u}_{\alpha, \beta, \gamma})  + \bar{u}_{\alpha, \beta, \gamma} 
 = 1 \ {\rm as} \ \e \to 0.
\end{align}
while, 
\begin{align}\label{honipany}
E_{\e, \DD, 13} / \tilde{v}_{\e, 1} & \to  \bar{E}_{0, \DD, 13}(1 - \bar{u}_{\alpha, \beta, \gamma})^{-1} \vspace{2mm} \nonumber \\
& = (1 - \bar{u}_{\alpha, \beta, \gamma})^{-1}  \ {\rm as} \ \e \to 0.
\end{align}

If $\beta + \gamma < 1 <  \alpha + \gamma$, then $\bar{u}_{\alpha, \beta, \gamma} = 0$. In this case,  $\bar{E}_{0, \DD, 13} = 
 \tilde{e}_{0, 13} \tilde{p}_{0, 13} =  E_{0, \DD, 13} = 1$ and  $\bar{E}_{0, \DD, 13}(1 - \bar{u}_{\alpha, \beta, \gamma})^{-1} = 
E_{0, \DD, 13} = 1$. 

Therefore, in the above relation (\ref{honipany}), the expectations of hitting times $E_{\e, \DD, 13}$ normalised by  function 
$\tilde{v}_{\e, 1}$ (used in the corresponding weak convergence relation for distributions of hitting times, for the case where 
$\tilde{p}_{0, 12} = 0$), converge to the first moment $E_{0, \DD, 13}$ for the corresponding limiting distribution $G_{0, \DD, 13}(\cdot)$.

If $\beta + \gamma  <  \alpha + \gamma = 1$, then $\bar{u}_{\alpha, \beta, \gamma} = \frac{1}{3}$. In this case,  $\bar{E}_{0, \DD, 13} = 
 \frac{2}{3}\tilde{e}_{0, 13} \tilde{p}_{0, 13} + \frac{1}{3} =  \frac{2}{3}E_{0, \DD, 13} + \frac{1}{3} = 1$, while  
 $\bar{E}_{0, \DD, 13}(1 - \bar{u}_{\alpha, \beta, \gamma})^{-1} = \tilde{e}_{0, 13} \tilde{p}_{0, 13} + \frac{1}{2} = 
 E_{0, \DD, 13} + \frac{1}{2} = \frac{3}{2}$. 
 
 Therefore, in the above relation (\ref{honipany}), the expectations of hitting times $E_{\e, \DD, 13}$ normalised by  function 
$\tilde{v}_{\e, 1}$ (used in the corresponding weak convergence relation for distributions of hitting times, for the case where 
$\tilde{p}_{0, 12} = 0$) converge to the constant, which differs of first moment $E_{0, \DD, 13}$ for the corresponding limiting distribution $G_{0, \DD, 13}(\cdot)$.

If $\beta + \gamma  <  \alpha + \gamma < 1$, then $\bar{u}_{\alpha, \beta, \gamma} = 1$. In this case, $\bar{E}_{0, \DD, 13} = 1$, 
while $\bar{E}_{0, \DD, 13}(1 - \bar{u}_{\alpha, \beta, \gamma})^{-1} = \infty$. 

 Therefore, in the above relation (\ref{honipany}), the expectations of hitting times $E_{\e, \DD, 13}$ normalised by  function 
$\tilde{v}_{\e, 1}$ (used in the corresponding weak convergence relation for distributions of hitting times for the case where 
$\tilde{p}_{0, 12} = 0$) converge to $\infty$. 

It also worth to note that finding of Laplace transforms for the corresponding limiting distributions in the weak convergence relations (\ref{hopceb}), (\ref{hopcena}), and (\ref{hopcenanu}) as well as the limits for expectations in asymptotic relations (\ref{hopcebnip}), (\ref{hopcebnipa}), and (\ref{honipa}) do require only performing of some recurrent rational algebraic transformations on the initial transition probabilities of  embedded Markov chains  and Laplace transforms of the transition times and their expectations.  

Analogous computations can be made in  the cases {\bf (ii)} $\delta = 1$ and {\bf (iii)} $\delta \in (1, \infty)$. In the case $\delta = 1$, states $1$ and $2$ are asymptotically equivalently  absorbing. Each state $1$ or $2$ can be used in the procedure of one-state exclusion  from phase space at the step 1 of the algorithm described in Section 4. In the case $\delta \in (1, \infty)$,  state $2$ is asymptotically less  absorbing than state $1$.  Thus, state $2$ should be used in the procedure of one-state exclusion  from phase space at the step 1 of the algorithm described in Section 4. \\

{\bf Appendix A: Families of Asymptotically Comparable Functions} \\

In Appendix A, we present results and comments concerning families of asymptotically comparable functions and discuss connection of these results with asymptotic recurrent algorithms of phase space reductions presented in the  paper.  \vspace{1mm} 

{\bf A1 Families of asymptotically comparable functions}.  In Appendix A,  we consider functions $h(\e)$ defined on  interval $(0, 1]$ and taking values in the interval  $(0, \infty)$. 

Let ${\cal H} =  \{  h(\e) \}$ be a non-empty  family of  such functions. \vspace{1mm}

{\bf Definition 1}. {\em ${\cal H} = \{ h(\cdot) \}$ is a  complete family of asymptotically comparable functions if: {\bf (1)} it is closed with respect to operations of summation, multiplication and division and {\bf (2)} there exists $\lim_{\e \to 0} h(\e) = a[h(\cdot)] \in [0, \infty]$, for any function $h(\cdot) \in {\cal H}$.} \vspace{1mm}

{\bf Definition 2}. {\em $\cal{H}'$ is a family of asymptotically comparable functions if it is a sub-family of 
some complete family of asymptotically comparable functions  ${\cal H}$}. \vspace{1mm}

Let functions $h_k(\e) \in {\cal H}', k = 1, 2, \ldots$, indices $\ell_k = \pm 1, k = 1, 2, \ldots$. Then, obviously,
the following limit exists, for any integer $n \geq 0$,
 \begin{equation}\label{limit}
 \lim_{\e \to 0} \prod_{k = 1}^n  h_{k}(\e)^{\ell_k}  = a[h_{1}(\cdot), \ell_1, \ldots, h_{n}(\cdot), \ell_n] \in [0, \infty].  
 \end{equation}
 
Indeed, according to  the obvious induction, functions $\prod_{k = 1}^n  h_{k}(\e)^{\ell_k}$ belongs to family ${\cal H}$.

The asymptotic relation (\ref{limit}) let one  give alternative variants of Definitions  1 and 2. 

 Let $\cal{H}'$ be a family of functions defined on  interval $(0, 1]$, taking values in the interval  $(0, \infty)$, for which the asymptotic relation (\ref{limit}) holds.
 
Let now define the family $[{\cal H}']$ of all functions given by the following relation, 

\begin{equation}\label{self}
h(\e) = \frac{\sum_{r' = 1}^{n'}  \prod_{k' = 1}^{n'_{r'}}  h'_{r', k'}(\e)^{\ell'_{r', k'}}}{\sum_{r'' = 1}^{n''}  
\prod_{k'' = 1}^{n''_{r''}}  h''_{r'', k''}(\e)^{\ell''_{r'', k''}}}, 
\end{equation}

\vspace{2mm}

\noindent where: (a) functions $h'_{r', k'}(\e) \in {\cal H}, 1 \leq k' \leq  n'_{r'} < \infty, 1 \leq r' \leq n' < \infty$; (b) functions $h''_{r'', k''}(\e) \in {\cal H}, 
1 \leq k'' \leq n''_{r''} < \infty, 1 \leq r'' \leq n'' < \infty$; (c) indices $\ell'_{r', k'} = \pm 1, 1 \leq k' \leq  n'_{r'} < \infty, 1 \leq r' \leq n' < \infty$, 
 (d) indices $\ell''_{r'', k''} = \pm 1, 1 \leq k'' \leq  n''_{r''} < \infty, 1 \leq r'' \leq n'' < \infty$.   \vspace{1mm}

Obviously,  ${\cal H}' \subseteq  [{\cal H}']$. Indeed, any function $h(\e) \in {\cal H}'$ can be represented in the form, 
$h(\e) = h(\e)^2 / h(\e) \in [{\cal H}]$.

Also constant $1 = 1(\e), \e \in (0, 1]$ belongs to family $[{\cal H}']$. Indeed, $1 = 1(\e)$ can be represented in the form, 
$1(\e) = h(\e) / h(\e) \in [{\cal H}']$.

We refer to family $[{\cal H}']$ as the closure of the family ${\cal H}'$. 

The following lemma presents the basic properties of  $[{\cal H}']$. \vspace{1mm}

{\bf Lemma A1}. {\em Let ${\cal H}'$ be a family of functions, for which the asymptotic relation {\rm (\ref{limit})} holds. Then its closure  $[{\cal H}']$ is a complete family of asymptotically comparable functions}. \vspace{1mm}

{\bf Proof}. Holding of assumption {\bf (a)} given in Definition 1  for family $[{\cal H}']$ and  the corresponding operational formulas for computing sums, products and quotients of functions from the family $[{\cal H}']$ are obvious. 

Holding of assumption {\bf (b)} given in Definition 1  for family $[{\cal H}']$ follows from the following relation, which holds for any function $h(\e)  \in [{\cal H}']$ given by relation ({\ref{self}),
\begin{align}\label{fewar}
& \ \ \lim_{\e \to 0}h(\e) = \lim_{\e \to 0} \frac{ \sum_{r' = 1}^{n'}  \prod_{k' = 1}^{n'_{r'}}  h'_{r', k'}(\e)^{\ell'_{r', k'}}}{\sum_{r'' = 1}^{n''}  
\prod_{k'' = 1}^{n''_{r''}}  h''_{r'', k''}(\e)^{\ell''_{r'', k''}}} \makebox[30mm]{}
\vspace{2mm} \nonumber \\
%\end{align*}
%\begin{align}
& \quad  \quad \quad \quad = \lim_{\e \to 0} \sum_{r' = 1}^{n'} \Big( \sum_{r'' = 1}^{n''}  
\prod_{k'' = 1}^{n''_{r''}}  h''_{r'', k''}(\e)^{\ell''_{r'', k''}} \prod_{k' = 1}^{n'_{r'}}  h'_{r', k'}(\e)^{- \ell'_{r', k'}} \Big)^{-1} \nonumber \\
& \quad \quad \quad  \quad = \sum_{r' = 1}^{n'} \Big( \sum_{r'' = 1}^{n''}  a[h''_{r'', k''}(\cdot), \ell''_{r'', k''}, k'' = 1, \ldots, n''_{r''}, \nonumber \\
& \quad \quad \quad \quad \quad \quad \quad \quad \quad \quad \ h'_{r', k'}(\cdot), -\ell'_{r', k'}, k' = 1, \ldots, n'_{r'}]\Big)^{-1} \in  [0, \infty]. 
\end{align}

Thus, $[{\cal H}']$ is, indeed, a complete family of asymptotically comparable functions. $\Box$ \vspace{1mm}

Since, ${\cal H}' \subseteq  [{\cal H}']$, we can say that ${\cal H}'$ is a family of asymptotically comparable functions if 
the asymptotic relation (\ref{limit})  holds for family  ${\cal H}'$.

Moreover,  the family of asymptotically comparable functions ${\cal H}'$ is complete if ${\cal H}' = [{\cal H}']$.} \vspace{1mm}

{\bf A2 Examples}. Let us present several examples of complete families of asymptotically comparable functions. \vspace{1mm}

{\bf Example 1}. The simplest example of a complete family of asymptotically comparable functions is the family of functions  ${\cal H}_1 =  \{  h(\cdot) \}$ such that, for any function  $h(\cdot) \in {\cal H}_1$, there exist constants $a_h > 0$ and $b_h \in (- \infty, \infty)$ such that the following asymptotic relation holds,
 \begin{equation}\label{repram}
 \frac{h(\e)}{a_h \e^{b_h}} \to 1 \ {\rm as} \ \e \to 0. 
 \end{equation}
 
 Let us  check that the family of functions ${\cal H}_1$ satisfies the assumptions formulated in Definition 1.
 
 Let $h_i(\cdot), i = 1, 2$ be two functions from the family ${\cal H}_1$, for which the asymptotic relation (\ref{repram}) holds, i.e.,
 $h_i(\e) / a_{i} \e^{b_{i}} \to 1$ as $\e \to 0$, for $i = 1, 2$, where  $a_{i} > 0, b_{i} \in (-\infty, \infty), i = 1, 2$. 
 
 In this case,  relation (\ref{repram}) holds for function $h(\e) = h_1(\e) + h_2(\e)$, with parameters $a = a_{1}{\rm I}(b_1 < b_2) +  (a_1 + a_2){\rm I}(b_1 = b_2)  + a_{2}{\rm I}(b_1 > b_2)$ and $b = b_{1} \wedge b_{2}$. 
 
 First, let us assume that  $b_1 < b_2$. Then,
 \begin{align}\label{sunm}
\frac{h(\e)}{a \e^{b}}  & = \frac{h_1(\e)}{a_{1} \e^{b_{1}}}  \frac{a_{1}}{a} \e^{b_{1} - b} + \frac{h_2(\e)}{a_{2}\e^{b_{2}}}  \frac{a_{2}}{a} \e^{b_{2} - b}  
\vspace{2mm} \nonumber \\
& = \frac{h_1(\e)}{a_{1} \e^{b_{1}}} + \frac{h_2(\e)}{a_{2}\e^{b_{2}}}  
 \frac{a_{2}}{a_1} \e^{b_{2} - b_1}
\to 1 \ {\rm as} \ \e \to 0.
 \end{align}
 
 Second, let $b_1 = b_2$. Then,
  \begin{align}\label{sunmop}
\frac{h(\e)}{a \e^{b}}  & = \frac{h_1(\e)}{a_{1} \e^{b_{1}}}  \frac{a_{1}}{a} \e^{b_{1} - b} + \frac{h_2(\e)}{a_{2}\e^{b_{2}}}  \frac{a_{2}}{a} \e^{b_{2} - b}  
\vspace{2mm} \nonumber \\
& = \frac{h_1(\e)}{a_{1} \e^{b_{1}}} \frac{a_1}{a_1 + a_2} 
+ \frac{h_2(\e)}{a_{2}\e^{b_{2}}}\frac{a_2}{a_1 + a_2}   
\to 1 \ {\rm as} \ \e \to 0.
 \end{align}
 
 The third case, where $b_1 > b_2$ is analogous to the first one considered above.  
 
 Thus, family ${\cal H}_1$ is closed with respect to the operation of summation.
 
Also, relation (\ref{repram}) holds for function $h_\ell(\e) = h_1(\e)  h_2(\e)^{\ell}$, with parameters $h_\ell = a_1 a_2^{\ell}$ and
$b_{\ell} = b_1 + b_2\ell$, for $\ell = \pm 1$. Indeed,
\begin{align}\label{sunma}
\frac{h(\e)}{a_\ell \e^{b_\ell}}  = \frac{h_1(\e)h_2(\e)^{\ell}}{a_{1} a_2^{\ell} \e^{b_{1} + b_2 \ell}} 
= \frac{h_1(\e)}{a_{1}  \e^{b_{1}}}  \big(\frac{h_2(\e)}{a_2 \e^{b_2}}\big)^{\ell}  \to   1 \ {\rm as} \ \e \to 0.
 \end{align}

Thus, family ${\cal H}_1$ is closed with respect to the operations of multiplication and division.

Finally, relation (\ref{repram}) obviously implies that, for any function $h(\cdot) \in {\cal H}_1$, for which the asymptotic relation (\ref{repram}) holds, i.e., $h(\e) / a_h \e^{b_h} \to 1$ as $\e \to 0$, for some $a_h > 0, b_h \in (-\infty, \infty)$, there exist the limit, 
 \begin{equation}\label{repramol}
  \lim_{\e \to 0}h(\e) =   \lim_{\e \to 0}a_h \e^{b_h} = a[h(\cdot)] = \left\{
 \begin{array}{ll}
 0 & \ \text{if} \ b_h > 0, \\
 a_h & \ \text{if} \ b_h = 0, \\
 \infty & \ \text{if} \ b_h < 0.  \\
 \end{array}
 \right.
  \end{equation}

Therefore, ${\cal H}_1$ is a complete family of asymptotically comparable functions.
 
Let ${\cal H}'_1$ b
non-empty sub-family of family ${\cal H}_1$. Then, for any functions $h_k(\e)  \in {\cal H}'_1$, for which relation (\ref{repram}) holds with parameters  $a_k > 0, b_k \in (-\infty, \infty)$, and  indices $\ell_k = \pm 1$, for $k = 1, 2, \ldots$,  there exists 
the following limit, for any $n \geq 1$, 
\begin{equation*}
 \lim_{\e \to 0} \prod_{k = 1}^n  h_{k}(\e)^{\ell_k}  =  \lim_{\e \to 0} \prod_{k = 1}^n a_k^{\ell_k}  
 \lim_{\e \to 0} \e^{\sum_{k = 1}^n \ell_k b_{k}} \makebox[34mm]{} 
\end {equation*}
\begin{equation}\label{compar}
\makebox[4mm]{}  = \left \{
\begin{array}{cll}
0 & \text{if} \  \sum_{k = 1}^n \ell_k b_{k} > 0, \vspace{2mm} \\
\prod_{k = 1}^n a^{\ell_k}_{k}  & \text{if} \  \sum_{k = 1}^n \ell_k b_{k} = 0, \vspace{2mm} \\
\infty & \text{if} \  \sum_{k = 1}^n \ell_k b_{k} < 0.
\end{array}
\right.
\end{equation}
\vspace{1mm}

{\bf Example 2}. Another example of a complete family of asymptotically comparable functions is the family of functions  ${\cal H}_2 =  
\{  h(\cdot) \}$ such that, for any function  $h(\cdot) \in {\cal H}_2$, there exists constants $a_h > 0$ and $b_h, c_h \in (- \infty, \infty)$ such that
the following asymptotic relation holds,
 \begin{equation}\label{repramno}
 \frac{h(\e)}{a_h \e^{b_h} e^{-c_h \e^{-1}}} \to 1 \ {\rm as} \ \e \to 0. 
 \end{equation}
 
 Let us  check that the family of functions ${\cal H}_2$ satisfies the assumptions formulated in Definition 1.
 
 Let $h_i(\cdot), i = 1, 2$ be two functions from the family ${\cal H}_1$, for which the asymptotic relation (\ref{repramno}) holds, i.e.,
 $h_i(\e) / a_{i} \e^{b_{i}}e^{-c_i \e^{-1}} \to 1$ as $\e \to 0$, for $i = 1, 2$, where  $a_{i} > 0, b_{i}, c_i \in (-\infty, \infty), i = 1, 2$. 
 
 In this case,  relation (\ref{repramno}) holds for function $h(\e) = h_1(\e) + h_2(\e)$, with parameters $a = a_{1}({\rm I}(c_1 < c_2)  + {\rm I}(c_1 = c_2, b_1 < b_2)) 
+  (a_1 + a_2){\rm I}(c_1 = c_2, b_1 = b_2) + a_{2}({\rm I}(c_1 > c_2)  + {\rm I}(c_1 = c_2, b_1 > b_2))$, $b = b_{1}{\rm I}(c_1 < c_2) +  
(b_1 \wedge b_2){\rm I}(c_1 = c_2) + b_2 {\rm I}(c_1 > c_2) $, and  $c = c_{1} \wedge c_{2}$. 

First, let us assume that $c_1 < c_2$ or $c_1 = c_2, b_1 < b_2$. Then,
 \begin{align}\label{sunmas}
\frac{h(\e)}{a \e^{b}e^{- c \e^{-1}}}  &  = 
\frac{h_1(\e)}{a_{1} \e^{b_{1}} e^{- c_1 \e^{-1}}}  
\frac{a_{1}}{a} \e^{b_{1} - b}e^{- (c_1 - c) \e^{-1}} \vspace{2mm} \nonumber \\
& \quad + \frac{h_2(\e)}{a_{2} \e^{b_{2}} e^{- c_2 \e^{-1}}}  
\frac{a_{2}}{a} \e^{b_{2} - b}e^{- (c_2 - c) \e^{-1}} \vspace{2mm} \nonumber \\
& = \frac{h_1(\e)}{a_{1} \e^{b_{1}} e^{- c_1 \e^{-1}}}  
+ \frac{h_2(\e)}{a_{2} \e^{b_{2}} e^{- c_2 \e^{-1}}}  
\frac{a_{2}}{a_1} \e^{b_{2} - b_1}e^{- (c_2 - c_1) \e^{-1}} 
\vspace{2mm} \nonumber \\
& \to   1 \ {\rm as} \ \e \to 0.
 \end{align}
 
 Second, let $c_1 = c_2, b_1 = b_2$. Then,
 \begin{align*}
\frac{h(\e)}{a \e^{b}e^{- c \e^{-1}}}  &  = 
\frac{h_1(\e)}{a_{1} \e^{b_{1}} e^{- c_1 \e^{-1}}}  
\frac{a_{1}}{a} \e^{b_{1} - b}e^{- (c_1 - c) \e^{-1}} \vspace{2mm} \nonumber \\
& \quad + \frac{h_2(\e)}{a_{2} \e^{b_{2}} e^{- c_2 \e^{-1}}}  
\frac{a_{2}}{a} \e^{b_{2} - b}e^{- (c_2 - c) \e^{-1}} \makebox[30mm]{} \vspace{2mm} \nonumber \\
\end{align*}
 \begin{align}\label{sunmasta}
& = \frac{h_1(\e)}{a_{1} \e^{b_{1}} e^{- c_1 \e^{-1}}} \frac{a_1}{a_1 + a_2}  
+ \frac{h_2(\e)}{a_{2} \e^{b_{2}} e^{- c_2 \e^{-1}}}  
\frac{a_{2}}{a_1 + a_2}  \vspace{2mm} \nonumber \\
& \to   1 \ {\rm as} \ \e \to 0.
 \end{align}

The third case, where $c_1 > c_2$ or $c_1 = c_2, b_1 > b_2$. is analogous to the first one considered above.  

Thus, family ${\cal H}_2$ is closed with respect to the operation of summation.
 
 Also, relation (\ref{repramno}) holds for function $h_\ell(\e) = h_1(\e)  h_2(\e)^{\ell}$, with parameters $h_\ell = a_1 a_2^{\ell}$, $b_{\ell} = b_1 + b_2\ell$, and $c_{\ell} = c_1 + c_2\ell$ for $\ell = \pm 1$. Indeed,
\begin{align}\label{sunmans}
\frac{h(\e)}{a_\ell \e^{b_\ell} e^{- c_\ell \e^{-1}}}  & = 
\frac{h_1(\e)h_2(\e)^{\ell}}{a_{1} a_2^{\ell} \e^{b_{1} + b_2 \ell} e^{- (c_1 + c_2\ell)\e^{-1}}} 
\vspace{2mm} \nonumber \\
& = \frac{h_1(\e)}{a_{1}  \e^{b_{1}} e^{- c_1 \e^{-1}}}  
\big(\frac{h_2(\e)}{a_2 \e^{b_2}e^{- c_2 \e^{-1}}}\big)^{\ell}  \to   1 \ {\rm as} \ \e \to 0.
 \end{align}
 
 Thus, family ${\cal H}_2$ is closed with respect to the operations of multiplication and division.

Finally, relation (\ref{repramno}) obviously implies that, for any function $h(\cdot) \in {\cal H}_2$, for which the asymptotic relation (\ref{repramno}) holds, i.e., $h(\e) / a_h \e^{b_h} e^{- c_h\e^{-1}} \to 1$ as $\e \to 0$, for some $a_h > 0, b_h, c_h \in (-\infty, \infty)$, there exist the limit, 
 \begin{equation*}
  \lim_{\e \to 0}h(\e) =   \lim_{\e \to 0}a_h \e^{b_h} e^{- c_h\e^{-1}}  \makebox[66mm]{} 
  \end{equation*}
 \begin{equation}\label{repramolnsmo}
 = a[h(\cdot)] = \left\{
 \begin{array}{ll}
 0 & \ \text{if} \ c_h > 0 \ \text{or} \ c_h = 0, \, b_h > 0, \vspace{2mm} \\
 a_h & \ \text{if} \ c_h = 0, \, b_h = 0, \vspace{2mm} \\
 \infty & \ \text{if} \ c_h < 0 \ \text{or} \  c_h = 0, \, b_h < 0.  \\
 \end{array}
 \right.
  \end{equation}
  
  Therefore, ${\cal H}_2$ is a complete family of asymptotically comparable functions.
 
Let ${\cal H}'_2$ be some non-empty sub-family of family ${\cal H}_2$. Then, for any functions $h_k(\e)  \in {\cal H}'_2$, for which relation (\ref{repramno}) holds with parameters  $a_k > 0, b_k, c_k \in (-\infty, \infty)$, and  indices $\ell_k = \pm 1$, for $k = 1, 2, \ldots$,  there exists 
the following limit, for any $n \geq 1$, 
\begin{align*}
 \lim_{\e \to 0} \prod_{k = 1}^n  h_{k}(\e)^{\ell_k}  & = \prod_{k = 1}^n a^{\ell_k}_{k}  
 \lim_{\e \to 0} \e^{ \sum_{k = 1}^n \ell_k b_{k}} e^{- \sum_{k =1}^n  \ell_k c_{k} \e^{-1}} 
 \makebox[25mm]{} 
\end {align*}
\begin{equation}\label{compart}
\makebox[21mm]{} = \left \{
\begin{array}{cll}
0 & \text{if} \  \sum_{k = 1}^n \ell_k c_{k} > 0 \vspace{2mm}  \\
& \quad \text{or} \ \sum_{k = 1}^n \ell_k c_{k} = 0, \ 
\sum_{k = 1}^n \ell_k b_{k} > 0, \vspace{2mm} \\
\prod_{k = 1}^n a^{\ell_k}_{k}  & \text{if} \ \sum_{k = 1}^n \ell_k c_{k} = 0,  \ 
\sum_{k = 1}^n \ell_k b_{k} = 0, \vspace{2mm}  \\
\infty & \text{if} \ \sum_{k = 1}^n \ell_k c_{k} < 0 \vspace{2mm}  \\
& \quad  \text{or} \ \sum_{k = 1}^n \ell_k c_{_k} = 0, \ 
\sum_{k = 1}^n \ell_k b_{k} < 0.
\end{array}
\right.
\end{equation}
\vspace{1mm}

{\bf Example 3}. One more  example of a complete family of asymptotically comparable functions is the family of functions  ${\cal H}_3 =  \{  h(\cdot) \}$ such that, for any function  $h(\cdot) \in {\cal H}_3$, there exists constants $a_h > 0$ and $b_h, d_h \in (- \infty, \infty)$ such that
the following asymptotic relation holds,
 \begin{equation}\label{repramnona}
 \frac{h(\e)}{a_h \e^{b_h} (1 + \ln \e^{-1})^{- d_h}} \to 1 \ {\rm as} \ \e \to 0. 
 \end{equation}
 
 Let us  check that the family of functions ${\cal H}_3$ satisfies the assumptions formulated in Definition 1.

 Let $h_i(\cdot), i = 1, 2$ be two functions from the family ${\cal H}_3$, for which the asymptotic relation (\ref{repramno}) holds, i.e.,
 $h_i(\e) / a_{i} \e^{b_{i}} (1 + \ln \e^{-1})^{- d_i} \to 1$ as $\e \to 0$, for $i = 1, 2$, where  $a_{i} > 0, b_{i}, d_i \in (-\infty, \infty), i = 1, 2$. 
 
 In this case,  relation (\ref{repramnona}) holds for function $h(\e) = h_1(\e) + h_2(\e)$, with parameters $a = a_{1}({\rm I}(b_1 < b_2)  + {\rm I}(b_1 = b_2, d_1 < d_2)) 
+  (a_1 + a_2){\rm I}(b_1 = b_2, d_1 = d_2) + a_{2}({\rm I}(b_1 > b_2)  + {\rm I}(b_1 = b_2, d_1 > d_2))$, 
$b = b_{1} \wedge b_{2}$, and $d = d_{1}{\rm I}(b_1 < b_2) +  
(d_1 \wedge d_2){\rm I}(b_1 = b_2) + d_2 {\rm I}(b_1 > b_2) $.  

First, let us assume that $b_1 < b_2$ or $b_1 = b_2, d_1 < d_2$. Then,
 \begin{align}\label{sunmasmom}
\frac{h(\e)}{a \e^{b}(1 + \ln \e^{-1})^{- d}}  &  = 
\frac{h_1(\e)}{a_{1} \e^{b_{1}} (1 + \ln \e^{-1})^{- d_1}}  
\frac{a_{1}}{a} \e^{b_{1} - b} (1 + \ln \e^{-1})^{- (d_1 - d)} \vspace{2mm} \nonumber \\
& \quad + \frac{h_2(\e)}{a_{2} \e^{b_{2}} (1 + \ln \e^{-1})^{- d_2}}  
\frac{a_{2}}{a} \e^{b_{2} - b}(1 + \ln \e^{-1})^{- (d_2 - d)}  \vspace{2mm} \nonumber \\
& = \frac{h_1(\e)}{a_{1} \e^{b_{1}} (1 + \ln \e^{-1})^{- d_1}}   \vspace{2mm} \nonumber \\
%\end{align*}
%\begin{align}
& \quad + \frac{h_2(\e)}{a_{2} \e^{b_{2}}(1 + \ln \e^{-1})^{- d_2}}   
\frac{a_{2}}{a_1}  \e^{b_{2} - b_{1}}(1 + \ln \e^{-1})^{- (d_2 - d_1)} 
\vspace{3mm} \nonumber \\
& \to   1 \ {\rm as} \ \e \to 0.
 \end{align}

Second, let $b_1 = b_2, d_1 = d_2$. Then,
 \begin{align}\label{sunmasmomu}
\frac{h(\e)}{a \e^{b}(1 + \ln \e^{-1})^{- d}}  &  = 
\frac{h_1(\e)}{a_{1} \e^{b_{1}} (1 + \ln \e^{-1})^{- d_1}}  
\frac{a_{1}}{a} \e^{b_{1} - b} (1 + \ln \e^{-1})^{- (d_1 - d)} \vspace{2mm} \nonumber \\
& \quad + \frac{h_2(\e)}{a_{2} \e^{b_{2}} (1 + \ln \e^{-1})^{- d_2}}  
\frac{a_{2}}{a} \e^{b_{2} - b}(1 + \ln \e^{-1})^{- (d_2 - d)}  \vspace{2mm} \nonumber \\
& = \frac{h_1(\e)}{a_{1} \e^{b_{1}} (1 + \ln \e^{-1})^{- d_1}} \frac{a_1}{a_1 + a_2}  \vspace{2mm} \nonumber \\
%\end{align*}
%\begin{align}
& \quad + \frac{h_2(\e)}{a_{2} \e^{b_{2}}(1 + \ln \e^{-1})^{- d_2}}   
\frac{a_{2}}{a_1 + a_2}  
\vspace{2mm} \nonumber \\
& \to   1 \ {\rm as} \ \e \to 0.
 \end{align}
 
 Thus, family ${\cal H}_3$ is closed with respect to the operation of summation.
 
 Also, relation (\ref{repramnona}) holds for function $h_\ell(\e) = h_1(\e)  h_2(\e)^{\ell}$, with parameters $h_\ell = a_1 a_2^{\ell}$, $b_{\ell} = b_1 + b_2\ell$, and $d_{\ell} = d_1 + d_2\ell$ for $\ell = \pm 1$. Indeed,
\begin{align}\label{sunmansby}
\frac{h(\e)}{a_\ell \e^{b_\ell} (1 + \ln \e^{-1})^{- d_\ell}}  & = 
\frac{h_1(\e)h_2(\e)^{\ell}}{a_{1} a_2^{\ell} \e^{b_{1} + b_2 \ell} (1 + \ln \e^{-1})^{- (d_1 + d_2 \ell)}}  
\vspace{2mm} \nonumber \\
%\end{align*}
%\begin{align}
& = \frac{h_1(\e)}{a_{1}  \e^{b_{1}} (1 + \ln \e^{-1})^{- d_1}}  
\big(\frac{h_2(\e)}{a_2 \e^{b_2}(1 + \ln \e^{-1})^{- d_2}} \big)^{\ell}  \vspace{2mm} \nonumber \\ 
&   \to   1 \ {\rm as} \ \e \to 0.
 \end{align}
 
 Thus, family ${\cal H}_3$ is closed with respect to the operations of multiplication and division.

Finally, relation (\ref{repramnona}) obviously implies that, for any function $h(\cdot) \in {\cal H}_3$, for which the asymptotic relation (\ref{repramnona}) holds, i.e., 
$h(\e) / a_h \e^{b_h} (1 + \ln \e^{-1})^{- d_h} \to 1$ as $\e \to 0$, for some $a_h > 0, b_h, d_h \in (-\infty, \infty)$, there exist the limit, 
 \begin{equation*}
  \lim_{\e \to 0}h(\e) =   \lim_{\e \to 0}a_h \e^{b_h} (1 + \ln \e^{-1})^{- d_h}   \makebox[50mm]{} 
  \end{equation*}
 \begin{equation}\label{reprmo}
\makebox[5mm]{} = a[h(\cdot)] = \left\{
 \begin{array}{ll}
 0 & \ \text{if} \ b_h > 0 \ \text{or} \ b_h = 0, \, d_h > 0, \\
 a_h & \ \text{if} \ b_h = 0, \, d_h = 0, \\
 \infty & \ \text{if} \ b_h < 0 \ \text{or} \  b_h = 0, \, d_h < 0.  \\
 \end{array}
 \right.
  \end{equation}
  
  Therefore, ${\cal H}_3$ is a complete family of asymptotically comparable functions.
  
 Let ${\cal H}'_3$ be some non-empty sub-family of family ${\cal H}_3$. Then, for any functions $h_k(\e)  \in {\cal H}'_3$, for which relation (\ref{repramnona}) holds with parameters  $a_k > 0, b_k, d_k \in (-\infty, \infty)$, and  indices $\ell_k = \pm 1$, for $k = 1, 2, \ldots$,  there exists 
the following limit, for any $n \geq 1$, 
\begin{align*}
& \lim_{\e \to 0} \prod_{k = 1}^n  h_{k}(\e)^{\ell_k}  \vspace{2mm} \\
& \quad \quad = \prod_{k = 1}^n a^{\ell_k}_{k}  \lim_{\e \to 0} 
\e^{ \sum_{k = 1}^n \ell_k b_{k}} (1 + \ln \e^{-1})^{- \sum_{k =1}^n  \ell_k d_{k}}  \makebox[12mm]{} 
\end {align*}
\begin{equation}\label{comparutt}
\makebox[21mm]{} = \left \{
\begin{array}{cll}
0 & \text{if} \  \sum_{k = 1}^n \ell_k b_{k} > 0 \vspace{2mm}  \\
& \quad \text{or} \ \sum_{k = 1}^n \ell_k b_{k} = 0, \ 
\sum_{k = 1}^n \ell_k d_{k} > 0, \vspace{2mm} \\
\prod_{k = 1}^n a^{\ell_k}_{k}  & \text{if} \ \sum_{k = 1}^n \ell_k b_{k} = 0,  \ 
\sum_{k = 1}^n \ell_k d_{k} = 0, \vspace{2mm}  \\
\infty & \text{if} \ \sum_{k = 1}^n \ell_k b_{k} < 0 \vspace{2mm}  \\
& \quad  \text{or} \ \sum_{k = 1}^n \ell_k b_{_k} = 0, \ 
\sum_{k = 1}^n \ell_k d_{\theta_k} < 0.
\end{array}
\right.
\end{equation}

Readers can readily construct other examples based on asymptotic representations analogous to 
(\ref{repram}), (\ref{repramno}), and (\ref{repramnona}), which would involve and combine power, exponential, logarithmic and other types of functions. \vspace{1mm}

{\bf A3 Asymptotic comparability conditions ${\bf G}$ and ${\bf H}$}. Let remind that it is assumed that condition ${\bf B}$ holds.

Let us first comment condition ${\bf G}$, which requires that transition probabilities $p_{\cdot, ij}, j \in \YY_{1, i}, i \in \overline{\DD}$ belong to some complete family of asymptotically comparable functions. 

In this case ${\cal P} = \{ p_{\cdot, ij}, j \in \YY_{1, i}, i \in \overline{\DD} \}$ is a family of asymptotically comparable functions. Thus, asymptotic relation  (\ref{limit}) holds for family 
${\cal P}$. This implies, in an obvious way, that condition ${\bf C}_{n}$ holds for any $n \geq 0$,  and, thus, proves Lemma 25.

Let us assume that ${\cal H}_1$ is the complete family of asymptotically comparable functions appearing in condition ${\bf G}$, i.e., it is assumed that, for every $j \in \YY_{1, i}, i \in \overline{\DD}$, there exist constants $a_{ij} > 0$ and 
$b_{ij} \in [0,  \infty)$ such that the following asymptotic relation holds,
 \begin{equation}\label{repramnop}
 \frac{p_{\e, ij}}{a_{ij} \e^{b_{ij}}} \to 1 \ {\rm as} \ \e \to 0. 
 \end{equation}
 
Note, that, in this case, relation (\ref{repram}) takes the form, where parameters $b_{ij} \in [0, \infty)$, since transition probabilities $p_{\e, ij} \in (0, 1], \e \in (0, 1]$, for $j \in \YY_{1, i}, i \in \overline{\DD}$.

Relation (\ref{repramnop}) implies that, for $j \in \YY_{1, i}, i \in \overline{\DD}$, 
\begin{equation}\label{repnop}
 \lim_{\e \to 0} p_{\e, ij} = p_{0, ij} = \left\{
 \begin{array}{ll}
 0 & \ \text{if} \ b_{ij} > 0, \\
 a_{ij} & \ \text{if} \ b_{ij} = 0.
 \end{array}
 \right.
 \end{equation}

Moreover, since matrix ${\mathbf P}_\e = \|  p_{\e, ij}  \|$ is stochastic, for every $\e \in (0, 1]$, the following relation holds, for $i \in \overline{\DD}$, 
\begin{equation}\label{repnope}
\sum_{j \in  \YY_{1, i}: \, b_{ij} = 0} a_{ij}  = 1. 
 \end{equation}
 
 An important role in the algorithms of phase space reduction is played by probabilities 
 $\bar{p}_{\e, ii} = 1 - p_{\e, ii} = \sum_{j \in  \YY_{1, i}, j \neq i} p_{\e, ij}$.  Condition ${\bf B}$ implies that $\YY_{1, i} \setminus \{ i \} \neq  \emptyset$  and, thus,  probabilities $\bar{p}_{\e, ii} > 0, \e \in (0, 1]$, for $i \in \overline{\DD}$. Moreover, since   ${\cal H}_1$ is closed with respect to summation operation,
 functions $\bar{p}_{\cdot, ii}, i \in \overline{\DD}$  belong to ${\cal H}_1$.
 
 Relation (\ref{repram}) takes for the above functions the following form, for $ i \in \overline{\DD}$
 \begin{equation}\label{asert}
\frac{\bar{p}_{\e, ii}}{\bar{a}_{ii} \e^{\bar{b}_{ii}}}   \to 1 \ {\rm as} \ \e \to 0,
 \end{equation}
 where $\bar{a}_{ii} = \sum_{j \in \YY_{1, i}: \,  j \neq i, b_{ij} = \bar{b}_{ii}} a_{ij}$ and $\bar{b}_{ii} = \wedge_{j \in \YY_{1, i}: \,  j \neq i} \, b_{ij}$. Note that  $\bar{a}_{ii}  > 0, \bar{b}_{ii} \geq 0$,  for $i \in \overline{\DD}$.  

 Relation (\ref{asert}) implies that, for $i \in \overline{\DD}$, 
\begin{equation}\label{repnoptyr}
 \lim_{\e \to 0} \bar{p}_{\e, ii} = \bar{p}_{0, ii} = \left\{
 \begin{array}{ll}
 0 & \ \text{if} \ \bar{b}_{ii} > 0, \\
 \bar{a}_{ii} & \ \text{if} \ \bar{b}_{ii} = 0.
 \end{array}
 \right.
 \end{equation}
 
 It is worth to note that quantities, $\bar{p}_{\e, ii}, i \in \overline{\DD}$, \,  $\tilde{p}_{\e, ij} = \frac{p_{\e, ij}}{1 - p_{\e, ii}}, j \in \YY_{1, i}, j \neq i, i \in \overline{\DD}$, \, $_kp_{\e, ij} = \tilde{p}_{\e, ij}  + \tilde{p}_{\e, ik} \tilde{p}_{\e, kj}, j \in \, _k\YY_{1, i}, i \in \, _k\overline{\DD}$, \, $\hat{q}_{\e}[ij] = \frac{\tilde{p}_{\e, ij}}{_kp_{\e, ij}}, j \in \, _k\YY_{1, i}, i \in \, _k\overline{\DD}$,  and $_k\tilde{p}_{\e, ij} = \frac{_kp_{\e, ij}}{1 - \, _kp_{\e, ii}}, j \in \, _k\YY_{1, i}, j \neq i, i \in \, _k\overline{\DD}$, all belong to ${\cal H}_1$, since this family of functions  is closed with respect to summation, multiplication and division operations.

The above quantities are involved in the procedures of removing of virtual transition and one-state reduction of phase space described in Sections 3 -- 5.

Therefore, computing of the  parameters $a$ and $b$ in the asymptotic relation (\ref{repram}) for all the above quantities, do require to apply  the corresponding arithmetic operations based or relations (\ref{sunm}) -- (\ref{sunma}).  

The same is related to analogues of the above quantities in the following recurrent steps of application the  above procedures, which are described in Section 6.
 
 Let us also comment condition ${\bf H}$, which requires that transition probabilities $p_{\cdot, ij}, j \in \YY_{1, i}, i \in \overline{\DD}$ and the initial normalisation functions 
 $v_{\cdot, i},  i \in \overline{\DD}$ belong to some complete family of asymptotically comparable functions. 

In this case ${\cal PV} = \{ p_{\cdot, ij}, v_{\cdot, i}, j \in \YY_{1, i}, i \in \overline{\DD} \}$ is a family of asymptotically comparable functions. Thus, asymptotic relation  (\ref{limit}) holds for family 
${\cal PV}$. This implies, in an obvious way, that condition ${\bf F}_{n}$ holds for any $n \geq 0$ and, thus, proves Lemma 27.

Let us assume that ${\cal H}_1$ is the complete family of asymptotically comparable functions appearing in condition ${\bf H}$, i.e., it is assumed that the asymptotic relation (\ref{repramnop}) holds for transition probabilities $p_{\cdot, ij}, j \in \YY_{1, i}, i \in \overline{\DD}$ and, also, for every $i \in \overline{\DD}$, there exist constants $a_{i} \geq 1$ and 
$b_{i} \in (- \infty, 0]$ such that the following asymptotic relation holds,
 \begin{equation}\label{repramnopnu}
 \frac{v_{\e, i}}{a_{i} \e^{b_{i}}} \to 1 \ {\rm as} \ \e \to 0, 
 \end{equation}

Note, that, in this case, relation (\ref{repram}) takes the form, where parameters $a_i \geq 1, b_{i} \in 
(- \infty, 0]$, since it was assumed that the initial normalisation functions  $v_{\e, i} \in [1, \infty), \, \e \in (0, 1]$, for $i \in \overline{\DD}$.

In this case, the normalisation functions $\tilde{v}_{\e, i} = (1 - p_{\e, ii})^{-1}v_{\e, i}, i \in \overline{\DD}$, 
$\tilde{w}_{\e, ji}  = \frac{\tilde{v}_{\e, j}}{\tilde{v}_{\e, i}} =  
\frac{(1 - p_{\e, jj})^{-1} v_{\e, j}}{(1 - p_{\e, ii})^{-1} v_{\e, i}}, i, j \in \overline{\DD}$ and  $_k\tilde{v}_{\e,  i} = (1 - \, _kp_{\e, ii})^{-1} \, _kv_{\e,  i} = (1 - \, _kp_{\e, ii})^{-1}(1 - p_{\e, ii})^{-1} v_{\e, i},  i \in \, _k\overline{\DD}$  also belong to ${\cal H}_1$, since this family of functions  is closed with respect to summation, multiplication and division operations.

The above quantities are used for computing the corresponding normalisation functions involved in the procedures of removing of virtual transition and one-state reduction of phase space described in Sections 3 -- 5.

Therefore, computing of the  parameters $a$ and $b$ in the asymptotic relation (\ref{repram}) for all the above quantities, do require to apply of the corresponding arithmetic operations based or relations (\ref{sunm}) -- (\ref{sunma}).  

The same is related to analogues of the above quantities in the following recurrent steps of application the  above procedures, which are described in Section 6, as well as to the corresponding final normalisation functions $_{\bar{k}_{\bar{m}}}\check{v}_{\e, i}, i \in \overline{\DD}$, their quotients $_{\bar{k}_{\bar{m}}}\check{v}_{\e, j}/ \, _{\bar{k}_{\bar{m}}}\check{v}_{\e, i}, i, j \in \overline{\DD}$, etc., which are involved in the weak convergence relations for hitting times given in theorems presented in Sections 7 -- 9.

Analogous comments can be made in the cases, where conditions ${\bf G}$ and ${\bf H}$ are based on the complete families of asymptotically comparable functions  ${\cal H}_2$ or ${\cal H}_3$. 

In conclusion, we would like to note that, in the case where the comparability condition ${\bf H}$ is based on one of the families, ${\cal H}_1$, ${\cal H}_2$, or ${\cal H}_3$,   computing of the limiting Laplace transforms and expectations  for hitting times do require only recurrent application of some rational transformations for the initial limiting transition probabilities $p_{0, ij}$, Laplace transforms $\phi_{0, ij}(s)$ and expectations $e_{0,ij}$. Coefficients of these rational transformations can be computed with the use of recurrent rational formulas based on operational rules (for computing limits for sum, products and quotients) for functions from the corresponding complete family of asymptotically 
comparable functions. The total number of required operations is  of the order $O(\bar{m}^3)$. 

The above estimate could be expected,  since, the asymptotic phase space reduction algorithm presenting on the present paper can be considered as some kind  of asymptotic stochastic variant of the Gauss elimination method.\\

{\bf Conclusion} \\

 This paper presents results of the complete asymptotic analysis related to distributions and expectations of hitting times for  singularly perturbed ergodic type finite semi-Markov processes. As we think,  the new asymptotic recurrent  algorithms of phase space reduction  presented in the paper have their own value.
 
 At the same time, we consider the results of this paper as intermediate ones for getting ergodic and quasi-ergodic theorems for perturbed  regenerative  processes modulated by singularly perturbed semi-Markov processes. Asymptotics of distributions and expectations of hitting times (in particular, return times) for modulating semi-Markov processes plays the key role in such theorems.
 
That is why we tried and, as we think, succeeded to get conditions of convergence for hitting times, which are based on minimal convergence conditions for transition characteristics of perturbed semi-Markov processes typical for ergodic theorems,  as well as to effectively compute normalisation functions, which play the role of asymptotic time compression factors,  and ``switching''  parameters determining  forms of the corresponding stationary and quasi-stationary distributions in ergodic theorems. 
 
 In the recent paper [75], the complete classification of ergodic theorems, based on about 20  short, long and super-long time ergodic theorems for perturbed  alternating regenerative  processes modulated by regularly, singularly and super-singularly  perturbed two-states semi-Markov  processes, has been given.
 
Our conjecture is that the model of perturbed regenerative processes modulated by  regularly and singularly perturbed multi-states semi-Markov processes can be reduced to the former one by aggregating regenerative periods. New regeneration moments can be defined as sequential hitting times  into the special two-states set by the modulating semi-Markov process. This two-states set should includes two the most asymptotically absorbing states for the modulating semi-Markov process. In this way, as we hope,  the results of the present paper let us expand ergodic theorems given in paper [75] to the model of perturbed regenerative processes modulated by  singularly perturbed multi-states semi-Markov processes. We plan to present the corresponding results in the near future.

\end{document}